\documentclass[10pt]{amsart}
\usepackage[margin=2.5cm]{geometry}
\usepackage{tikz-cd}
\usepackage{graphicx}
\usepackage{mathdots}
\usepackage{amssymb}
\usepackage{amsmath}
\usepackage{bbm}
\usepackage{relsize}
\usepackage{subfig, graphicx}
\usepackage{mathrsfs}
\usepackage{amssymb, amsthm, indentfirst}
\usepackage{relsize}
\usepackage{hyperref}
\usepackage{stmaryrd}
\usepackage{lineno}
\usepackage{mathabx}
\numberwithin{equation}{section}
\usepackage{color}
\theoremstyle{plain}
\newtheorem{thm}{Theorem}[section]

\newtheorem{conj}[thm]{Conjecture}
\newtheorem{lemma}[thm]{Lemma}
\newtheorem{prop}[thm]{Proposition}

\newtheorem*{conja}{Conjecture A}
\newtheorem*{conjb}{Conjecture B}
\newtheorem*{thma}{Theorem A}
\newtheorem*{thmb}{Theorem B}
\theoremstyle{definition}
\newtheorem{definition}[thm]{Definition}
\newtheorem*{example*}{Example}
\newtheorem{rmk}[thm]{Remark}
\def\Gal{\operatorname{Gal}}

\newtheorem*{hypothesis*}{Hypothesis}

\author{Shih-Yu Chen}
\address{Department of Mathematics, National Tsing Hua University, 101, Section 2, Kuang-Fu Road, Hsinchu 300, Taiwan, Republic of China}
\email{sychen.math@gmail.com}

\def\GL{{\rm{GL}}}
\def\GSp{{\rm GSp}}

\def\Sp{{\rm Sp}}

\def\A{{\mathbb A}}
\def\C{{\mathbb C}}
\def\E{{\mathbb E}}
\def\F{{\mathbb F}}
\def\K{{\mathbb K}}

\def\R{{\mathbb R}}
\def\Q{{\mathbb Q}}
\def\Z{{\mathbb Z}}

\def\<{\langle}
\def\>{\rangle}

\def\bp{\begin{pmatrix}}
\def\ep{\end{pmatrix}}

\def\<{\langle}
\def\>{\rangle}

\def\GL{\operatorname{GL}}

\def\GSp{\operatorname{GSp}}

\def\Sp{\operatorname{Sp}}

\def\1{\mathbf{1}}

\def\itPi{\mathit{\Pi}}
\def\itPsi{\mathit{\Psi}}
\def\itSigma{\mathit{\Sigma}}

\makeatletter
\newcommand{\extp}{\@ifnextchar^\@extp{\@extp^{\,}}}
\def\@extp^#1{\mathop{\bigwedge\nolimits^{\!#1}}}
\makeatother

\makeatletter
\newcommand{\exterior}[1]{\mathop{\mathpalette\exterior@{#1}}}
\newcommand{\exterior@}[2]{%
  \raisebox{\depth}{%
  \fontsize{\sf@size}{0}%
  \m@th
  $\ifx#1\displaystyle\textstyle\else#1\fi\bigwedge$}%
  ^{\mspace{-2mu}#2}%
  \kern-\scriptspace
}
\makeatother

\title{Algebraicity of ratios of Rankin-Selberg L-functions and applications to Deligne's conjecture}
\begin{document}

\begin{abstract}
In this paper, we prove Deligne's conjecture on the algebraicity of the critical values of symmetric power $L$-functions associated with modular forms of weight at least 5. 
We also establish new cases of Blasius' conjecture on the algebraicity of the critical values of tensor product $L$-functions associated with modular forms. Additionally, we prove an algebraicity result for the critical values of Rankin--Selberg $L$-functions for $\GL_n \times \GL_2$ in the unbalanced case, which extends the previous results of Furusawa and Morimoto for ${\rm SO}(V) \times \GL_2$.
These results are applications of our main theorem on the algebraicity of cross ratios of special values of Rankin--Selberg $L$-functions.
%In this paper, we prove Deligne's conjecture on the algebraicity of critical values of symmetric power $L$-functions associated to modular forms of weight at least $5$.
%We also prove new cases of Blasius' conjecture on the algebraicity of critical values of tensor product $L$-functions associated to modular forms, and an algebraicity result on critical values of Rankin--Selberg $L$-functions for $\GL_n \times \GL_2$ in the unbalanced case, which extends the previous results of Furusawa and Morimoto for ${\rm SO}(V) \times \GL_2$.
%These are applications of our main result on the algebraicity of ratios of special values of Rankin--Selberg $L$-functions. 
%\begin{CJK}{UTF8}{bkai}

%本文章的目的為

%\end{CJK}
%\begin{CJK}{UTF8}{bsmi}

%本文章的目的為

%\end{CJK}
\end{abstract}

%% \tableofcontents %% Just for papers exceeding 50 pages.
\maketitle
%\pagewiselinenumbers
\tableofcontents

\section{Introduction}
%In \cite{Deligne1979}, Deligne proposed a remarkable conjecture on the algebraicity of critical values of $L$-functions of motives, in terms of the periods obtained by comparing the Betti and de Rham realizations of the motives.
%In the introduction, we recall two families of motivic $L$-functions associated to elliptic modular forms, and state our results toward Deligne's conjecture for these $L$-functions. 
%We begin with the symmetric power $L$-functions of modular forms.
In \cite{Deligne1979}, Deligne proposed a remarkable conjecture on the algebraicity of critical values of $L$-functions of motives, in terms of the periods obtained by comparing the Betti and de Rham realizations of the motives.
In the introduction, we recall two families of motivic $L$-functions associated with elliptic modular forms and state our results toward Deligne's conjecture for these $L$-functions. Our discussion begins with the symmetric power $L$-functions of modular forms.
Let 
\[
f(\tau) = \sum_{n=1}^\infty a_{f}(n) q^n \in S_{\kappa}(N,\omega),\quad q=e^{2\pi\sqrt{-1}
\,\tau}
\]
be a normalized elliptic modular newform of weight $\kappa \geq 2$, level $N$, and with nebentypus $\omega$.
For each prime $p \nmid N$, denote by $\alpha_{p},\beta_{p}$ the Satake parameters of $f$ at $p$ and put %$A_{i,p} = {\rm diag}(\alpha_{i,p},\beta_{i,p}) \in \GL_2(\C)$. 
\[
A_{p} = \bp \alpha_{p} & 0 \\ 0 & \beta_{p}\ep.
\]
Recall that $\alpha_{p},\beta_{p}$ are the roots of the Hecke polynomial $X^2-a_{f}(p)X+p^{\kappa-1}\omega(p)$.
For $n \geq 1$, the symmetric $n$-th power $L$-function $L(s,{\rm Sym}^n(f))$ is defined by an Euler product
\[
L(s,{\rm Sym}^n(f)) = \prod_{p}L_p(s,{\rm Sym}^n(f)),\quad {\rm Re}(s)> 1+\tfrac{n(\kappa-1)}{2}.
\]
Here the Euler factors are defined by
\[
L_p(s,{\rm Sym}^n(f)) = \det \left( {\bf 1}_{n+1}-{\rm Sym}^n(A_p)\cdot p^{-s}\right)^{-1}
\]
for $p\nmid N$, where ${\rm Sym}^n : \GL_2(\C) \rightarrow \GL_{n+1}(\C)$ is the symmetric $n$-th power representation.
By the result of Barnet-Lamb, Geraghty, Harris, and Taylor \cite[Theorem B]{BLGHT2011}, the symmetric power $L$-functions admit meromorphic continuation to the whole complex plane and satisfy functional equations relating $L(s,{\rm Sym}^n(f))$ to $L(1+n(\kappa-1)-s,{\rm Sym}^n(f^\vee))$, where $f^\vee \in S_{\kappa}(N,\omega^{-1})$ is the normalized newform dual to $f$.
Associated to $f$, we have a pure motive $M_f$ over $\Q$ of rank $2$ with coefficients in the Hecke field $\Q(f)$ of $f$, which was constructed by Deligne \cite{Deligne1971} and Scholl \cite{Scholl1990}, such that
\[
L(M_f,s) = \left( L(s,{}^\sigma\!f) \right)_{\sigma : \Q(f) \rightarrow \C}.
\]
Denote by $\delta(M_f),\,c^\pm(M_f) \in (\Q(f)\otimes_\Q\C)^\times / \Q(f)^\times$ the Deligne's periods of $M_f$ (cf.\,\cite[(1.7.2) and (1.7.3)]{Deligne1979}). 
Under the canonical isomorphism $\Q(f)\otimes_\Q\C \cong \prod_{\sigma : \Q(f)\rightarrow\C}\C$, we write
\[
\delta(M_f) =  \left( \delta({}^\sigma\!f)\right)_{\sigma : \Q(f)\rightarrow\C},\quad c^\pm(M_f) = \left( c^\pm({}^\sigma\!f)\right)_{\sigma : \Q(f)\rightarrow\C}.
\]
In \cite[Proposition 7.7]{Deligne1979}, Deligne explicitly computed the periods $c^\pm({\rm Sym}^n M_f)$ in terms of $\delta(M_f)$ and $c^\pm(M_f)$.
In particular, we have the following conjecture on the algebraicity of critical values of symmetric power $L$-functions.

\begin{conja}[Deligne, Conjecture \ref{C:Deligne Sym}]\label{C:Deligne Sym introduction}
For a critical point $m \in \Z$ for $L(s,{\rm Sym}^n(f))$, we have
\begin{align*}
&\sigma \left( \frac{L(m,{\rm Sym}^n(f))}{(2\pi\sqrt{-1})^{d_n^{(-1)^m} m}\cdot c^{(-1)^m}({\rm Sym}^n(f))}\right) = \frac{L(m,{\rm Sym}^n({}^\sigma\!f))}{(2\pi\sqrt{-1})^{d_n^{(-1)^m} m}\cdot c^{(-1)^m}({\rm Sym}^n({}^\sigma\!f))}, \quad  \sigma\in{\rm Aut}(\C).
\end{align*}
Here $d_n^+ = \tfrac{n}{2}+1,\,d_n^- = \tfrac{n}{2}$ (resp.\,$d_n^\pm = \tfrac{n+1}{2}$) if $n$ is even (resp.\,odd) and
\[
c^\pm({\rm Sym}^n(f)) = \begin{cases}
\delta(f)^{r(r\pm1)/2}\cdot(c^+(f)\cdot c^-(f))^{r(r+1)/2} & \mbox{ if $n=2r$},\\
\delta(f)^{r(r+1)/2}\cdot c^\pm(f)^{(r+1)(r+2)/2}\cdot c^\mp(f)^{r(r+1)/2} & \mbox{ if $n=2r+1$}.
\end{cases}
\]
%\[
%c^\pm({\rm Sym}^{2r+1}(f)) = (2\pi\sqrt{-1})^{r(r+1)(1-\kappa)/2}\cdot G(\omega)^{r(r+1)/2}\cdot c^\pm(f)^{(r+1)(r+2)/2}\cdot c^\mp(f)^{r(r+1)/2}.
%\]
\end{conja}

%For $n=1$, the conjecture is known, as explained by Deligne in \cite[\S\,7]{Deligne1979}.
%The conjecture was proved by Sturm  \cite{Sturm1980}, \cite{Sturm1989} for $n=2$, and by Garrett and Harris \cite{GH1993} for $n=3$ assuming $\kappa \geq 5$. We extend the result of Garrett and Harris to $\kappa=3,4$ in \cite{Chen2021d}. 
%For $n=4,6$, Morimoto proves the algebraicity of the ratio in \cite{Morimoto2021} assuming $\kappa \geq 6$ and $\omega$ is trivial.
%Based on similar idea as in $loc.$ $cit.$, we prove the conjecture for $n=4$ in \cite{Chen2021} assuming $\kappa \geq 3$.
%For $n=5,6$, we prove the conjecture assuming $\kappa \geq 6$ in \cite{Chen2021g} and \cite{Chen2021f}.
%We prove in this paper the following theorem regarding Deligne's conjecture for arbitrary $n$:

For $n=1$, the conjecture is well-known, as explained by Deligne in \cite[\S\,7]{Deligne1979}. The case $n=2$ was proved by Zagier in \cite{Zagier1977} and Sturm in \cite{Sturm1980} and \cite{Sturm1989}, and for $n=3$, Garrett and Harris established the result under the assumption $\kappa \geq 5$ in \cite{GH1993}. In our previous work \cite{Chen2021d}, we extended the results of Garrett and Harris to the cases $\kappa = 3$ or $\kappa = 4$.
For $n = 4$ or $n = 6$, Morimoto proved the algebraicity of the ratio in \cite{Morimoto2021}, assuming $\kappa \geq 6$ and that $\omega$ is trivial. Building on a similar approach as in \textit{loc.\ cit.}, we proved the conjecture for $n = 4$ under the weaker assumption $\kappa \geq 3$ in \cite{Chen2021}. For $n = 5$ and $n = 6$, we established the conjecture assuming $\kappa \geq 6$ in \cite{Chen2021g} and \cite{Chen2021f}, respectively.
In this paper, we present the following theorem, which provides a proof of Deligne's conjecture for arbitrary $n$:

\begin{thma}[Theorem \ref{T:Sym odd}]\label{T:Sym odd introduction}
Conjecture A holds when $\kappa \geq 5$.
\end{thma}

As another example of Deligne's conjecture on critical values of motivic $L$-functions, we consider tensor product $L$-functions of modular forms. 
Let $n$ be a positive integer and 
$f_i \in S_{\kappa_i}(N_i,\omega_i)$
a normalized elliptic modular newform for $1 \leq i \leq n$.
Define the tensor product $L$-function $L(s,f_1 \otimes \cdots \otimes f_n)$ by an Euler product
\[
L(s,f_1 \otimes \cdots \otimes f_n) = \prod_{p} L_p(s,f_1\otimes\cdots\otimes f_n), \quad {\rm Re}(s)>1+\sum_{i=1}^n\tfrac{\kappa_i-1}{2}.
\]
Here the Euler factors are given by
\[
L_p(s,f_1\otimes\cdots\otimes f_n) = \det \left( {\bf 1}_{2^n}-\otimes_{i=1}^n A_{i,p}\cdot p^{-s}\right)^{-1}
\]
for $p \nmid N_1\cdots N_n$, where $A_{i,p}$ is a $2$ by $2$ semisimple matrix whose eigenvalues are Satake parameters of $f_i$ at $p$.
Conjecturally, it admits meromorphic continuation to the whole complex plane and satisfies a functional equation relating $L(s,f_1 \otimes \cdots \otimes f_n)$ to $L(1+\sum_{i=1}^n(\kappa_i-1)-s,f_1^\vee \otimes \cdots \otimes f_n^\vee)$.
In \cite{Blasius1987}, Blasius proposed the following conjecture on the algebraicity of special values of the tensor product $L$-function, which is a refinement of Deligne's conjecture \cite{Deligne1979} for the tensor product motive associated to $f_1\otimes\cdots\otimes f_n$.
For $1 \leq i \leq n$, let $G(\omega_i)$ be the Gauss sum of $\omega_i$ and $\< f_i ,f_i\>$ be the Petersson norm of $f_i$ defined by
\[
\< f_i ,f_i\> = {\rm vol}(\Gamma_0(N_i)\backslash\frak{H})^{-1}\int_{\Gamma_0(N_i)\backslash\frak{H}}|f_i(\tau)|^2y^{\kappa_i-2}\,d\tau.
\]

\begin{conjb}[Blasius, Conjecture \ref{C:Blasius}]\label{C:Blasius introduction}
Assume $n \geq 2$.
Let $m \in \Z$ be a critical point for $L(s,f_1 \otimes \cdots \otimes f_n)$. Then the tensor product $L$-function is holomorphic at $s=m$ and we have
\[
\sigma \left( \frac{L(m,f_1\otimes\cdots\otimes f_n)}{(2\pi\sqrt{-1})^{2^{n-1}m}\cdot c(f_1\otimes \cdots\otimes f_n)} \right) = \frac{L(m,{}^\sigma\!f_1\otimes\cdots\otimes {}^\sigma\!f_n)}{(2\pi\sqrt{-1})^{2^{n-1}m}\cdot c({}^\sigma\!f_1\otimes \cdots\otimes {}^\sigma\!f_n)},\quad \sigma \in {\rm Aut}(\C).
\]
Here 
\[
c(f_1 \otimes \cdots \otimes f_n) = (2\pi\sqrt{-1})^{2^{n-2}\sum_{i=1}^n(1-\kappa_i)}\cdot \prod_{i=1}^n G(\omega_i)^{2^{n-2}}\cdot (\pi^{\kappa_i}\cdot \< f_i ,f_i\>)^{2^{n-2}-t_i}
\]
with $t_i$ equal to the cardinality of the set
\[
\left\{(\varepsilon_1,...,\varepsilon_n)\in \{\pm1\}^n\,\left\vert\, 2(\kappa_i-1)<\sum_{j=1}^n \varepsilon_j(\kappa_j-1),\quad \varepsilon_i=1\right\}\right. ,\quad 1 \leq i \leq n.
\]
%for $1 \leq i \leq n$.
\end{conjb}

For $n=2$, we have the pioneering result of Shimura in \cite{Shimura1976}.
For $n = 3$, the conjecture has been considered by various authors and partially proved. Notable contributions include the works of Harris–Kudla \cite{HK1991}, Garrett–Harris \cite{GH1993}, and Furusawa–Morimoto \cite{FM2014}, \cite{FM2016}. 
In this paper, we prove the following theorem regarding Blasius' conjecture for arbitrary $n$:

%For $n=2$, we have the pioneering result of Shimura in \cite{Shimura1976}.
%For $n=3$, the conjecture was considered by various authors and is partially proved, notably we have the results of Harris--Kudla \cite{HK1991}, Garrett--Harris \cite{GH1993}, and Furusawa--Morimoto \cite{FM2014}, \cite{FM2016}.
%We prove in this paper the following theorem regarding Blasius' conjecture for arbitrary $n$:

\begin{thmb}[Theorem \ref{T:tensor product for GL_2}]\label{T:Blasius introduction}
Conjecture B holds if the following conditions are satisfied:
\begin{itemize}
\item[(1)]
$|\sum_{i=1}^n\varepsilon_i(\kappa_i-1)| \geq 6$ and $|\sum_{i=1}^n(\varepsilon_i-\varepsilon_i')(\kappa_i-1)| \geq 6$ for all $(\varepsilon_1,...,\varepsilon_n)\neq(\varepsilon_1',...,\varepsilon_n')$ in $\{\pm1\}^n$.\\
In particular, $\kappa_i \geq 4$ and $|\kappa_i-\kappa_j|\geq 3$ for all $1 \leq i<j\leq n$.
 %with $(\varepsilon_1,\cdots,\varepsilon_n)\neq \pm(\varepsilon_1',\cdots,\varepsilon_n')$.
\item[(2)] If $n \geq 5$, then $N_1=\cdots=N_{n-2}=1$.
\end{itemize}
\end{thmb}

%Besides Theorems A and B, in Theorem \ref{T:RS for GL_n x GL_2} below we also prove an automorphic analogue of Deligne's conjecture \cite{Deligne1979} on critical values of Rankin--Selberg $L$-functions of essentially self-dual cuspidal automorphic representations of $\GL_n \times \GL_2$, which extends the previous results of Furusawa and Morimoto \cite{FM2014}, \cite{FM2016} for ${\rm SO}(V) \times \GL_2$.
%Theorems A and B are proved by mathematical induction. In the induction step, the key ingredient is our main result Theorem \ref{T:main 1} which we shall now introduce.

Besides Theorems A and B, in Theorem \ref{T:RS for GL_n x GL_2} below, we also establish an automorphic analogue of Deligne's conjecture \cite{Deligne1979} on the critical values of Rankin--Selberg $L$-functions for essentially self-dual cuspidal automorphic representations of $\GL_n \times \GL_2$. This result generalizes the prior work of Furusawa and Morimoto \cite{FM2014}, \cite{FM2016} for ${\rm SO}(V) \times \GL_2$.
Theorems A and B are proved via mathematical induction. In the induction step, the key ingredient is our main result, Theorem \ref{T:main 1}, which we shall now introduce.

\subsection{Cross-ratio of critical Rankin--Selberg $L$-values}

%Denote by $\A$ the ring of adeles of $\Q$.
%For isobaric automorphic representations $\itSigma$ and $\itPi$ of $\GL_n(\A)$ and $\GL_{n'}(\A)$ respectively, let $L(s,\itSigma \times \itPi)$ be the Rankin--Selberg $L$-function of $\itSigma \times \itPi$.
%These $L$-functions play important roles in number theory.
%The analytic properties were studied by various authors (cf.\,\cite{JS1981}, \cite{JS1981b}, \cite{JPSS1983}, \cite{Shahidi1981}, \cite{MW1989}, and \cite{CP2004}). In particular, $L(s,\itSigma \times \itPi)$ admits meromorphic continuation to the whole complex plane and satisfies a functional equation relating $L(s,\itSigma \times \itPi)$ to $L(1-s,\itSigma^\vee \times \itPi^\vee)$.
%Besides its analytic properties, the algebraic aspect is also an important subject of study.
%In this paper, we investigate the algebraicity of ratios of special values of Rankin--Selberg $L$-functions.
%Following Clozel in \cite{Clozel1990}, we consider $algebraic$ automorphic representations of general linear groups.
%For an algebraic automorphic representation $\itPi$, we denote by ${}^\sigma\!\itPi = \itPi_\infty \otimes {}^\sigma\!\itPi_f$ its $\sigma$-conjugate for $\sigma \in {\rm Aut}(\C)$.
%It is conjectured that these conjugates are automorphic as well, which was proved in $loc.$ $cit.$ for regular algebraic automorphic representations.
Let $\A$ denote the ring of adeles of $\Q$. For isobaric automorphic representations $\itSigma$ and $\itPi$ of $\GL_n(\A)$ and $\GL_{n'}(\A)$, respectively, let $L(s, \itSigma \times \itPi)$ denote the Rankin--Selberg $L$-function associated with $\itSigma$ and $\itPi$. These $L$-functions play an important role in number theory.
The fundamental analytic properties of these $L$-functions have been extensively studied by various authors (cf.\,\cite{JS1981}, \cite{JS1981b}, \cite{JPSS1983}, \cite{Shahidi1981}, \cite{MW1989}, and \cite{CP2004}). In particular, it is known that $L(s, \itSigma \times \itPi)$ admits a meromorphic continuation to the entire complex plane and satisfies a functional equation that relates $L(s, \itSigma \times \itPi)$ to $L(1-s, \itSigma^\vee \times \itPi^\vee)$.
Beyond their analytic properties, the algebraic aspects of $L$-functions also constitute a significant area of research. In this paper, we focus on the algebraicity of ratios of special values of Rankin--Selberg $L$-functions. Following Clozel's framework in \cite{Clozel1990}, we consider \textit{algebraic} automorphic representations of general linear groups.
For an algebraic automorphic representation $\itPi$, we denote its $\sigma$-conjugate by ${}^\sigma\!\itPi = \itPi_\infty \otimes {}^\sigma\!\itPi_f$ for $\sigma \in \mathrm{Aut}(\C)$. It is believed that $\itPi_f$ is defined over a number field and these conjugates are themselves automorphic, a result that was proved in \textit{loc.\ cit.} for regular algebraic automorphic representations.
%In \cite{HR2020}, Harder and Raghuram prove a striking result on the algebraicity of ratios of critical values of a fixed Rankin--Selberg $L$-function. More precisely, let $\itSigma$ and $\itPi$ be regular algebraic cuspidal automorphic representations of $\GL_n(\A)$ and $\GL_{n'}(\A)$, respectively. Assume $n$ and $n'$ are both even, they prove that 
%\[
%\sigma\left( \frac{L(m_0,\itSigma\times\itPi)}{L(m_0+1,\itSigma\times\itPi)}\right) = \frac{L(m_0,{}^\sigma\!\itSigma\times{}^\sigma\!\itPi)}{L(m_0+1,{}^\sigma\!\itSigma\times{}^\sigma\!\itPi)},\quad \sigma \in {\rm Aut}(\C)
%\]
%for all consecutive critical points $m_0,m_0+1$ for $L(s,\itSigma\times\itPi)$ with $L(m_0+1,\itSigma\times\itPi) \neq 0$.
%If $n$ is even and $n'$ is odd, then the ratios are expressed in terms of the relative periods of $\itSigma$ (cf.\,\cite[Definition 5.3]{HR2020}). 
%The result is compatible with Deligne's conjecture \cite{Deligne1979} on critical values of motivic $L$-functions, and in practice
%it reduces the study of algebraicity of critical values to that of the rightmost and leftmost critical values.
In \cite{HR2020}, Harder and Raghuram establish a striking result concerning the algebraicity of ratios of critical values for a fixed Rankin--Selberg $L$-function. More precisely, let $\itSigma$ and $\itPi$ be regular algebraic cuspidal automorphic representations of $\GL_n(\A)$ and $\GL_{n'}(\A)$, respectively. Assuming that $n$ and $n'$ are both even, they prove that  
\[
\sigma\left( \frac{L(m_0, \itSigma \times \itPi)}{L(m_0+1, \itSigma \times \itPi)} \right) = \frac{L(m_0, {}^\sigma\!\itSigma \times {}^\sigma\!\itPi)}{L(m_0+1, {}^\sigma\!\itSigma \times {}^\sigma\!\itPi)}, \quad \sigma \in \mathrm{Aut}(\C)
\]
for any consecutive critical points $m_0, m_0+1$ for $L(s, \itSigma \times \itPi)$, provided that $L(m_0+1, \itSigma \times \itPi) \neq 0$. 
When $n$ is even and $n'$ is odd, the ratios are expressed in terms of the relative periods of $\itSigma$ (cf.\,\cite[Definition 5.3]{HR2020}). This result is compatible with Deligne's conjecture \cite{Deligne1979} on the critical values of motivic $L$-functions. In practice, it reduces the study of the algebraicity of critical values to that of the rightmost or leftmost critical values.
%As a different aspect of ratios of critical values, we consider ratios of product of different Rankin--Selberg $L$-functions at a fixed critical point. 
%Under reasonable assumptions, we expect these ratios to be algebraic and propose the following:
As another perspective on the study of ratios of critical values, we consider cross ratios of Rankin--Selberg $L$-functions evaluated at a fixed critical point. Under reasonable assumptions, these cross ratios are expected to be algebraic and we propose the following:

\begin{conj}[Cross-Ratio Formula]\label{C:main 1}
Let $\itSigma, \itSigma'$ (resp.\,$\itPi,\itPi'$) be algebraic automorphic representations of $\GL_n(\A)$ (resp.\,$\GL_{n'}(\A)$) such that %the archimedean components are essentially tempered with
\[
\itSigma_\infty = \itSigma_\infty',\quad \itPi_\infty = \itPi_\infty'.
\]
Let $m_0 \in \Z+\tfrac{n+n'}{2}$ be a critical point for $L(s,\itSigma \times \itPi)$ such that
$L(m_0,\itSigma \times \itPi')\cdot L(m_0,\itSigma' \times \itPi) \neq 0$.
Then we have
\[
\sigma \left( \frac{L(m_0,\itSigma \times \itPi)\cdot L(m_0,\itSigma' \times \itPi')}{L(m_0,\itSigma \times \itPi')\cdot L(m_0,\itSigma' \times \itPi)}\right) = \frac{L(m_0,{}^\sigma\!\itSigma \times {}^\sigma\!\itPi)\cdot L(m_0,{}^\sigma\!\itSigma' \times {}^\sigma\!\itPi')}{L(m_0,{}^\sigma\!\itSigma \times {}^\sigma\!\itPi')\cdot L(m_0,{}^\sigma\!\itSigma' \times {}^\sigma\!\itPi)},\quad \sigma \in {\rm Aut}(\C).
\]
%In particular, we have
%\[
% \frac{L(m_0,\itSigma \times \itPi)\cdot L(m_0,\itSigma' \times \itPi')}{L(m_0,\itSigma \times \itPi')\cdot L(m_0,\itSigma' \times \itPi)} \in  \Q(\itSigma)\cdot \Q(\itSigma')\cdot\Q(\itPi)\cdot\Q(\itPi').
%\]
\end{conj}

%The conjecture is also compatible with Deligne's conjecture, as will be explained in \S\,\ref{SS:DC} below. 
%When the automorphic representations are regular algebraic and cuspidal, we listed below some special cases for which the conjecture can be directly verified using the algebraicity result in each case:
The conjecture is also compatible with Deligne's conjecture, as will be explained in \S\ref{SS:DC} below. When the automorphic representations are regular algebraic and cuspidal, we listed below some special cases where the conjecture can be directly verified using known algebraicity results for each case:
\begin{itemize}
\item $n=2$ and $n' \leq 2$, see \cite{Shimura1976} and \cite{Shimura1977}.
\item $n$ is even, $n'=1$, and $\itSigma$, $\itSigma'$ are essentially self-dual and symplectic, see \cite{GR2014}.
\item $n$ is odd, $n'=1$, and $\itSigma$, $\itSigma'$ are self-dual and orthogonal, see \cite{Shimura2000}, \cite{HPSS2021}, and \cite{Liu2019b}.
\item $n'=n-1$ and $(\itSigma_\infty,\itPi_\infty)$ is balanced/interlaced, see \cite{Raghuram2009} and \cite{Raghuram2016}.
\end{itemize}

In practice, compared with known methods in the literature, Conjecture \ref{C:main 1} offers an alternative approach to establishing algebraicity results for automorphic or motivic $L$-functions. For instance, in \cite{Chen2021g}, we applied Conjecture \ref{C:main 1} to $\GL_4 \times \GL_3$ to prove Deligne's conjecture for the symmetric fifth $L$-functions of modular forms. Further applications to Deligne's conjecture for motivic $L$-functions of higher degrees are discussed in \S\,\ref{S:applications} below.
The following theorem constitutes the main result of this paper. We prove Conjecture \ref{C:main 1} when $nn'$ is even, under certain regularity assumptions on $\itSigma_\infty$ and $\itPi_\infty$. Please refer to (\ref{E:regularity 1}) and (\ref{E:regularity 2}) for regularity conditions which are more refined than condition (3) below, and Remark \ref{R:regularity} for the interpretation of these conditions.

\begin{thm}[Theorem \ref{T:Main}]\label{T:main 1}
{Conjecture \ref{C:main 1} holds if the following conditions are satisfied:
\begin{itemize}
\item[(1)] $nn'$ is even.
\item[(2)] $\itSigma_\infty$ and $\itPi_\infty$ are essentially tempered and cohomological.
\item[(3)]
%Let $\underline{\kappa} = (\kappa_1,\cdots,\kappa_{\lfloor \frac{n}{2}\rfloor}) \in \Z^{\lfloor \frac{n}{2}\rfloor}$ and $\underline{\ell} = (\ell_1,\cdots,\ell_{\lfloor \frac{n'}{2}\rfloor}) \in \Z^{\lfloor \frac{n'}{2}\rfloor}$ be the highest weights of the minimal ${\rm O}_n(\R)$ and ${\rm O}_{n'}(\R)$ types of $\itSigma_\infty$ and $\itPi_\infty$ respectively.
We have $\min\{\kappa_j,\ell_j\} \geq 5$ and $\underline{\kappa}\sqcup\underline{\ell}$ is $6$-regular (see \S\,\ref{SS:notation} for this notion), where $\underline{\kappa}\in \mathbb{N}^{\lfloor \frac{n}{2}\rfloor}$ and $\underline{\ell}\in \mathbb{N}^{\lfloor \frac{n'}{2}\rfloor}$ are the highest weights of the minimal ${\rm O}_n(\R)$ and ${\rm O}_{n'}(\R)$ types of $\itSigma_\infty$ and $\itPi_\infty$ respectively. 
\end{itemize}}
\end{thm}

\subsection{Outline of the proof}

We refer to the beginning of \S\,\ref{S:applications} for the strategy on proving Theorems A and B.
For Theorem \ref{T:main 1}, the central idea behind the proof is to separate the dependence of the algebraicity of the critical values of $L(s,\itSigma \times \itPi)$ on $\itSigma$ and $\itPi$.  Specifically, for a critical point $m_0$ for $L(s,\itSigma \times \itPi)$, we show that there exist $q(m_0,\itSigma),\,q(m_0,\itPi) \in \C^\times$ depending on $\itSigma$ and $\itPi$ respectively and also on $m_0$, and an archimedean constant $p_\infty(m_0) \in \C^\times$ depending on $\itSigma_\infty, \,\itPi_\infty,\,m_0$ such that
\begin{align}\label{E:sketch 1}
L(m_0,\itSigma \times \itPi) \sim p_\infty(m_0)\cdot q(m_0,\itSigma)\cdot q(m_0,\itPi).
\end{align}
Subject to the existence of pure motives associated to $\itSigma$ and $\itPi$, (\ref{E:sketch 1}) is actually a consequence of Deligne's conjecture \cite{Deligne1979} and Yoshida's computation \cite{Yoshida2001} on Deligne's periods of tensor product motives (cf.\,\S\,\ref{SS:DC}).
Under the assumption $\itSigma_\infty = \itSigma_\infty'$ and $\itPi_\infty = \itPi_\infty'$, (\ref{E:sketch 1}) also holds for $L(s,\itSigma \times \itPi')$, $L(s,\itSigma' \times \itPi)$, and $L(s,\itSigma' \times \itPi')$ with $q(m_0,\itSigma)$ or $q(m_0,\itPi)$ replaced by $q(m_0,\itSigma')$ or $q(m_0,\itPi')$ in each case, which is also compatible with \cite{Deligne1979} and \cite{Yoshida2001} mentioned above. Conjecture \ref{C:main 1} for $m_0$ then follows immediately from (\ref{E:sketch 1}) for these Rankin--Selberg $L$-functions.

By assumption (1), we may assume $n$ is even.
Note that $\min\{\kappa_j,\ell_j\} \geq 3$ and $\underline{\kappa}\sqcup\underline{\ell}$ is $2$-regular imply that $L(s,\itSigma \times \itPi)$ admits at least two critical points (cf.\,\S\,\ref{SS:RS}).
To prove Conjecture \ref{C:main 1}, it then suffices to show that (\ref{E:sketch 1}) holds for an arbitrary non-central critical point by applying the results of Harder--Raghuram \cite[Theorem 7.21]{HR2020} and the author \cite[Theorem 3.2]{Chen2023}.
We assume $m_0$ is slightly on the right of the central strip as in (\ref{E:main proof 0}).
Consider an algebraic automorphic representation $\tau$ of $\GL_{n+n'}(\A)$ defined by
\[
\tau = (\itSigma \otimes |\mbox{ }|_\A^{(n+n')/2}) \boxplus (\itPi^\vee \otimes |\mbox{ }|_\A^{1-m_0+(n+n')/2}).
\]
Based on assumptions (1) and (2), applying Lemma \ref{L:combinatorial} we see that $\tau_\infty$ contributes non-trivially to the bottom degree relative Lie algebra cohomology. With respect to a canonical choice of generator in the relative Lie algebra cohomology in \S\,\ref{SS:generator}, we construct in \S\,\ref{S:Eisenstein} certain Eisenstein cohomology classes associated to $\tau$ in the cohomology of the locally symmetric space for $\GL_{n+n'}$ in the bottom degree. By restricting to the boundary of the Borel--Serre compactification of the locally symmetric space, we prove in Theorem \ref{T:period relation} that the algebraicity of these cohomology classes can be expressed in terms of products of the Betti--Whittaker periods (cf.\,\S\,\ref{SS:BW periods}) of the cuspidal summands of $\tau$.
The assumption (3) ensures the existence of an auxiliary regular algebraic cuspidal automorphic representation $\itPsi$ of $\GL_{n+n'+1}(\A)$ together with a half-integer $M_0$ such that (i) $M_0$ is critical for $L(s,\itPsi \times \tau)$; (ii) $(\itPsi_\infty,\tau_\infty)$ is balanced (cf.\,(\ref{E:balanced})); (iii) $M_0$ shifted to non-central critical points for the two Rankin--Selberg $L$-functions that appear on the right-hand side of the decomposition 
\[
L(s,\itPsi \times \tau) = L(s+\tfrac{n+n'}{2},\itPsi\times\itSigma)\cdot L(s+1-m_0+\tfrac{n+n'}{2},\itPsi\times\itPi^\vee).
\]
In particular, $L(M_0,\itPsi \times \tau)$ is non-zero.
%More precisely, assumption (3) implies the existence of $M_0$ such that $M_0+\tfrac{n+n'}{2}$ and $M_0+1-m_0+\tfrac{n+n'}{2}$ are non-central critical points for $L(s,\itPsi \times \itSigma)$ and $L(s,\itPsi \times \itPi^\vee)$ respectively.
For the Rankin--Selberg $L$-function $L(s,\itPsi \times \tau)$, by following the methods of Mahnkopf \cite{Mahnkopf2005} and Raghuram \cite{Raghuram2009} together with the algebraicity of Eisenstein cohomology classes associated to $\tau$, we prove in Theorem \ref{T:Raghuram} that the algebraicity of the critical value $L(M_0,\itPsi\times\tau)$ can be expressed in terms of a product of the Betti--Whittaker periods (cf.\,\S\,\ref{SS:BW periods 2}) of $\itPsi$ and $\tau$ together with an archimedean constant  $Z_\infty(M_0)$ depending on $\itPsi_\infty$, $\tau_\infty$, and $M_0$. More precisely, the integral representation of $L(s,\itPsi \times \tau)$ in \cite{JS1981} and \cite{JS1981b} can be cohomologically interpreted as a cup product between bottom degree cohomology of the locally symmetric spaces for $\GL_{n+n'+1}$ and $\GL_{n+n'}$. The cup product then induces a pairing between relative Lie algebra cohomology groups and defines a non-zero archimedean constant $Z_\infty(M_0)$ by evaluation of the pairing at the canonical generators. 
%We prove in appendix \ref{S:appendix} that $Z_\infty(M_0)$ is non-zero, which extends the result of Sun \cite{Sun2017} for the case $\tau_\infty$ is essentially tempered.
We then take 
\[
p_\infty(m_0) = Z_\infty(M_0).
\]
Note that $\tau_\infty$ is not essentially tempered when $L(s,\itSigma \times \itPi)$ admits central critical point.
Finally, we prove in Lemma \ref{L:period relation} the following period relation between Betti--Whittaker periods, which is an easy consequence of Definition \ref{D:Betti-Whittaker}:
\[
p(\tau,\varepsilon) \sim G(\omega_\itPi)^{-n}\cdot L^{(\infty)}(m_0, \itSigma \times \itPi)\cdot p(\itSigma \otimes |\mbox{ }|_\A^{(n+2n')/2},\varepsilon)\cdot p(\itPi^\vee\otimes|\mbox{ }|_\A^{1-m_0+n'/2},\varepsilon),
\]
where $L^{(\infty)}(s, \itSigma \times \itPi)$ is the finite part of the Rankin--Selberg $L$-function of $\itSigma \times \itPi$.
Therefore, (\ref{E:sketch 1}) for $m_0$ then follows once we take $q(m_0,\itSigma)$ and $q(m_0,\itPi)$ to be ratios of the critical values $L(M_0+\tfrac{n+n'}{2},\itPsi\times\itSigma)$ and $L(M_0+1-m_0+\tfrac{n+n'}{2},\itPsi\times\itPi^\vee)$ by the Betti--Whittaker periods of $\itSigma \otimes |\mbox{ }|_\A^{(n+2n')/2}$ and $\itPi^\vee\otimes|\mbox{ }|_\A^{1-m_0+n'/2}$ respectively (with respect to suitable signature $\varepsilon$), and absorb the Betti--Whittaker period of $\itPsi$ into either of the ratios. 

This paper is organized as follows. 
In \S\,\ref{S:BW periods}, we define the Betti–Whittaker periods of a regular algebraic automorphic representation whose archimedean component is cohomological of a particular type, which may be non-tempered. 
Associated with such an automorphic representation, in \S\,\ref{S:Eisenstein}, we construct Eisenstein cohomology classes in the bottom degree and prove their algebraicity in terms of Betti–Whittaker periods in Theorem \ref{T:period relation}.
In \S\,\ref{S:proof main}, we prove our main result, Theorem \ref{T:Main}. A key ingredient in this proof is an algebraicity result for \( \GL_n \times \GL_{n-1} \), which is established in \S\,\ref{SS:Raghuram}. The proof that the archimedean constant appearing in this algebraicity result is non-zero is deferred to Appendix \ref{S:appendix}.
In \S\,\ref{S:applications}, we present applications of Theorem \ref{T:Main} to Deligne's conjecture on motivic \( L \)-functions. This section is logically independent of \S\,\ref{S:BW periods}–\S\,\ref{S:proof main}, except for some notation introduced in \S\,\ref{SS:cohomological} and preliminaries on Deligne's conjecture in \S\,\ref{SS:DC}. 
Readers primarily interested in the applications may start with \S\,\ref{S:applications}.

%This paper is organized as follows. 
%In \S\,\ref{S:BW periods}, we define the Betti--Whittaker periods of a regular algebraic automorphic representation whose archimedean component is cohomological of a particular type which might be non-tempered. 
%Associated to such an automorphic representation, in \S\,\ref{S:Eisenstein} we construct Eisenstein cohomology classes in the bottom degree and prove in Theorem \ref{T:period relation} its algebraicity in terms of Betti--Whittaker periods.
%In \S\,\ref{S:proof main}, we prove our main result Theorem \ref{T:Main}. One of the key ingredients is an algebraicity result for $\GL_n \times \GL_{n-1}$ proved in \S\,\ref{SS:Raghuram}. The archimedean constant appearing in the algebraicity result is non-zero and we postpone its proof to appendix \ref{S:appendix}.
%In \S\,\ref{S:applications}, we give applications of Theorem \ref{T:Main} to Deligne's conjecture on motivic $L$-functions. This section is logically independent of \S\,\ref{S:BW periods}-\S\,\ref{S:proof main} except for some notation in \S\,\ref{SS:cohomological} and preliminaries on Deligne's conjecture in \S\,\ref{SS:DC}.
%Readers who are interested in the applications can begin with \S\,\ref{S:applications}.

\subsection{Notations}\label{SS:notation}

Let $T_n$ and $U_n$ be the maximal torus and maximal unipotent subgroup of $\GL_n$ consisting of diagonal matrices and upper triangular unipotent matrices respectively. Let $X^+(T_n)$ be the set of dominant integral weights of $T_n$ with respect to $U_n$.
We identify $X^+(T_n)$ with the set of $n$-tuples ${\mu} = (\mu_1,...,\mu_n) \in \Z^n$ such that $\mu_1\geq \cdots\geq \mu_n$. We say $\mu \in X^+(T_n)$ is $N$-regular if $|\mu_i-\mu_j| \geq N$ for all $i \neq j$.
We denote by $W_n$ the Weyl group of $\GL_n$ with respect to $T_n$.
We will take permutation matrices in $\GL_n$ as representatives of elements in $W_n$.
For an ordered partition $(n_1,...,n_k)$ of $n$, let $P_{(n_1,...,n_k)}$ be the standard parabolic subgroup of $\GL_n$ defined by
\[
P_{(n_1,...,n_k)} = \left.\left \{ \bp  g_{11} & g_{12} & \cdots & g_{1k}\\
0 & g_{22}& \cdots & g_{2k} \\
0&0&\ddots&\vdots\\
0&0&0&g_{kk}
\ep \,\right\vert\, g_{ii} \in \GL_{n_i},\,g_{ij} \in {\rm M}_{n_i,n_j}\mbox{ for $1\leq i < j \leq k$} \right\}.
\]
For a standard parabolic subgroup $P$ of $\GL_n$, we denote by $U_P$ and $M_P$ the unipotent radical and standard Levi subgroup consisting of block-diagonal matrices of $P$, respectively.
Let $W_n^P$ be the minimal coset representatives of $W_{M_P}\backslash W_n$, where $W_{M_P}$ is the Weyl group of $M_P$.
Note that $w \in W_n^P$ if and only if $w(X^+(T_n))\subseteq \prod_{i=1}^kX^+(T_{n_i})$, where $P=P_{(n_1,...,n_k)}$.

Let $\A$ be the ring of adeles of $\Q$. Let $\A_f$ be the finite part of $\A$, and $\widehat{\Z} = \prod_p \Z_p$ be its maximal compact subring. 
For each place $v$ of $\Q$, let $|\mbox{ }|_v$ be the absolute value on $\Q_v$ normalized so that $|p|_p=p^{-1}$ if $v=p$ is a prime number and $|\mbox{ }|_\infty = |\mbox{ }|$ is the ordinary absolute value on $\R$ if $v = \infty$. 
Let $|\mbox{ }|_\A = \prod_v |\mbox{ }|_v$ be the normalized absolute value on $\A$.
Let $\psi=\otimes_v\psi_v$ be the standard additive character of $\Q\backslash \A$ defined so that
\begin{align*}
\psi_p(x) & = e^{-2\pi \sqrt{-1}x} \mbox{ for }x \in \Z[p^{-1}],\\
\psi_\infty(x) & = e^{2\pi \sqrt{-1}x} \mbox{ for }x \in \R.
\end{align*}
Let $\psi_{U_n} : U_n(\Q)\backslash U_n(\A) \rightarrow \C^\times$ be the additive character defined by
\[
\psi_{U_n}(u) = \psi(u_{12}+u_{23}+\cdots+u_{n-1,n}),\quad u=(u_{ij}) \in U_n(\A).
\]
For each place $v$ of $\Q$, let $\psi_{U_n,v}$ be the local component of $\psi_{U_n}$ at $v$. 
Let $\psi_{U_n,f} = \otimes_{p} \psi_{U_n,p}$.

We denote by $\mathcal{A}(\GL_n(\A))$ the space of automorphic forms on $\GL_n(\Q)\backslash\GL_n(\A)$, that is, smooth, ${\rm O}_n(\R)$-finite, and $\frak{z}_n$-finite functions of moderate growth, where $\frak{z}_n$ is the center of the universal enveloping algebra of the complexified Lie algebra of $\GL_n(\R)$. Let $\mathcal{A}_0(\GL_n(\A)) \subset \mathcal{A}(\GL_n(\A))$ be the space of cusp forms on $\GL_n(\A)$. For a cuspidal automorphic representation of $\GL_n(\A)$, we take the automorphic realization in $\mathcal{A}_0(\GL_n(\A))$ as its underlying space.

Let $\chi$ be an algebraic Hecke character of $\A^\times$. We denote by $G(\chi)$ the Gauss sum of $\chi$ defined by
\[
G(\chi) = \prod_p\varepsilon(0,\chi_p,\psi_p),
\]
where $\varepsilon(s,\chi_p,\psi_p)$ is the $\varepsilon$-factor of $\chi_p$ with respect to $\psi_p$ defined in \cite{Tate1979}.
For $\sigma \in {\rm Aut}(\C)$, let ${}^\sigma\!\chi$ the unique algebraic Hecke character of $\A^\times$ such that ${}^\sigma\!\chi(a) = \sigma (\chi(a))$ for $a \in \A_f^\times$.
It is easy to verify that 
\begin{align}\label{E:Galois Gauss sum}
\begin{split}
\sigma(G(\chi)) = {}^\sigma\!\chi(u_\sigma)G({}^\sigma\!\chi),
\end{split}
\end{align} 
where $u_\sigma \in \widehat{\Z}^\times$ is the unique element such that $\sigma(\psi(x)) = \psi(u_\sigma x)$ for $x \in \A_f$.

The following notation will appear only in the proof of algebraicity results:
Let $(a_\sigma)_{\sigma \in {\rm Aut}(\C)}$ and $(b_\sigma)_{\sigma \in {\rm Aut}(\C)}$ be families of complex numbers indexed by ${\rm Aut}(\C)$ such that $b_\sigma \neq 0$ for all $\sigma$. We write $a_{\rm id} \sim b_{\rm id}$ if 
\[
\sigma \left ( \frac{a_{\rm id}}{b_{\rm id}} \right) = \frac{a_\sigma}{b_\sigma},\quad \sigma \in {\rm Aut}(\C).
\]

\subsubsection*{Acknowledgements}
The author would like to thank Don Blasius, Yao Cheng, Laurent Clozel, Harald Grobner, Ming-Lun Hsieh, Atsushi Ichino, Yubo Jin, Hidenori Katsurada, Wen-Wei Li, Jie Lin, Dongwen Liu, Kazuki Morimoto, A. Raghuram, Binyong Sun, Jacques Tilouine, and Fu-Tsun Wei for their helpful comments and encouragement.
The author is also grateful to Michael Harris for presenting this work in his Bourbaki Seminar lecture \cite{Harris2025}.
The author would also like to thank the anonymous referees for their careful reading and many valuable comments and suggestions, which have greatly improved the exposition of the paper.
Some ideas of this work grew during the scientific activities supported by the National Center for Theoretical Sciences (NCTS).
The author gratefully acknowledges the continuous support of NCTS throughout his academic career, from graduate study onward.
This work was supported by the Japan Society for the Promotion of Science (JSPS) through JSPS KAKENHI Grant Number 22F22316, and by the National Science and Technology Council (NSTC) of R.O.C. under Grant Number 110-2628-M-001-004 and 113-2115-M-007-002-MY3.

\section{Betti--Whittaker periods for $\GL_n$}\label{S:BW periods}

The purpose of this section is to define the Betti--Whittaker periods of a regular algebraic automorphic representation whose archimedean component is cohomological of a particular type introduced in \S\,\ref{SS:cohomological} below. 
Any regular algebraic automorphic representation appears in this paper will have archimedean component of this particular type.
The definition of Betti--Whittaker periods depends on choice of generators in relative Lie algebraic cohomology groups. We specify a canonical choice in \S\,\ref{SS:generator}.

\subsection{Cohomological representations}\label{SS:cohomological}

Let $K_n^\circ$ and $K_n$ be closed subgroups of $\GL_n(\R)$ defined by
\[
K_n^\circ = \R_+\cdot{\rm SO}_n(\R),\quad K_n = \R_+\cdot{\rm O}_n(\R).
\]
Here we regard the set of positive real numbers $\R_+$ as the topological connected component of the center of $\GL_n(\R)$.
We denote by $\frak{g}_n$ and $\frak{k}_n$ the complexified Lie algebras of $\GL_n(\R)$ and $K_n$, respectively.
For an irreducible admissible $(\frak{g}_n,{\rm O}_n(\R))$-module $\pi$, we denote by $\pi^\infty$ the Casselman--Wallach smooth globalization of $\pi$.
Recall $\pi$ is cohomological if there exists some irreducible algebraic representation $M_\C$ of $\GL_n(\C)$ such that the $(\frak{g}_n,K_n^\circ)$-cohomology of $\pi\otimes M_\C$ is non-zero:
\[
H^\bullet(\frak{g}_n,K_n^\circ;\pi\otimes M_\C)\neq 0.
\]
In this case, we say $\pi$ is cohomological with coefficients in $M_\C$. 
%Moreover, if we assume further that $\itPi$ is essentially tempered, then $M_\C$ is uniquely determined.

Let $\mu \in X^+(T_n)$. 
We denote by $M_\mu$ the irreducible algebraic representation of $\GL_n(\Q)$ with highest weight $\mu$ and 
$
M_{\mu,\C} = M_\mu \otimes_\Q \C
$
be the base change to a representation of $\GL_n(\C)$.
%In the following discussion, we associate to each $\mu \in X^+(T_n)$ a set $\Omega_\mu$ consisting of cohomological irreducible admissible $(\frak{g}_n,{\rm O}_n(\R))$-modules with coefficients in $M_{\mu,\C}$.
%Let $\mu \in X^+(T_n)$. 
We define two tuples of integers $\underline{\kappa} \in \Z^{\lfloor \frac{n}{2} \rfloor}$ and $\underline{\sf w} \in \Z^{\lfloor \frac{n+1}{2} \rfloor}$ by the following formula:
\begin{align}\label{E:parameter}
\mu = -\rho_n+
\begin{cases}
(\tfrac{\kappa_1-1-{\sf w}_1}{2},...,\tfrac{\kappa_r-1-{\sf w}_r}{2},\tfrac{1-\kappa_r-{\sf w}_r}{2},...,\tfrac{1-\kappa_1-{\sf w}_1}{2}) & \mbox{ if $n=2r$},\\
(\tfrac{\kappa_1-1-{\sf w}_1}{2},...,\tfrac{\kappa_r-1-{\sf w}_r}{2},-\tfrac{{\sf w}_{r+1}}{2},\tfrac{1-\kappa_r-{\sf w}_r}{2},...,\tfrac{1-\kappa_1-{\sf w}_1}{2}) & \mbox{ if $n=2r+1$}.
\end{cases}
\end{align}
Here $\rho_n = (\tfrac{n-1}{2},\tfrac{n-3}{2},...,\tfrac{1-n}{2})$ is half the sum of positive roots.
Note that $\kappa_{1} > \cdots > \kappa_{r} \geq 2$, ${\sf w}_{r+1}$ is even if $n=2r+1$, and
\begin{align}\label{E:parity}
\kappa_i \equiv {\sf w}_i+n\,({\rm mod}\,2),\quad 1 \leq i \leq {\lfloor \tfrac{n}{2} \rfloor}.
\end{align}
%Conversely, any two tuples of integers $\underline{\kappa}$ and $\underline{\sf w}$ that satisfy these conditions would uniquely determine a dominant integral weight $\mu$ by formula (\ref{E:parameter}).
Let $\pi_\mu$ be the $(\frak{g}_n,{\rm O}_n(\R))$-module realized as the space of ${\rm O}_n(\R)$-finite vectors of the following induced representation of $\GL_n(\R)$: 
\begin{align}\label{E:cohomological}
\begin{cases}
{\rm Ind}_{P_{(2,...,2)}(\R)}^{\GL_n(\R)}\left((D_{\kappa_{1}}{\otimes} |\mbox{ }|^{{\sf w}_1/2})\,\widehat{\otimes}\cdots\widehat{\otimes}\, (D_{\kappa_{r}}\otimes |\mbox{ }|^{{\sf w}_r/2})\right) & \mbox{ if $n=2r$},\\
{\rm Ind}_{P_{(2,...,2,1)}(\R)}^{\GL_n(\R)}\left((D_{\kappa_{1}}\otimes |\mbox{ }|^{{\sf w}_1/2})\,\widehat{\otimes}\cdots\widehat{\otimes} \,(D_{\kappa_{r}}\otimes |\mbox{ }|^{{\sf w}_r/2}) \,\widehat{\otimes}\, |\mbox{ }|^{{\sf w}_{r+1}/2}\right) & \mbox{ if $n=2r+1$}.
\end{cases}
\end{align}
%for some $\kappa_{1} > \cdots > \kappa_{r} \geq 2$ and $\delta \in \{0,1\}$ if $n$ is odd, such that
%\begin{align}\label{E:parity}
%\begin{cases}
%\kappa_{1} \equiv \cdots \equiv \kappa_{r} \equiv {\sf w}\,({\rm mod}\,2) & \mbox{ if $n=2r$},\\
%\kappa_1 \equiv \cdots \equiv \kappa_{r}\equiv {\sf w}+1\equiv 1\,({\rm mod}\,2) & \mbox{ if $n=2r+1$}.
%\end{cases}
%\end{align}
Here $\widehat{\otimes}$ refers to completed projective tensor product, and $D_\kappa$ is the discrete series representation of $\GL_2(\R)$ with minimal weight $\kappa \geq 2$ normalized so that its central character is trivial on $\R_+$.
%We write $\underline{\kappa} = (\kappa_1,\cdots,\kappa_r)$ and call $(\underline{\kappa};\,{\sf w})$ the infinity type of $\itPi$. 
%By the result of Shahidi on local coefficients in \cite[Theorem 3.1]{Shahidi1985} and 
By the irreducibility criterion of Speh in \cite[Theorem 10b]{Moeglin1996}, the condition that $\mu$ is dominant implies $\pi_\mu$ is irreducible. 
Also it is clear that the infinitesimal character of $\pi_\mu$ is equal to $\mu^\vee + \rho_n$, where $\mu^\vee = (-\mu_n,...,-\mu_1)$.
Define a set $\Omega_\mu$ of irreducible dmissible $(\frak{g}_n,{\rm O}_n(\R))$-modules by
\[
\Omega_\mu = \begin{cases}
\left\{\pi_\mu\right\} & \mbox{ if $n$ is even},\\
\left\{\pi_\mu,\pi_\mu\otimes{\rm sgn}\right\} & \mbox{ if $n$ is odd}.
\end{cases}
\]
Let $\pi \in \Omega_\mu$. We call 
\[
(\underline{\kappa};\,\underline{\sf w})\in \Z^{{\lfloor \frac{n}{2} \rfloor}} \times \Z^{{\lfloor \frac{n+1}{2} \rfloor}}
\] 
the infinity type of $\pi$. Note that $\underline{\kappa}$ is the highest weight of a minimal ${\rm SO}_n(\R)$-type of $\pi$ with respect the Borel subalgebra $\frak{b}_n^c$ of ${\rm Lie}({\rm SO}_n(\R))_\C$ in appendix \ref{S:appendix}.
Also note that $\pi$ is essentially tempered if and only if $\underline{{\sf w}} = ({\sf w},...,{\sf w})$ for some ${\sf w}\in\Z$. In this case we also write $(\underline{\kappa};\,{\sf w})$ for the infinity type of $\pi$.
When $n=2r+1$ is odd, we define the signature $\varepsilon(\pi)$ of $\pi$ by
\begin{align}\label{E:signature}
\varepsilon(\pi) = \begin{cases}
(-1)^{r+{\sf w}_{r+1}/2} & \mbox{ if $\pi = \pi_\mu$},\\
(-1)^{r+1+{\sf w}_{r+1}/2} & \mbox{ if $\pi = \pi_\mu\otimes{\rm sgn}$}.
\end{cases}
\end{align}
By arguments as in \cite[Lemme 3.14]{Clozel1990}, $\pi$ is cohomological with coefficients in $M_{\mu,\C}$.
Moreover, 
\[
b_n := \lfloor \tfrac{n^2}{4} \rfloor
\]
is the least non-vanishing degree (bottom degree) of the $(\frak{g}_n,K_n^\circ)$-cohomology of $\pi\otimes M_{\mu,\C}$
%We have
%\[
%H^{b_n}(\frak{g}_n,K_n^\circ;\itPi\otimes M_{\mu,\C}) = \left(  \itPi \otimes \bigwedge^{b_n}(\frak{g}_{n,\C}/\frak{k}_{n,\C})^*\otimes M_{\mu,\C}\right)^{K_n^\circ}
%\]
and
\[
{\rm dim}_\C \,H^{b_n}(\frak{g}_n,K_n^\circ;\pi\otimes M_{\mu,\C}) = \begin{cases}
2 & \mbox{ if $n$ is even},\\
1 & \mbox{ if $n$ is odd}.
\end{cases}
\]
Note that $\pi_0(\GL_n(\R)) = \GL_n(\R)/\GL_n(\R)^\circ\cong {\rm O}_n(\R)/{\rm SO}_n(\R)$ naturally acts on the $(\frak{g}_n,K_n^\circ)$-cohomology groups. For $\varepsilon \in \{\pm1\}$, we denote by 
\[
H^{b_n}(\frak{g}_n,K_n^\circ;\pi\otimes M_{\mu,\C})[\varepsilon]
\]
the $\varepsilon$-isotypic space under the action of $\pi_0(\GL_n(\R))$.
The $\varepsilon$-isotypic space is one-dimensional when $n$ is even, and is non-zero for $\varepsilon = \varepsilon(\pi)$ when $n$ is odd.
In \S\,\ref{SS:generator} below, we will specify a canonical choice of generator of the $\varepsilon$-isotypic space under the Whittaker realization of $\pi$.

For $n \geq 1$, we write 
\[
\Omega(n) = \bigcup_{\mu \in X^+(T_n)}\Omega_\mu. 
\]
Let $P=P_{(n_1,...,n_k)}$ be a standard parabolic subgroup of $\GL_n$ with type $(n_1,...,n_k)$. 
For $1 \leq i \leq k$, let $\pi_i$ be an irreducible admissible $(\frak{g}_{n_i},{\rm O}_{n_i}(\R))$-module.
Let $\delta_{P}$ be the square-root of the modulus character of $P(\R)$ and $\delta_i = \delta_P\vert_{\GL_{n_i}(\R)}$ be the restriction of $\delta_{P}$ to $\GL_{n_i}(\R)$ as the $i$-th component of $M_P(\R)$ for $1 \leq i \leq k$.
Note that
\[
\delta_i = |\det|^{(n-n_i)/2 - (n_1+\cdots+n_{i-1})}.
\]
We write
\[
\pi = {\rm Ind}_{{P}(\R)}^{\GL_n(\R)}(\pi_1 \otimes\cdots \otimes \pi_k)
\]
for the $(\frak{g}_n,{\rm O}_n(\R))$-module realized as the space of ${\rm O}_n(\R)$-finite vectors of the induced representation 
\[
{\rm Ind}_{P(\R)}^{\GL_n(\R)}(\pi_1^\infty\, \widehat{\otimes}\cdots \widehat{\otimes} \,\pi_k^\infty ).
\]
In the following lemma, we give sufficient conditions under which $\pi$ belongs to $\Omega(n)$.
In the special case when $P$ is maximal and $\pi_1, \pi_2$ are essentially tempered, the assertion follows from the combinatorial lemma \cite[Lemma 7.14]{HR2020} which was conjectured by Harder--Raghuram and proved by Weselmann.
We prove the following lemma by generalizing the arguments in \textit{loc.\ cit.}.
%{\color{red}The idea for the proof is similar to the one for $loc.$ $cit.$ and we give the details for the reader's convenience.}

\begin{lemma}\label{L:combinatorial}
Let 
$\pi = {\rm Ind}_{P_{(n_1,...,n_k)}(\R)}^{\GL_n(\R)}(\pi_1 \otimes\cdots \otimes \pi_k)$
for some irreducible admissible $(\frak{g}_{n_i},{\rm O}_{n_i}(\R))$-module $\pi_i$ for $1 \leq i \leq k$.
Assume the following conditions are satisfied:
\begin{itemize}
\item[(1)] At most one of the $n_i$'s is odd.
\item[(2)] $\pi_i \otimes \delta_i \in \Omega_{\mu^{(i)}}$ for $1 \leq i \leq k$.
\item[(3)] The $L$-factors $L(s,\pi_i \times \pi_j^\vee)$ and $L(s,\pi_i^\vee \times \pi_j)$ are holomorphic at $s=0$ for $1 \leq i < j \leq k$.
\end{itemize}
Then there exists a unique Weyl element $w \in W_n^P$ such that $\ell(w) = b_n-\sum_{i=1}^kb_{n_i}$ and 
\[
w^{-1}\left( (\mu^{(1)},...,\mu^{(k)})+\rho_n \right)-\rho_n \in X^+(T_n).
\]
Moreover,
$\pi \in \Omega_\mu$ with 
\[
\mu = w^{-1}\left( (\mu^{(1)},...,\mu^{(k)})+\rho_n \right)-\rho_n.
\]
\end{lemma}

\begin{proof}
For $1 \leq i \leq k$, let $r_i = {\lfloor \tfrac{n_i}{2} \rfloor}$, $a_i\in \tfrac{n-n_i}{2}+\Z$ such that $\delta_i = |\det|^{a_i}$, and $(\underline{\kappa}^{(i)};\,\underline{\sf w}^{(i)})$ be the infinity type of $\pi_i \otimes \delta_i$. 
It is easy to verify that
\begin{align}\label{E:combinatorial proof 0}
(\mu^{(1)},...,\mu^{(k)})+\rho_n = ( \mu^{(1)}+\rho_{n_1}+\underline{a_1},...,\mu^{(k)}+\rho_{n_k}+\underline{a_k}),
\end{align}
where $\underline{a_i} = (a_i,...,a_i) \in \Q^{n_i}$ for $1 \leq i \leq k$.
By definition of the infinity type, we have
\[
\mu^{(i)}+\rho_{n_i} = 
\begin{cases}
\left(\tfrac{\kappa_1^{(i)}-1-{\sf w}_1^{(i)}}{2},...,\tfrac{\kappa_{r_i}^{(i)}-1-{\sf w}_{r_i}^{(i)}}{2},\tfrac{1-\kappa_{r_i}^{(i)}-{\sf w}_{r_i}^{(i)}}{2},...,\tfrac{1-\kappa_1^{(i)}-{\sf w}_1^{(i)}}{2}\right) & \mbox{ if $n_i=2r_i$},\\
\left(\tfrac{\kappa_1^{(i)}-1-{\sf w}_1^{(i)}}{2},...,\tfrac{\kappa_{r_i}^{(i)}-1-{\sf w}_{r_i}^{(i)}}{2},-\tfrac{{\sf w}_{r_i+1}^{(i)}}{2},\tfrac{1-\kappa_{r_i}^{(i)}-{\sf w}_{r_i}^{(i)}}{2},..., \tfrac{1-\kappa_1^{(i)}-{\sf w}_1^{(i)}}{2}\right) & \mbox{ if $n_i=2r_i+1$}.
\end{cases}
\]
For $1 \leq i < j \leq k$, by condition (1) and (\ref{E:parity}), $L(s,\pi_i \times \pi_j^\vee)$ is a finite products of some $\Gamma_\C(s+m)$ for $m \in \Z$.
More precisely, for $1 \leq \alpha \leq r_i$ and $1 \leq \beta \leq r_j$, 
\[
\Gamma_\C\left(s+\tfrac{|\kappa_\alpha^{(i)} - \kappa_\beta^{(j)}|+{\sf w}_\alpha^{(i)}-{\sf w}_\beta^{(j)}}{2}-a_i+a_j\right)\Gamma_\C\left(s+\tfrac{\kappa_\alpha^{(i)} + \kappa_\beta^{(j)}+{\sf w}_\alpha^{(i)}-{\sf w}_\beta^{(j)}}{2}+a_i-a_j-1\right)
\]
appear as factors of $L(s,\pi_i \times \pi_j^\vee)$. 
If $n_i$ is odd or $n_j$ is odd, then either
\[
\Gamma_\C\left(s+\tfrac{\kappa_\beta^{(j)} - 1+{\sf w}_{r_i+1}^{(i)}-{\sf w}_{\beta}^{(j)}}{2}-a_i+a_j\right) 
\mbox{\, or \,}
\Gamma_\C\left(s+\tfrac{\kappa_\alpha^{(i)} - 1+{\sf w}_\alpha^{(i)}-{\sf w}_{r_j+1}^{(j)}}{2}-a_i+a_j\right)
\]
appears as a factor of $L(s,\pi_i \times \pi_j^\vee)$. 
Similarly for $L(s,\pi_i^\vee \times \pi_j)$. Therefore, if follows from condition (3) that
\begin{align}\label{E:combinatorial proof 1}
\begin{cases}
1-\tfrac{|\kappa_\alpha^{(i)} - \kappa_\beta^{(j)}|}{2} \leq \tfrac{{\sf w}_\alpha^{(i)}-{\sf w}_\beta^{(j)}}{2}-a_i+a_j \leq -1+\tfrac{|\kappa_\alpha^{(i)} - \kappa_\beta^{(j)}|}{2} & \\
\tfrac{3- \kappa_\beta^{(j)}}{2} \leq \tfrac{{\sf w}_{r_i+1}^{(i)}-{\sf w}_\beta^{(j)}}{2}-a_i+a_j \leq \tfrac{ \kappa_\beta^{(j)}-3}{2} & \mbox{ if $n_i$ is odd},\\
\tfrac{3- \kappa_\alpha^{(i)}}{2} \leq \tfrac{{\sf w}_{r_j+1}^{(j)}-{\sf w}_\alpha^{(i)}}{2}-a_j+a_i \leq \tfrac{ \kappa_\alpha^{(i)}-3}{2} & \mbox{ if $n_j$ is odd}.
\end{cases}
\end{align}
An immediate consequence of (\ref{E:combinatorial proof 1}) is that the entries of $(\mu^{(1)},...,\mu^{(k)})+\rho_n$ are all distinct.
Thus there exists a unique Weyl element $w \in W_n$ such that the entries of $w^{-1}\left( (\mu^{(1)},...,\mu^{(k)})+\rho_n \right)$ are strictly decreasing. As these entries belong to $\tfrac{n-1}{2}+\Z$, we conclude that
\[
\mu:=w^{-1}\left( (\mu^{(1)},...,\mu^{(k)})+\rho_n \right)-\rho_n \in X^+(T_n).
\]
%Since $w\cdot (\mu+\rho_n) = (\mu^{(1)},\cdots,\mu^{(k)})+\rho_n$ and the entries of $\mu^{(i)}+\rho_{n_i}$ are strictly decreasing for $1 \leq i \leq k$, we have $w \in W_n^P$. 
Since the entries of $\mu^{(i)}+\rho_{n_i}$ are strictly decreasing for $1 \leq i \leq k$, by (\ref{E:combinatorial proof 0}) we have $w \in W_n^P$.
Now we show that $\ell(w) = b_n-\sum_{i=1}^kb_{n_i}$. Write $(\mu^{(1)},...,\mu^{(k)})+\rho_n = (\lambda_1,...,\lambda_n)$.
Then 
\[
\ell(w) = {}^\sharp\!\left\{ (\alpha,\beta)\,\vert\,1 \leq \alpha <\beta \leq n,\,\lambda_\alpha<\lambda_\beta \right\}.
\]
Put $m_1=0$ and $m_i = n_1+\cdots+n_{i-1}$ for $2 \leq i \leq k$.
Since $\lambda_{m_i+1}>\cdots>\lambda_{m_i+n_i}$ for $1 \leq i \leq k$, we have 
\[
\ell(w) = \sum_{i<j}{}^\sharp S_{i,j,<}
\]
 with
\[
S_{i,j,<} = \{(\alpha,\beta)\,\vert\, 1\leq \alpha \leq n_i,\,1\leq \beta \leq n_j,\, \lambda_{m_i+\alpha}<\lambda_{m_j+\beta} \},\quad 1 \leq i < j \leq k.
\]
Let $1 \leq i < j \leq k$, $1 \leq \alpha \leq r_i$, and $1 \leq \beta \leq r_j$. Then (\ref{E:combinatorial proof 1}) implies that 
\begin{align}\label{E:combinatorial proof 2}
\begin{cases}
\lambda_{m_i+\alpha}>\lambda_{m_j+\beta}> \lambda_{m_j+n_j+1-\beta} > \lambda_{m_i+n_i+1-\alpha}& \mbox{ if $\kappa_\alpha^{(i)}>\kappa_\beta^{(j)}$},\\
\lambda_{m_j+\beta}>\lambda_{m_i+\alpha}> \lambda_{m_i+n_i+1-\alpha}>\lambda_{m_j+n_j+1-\beta}  & \mbox{ if $\kappa_\beta^{(j)}>\kappa_\alpha^{(i)}$},\\
\lambda_{m_j+\beta}> \lambda_{m_i+r_i+1} >\lambda_{m_j+n_j+1-\beta} & \mbox{ if $n_i$ is odd},\\
\lambda_{m_i+\alpha}> \lambda_{m_j+r_j+1} >\lambda_{m_i+n_i+1-\alpha} & \mbox{ if $n_j$ is odd}.
\end{cases}
\end{align}
Note that 
\[
\{\lambda_{m_i+1},...,\lambda_{m_i+n_i}\} = \bigcup_{1 \leq \alpha \leq r_i}\{\lambda_{m_i+\alpha},\lambda_{m_i+n_i+1-\alpha}\} \cup \begin{cases}
\emptyset & \mbox{ if $n_i=2r_i$},\\
\{\lambda_{m_i+r_i+1}\} & \mbox{ if $n_i=2r_i+1$}.
\end{cases}
\]
We thus deduce from (\ref{E:combinatorial proof 2}) that ${}^\sharp S_{i,j,<} = \tfrac{1}{2}n_in_j$. Therefore, 
\[
\ell(w) = \tfrac{1}{2}\sum_{i<j}n_in_j = b_n-\sum_{i=1}^kb_{n_i}.
\]
Finally, it is clear that $\pi \in \Omega_\mu$ by writing $\pi_i$'s as induced representations of the form (\ref{E:cohomological}) and using (\ref{E:combinatorial proof 0}).
Note that the inducing data of $\pi$ from that of the $\pi_i$'s can be rearranged by the result of Shahidi on local coefficients in \cite[Theorem 3.1]{Shahidi1985} and the irreducibility criterion of Speh in \cite[Theorem 10b]{Moeglin1996}.
This completes the proof.
\end{proof}

\begin{rmk}\label{R:Kostant}
In \cite[Definition 5.9]{HR2020}, a Kostant representative $w \in W_n^P$ is called balanced if $\ell(w) = \tfrac{1}{2}{\rm dim}(U_P)$. In general, we have $0 \leq \ell(w) \leq {\rm dim}(U_P)$. Thus the unique $w \in W_n^P$ in Lemma \ref{L:combinatorial} is an example of a balanced Kostant representative.
\end{rmk}

%\subsection{Tamely isobaric automorphic representations}\label{SS:isobaric}

\subsection{Betti--Whittaker periods (I)}\label{SS:BW periods}

Let $\itPi \cong \otimes_v \itPi_v$ be a regular algebraic cuspidal automorphic representation of $\GL_n(\A)$.
By the results of Clozel \cite{Clozel1990}, under the cuspidality condition, $\itPi$ is regular algebraic if and only if $\itPi_\infty \in \Omega(n)$.
Moreover, in this case, $\itPi_\infty$ is essentially tempered.
For $\sigma \in {\rm Aut}(\C)$, let ${}^\sigma\!\itPi$ be the irreducible admissible $((\frak{g}_n,{\rm O}_n(\R)) \times \GL_n(\A_f))$-module defined by
\[
{}^\sigma\!\itPi = \itPi_\infty \otimes {}^\sigma\!\itPi_f,
\]
where ${}^\sigma\!\itPi_f$ is the $\sigma$-conjugate of $\itPi_f = \otimes_p\itPi_p$.
Then ${}^\sigma\!\itPi$ is cuspidal automorphic by the result of Clozel \cite[Th\'eor\`eme 3.19]{Clozel1990}.
Let $\Q(\itPi)$ be the fixed field of $\{\sigma \in {\rm Aut}(\C)\,\vert\,{}^\sigma\!\itPi \cong \itPi\}$. It is a number field (cf.\,\textit{loc.\ cit.}) and we call it the rationality field of $\itPi$.
In this section, we recall in Lemma \ref{L:BW period} the definition of the Betti--Whittaker periods of $\itPi$ introduced by Mahnkopf \cite{Mahnkopf2005} with respect to a fixed choice of generators for the bottom degree relative Lie algebra cohomology, obtained by comparing two $\Q(\itPi)$-rational structures on $\itPi_f$ recalled in \S\,\ref{SS:RS Whittaker} and \S\,\ref{SS:rational sheaf} (see also \cite{Harder1983}, \cite{Hida1994} for $n=2$ and \cite{RS2008} for general number field).

\subsubsection{Rational structure via Whittaker model}\label{SS:RS Whittaker}

For a cusp form $\varphi \in \itPi$, the Whittaker function $W(\varphi)$ with respect to $\psi_{U_n}$ is defined by
\begin{align}\label{E:Whittaker functional}
W(g,\varphi) = \int_{U_n(\Q)\backslash U_n(\A)}\varphi(ug)\overline{\psi_{U_n}(u)}\,du^{\rm Tam},\quad g \in \GL_n(\A).
\end{align}
Here $du^{\rm Tam}$ is the Tamagawa measure on $U_n(\A)$.
Let $\mathcal{W}(\itPi)$ be the space of Whittaker functions of $\itPi$ with respect to $\psi_{U_n}$.
It is well-known that $\itPi$ is globally generic, that is, $\mathcal{W}(\itPi)$ is non-zero.
For each place $v$ of $\Q$, let $\mathcal{W}(\itPi_v)$ be the space of Whittaker functions of $\itPi_v$ with respect to $\psi_{U_n,v}$. 
When $v=p$ is a finite place such that $\itPi_p$ is unramified, let $W_{\itPi_p}^\circ \in \mathcal{W}(\itPi_p)$ be the $\GL_n(\Z_p)$-invariant Whittaker function normalized so that $W_{\itPi_p}^\circ({\bf 1}_n)=1$.
Let $\mathcal{W}(\itPi_f) = \prod_p \mathcal{W}(\itPi_p)$ be the restricted product with respect to $W_{\itPi_p}^\circ$ for unramified $p$.
Then 
\[
\mathcal{W}(\itPi) = \prod_v \mathcal{W}(\itPi_v)
\]
by the multiplicity one theorem for Whittaker functionals of $\GL_n$.
For $W_\infty \in \mathcal{W}(\itPi_\infty)$ and $W_f \in \mathcal{W}(\itPi_f)$, there exists a unique cusp form $\varphi \in \itPi$ such that
\[
W(\varphi) = W_\infty\cdot W_f.
\]
We obtain an $((\frak{g}_n,{\rm O}_n(\R))\times \GL_n(\A_f))$-equivariant isomorphism
\begin{align}\label{E:Whittaker isomorphism}
\Upsilon_\itPi:\mathcal{W}(\itPi) \longrightarrow \itPi.
\end{align}
For $\sigma \in {\rm Aut}(\C)$, let $t_\sigma : \mathcal{W}(\itPi_f) \rightarrow \mathcal{W}({}^\sigma\!\itPi_f)$ be the $\sigma$-linear $\GL_n(\A_f)$-equivariant isomorphism defined by
\begin{align}\label{E:sigma-linear Whittaker}
t_\sigma W(g) = \sigma\left(W({\rm diag}(u_\sigma^{-n+1},u_\sigma^{-n+2},...,1)g)\right),\quad g \in \GL_n(\A_f).
\end{align}
Here $u_\sigma \in \widehat{\Z}^\times$ is the unique element such that $\sigma(\psi(x)) = \psi(u_\sigma x)$ for all $x \in \A_f$.
We thus obtain a $\Q(\itPi)$-rational structure on $\mathcal{W}(\itPi_f)$ given by taking the Galois invariants over $\Q(\itPi)$:
\begin{align}\label{E:rational structure 1}
\mathcal{W}(\itPi_f)^{{\rm Aut}(\C/\Q(\itPi))} = \left.\left\{W\in \mathcal{W}(\itPi_f)\,\right\vert\,t_\sigma W=W\mbox{ for $\sigma \in {\rm Aut}(\C/\Q(\itPi))$}\right\}.
\end{align}
Note that $\mathcal{W}(\itPi_f)^{{\rm Aut}(\C/\Q(\itPi))}$ is non-zero by the newform theory for $\GL_n(\A_f)$ \cite{JPSS1981}, the existence of the Kirillov model for $\itPi_f$ \cite[Theorem 5.20]{BZ1976}, and \cite[Lemme I.1]{Wald1985B}.

\subsubsection{Rational structure via sheaf cohomology}\label{SS:rational sheaf}
By the results \cite[Lemme 3.14 and Lemme 4.9]{Clozel1990} due to Clozel, we have $\itPi_\infty \in \Omega_\mu$ for some uniquely determined $\mu \in X^+(T_n)$ which is pure, that is, $\mu_i+\mu_{n-i+1} = \mu_1+\mu_n$ for $1 \leq i \leq n$.
In particular, $\itPi_\infty$ is essentially tempered and we have
\[
H^{b_n}(\frak{g}_n,K_n^\circ;\itPi_\infty\otimes M_{\mu,\C}) \neq 0.
\]
Consider the topological space
\begin{align}\label{E:orbifold}
\mathcal{S}_n = \GL_n(\Q)\backslash \GL_n(\A) /K_n^\circ. 
\end{align}
The algebraic representation $M_{\mu}$ defines a sheaf
$\mathcal{M}_{\mu}$ of $\Q$-vector spaces on $\mathcal{S}_n$ (cf.\,\cite[\S\,2.2.8]{HR2020}).
We denote by
\[
H^\bullet(\mathcal{S}_n,\mathcal{M}_{\mu})
\]
the sheaf cohomology groups of $\mathcal{M}_{\mu}$.
It is naturally endowed with an action of $\pi_0(\GL_n(\R))\times \GL_n(\A_f)$ (cf.\,\cite[\S\,2.3]{HR2020}).
When we go to the transcendental level, the base change ${M}_{\mu,\C} = {M}_\mu\otimes_\Q\C$ defines a sheaf $\mathcal{M}_{\mu,\C}$ of $\C$-vector spaces on $\mathcal{S}_n$ and we have
\begin{align}\label{E:transcendental}
H^\bullet(\mathcal{S}_n,\mathcal{M}_{\mu,\C}) = H^\bullet(\mathcal{S}_n,\mathcal{M}_\mu) \otimes_\Q \C.
\end{align}
At the transcendental level, via the de Rham complex, the cohomology groups $H^\bullet(\mathcal{S}_n,\mathcal{M}_{\mu,\C})$ are canonically isomorphic to the $(\frak{g}_n,K_n^\circ)$-cohomology groups of $C^\infty(\GL_n(\Q)\backslash\GL_n(\A),\xi_\mu) \otimes M_{\mu,\C}$.
Here $\xi_\mu$ is the algebraic character of $\R_+$ such that $\R_+$ acts on $M_{\mu,\C}$ by $\xi_\mu^{-1}$, and $C^\infty(\GL_n(\Q)\backslash\GL_n(\A),\xi_\mu)$ is the space of smooth functions on $\GL_n(\Q)\backslash\GL_n(\A)$ such that the action of right multiplication by $\R_+$ is given by the character $\xi_\mu$.  
%Moreover, we have the celebrated result of Franke \cite[Theorem 18]{Franke1998} which shows that we can replace the space of smooth functions on $\GL_n(\Q)\backslash\GL_n(\A)$ by the space $\mathcal{A}(\GL_n)$ of automorphic forms on $\GL_n(\A)$.
Let $H^\bullet_{\rm cusp}(\mathcal{S}_n,\mathcal{M}_{\mu,\C})$ be the cuspidal cohomology of $\GL_n$ with coefficients in $M_{\mu,\C}$ defined by
\[
H^\bullet_{\rm cusp}(\mathcal{S}_n,\mathcal{M}_{\mu,\C}) = H^\bullet(\frak{g}_n,K_n^\circ;\mathcal{A}_0(\GL_n(\A),\xi_\mu) \otimes M_{\mu,\C}),
\]
where $\mathcal{A}_0(\GL_n(\A),\xi_\mu) = \mathcal{A}_0(\GL_n(\A))\cap C^\infty(\GL_n(\Q)\backslash\GL_n(\A),\xi_\mu)$. %is the space of cusp forms on $\GL_n(\A)$ such that the action of right multiplication by $\R_+$ is $\xi_\mu$.
The natural inclusion 
\[
\iota : \mathcal{A}_0(\GL_n(\A),\xi_\mu) \hookrightarrow C^\infty(\GL_n(\Q)\backslash\GL_n(\A),\xi_\mu)
\] 
then induces an $(\pi_0(\GL_n(\R))\times \GL_n(\A_f))$-equivariant homomorphism 
\begin{align}\label{E:Psi}
\iota^\bullet : H^\bullet_{\rm cusp}(\mathcal{S}_n,\mathcal{M}_{\mu,\C}) \longrightarrow H^\bullet(\mathcal{S}_n,\mathcal{M}_{\mu,\C}).
\end{align}
%It is well-known that $\iota^\bullet$ is injective.
%We mention that the image of $\Psi_\itPi$ lies in the Eisenstein cohomology (resp.\,cuspidal cohomology) if $k \geq 2$ (resp.\,$k=1$).
For $\sigma \in {\rm Aut}(\C)$, we have the $\sigma$-linear isomorphism $M_{\mu,\C} \rightarrow M_{\mu,\C}$ defined by
\[
{\bf v}\otimes z \longmapsto {\bf v} \otimes \sigma(z)
\]
for ${\bf v}\otimes z \in M_{\mu,\C} = M_{\mu}\otimes_\Q\C$.
It naturally induces a $\sigma$-linear $(\pi_0(\GL_n(\R))\times \GL_n(\A_f))$-equivariant isomorphism
\[
\sigma^\bullet : H^\bullet(\mathcal{S}_n,\mathcal{M}_{\mu,\C}) \longrightarrow H^\bullet(\mathcal{S}_n,\mathcal{M}_{\mu,\C}).
\]
%Following the proof of \cite[Propositions 1.6 and 1.7]{Grobner2018b} (based on the result of Grobner \cite[Theorem 4]{Grobner2013} which is a refinement of the result of Franke \cite[Theorem 14]{Franke1998}), %(cf.\,\cite[Theorem 7.23]{GR2014b} and \cite[Proposition 4.3]{LLS2021}), we deduce that $\Psi_\itPi^{b_n}$ is injective and the image of $\sigma^*\circ \Psi_\itPi^{b_n}$ is equal to the image of $\Psi_{{}^\sigma\!\itPi}^{b_n}$. 
%We will write $\Psi_\itPi^{b_n} = \Psi_\itPi$ for brevity.
By \cite[Th\'{e}or\`{e}me 3.19]{Clozel1990}, we have
\begin{align}\label{E:Clozel}
\sigma^\bullet\circ\iota^\bullet(H^\bullet_{\rm cusp}(\mathcal{S}_n,\mathcal{M}_{\mu,\C})) = \iota^\bullet(H^\bullet_{\rm cusp}(\mathcal{S}_n,\mathcal{M}_{\mu,\C})).
\end{align}
%This induces a $\sigma$-linear $\pi_0(\GL_n(\R))\times \GL_n(\A_f)$-equivariant isomorphism
%\[
%\sigma^* : H^{b_n}(\frak{g}_n,K_n^\circ;\itPi \otimes M_{\mu,\C}) \longrightarrow H^{b_n}(\frak{g}_n,K_n^\circ;{}^\sigma\!\itPi \otimes M_{\mu,\C}).
%\]
By the strong multiplicity one theorem for $\GL_n$, the $\itPi_f$-isotypic component of the cuspidal cohomology is equal to $H^\bullet(\frak{g}_n,K_n^\circ;\itPi \otimes M_{\mu,\C})$. %We denote by $\iota_\itPi^\bullet$ the restriction of $\iota^\bullet$ to the $\itPi_f$-isotypic component.
Let $\varepsilon \in \{\pm1\}$ if $n$ is even, and $\varepsilon = \varepsilon(\itPi_\infty)$ if $n$ is odd.
Let
\[
H^\bullet(\mathcal{S}_n,\mathcal{M}_{\mu,\C})[\varepsilon\times\itPi_f] := \iota^\bullet(H^\bullet(\frak{g}_n,K_n^\circ;\itPi \otimes M_{\mu,\C})[\varepsilon]).
\]
Since $\iota^\bullet$ is injective (cf.\,\cite[Corollary 5.5]{Borel1981} and \cite[Lemme 3.15]{Clozel1990}) and $\iota^\bullet,\,\sigma^\bullet$ are equivariant, by (\ref{E:Clozel}) we have 
\[
\sigma^\bullet(H^\bullet(\mathcal{S}_n,\mathcal{M}_{\mu,\C})[\varepsilon\times\itPi_f]) = H^\bullet(\mathcal{S}_n,\mathcal{M}_{\mu,\C})[\varepsilon\times{}^\sigma\!\itPi_f].
\]
%When $\bullet = b_n$, we write $\iota^* = \iota^{b_n}$ and $\sigma^* = \sigma^{b_n}$.
Since $H^{b_n}(\frak{g}_n,K_n^\circ;\itPi_\infty \otimes M_{\mu,\C})[\varepsilon]$ is one-dimensional, by taking the Galois invariants over $\Q(\itPi)$, we obtain a $\Q(\itPi)$-rational structure on $H^{b_n}(\mathcal{S}_n,\mathcal{M}_{\mu,\C})[\varepsilon\times\itPi_f]\cong\itPi_f$:
\begin{align}\label{E:rational structure 2}
\begin{split}
&H^{b_n}(\mathcal{S}_n,\mathcal{M}_{\mu,\C})[\varepsilon\times\itPi_f]^{{\rm Aut}(\C/\Q(\itPi))}\\
& = \left.\left\{\omega\in H^{b_n}(\mathcal{S}_n,\mathcal{M}_{\mu,\C})[\varepsilon\times\itPi_f]\,\right\vert\,\sigma^{b_n} \omega=\omega\mbox{ for $\sigma \in {\rm Aut}(\C/\Q(\itPi))$}\right\}.
\end{split}
\end{align}
Now we fix a generator
\[
[\itPi_\infty]^\varepsilon \in H^{b_n}(\frak{g}_n,K_n^\circ; \mathcal{W}(\itPi_\infty) \otimes M_{\mu,\C})[\varepsilon].
\]
Let 
\[
\Phi_\itPi^\varepsilon : \mathcal{W}(\itPi_f) \longrightarrow H^{b_n}(\mathcal{S}_n,\mathcal{M}_{\mu,\C})[\varepsilon\times\itPi_f]
\]
be the $\GL_n(\A_f)$-equivariant isomorphism defined by
\begin{align}\label{E:Betti-Whittaker}
\Phi_\itPi^\varepsilon = \iota^{b_n}\circ\Upsilon_\itPi^{b_n}\circ([\itPi_\infty]^\varepsilon\otimes\,\cdot\,).
\end{align}
Here $\Upsilon_\itPi^{b_n} : H^{b_n}(\frak{g}_n,K_n^\circ; \mathcal{W}(\itPi) \otimes M_{\mu,\C}) \rightarrow H^{b_n}(\frak{g}_n,K_n^\circ; \itPi \otimes M_{\mu,\C})$ is the $(\pi_0(\GL_n(\R)) \times \GL_n(\A_f))$-equivariant isomorphism induced by $\Upsilon_\itPi$ in (\ref{E:Whittaker isomorphism}).
%as follows: For $W \in \mathcal{W}(\itPi_f)$, we have
%\[
%[\itPi_\infty]^\varepsilon\otimes W \in H^{b_n}(\frak{g}_n,K_n^\circ; \mathcal{W}(\itPi_\infty) \otimes \mathcal{W}(\itPi_f)\otimes M_{\mu,\C})[\varepsilon].
%\] 
%Then $\Phi_\itPi^\varepsilon(W)$ is the image of $[\itPi_\infty]^\varepsilon\otimes W$ under the $\GL_n(\A_f)$-equivariant isomorphism induced by the isomorphism (\ref{E:Whittaker isomorphism}).
Comparing the $\Q(\itPi)$-rational structures given by (\ref{E:rational structure 1}) and (\ref{E:rational structure 2}), we have the following lemma/definition for the Betti--Whittaker periods of $\itPi$.

\begin{lemma}\label{L:BW period}
Let $\varepsilon \in \{\pm1\}$ if $n$ is even, and $\varepsilon = \varepsilon(\itPi_\infty)$ if $n$ is odd.
There exists $p(\itPi,\varepsilon) \in \C^\times$, unique up to $\Q(\itPi)^\times$, such that
\[
\frac{\Phi_\itPi^\varepsilon(\mathcal{W}(\itPi_f)^{{\rm Aut}(\C/\Q(\itPi))})}{p(\itPi,\varepsilon)} = H^{b_n}(\mathcal{S}_n,\mathcal{M}_{\mu,\C})[\varepsilon\times\itPi_f]^{{\rm Aut}(\C/\Q(\itPi))}.
\]
Moreover, we can normalize the periods so that
\begin{align*}
\sigma^{b_n}\circ\left( \frac{\Phi_\itPi^\varepsilon}{p(\itPi,\varepsilon)} \right) = \left(\frac{\Phi_{{}^\sigma\!\itPi}^\varepsilon}{p({}^\sigma\!\itPi,\varepsilon)}\right)\circ t_\sigma,\quad \sigma \in {\rm Aut}(\C).
%& = \frac{p({}^\sigma\!\itPi,\varepsilon)}{\prod_{1 \leq i < j \leq k}G({}^\sigma\!\omega_{\itPi_j\rho_j})^{n_i}L^{(\infty)}(1,{}^\sigma\!(\itPi_i\rho_i)\rho_i^{-1} \times {}^\sigma\!(\itPi_j^\vee\rho_j^{-1})\rho_j)\cdot \prod_{i=1}^{k} p({}^\sigma\!(\itPi_i\rho_i),\varepsilon)}.
\end{align*} 
\end{lemma}

It is clear that the Betti--Whittaker period $p(\itPi,\varepsilon)$ of $\itPi$ depends on the choice of $[\itPi_\infty]^\varepsilon$. In \S\,\ref{SS:generator} below, we will specify a canonical choice of generator.
%When $n$ is odd, for brevity we will also write 
%\begin{align}\label{E:convention}
%p(\itPi) = p(\itPi,\varepsilon).
%\end{align}

%In the following theorem, we recall the result of Raghuram and Shahidi \cite{RS2008} on the period relations for Betti--Whittaker periods upon twisting by algebraic Hecke characters. 

%\begin{thm}[Raghuram--Shahidi]\label{T:twisted period relation}
%Let $\itPi$ be a cohomological tamely isobaric automorphic representation of $\GL_n(\A)$ and $\chi$ an algebraic Hecke character of $\A^\times$.
%Let $\varepsilon \in \{\pm1\}$ if $n$ is even, and $\varepsilon = \varepsilon(\itPi_\infty)$ if $n$ is odd.
%Assume the generators $[\itPi_\infty]^\varepsilon$ and $[\itPi_\infty\otimes\chi_\infty]^{\varepsilon\cdot\varepsilon(\chi_\infty)}$ satisfy the compatibility relation in Lemma \ref{L:compatibility} below.
%For $\sigma \in {\rm Aut}(\C)$, we have
%\[
%\sigma \left( \frac{p(\itPi \otimes \chi,\varepsilon\cdot\varepsilon(\chi_\infty))}{G(\chi)^{n(n-1)/2}\cdot p(\itPi,\varepsilon)} \right) = \frac{p({}^\sigma\!\itPi \otimes {}^\sigma\!\chi,\varepsilon\cdot \varepsilon(\chi_\infty))}{G({}^\sigma\!\chi)^{n(n-1)/2}\cdot p({}^\sigma\!\itPi,\varepsilon)}.
%\]
%\end{thm}

%\begin{rmk}
%The result of Raghuram and Shahidi is stated for cuspidal $\itPi$. Nonetheless, the proof goes without change for general case. One can also deduce the general case from the cuspidal one by using the period relation in Theorem \ref{T:period relation} below.
%\end{rmk}

\subsection{Betti--Whittaker periods (II)}\label{SS:BW periods 2}

Let $\itPi$ be a regular algebraic automorphic representation of $\GL_n(\A)$ (cf.\,\cite[Definition 1.8]{Clozel1990}). It means that $\itPi$ is isobaric and the infinitesimal character of $\itPi_\infty$ is regular and belongs to $(\Z+\tfrac{n-1}{2})^n$. 
We have
\[
\itPi = \itPi_1 \boxplus \cdots \boxplus \itPi_k
\]
for some cuspidal automorphic representation $\itPi_i$ of $\GL_{n_i}(\A)$ for $1 \leq i \leq k$. 
%Let $P=P_{(n_1,\cdots,n_k)}$ be a standard parabolic subgroup of $\GL_n$ of type $(n_1,\cdots,n_k)$.
%For $1\leq i \leq k$, let $\delta_i = \delta_P \vert_{\GL_{n_i}(\A)}$, where $\delta_P$ is the square-root of the modulus character of $P(\A)$.
It is clear that $\itPi_i\otimes|\mbox{ }|_\A^{(n-n_i)/2}$ is regular algebraic for $1 \leq i \leq k$.
Therefore, for any $\sigma \in {\rm Aut}(\C)$, ${}^\sigma\!\itPi:=\itPi_\infty\otimes {}^\sigma\!\itPi_f$ is automorphic and
\[
{}^\sigma\!\itPi = \bigboxplus_{i=1}^k{}^\sigma\!(\itPi_i\otimes|\mbox{ }|_\A^{(n-n_i)/2})\otimes|\mbox{ }|_\A^{(n_i-n)/2}.
\]
Let $P=P_{(n_1,...,n_k)}$ be the standard parabolic subgroup of $\GL_n$ of type $(n_1,...,n_k)$ and $\delta_P$ the square-root of the modulus character of $P(\A)$.
For $1\leq i \leq k$, let $\delta_i = \delta_P \vert_{\GL_{n_i}(\A)}$ and we write $\itPi_i\delta_i:=\itPi_i\otimes\delta_i$. Note that $\itPi_i\delta_i$ is regular algebraic, since $\delta_i$ differs from $|\mbox{ }|_\A^{(n-n_i)/2}$ by an integral power of $|\mbox{ }|_\A$.
We give the following definition of the Betti--Whittaker periods of $\itPi$, which is inspired by the result of Grobner and Lin \cite[Theorem 2.6]{GL2021} where the base field is a CM-field and $\itPi_\infty$ is essentially tempered.
\begin{definition}\label{D:Betti-Whittaker}
Assume the following conditions are satisfied:
\begin{itemize}
\item[(1)] At most one of the $n_i$'s is odd.
\item[(2)] $L(s,\itPi_{i,\infty}\times \itPi^\vee_{j,\infty})$ and $L(s,\itPi_{i,\infty}^\vee\times \itPi_{j,\infty})$ are holomorphic at $s=0$ for $1 \leq i<j \leq k$.
\item[(3)] $e(\itPi_1)\geq \cdots\geq e(\itPi_k)$, where $e(\itPi_i)$ is the unique real number such that $\itPi_i \otimes |\mbox{ }|_\A^{-e(\itPi_i)}$ is unitary.
\end{itemize}
Let $\varepsilon \in \{\pm1\}$ if $n$ is even, and $\varepsilon = \varepsilon(\itPi_\infty)$ if $n$ is odd.
We define $p(\itPi,\varepsilon) \in \C^\times / \Q(\itPi)^\times$ by
\[
p(\itPi,\varepsilon) = \prod_{1 \leq i < j \leq k}G(\omega_{\itPi_j\delta_j})^{n_i}L^{(\infty)}(1,\itPi_i \times \itPi_j^\vee)\cdot \prod_{i=1}^k p(\itPi_i\delta_i,\varepsilon),
\]
where $L^{(\infty)}(s,\itPi_i \times \itPi_j^\vee)$ is finite part of the Rankin--Selberg $L$-function of $\itPi_i \times \itPi_j^\vee$.
\end{definition}

\begin{rmk}\label{R:Betti-Whittaker}
Conditions (1) and (2) imply that $\itPi_{i,\infty}$ and $\itPi_{j,\infty} \otimes |\mbox{ }|^t$ are inequivalent for any $t \in \C$ and $i \neq j$.
In particular, the rationality field $\Q(\itPi)$ of $\itPi$ is equal to the composite of the rationality field of $\itPi_i\otimes|\mbox{ }|_\A^{(n-n_i)/2}$ for $1 \leq i \leq k$.
Also the $L$-function $L(s,\itPi_i\times\itPi_j^\vee)$ is entire and non-vanishing for ${\rm Re}(s)\geq 1-e(\itPi_i)+e(\itPi_j)$. Condition (3) then implies that the $L(s,\itPi_i\times\itPi_j^\vee)$ is non-vanishing at $s=1$ for $i < j$.
\end{rmk}

\begin{rmk}
The conditions in Definition \ref{D:Betti-Whittaker} are satisfied when $\itPi_\infty$ is essentially tempered.
\end{rmk}

In the following lemma, we show that the Betti--Whittaker periods are well-defined.
\begin{lemma}\label{L:well-defined}
The Betti--Whittaker periods of $\itPi$ are independent of the ordering of the cuspidal summands subject to the conditions in Definition \ref{D:Betti-Whittaker}.
\end{lemma}

\begin{proof}
It is clear that we may assume $k=2$ with $e(\itPi_1)=e(\itPi_2)$. 
Let $\delta_1,\delta_2$ (resp.\,$\delta_1',\delta_2'$) be the characters defined with respect to $P_{(n_1,n_2)}$ (resp.\,$P_{(n_2,n_1)}$).
We have to show that
\begin{align}\label{E:well-defined pf 0}
\begin{split}
&G(\omega_{\itPi_2\delta_2})^{n_1}\cdot L^{(\infty)}(1,\itPi_1 \times \itPi_2^\vee)\cdot p(\itPi_1\delta_1,\varepsilon)\cdot p(\itPi_2\delta_2,\varepsilon)\\
& \sim G(\omega_{\itPi_1\delta_2'})^{n_2}\cdot L^{(\infty)}(1,\itPi_1^\vee \times \itPi_2)\cdot p(\itPi_1\delta_2',\varepsilon)\cdot p(\itPi_2\delta_1',\varepsilon).
\end{split}
\end{align}
Note that $\delta_1 = |\mbox{ }|_\A^{n_2/2}$, $\delta_2 = |\mbox{ }|_\A^{-n_1/2}$, $\delta_1' = |\mbox{ }|_\A^{n_1/2}$, $\delta_2' = |\mbox{ }|_\A^{-n_2/2}$. 
By the period relation proved by Raghuram and Shahidi \cite{RS2008}, we have
\begin{align}\label{E:well-defined pf 1}
p(\itPi_1\delta_1,\varepsilon) \sim p(\itPi_1\delta_2',(-1)^{n_2}\varepsilon), \quad p(\itPi_2\delta_2,\varepsilon) \sim p(\itPi_2\delta_1',(-1)^{n_1}\varepsilon).
\end{align}
By condition (2), $L(s, \itPi_1\delta_1 \times (\itPi_2\delta_2)^\vee)$ has consecutive critical points $m_0$ and $m_0+1$, where $m_0 = \tfrac{-n_1-n_2}{2}$ (see \S\,\ref{SS:RS} below). By condition (1), we may assume further that $n_1$ is even. In \cite[Theorem 7.21]{HR2020}, Harder and Raghuram prove the following result on the algebraicity of the ratio:
\begin{align}\label{E:well-defined pf 2}
\begin{split}
\frac{L^{(\infty)}(m_0,\itPi_1\delta_1 \times (\itPi_2\delta_2)^\vee)}{L^{(\infty)}(m_0+1,\itPi_1\delta_1 \times (\itPi_2\delta_2)^\vee)} \sim (2\pi\sqrt{-1})^{-n_1n_2/2}\cdot \begin{cases}
1 & \mbox{ if $n_2$ is even},\\
\displaystyle{\frac{p(\itPi_1\delta_1,\varepsilon)}{p(\itPi_1\delta_1,-\varepsilon)}} & \mbox{ if $n_2$ is odd}.
\end{cases}
\end{split}
\end{align}
On the other hand, we prove in \cite[Theorem 3.2]{Chen2023} the algebraicity of critical values under duality:
\begin{align}\label{E:well-defined pf 3}
\begin{split}
&L^{(\infty)}(m_0,\itPi_1\delta_1 \times (\itPi_2\delta_2)^\vee) \\
&\sim \gamma(m_0,(\itPi_1\delta_1)_\infty \times (\itPi_2\delta_2)_\infty^\vee,\psi_\infty)\cdot G(\omega_{\itPi_1\delta_1})^{n_2}G(\omega_{\itPi_2\delta_2})^{-n_1}\cdot L^{(\infty)}(1-m_0,(\itPi_1\delta_1)^\vee \times \itPi_2\delta_2),
\end{split}
\end{align}
where $\gamma(s,(\itPi_1\delta_1)_\infty \times (\itPi_2\delta_2)_\infty^\vee,\psi_\infty)$ is the $\gamma$-factor of $(\itPi_1\delta_1)_\infty \times (\itPi_2\delta_2)_\infty^\vee$ with respect to $\psi_\infty$.
Therefore, (\ref{E:well-defined pf 0}) follows immediately from (\ref{E:well-defined pf 1})-(\ref{E:well-defined pf 3}) and the following relation:
\begin{align}\label{E:well-defined pf 4}
\gamma(m_0,(\itPi_1\delta_1)_\infty \times (\itPi_2\delta_2)_\infty^\vee,\psi_\infty) \in (2\pi\sqrt{-1})^{-n_1n_2/2}\cdot \Q^\times.
\end{align}
It remains to prove (\ref{E:well-defined pf 4}).
Let $r_1 = \lfloor \tfrac{n_1}{2}\rfloor$ and $r_2 = \lfloor \tfrac{n_2}{2}\rfloor$. Let $(\underline{\kappa};\,{\sf w})$ and $(\underline{\ell};\,{\sf u})$ be the infinity types of $\itPi_1\delta_1$ and $\itPi_2\delta_2$ respectively. 
%By condition (1) we may assume $n_1$ is even.
Note that
\begin{align*}
\gamma(m_0,(\itPi_1\delta_1)_\infty \times (\itPi_2\delta_2)_\infty^\vee,\psi_\infty) & = \varepsilon(m_0,(\itPi_1\delta_1)_\infty \times (\itPi_2\delta_2)_\infty^\vee,\psi_\infty)\cdot \frac{L(1-m_0,(\itPi_1\delta_1)_\infty ^\vee\times (\itPi_2\delta_2)_\infty)}{L(m_0,(\itPi_1\delta_1)_\infty \times (\itPi_2\delta_2)_\infty^\vee)}.
\end{align*}
Note that we have assume $n_1$ is even.
By the condition $e(\itPi_1) = e(\itPi_2)$ and formulas of archimedean local factors (cf.\,\cite[(3.6) and (3.7)]{Knapp1994}), the $\varepsilon$-facotr is equal to
\[
\prod_{1 \leq i \leq r_1,1\leq j\leq r_2}(\sqrt{-1})^{\kappa_i+\ell_j+|\kappa_i-\ell_j|}\cdot
\begin{cases}
1 & \mbox{ if $n_2$ is even},\\
\prod_{i=1}^{r_1}(\sqrt{-1})^{\kappa_i} & \mbox{ if $n_2$ is odd},
\end{cases}
\]
and the ratio of $L$-factors is equal to
\begin{align*}
\prod_{1 \leq i \leq r_1,1\leq j\leq r_2} \frac{\Gamma_\C\left(1+\tfrac{\kappa_i+\ell_j-2}{2}\right)\Gamma_\C\left(1+\tfrac{|\kappa_i-\ell_j|}{2}\right)}{\Gamma_\C\left(\tfrac{\kappa_i+\ell_j-2}{2}\right)\Gamma_\C\left(\tfrac{|\kappa_i-\ell_j|}{2}\right)}\cdot\
\begin{cases}
1 & \mbox{ if $n_2$ is even},\\
\prod_{i=1}^{r_1} \frac{\Gamma_\C\left( 1+\tfrac{\kappa_i-1}{2} \right)}{\Gamma_\C\left( \tfrac{\kappa_i-1}{2} \right)} & \mbox{ if $n_2$ is odd}.
\end{cases}
\end{align*}
By (\ref{E:parity}), ${\sf u}$ is even if $n_2$ is odd, and
\[
\kappa_i \equiv {\sf w}\,({\rm mod}\,2),\quad \ell_j \equiv {\sf u}+n_2\,({\rm mod}\,2).
\]
Note that the condition $e(\itPi_1) = e(\itPi_2)$ implies that ${\sf w}+{\sf u} \equiv n_2\,({\rm mod}\,2)$.
In particular, $\kappa_i \equiv \ell_j\,({\rm mod}\,2)$ and $\kappa_i$ is odd if $n_2$ is odd. Therefore, up to $\Q^\times$, the ratio of $L$-factors is equal to $(2\pi)^{-2r_1r_2}$ if $n_2$ is even and $(2\pi)^{-2r_1r_2-r_1}$ if $n_2$ is odd. In both cases, it is equal to $(2\pi)^{-n_1n_2/2}$.
Also, up to $\pm 1$, the $\varepsilon$-factor is equal to 1 if $n_2$ is even and $(\sqrt{-1})^{r_1}$ if $n_2$ is odd. 
In both cases, we can replace it by $(\sqrt{-1})^{-n_1n_2/2}$.
This completes the proof.
\end{proof}

We rewrite the Betti--Whittaker period in the following lemma in a form which is more convenient to apply in the proof of Theorem \ref{T:Main}.

\begin{lemma}\label{L:period relation}

Let $\itSigma = \itSigma_1 \boxplus \cdots \boxplus \itSigma_l$ and $\itPi = \itPi_1 \boxplus \cdots \boxplus \itPi_k$ be regular algebraic automorphic representations of $\GL_n(\A)$ and $\GL_{n'}(\A)$ respectively with $n$ even. 
Let $\delta \in \{0,1\}$ so that $\delta \equiv n'\,({\rm mod}\,2)$.
Assume the algebraic automorphic representation
\[
(\itSigma\otimes|\mbox{ }|_\A^{-\delta/2}) \boxplus \itPi := (\itSigma_1\otimes|\mbox{ }|_\A^{-\delta/2}) \boxplus\cdots\boxplus(\itSigma_l\otimes|\mbox{ }|_\A^{-\delta/2})\boxplus\itPi_1 \boxplus \cdots \boxplus \itPi_k
\]
is regular and satisfies the conditions in Definition \ref{D:Betti-Whittaker}.
Then we have
\[
p((\itSigma\otimes|\mbox{ }|_\A^{-\delta/2}) \boxplus \itPi,\varepsilon) \sim G(\omega_\itPi)^{n}\cdot L^{(\infty)}(1-\tfrac{\delta}{2}, \itSigma \times \itPi^\vee)\cdot p(\itSigma \otimes |\mbox{ }|_\A^{(n'-\delta)/2},\varepsilon)\cdot p(\itPi\otimes|\mbox{ }|_\A^{-n/2},\varepsilon).
\]
\end{lemma}

\begin{proof}
For $1 \leq i \leq l$ and $1 \leq j \leq k$, let $n_i$ and $n_j'$ be the degree of $\itSigma_i$ and $\itPi_j$.
Let $\delta_1,...,\delta_{l+k}$ be the characters defined with respect to $P_{(n_1,...,n_l,n_1',...,n_k')}$.
By definition, we have
\begin{align}\label{E:period relation pf 1}
\begin{split}
p((\itSigma\otimes|\mbox{ }|_\A^{-\delta/2}) \boxplus \itPi,\varepsilon) &= \prod_{1 \leq i < j \leq l}G\left(\omega_{\itSigma_j\delta_j\otimes |\mbox{ }|_\A^{-\delta/2}}\right)^{n_i}L^{(\infty)}(1,\itSigma_i \times \itSigma_j^\vee)\cdot \prod_{i=1}^l p(\itSigma_i\delta_i \otimes |\mbox{ }|_\A^{-\delta/2},\varepsilon)\\
&\times \prod_{1 \leq i < j \leq k}G(\omega_{\itPi_{j}\delta_{l+j}})^{n_{i}'}L^{(\infty)}(1,\itPi_{i} \times \itPi_{j}^\vee)\cdot \prod_{i=1}^k p(\itPi_{i}\delta_{l+i},\varepsilon)\\
&\times \prod_{1 \leq i \leq l \atop 1 \leq j \leq k} G(\omega_{\itPi_j\delta_{l+j}})^{n_i}L^{(\infty)}(1-\tfrac{\delta}{2},\itSigma_i \times \itPi_j^\vee).
\end{split}
\end{align}
Note that
\[
\prod_{1 \leq i \leq l \atop 1 \leq j \leq k} L(s,\itSigma_i \times \itPi_j^\vee) = L(s,\itSigma \times \itPi^\vee).
\]
Let $\delta_1',...,\delta_l'$ and $\delta_1'',...,\delta_k''$ be the characters defined with respect to $P_{(n_1,...,n_l)}$ and $P_{(n_1',...,n_k')}$ respectively.
It is easy to verify that for for $1 \leq i \leq l$ and $1 \leq j \leq k$, we have
\[
\delta_i = \delta_i'|\mbox{ }|_\A^{n'/2},\quad \delta_{l+j} = \delta_j''|\mbox{ }|_\A^{-n/2}.
\]
Hence the product in the first and second lines on the right-hand side of (\ref{E:period relation pf 1}) are equal to $p(\itSigma \otimes |\mbox{ }|_\A^{(n'-\delta)/2},\varepsilon)$ and $p(\itPi\otimes|\mbox{ }|_\A^{-n/2},\varepsilon)$ respectively.
Since $\omega_\itPi$ differs from $\prod_{1 \leq i \leq l } \omega_{\itPi_j\delta_{l+j}}$ by an integral power of $|\mbox{ }|_\A$, by (\ref{E:Galois Gauss sum}) we have
\[
\prod_{1 \leq i \leq l \atop 1 \leq j \leq k} G(\omega_{\itPi_j\delta_{l+j}})^{n_i} = \prod_{1 \leq i \leq l} G(\omega_{\itPi_j\delta_{l+j}})^{n} \sim G(\omega_\itPi)^n.
\]
This completes the proof.
\end{proof}

\subsection{Generators for relative Lie algebra cohomology}\label{SS:generator}

%Let $\itPi$ be an irreducible admissible essentially tempered $(\frak{g}_n,{\rm O}_n(\R))$-module. 
%Assume $\itPi$ is cohomological with coefficients in $M_{\mu,\C}$.
Let $\pi \in \Omega_\mu$ be a cohomological irreducible admissible $(\frak{g}_n,{\rm O}_n(\R))$-module.
Let $\mathcal{W}(\pi)$ be the space of Whittaker functions of $\pi$ with respect to $\psi_{U_n,\infty}$.
Let $\varepsilon \in \{\pm1\}$ if $n$ is even, and $\varepsilon = \varepsilon(\pi)$ if $n$ is odd.
The aim of this section is to fix a canonical choice of generator of $H^{b_n}(\frak{g}_n,K_n^\circ;\mathcal{W}(\pi)\otimes M_{\mu,\C})[\varepsilon]$ in Definition \ref{D:generator}.

We begin with the case $n=1$. Then $\pi = \chi = {\rm sgn}^\delta|\mbox{ }|^{{\sf w}/2}$ for some $\delta \in \{0,1\}$ and ${\sf w}\in2\Z$. In this case $\varepsilon = (-1)^{\delta+{\sf w}/2}$ and $M_\mu = \Q$ with 
$\GL_1(\Q)=\Q^\times$ acts by $(-\tfrac{{\sf w}}{2})$-th power multiplication. Note that $\mathcal{W}(\chi)$ consisting of functions $W: \R^\times \rightarrow \C$ such that $W(a) = \chi(a)W(1)$. Let $W_\chi \in \mathcal{W}(\chi)$ be the Whittaker function normalized so that $W_\chi(1)=1$.
Then we define the class
\[
[\chi]^\varepsilon \in H^0(\frak{g}_1,K_1^\circ;\mathcal{W}(\chi)\otimes M_{\mu,\C})[\varepsilon] = \left(\mathcal{W}(\chi) \otimes M_{\mu,\C}\right)^{K_1^\circ}
\]
by
\begin{align}\label{E:GL_1 generator}
[\chi]^\varepsilon = W_\chi \otimes 1.
\end{align}
%Here we identify $\bigwedge^0(\frak{g}_{1,\C}/\frak{k}_{1,\C})^*$ with $\C$.

Secondly, we consider the case $n=2$. Then $\pi = D_\kappa \otimes |\mbox{ }|^{{\sf w}/2}$ for some $\kappa \geq 2$ and ${\sf w}\in\Z$ such that $\kappa \equiv {\sf w}\,({\rm mod}\,2)$. In this case, $M_\mu$ can be realized as the $\Q$-vector space consisting of homogeneous polynomials over $\Q$ of degree $\kappa-2$ in variables $x,y$, and $\GL_2(\Q)$ acts on it by
\[
(g\cdot P)(x,y) = (\det g)^{(-\kappa+2-{\sf w})/2}\cdot P((x,y)g).
\]
Let $Y_\pm \in \frak{g}_{2,\C} / \frak{k}_{2,\C}$ defined by
\[
Y_\pm = \frac{1}{2}\bp -1 & 0\\ 0& 1\ep \otimes \sqrt{-1} \pm \frac{1}{2}\bp 0 & 1 \\ 1 & 0 \ep \otimes 1,
\]
and $\{Y_+^*,Y_-^*\} \subset (\frak{g}_{2,\C} / \frak{k}_{2,\C})^*$ be the corresponding dual basis.
Let $W_{\pi}^\pm \in \mathcal{W}(\pi)$ be the Whittaker function of weight $\pm \kappa$ under the action of ${\rm SO}_2(\R)$ and normalized so that
\[
W_\pi^\pm ({\rm diag}(a,1)) = |a|^{(\kappa+{\sf w})/2}e^{-2\pi|a|}\cdot\mathbb{I}_{\R_+}(\pm a).
\]
We define the class 
\[
[\pi]^\varepsilon\in  H^1(\frak{g}_2,K_2^\circ;\mathcal{W}(\pi) \otimes M_{\mu,\C})[\varepsilon] = \left(\mathcal{W}(\pi)\otimes (\frak{g}_{2}/\frak{k}_{2})^* \otimes M_{\mu,\C}\right)^{K_2^\circ}[\varepsilon]
\]
by
\begin{align}\label{E:GL_2 generator}
[\pi]^\varepsilon &= W_\pi^+ \otimes Y_+^*\otimes (\sqrt{-1}\,x+y)^{\kappa-2}+\varepsilon\cdot  (\sqrt{-1})^{{\sf w}}\cdot W_{\pi}^- \otimes Y_-^*\otimes (x+\sqrt{-1}\,y)^{\kappa-2}.
\end{align}
%Note that the normalization here differs from the one in \cite[(3.6) and (3.7)]{RT2011} by a rational number.

Now we consider general $n$. Recall we have the result of Borel and Wallach \cite[III, Theorem 3.3]{BW2000} (also known as Delorme's lemma) %in \cite[Theorem 9.2.1]{Harderbook})
 which expresses the relative Lie algebra cohomology of an induced representation in terms of the corresponding inducing datum. The result was proved by consider a particular spectral sequence (cf.\,\cite[I, Theorem 6.5]{BW2000}). For our purpose, in the following proposition we establish an explicit version of the theorem which describes the isomorphism in bottom degree explicitly. 
The construction will be crucial to our proof of Proposition \ref{L:cohomology} below.
\begin{prop}\label{L:Delorme}
Assume $\pi \cong {\rm Ind}_{P(\R)}^{\GL_n(\R)}(\pi_1\otimes\cdots\otimes\pi_k)$
for some irreducible admissible essentially unitary $(\frak{g}_{n_i},{\rm O}_{n_i}(\R))$-module $(\pi_i, V_{\pi_i})$ for $1 \leq i \leq k$ satisfying conditions (1)-(3) in Lemma \ref{L:combinatorial}, where $P = P_{(n_1,...,n_k)}$.
Fix an $M_P(\C)$-equivariant embedding
\[
\otimes_{i=1}^k M_{\mu^{(i)},\C} \longrightarrow {\rm Hom}\left(\extp^{b_n-\sum_{i=1}^kb_{n_i}}{\rm Lie}(U_P(\Q))_\C,\,M_{\mu,\C}\right)
\] 
defined over $\Q$.
%Then there exists a unique Weyl element $w \in W_n^P$ such that $\ell(w) = b_n-\sum_{i=1}^kb_{n_i}$ and $\pi \in \Omega_\mu$ with $\mu = w^{-1}\cdot\left( (\mu^{(1)},\cdots,\mu^{(k)})+\rho_n \right)-\rho_n$.
Then there exist canonical $\pi_0(\GL_n(\R))$-equivariant isomorphism $I_{\otimes_{i=1}^kV_{\pi_i}}$ and injective homomorphism $\iota_{\otimes_{i=1}^kV_{\pi_i}}$ such that the following diagram is commutative
\[
\begin{tikzcd}[row sep=normal, column sep=normal]
{\rm Hom}_{K_n^\circ}\left(\extp^{b_n}(\frak{g}_n/\frak{k}_n),\,{\rm Ind}_{P(\R)}^{\GL_n(\R)}(\otimes_{i=1}^k V_{\pi_i})\otimes M_{\mu,\C}\right)^{d=0} \arrow[twoheadrightarrow]{r} & H^{b_n}(\frak{g}_n,K_n^\circ;{\rm Ind}_{P(\R)}^{\GL_n(\R)}(\otimes_{i=1}^k V_{\pi_i})\otimes M_{\mu,\C})\arrow[ld, "I_{\otimes_{i=1}^kV_{\pi_i}}"]\\
\left( \bigotimes_{i=1}^k H^{b_{n_i}}(\frak{g}_{n_i},K_{n_i}^\circ;V_{\pi_i\otimes \delta_i} \otimes M_{\mu^{(i)},\C})\right)^{\pi_0(K_{n}^{M_P})}\arrow[hookrightarrow]{u}{\iota_{\otimes_{i=1}^kV_{\pi_i}}}&
\end{tikzcd}
\]
%\[
%\begin{tikzcd}[row sep=normal, column sep=normal]
%H^{b_n}(\frak{g}_n,K_n^\circ;{\rm Ind}_{P(\R)}^{\GL_n(\R)}(\otimes_{i=1}^k V_{\pi_i})\otimes M_{\mu,\C}) \arrow[r, "I_{\otimes_{i=1}^kV_{\pi_i}}"] & \left( \bigotimes_{i=1}^k H^{b_{n_i}}(\frak{g}_{n_i},K_{n_i}^\circ;V_{\pi_i\otimes \delta_i} \otimes M_{\mu^{(i)},\C})\right)^{\pi_0(K_{n}^{M_P})}\arrow[hookrightarrow]{dl}{\iota_{\otimes_{i=1}^kV_{\pi_i}}}\\
%{\rm Hom}_{K_n^\circ}\left(\extp^{b_n}(\frak{g}_n/\frak{k}_n),\,{\rm Ind}_{P(\R)}^{\GL_n(\R)}(\otimes_{i=1}^k V_{\pi_i})\otimes M_{\mu,\C}\right)^{d=0} \arrow[twoheadrightarrow]{u}
%\end{tikzcd}
%\]
%\begin{align}\label{E:BW}
%\begin{split}
%I_{\otimes_{i=1}^kV_{\pi_i}}& :  H^{b_n}(\frak{g}_n,K_n^\circ;{\rm Ind}_{P(\R)}^{\GL_n(\R)}(\otimes_{i=1}^k V_{\pi_i})\otimes M_{\mu,\C}) \\
%&\quad\quad\quad\quad\quad\quad\quad\quad\quad\quad\quad\quad\longrightarrow \left( \bigotimes_{i=1}^k H^{b_{n_i}}(\frak{g}_{n_i},K_{n_i}^\circ;V_{\pi_i\otimes \delta_i} \otimes M_{\mu^{(i)},\C})\right)^{\pi_0(K_{n}^{M_P})},
%\end{split}
%\end{align}
%\begin{align}\label{E:BW}
%\begin{split}
%I_{\otimes_{i=1}^k\mathcal{W}(\itPi_i)}& :  H^{b_n}(\frak{g}_n,K_n^\circ;{\rm Ind}_{P(\R)}^{\GL_n(\R)}(\otimes_{i=1}^k \mathcal{W}(\itPi_i))_{{\rm O}_n(\R)} \otimes M_{\mu,\C}) \\
%&\quad\quad\quad\quad\quad\quad\quad\quad\longrightarrow {\rm Ind}_{\pi_0(P(\R))}^{\pi_0(\GL_n(\R))}\left[ \left( \bigotimes_{i=1}^k H^{b_{n_i}}(\frak{g}_{n_i},K_{n_i}^\circ;\mathcal{W}(\itPi_i\rho_i) \otimes M_{\mu_i,\C})\right)^{\pi_0(K_{n}^{M_P})}\right].
%\end{split}
%\end{align}
Here $d$ is the differential for the $(\frak{g}_n,K_n^\circ)$-cohomology, $K_n^{M_P} = K_n^\circ \cap M_P(\R)$, and $V_{\pi_i\otimes\delta_i} = V_{\pi_i}$ is the underlying space of $\pi_i\otimes\delta_i$. % is the image of $K_n^\circ \cap P(\R)$ under the canonical map $P \rightarrow M_P$.
\end{prop}

\begin{proof}
We denote by ${}^a{\rm Ind}$ the algebraic induction (without the normalizing factor $\delta_P$) and let $V=\otimes_{i=1}^kV_{\pi_i\otimes\delta_i}$.
Let $K = {\rm SO}_{n_1}(\R)\times \cdots\times {\rm SO}_{n_k}(\R) \subset M_P(\R)$, and $\frak{m}_P$ and $\frak{u}_P$ be the complexified Lie algebras of $M_P(\R)$ and $U_P(\R)$ respectively.
Consider the $\pi_0(\GL_n(\R))$-equivariant homomorphism
\begin{align}\label{E:Delorme proof 1}
\begin{split}
{\rm Hom}_{K_n^\circ}\left(\extp^\bullet(\frak{g}_n/\frak{k}_n),\,{}^a{\rm Ind}_{P(\R)}^{\GL_n(\R)}(V)\otimes M_{\mu,\C}\right) &\longrightarrow {\rm Hom}_{K}\left( \extp^\bullet(\frak{g}_n/\frak{k}_n),\,V\otimes M_{\mu,\C}\right)^{\pi_0(K_n^{M_P})}\\
f & \longmapsto \left( X \mapsto f(X)(1)\right).
\end{split}
\end{align}
It is easy to check that this is an isomorphism. 
We will define $\iota_{\otimes_{i=1}^kV_{\pi_i}}$ to be the restriction of the inverse of this isomorphism to a particular summand of the right-hand side.
%We show that it defines the isomorphism (\ref{E:BW}) for $q=b_n$.
The natural inclusion $\frak{m}_P\oplus\frak{u}_P \subset \frak{g}_n$ induces an isomorphism $\left(\frak{m}_P / \frak{m}_P\cap\frak{k}_n\right)\oplus\frak{u}_P\cong\frak{g}_n/\frak{k}_n$. This defines a canonical isomorphism between the right-hand side of (\ref{E:Delorme proof 1}) and
\begin{align*}
\bigoplus_{a+b = \bullet}C^{a,b}(V\otimes M_{\mu,\C}),
\end{align*}
where the bigraded complex is defined by
\[
C^{a,b}(V\otimes M_{\mu,\C}) = {\rm Hom}_K\left( \extp^a(\frak{m}_P / \frak{m}_P\cap\frak{k}_n),\,V\otimes{\rm Hom}\left(\extp^b\frak{u}_P,M_{\mu,\C}\right)\right)^{\pi_0(K_n^{M_P})}.
\]
Let $d$ be the differential on the left-hand side of (\ref{E:Delorme proof 1}) for the $(\frak{g}_n,K_n^\circ)$-cohomology, and $d_1$ and $d_2$ be the differentials on $C^{a,b}$ for the $(\frak{m}_P,K_n^{M_P})$-cohomology and $\frak{u}_P$-cohomology respectively. Via the isomorphism (\ref{E:Delorme proof 1}), on $C^{a,b}$ we have
\begin{align}\label{E:Delorme proof 2}
d=(d_1\otimes 1) + (-1)^a\cdot (1 \otimes d_2).
\end{align}
Let $w \in W_n^P$ be the unique Kostant representative in Lemma \ref{L:combinatorial}. 
The $M_P(\C)$-representation $\otimes_{i=1}^k M_{\mu^{(i)},\C}$ has highest weight $w(\mu+\rho_n)-\rho_n$.
By the result of Kostant \cite[Corollary 5.7 and Theorem 5.14]{Kostant1961}, $\otimes_{i=1}^k M_{\mu^{(i)},\C}$ appears with multiplicity one in ${\rm Hom}\left(\extp^{\ell(w)}\frak{u}_P,M_{\mu,\C}\right)$. Moreover, its image belongs to the kernel of the differential $d_2$. We fix an embedding defined over $\Q$, where the $\Q$-rational structure of $\frak{u}_P$ is the Lie algebra of $U_P(\Q)$.
Let $\frak{a}$ and $\frak{a}_P$ be the complexified Lie algebras of $\R^\times$ and the center of $M_P(\R)$ respectively. Then we have a canonical isomorphism $\frak{m}_P / \frak{m}_P\cap\frak{k}_n \cong \frak{a}_P/\frak{a} \,\oplus \left(\frak{g}_{n_1}/\frak{k}_{n_1}\oplus\cdots\oplus\frak{g}_{n_1}/\frak{k}_{n_k}\right)$. The fixed embedding thus determines a summand of $C^{a,b}$ for $a = \sum_{i=1}^kb_{n_i}$ and $b =\ell(w)$ (thus $\bullet=b_n$):
\begin{align}\label{E:Delorme proof 3}
\left(\bigotimes_{i=1}^k {\rm Hom}_{K_{n_i}^\circ}\left(\extp^{b_{n_i}}(\frak{g}_{n_i}/\frak{k}_{n_i}),\,V_{\pi_i\otimes\delta_i}\otimes M_{\mu^{(i)},\C}\right) \right)^{\pi_0(K_n^{M_P})}.
\end{align}
Here we have replaced ${\rm SO}_{n_i}(\R)$ by $K_{n_i}^\circ$ since $\R_+$ acts trivially on $(\pi_i\otimes\delta_i)\otimes M_{\mu^{(i)},\C}$.
Since $\pi_i$ is essentially unitary and the infinitesimal characters of ${\pi_i\otimes\delta_i}$ and $M_{\mu^{(i)},\C}^\vee$ are equal, by \cite[II, Proposition 3.1-(b)]{BW2000} the complex defining the $(\frak{g}_{n_i},K_{n_i}^\circ)$-cohomology of $({\pi_i\otimes\delta_i})\otimes M_{\mu^{(i)},\C}$ is closed in all degrees. In particular, the summand $(\ref{E:Delorme proof 3})$ of $C^{a,b}$ is in the kernel of the differential $d_1$ and is equal to the target of $I_{\otimes_{i=1}^kV_{\pi_i}}$.
Therefore, by (\ref{E:Delorme proof 2}), the preimage of the summand under the isomorphism (\ref{E:Delorme proof 1}) belongs to the kernel of $d$ for $\bullet=b_n$. Since the dimensions of both sides of $I_{\otimes_{i=1}^kV_{\pi_i}}$ are equal, it remains to show that the intersection of the summand (\ref{E:Delorme proof 3}) and the coboundary is zero.
Suppose $f$ is in the intersection, then we have
\[
f= (d_1\otimes 1)f_1 + (1\otimes d_2)f_2
\]
for some $f_1 \in C^{a-1,b}$ and $f_2 \in C^{a,b-1}$. By \cite[(5.7.4)]{Kostant1961}, there exists a $M_P(\C)$-subrepresentation $M'$ of ${\rm Hom}\left(\extp^b\frak{u}_P,M_{\mu,\C}\right)$ containing the image of $d_2$ such that 
\[
{\rm Hom}\left(\extp^b\frak{u}_P,M_{\mu,\C}\right) = \otimes_{i=1}^k M_{\mu^{(i)},\C}\bigoplus M'.
\]
Then $f_1 = h_1 + h_1'$ for some $h_1$ and $h_1'$ in $C^{a-1,b}$ valued in $V\otimes (\otimes_{i=1}^k M_{\mu^{(i)},\C})$ and $V \otimes M'$ respectively.
Since $f$ and $(1\otimes d_2)f_2$ are valued in $V\otimes (\otimes_{i=1}^k M_{\mu^{(i)},\C})$ and $V\otimes M'$ respectively, we must have $(d_1\otimes 1)h_1'+(1\otimes d_2)f_2=0$. By the same reasoning as above using \cite[II, Proposition 3.1-(b)]{BW2000}, we have $(d_1\otimes 1)h_1=0$. Therefore, $f=0$. This completes the proof. 
\end{proof}
%Note that the tensor product on the right-hand side is naturally a $\pi_0(M_P(\R))$-module. After taking the $\pi_0(K_n^{M_P})$-invariant part, it is a $\pi_0(\GL_n(\R))$-module by the short exact sequence
%\[
%1 \longrightarrow \pi_0(K_n^{M_P}) \longrightarrow \pi_0(M_P(\R)) \longrightarrow \pi_0(\GL_n(\R)) \longrightarrow 1.
%\]
\begin{rmk}\label{R:Delorme}
Since $I_{\otimes_{i=1}^kV_{\pi_i}}$ is $\pi_0(\GL_n(\R))$-equivariant, by taking the $\varepsilon$-isotypic component on both sides of $I_{\otimes_{i=1}^kV_{\pi_i}}$, we obtain the isomorphism 
\begin{align}\label{E:Delorme}
\begin{split}
I_{\otimes_{i=1}^kV_{\pi_i}}^\varepsilon &: H^{b_n}(\frak{g}_n,K_n^\circ;{\rm Ind}_{P(\R)}^{\GL_n(\R)}(\otimes_{i=1}^k V_{\pi_i}) \otimes M_{\mu,\C})[\varepsilon]\\
& \quad\quad\quad\quad\quad\quad\quad\quad\quad\quad\quad\quad\longrightarrow  \bigotimes_{i=1}^k H^{b_{n_i}}(\frak{g}_{n_i},K_{n_i}^\circ;V_{\pi_i\otimes \delta_i} \otimes M_{\mu^{(i)},\C})[\varepsilon].
\end{split}
\end{align}
\end{rmk}

Keep the notation and assumptions in Proposition \ref{L:Delorme}. We take the underlying space $V_{\pi_i} = \mathcal{W}(\pi_i)$ to be the Whittaker model of $\pi_i$ for $1 \leq i \leq k$.
We have the equivariant isomorphism 
\[
{\mathbb W} : {\rm Ind}_{P(\R)}^{\GL_n(\R)}\left(\otimes_{i=1}^k\mathcal{W}(\pi_i)\right) \longrightarrow \mathcal{W}(\pi)
\]
defined by Jacquet's integral \cite{Jacquet1967} (cf.\,(\ref{E:Jacquet integral}) below).
This in turn defines an $\pi_0(\GL_n(\R))$-equivariant isomorphism
\begin{align*}
{\mathbb W}^{b_n} : H^{b_n}(\frak{g}_n,K_n^\circ;{\rm Ind}_{P(\R)}^{\GL_n(\R)}(\otimes_{i=1}^k \mathcal{W}(\pi_i)) \otimes M_{\mu,\C}) \longrightarrow H^{b_n}(\frak{g}_n,K_n^\circ;\mathcal{W}(\pi) \otimes M_{\mu,\C}).
\end{align*}

\begin{definition}\label{D:generator}
For $\pi \in \Omega_\mu \subset \Omega(n)$, let
\[
[\pi]^\varepsilon \in H^{b_n}(\frak{g}_n,K_n^\circ;\mathcal{W}(\pi) \otimes M_{\mu,\C})[\varepsilon]
\]
be the cohomology class defined by 
\[
[\pi]^\varepsilon = {\mathbb W}^{b_n}\circ \left(I_{\otimes_{i=1}^k\mathcal{W}(\pi_i)}^\varepsilon\right)^{-1}\left( \otimes_{i=1}^k[\pi_i\otimes\delta_i]^\varepsilon\right)
\]
for some specific choices of $\pi_i$ and $[\pi_i\otimes\delta_i]^\varepsilon$ described as follows.
Put $r = {\lfloor \tfrac{n}{2} \rfloor}$.
Let $k=r$ (resp.\,$k=r+1$) if $n$ is even (resp.\,odd).
Let $(\underline{\kappa};\,\underline{\sf w})$ be the infinity type of $\pi$.
If $n$ is odd, let $\delta \in \{0,1\}$ such that $\varepsilon(\pi) = (-1)^{r+\delta+{\sf w}_{r+1}/2}$.
Then $\pi_1,...,\pi_k$ are the representations uniquely determined by 
\[
\left\{\pi_1,...,\pi_k\right\} = \{D_{\kappa_1}\otimes |\mbox{ }|^{{\sf w}_1/2},...,D_{\kappa_r}\otimes |\mbox{ }|^{{\sf w}_r/2}\} \cup \begin{cases}
\emptyset & \mbox{ if $n$ is even},\\
\{{\rm sgn}^\delta|\mbox{ }|^{{\sf w}_{r+1}/2}\} & \mbox{ if $n$ is odd},
\end{cases}
\]
and the following conditions for $1 \leq i < j \leq k$:
\begin{itemize}
\item If $\pi_i = D_{\kappa_\alpha}\otimes|\mbox{ }|^{{\sf w}_\alpha/2}$ and $\pi_j = D_{\kappa_\beta}\otimes|\mbox{ }|^{{\sf w}_\beta/2}$, then either ${\sf w}_\alpha>{\sf w}_\beta$, or ${\sf w}_\alpha = {\sf w}_\beta$ and $\alpha<\beta$.
\item If $\pi_i = D_{\kappa_\alpha}\otimes|\mbox{ }|^{{\sf w}_\alpha/2}$ and $\pi_j = {\rm sgn}^\delta|\mbox{ }|^{{\sf w}_{r+1}/2}$, then ${\sf w}_\alpha \geq {\sf w}_{r+1}$.
\item If $\pi_i = {\rm sgn}^\delta|\mbox{ }|^{{\sf w}_{r+1}/2}$ and $\pi_j = D_{\kappa_\alpha}\otimes|\mbox{ }|^{{\sf w}_\alpha/2}$, then ${\sf w}_{r+1}>{\sf w}_\alpha$.
\end{itemize}
%Let 
%\begin{align*}
%P=
%\begin{cases}
%P_{(2,\cdots,2)} & \mbox{ if $n$ is even},\\
%P=P_{(2,\cdots,2,1)} & \mbox{ if $n$ is odd},
%\end{cases}\quad 
%k= \begin{cases}
%r & \mbox{ if $n$ is even},\\
%r+1 & \mbox{ if $n$ is odd}.
%\end{cases}
%\end{align*}
%Define $\itPi_i = D_{\kappa_i}\otimes |\mbox{ }|^{{\sf w}_i/2}$ for $1 \leq i \leq r$. %When $n$ is odd, let $\itPi_{r+1}$ be the algebraic character of $\R^\times$ such that ${\itPi_{r+1}}(-1)=(-1)^{r+{\sf w}_{r+1}/2}\varepsilon(\itPi)$ and $|{\itPi_{r+1}}| = |\mbox{ }|^{{\sf w}_{r+1}/2}$.
The generators $[\pi_1\otimes\delta_1]^\varepsilon,...,[\pi_k\otimes\delta_k]^\varepsilon$ are the ones fixed in (\ref{E:GL_1 generator}) and (\ref{E:GL_2 generator}) for $n=1,2$.
\end{definition}

The following lemma is one of the crucial ingredients in the proof of Theorem \ref{T:Raghuram}.
We show that, up to multiplication by $\Q^\times$, the generator in Definition \ref{D:generator} is independent of the choice of $\pi_1,...,\pi_k$ subject to the ordering with $e(\pi_1)\geq \cdots \geq e(\pi_k)$ and conditions (1)-(3) in Lemma \ref{L:combinatorial}.

\begin{lemma}\label{L:GL_n generator}
Assume $\pi \cong {\rm Ind}_{P(\R)}^{\GL_n(\R)}(\otimes_{i=1}^k \pi_i)$ with $\otimes_{i=1}^k \pi_i$ satisfying conditions (1)-(3) in Lemma \ref{L:combinatorial} and $\pi_1,...,\pi_{k-1}$ are essentially tempered with $e(\pi_1)\geq \cdots \geq e(\pi_k)$, where $e(\pi_i)$ is the unique real number such that $\pi_i\otimes |\mbox{ }|^{-e(\pi_i)}$ is tempered.
Then we have
\[
{\mathbb W}^{b_n}\circ \left(I^\varepsilon_{\otimes_{i=1}^k\mathcal{W}(\pi_i)}\right)^{-1}\left( \otimes_{i=1}^k[\pi_i\otimes\delta_i]^\varepsilon\right) = C\cdot [\pi]^\varepsilon
\]
for some constant $C\in \Q^\times$.
\end{lemma}

\begin{proof}
Note that the isomorphisms $I_{\otimes_{i=1}^k\mathcal{W}(\pi_i)}$ and ${\mathbb W}^{b_n}$ commute with parabolic induction in stages.
By the inductive definition of $[\pi]^\varepsilon,[\pi_1\otimes\delta_1]^\varepsilon,...,[\pi_k\otimes\delta_k]^\varepsilon$ and the assumption that $\pi_1,...,\pi_{k-1}$ are essentially tempered with $e(\pi_1)\geq \cdots \geq e(\pi_k)$, to prove the the lemma, we are reduced to prove the assertion for $k=2$ with $e(\pi_1)=e(\pi_2)$.
%: Let $\itSigma$ be an irreducible admissible $(\frak{g}_{N},{\rm O}_N(\R))$-module such that
%\[
%\itSigma \cong {\rm Ind}_{P_{(n_1,n_2)}(\R)}^{\GL_N(\R)}(\itSigma_1 \otimes \itSigma_2)
%\] 
%for some irreducible admissible $(\frak{g}_{n_i},{\rm O}_{n_i}(\R))$-module $\itSigma_i$ for $i=1,2$ such that
We recall the assumptions as follows:
\begin{itemize}
\item $n_1n_2$ is even,
\item $\pi_1\otimes|\mbox{ }|^{n_2/2} \in \Omega(n_1)$ and $\pi_2\otimes|\mbox{ }|^{n_1/2} \in \Omega(n_2)$,
\item $L(s,\pi_1 \times \pi_2^\vee)$ and $L(s,\pi_1^\vee \times \pi_2)$ are holomorphic at $s=0$,
\item $\pi_1$ and $\pi_2$ are essentially tempered with $e(\pi_1) = e(\pi_2)$.
\end{itemize}
%$n_1+n_2=N$ such that $n_1n_2$ is even, $L(s,\itSigma_1 \times \itSigma_2^\vee)$ and $L(s,\itSigma_1^\vee \times \itSigma_2)$ are holomorphic at $s=0$, and $e(\itSigma_1) = e(\itSigma_2)$. 
We need to show that
\begin{align}\label{E:GL_n generator proof 1}
\begin{split}
&{\mathbb W}^{b_n}\circ \left(I^{\varepsilon}_{\mathcal{W}(\pi_1)\otimes \mathcal{W}(\pi_2)}\right)^{-1}\left([\pi_1\otimes|\mbox{ }|^{n_2/2}]^{\varepsilon} \otimes [\pi_2\otimes|\mbox{ }|^{-n_1/2}]^{\varepsilon}\right)\\
& = C\cdot {\mathbb W}^{b_n}\circ \left(I^{\varepsilon}_{\mathcal{W}(\pi_2)\otimes \mathcal{W}(\pi_1)}\right)^{-1}\left([\pi_2\otimes|\mbox{ }|^{n_1/2}]^{\varepsilon} \otimes [\pi_1\otimes|\mbox{ }|^{-n_2/2}]^{\varepsilon}\right)
\end{split}
\end{align}
for some $C \in \Q^\times$. %Here $\varepsilon' \in \{\pm1\}$ if $N$ is even and $\varepsilon' = \varepsilon(\itSigma)$ if $N$ is odd. (It is clear that we can even assume $N \leq 4$.)
Write $P=P_{(n_1,n_2)}$ and $Q=P_{(n_2,n_1)}$.
To prove (\ref{E:GL_n generator proof 1}), we begin with the intertwining operator 
\[
M^Q : {\rm Ind}_{P(\R)}^{\GL_n(\R)}(\mathcal{W}(\pi_1) \otimes \mathcal{
W}(\pi_2)) \longrightarrow {\rm Ind}_{Q(\R)}^{\GL_n(\R)}(\mathcal{W}(\pi_2) \otimes \mathcal{
W}(\pi_1)),\quad f \longmapsto M^Q(f)
\]
so that $M^Q(f)(g)$ is the evaluation at $\underline{s}=0$ of the meromorphic function defined by the integral
\[
\int_{U_Q(\R)}f_{\underline{s}}\left(\bp  & {\bf 1}_{n_1} \\ {\bf 1}_{n_2} & \ep ug\right)\,du
\]
when ${\rm Re}(s_1-s_2)$ is sufficiently large, 
where $du$ is the product measure of the Lebesgue measures on $\R$.
Let $(M^Q)^{b_n}$ be the corresponding isomorphism for the relative Lie algebra cohomology group in bottom degree $b_n$.
Consider the following diagram of equivariant isomorphisms:
\[
\begin{tikzcd}[row sep=normal, column sep=huge]
{\rm Ind}_{P(\R)}^{\GL_n(\R)}(\mathcal{W}(\pi_1) \otimes \mathcal{
W}(\pi_2)) \arrow[r, "{\mathbb W}"] \arrow[d, "M^Q"]& \mathcal{W}(\pi)\\
{\rm Ind}_{Q(\R)}^{\GL_n(\R)}(\mathcal{W}(\pi_2) \otimes \mathcal{
W}(\pi_1))\arrow[ru, "{\mathbb W}"']&
\end{tikzcd}
\]
As a special case of the result of Shahidi on local coefficients in \cite[Theorem 3.1]{Shahidi1985}, we have
\[
{\mathbb W} = \epsilon\cdot \gamma(0,\pi_1 \times \pi_2^\vee,\psi_\infty) \cdot {\mathbb W}\circ M^Q
\]
for some $\epsilon \in \{\pm1\}$, where $\gamma(s,\pi_1 \times \pi_2^\vee,\psi_\infty)$ is the $\gamma$-factor of $\pi_1 \times \pi_2^\vee$ with respect to $\psi_\infty$.
We refer to \cite[Remark 5.1.3]{Shahidi2010} for the appearance of $\epsilon$.
This induces an isomorphism at the level of relative Lie algebra cohomology groups in degree $b_n$:
\[
{\mathbb W}^{b_n} = \epsilon\cdot\gamma(0,\pi_1 \times \pi_2^\vee,\psi_\infty) \cdot {\mathbb W}^{b_n}\circ (M^Q)^{b_n}.
\]
Therefore, (\ref{E:GL_n generator proof 1}) holds if and only if
\begin{align*}
&\gamma(0,\pi_1 \times \pi_2^\vee,\psi_\infty) \cdot (M^Q)^{b_n}\circ \left(I^{\varepsilon}_{\mathcal{W}(\pi_1)\otimes \mathcal{W}(\pi_2)}\right)^{-1}\left([\pi_1\otimes|\mbox{ }|^{n_2/2}]^{\varepsilon} \otimes [\pi_2\otimes|\mbox{ }|^{-n_1/2}]^{\varepsilon}\right)\\
& = C\cdot \left(I^{\varepsilon}_{\mathcal{W}(\pi_2)\otimes \mathcal{W}(\pi_1)}\right)^{-1}\left([\pi_2\otimes|\mbox{ }|^{n_1/2}]^{\varepsilon} \otimes [\pi_1\otimes|\mbox{ }|^{-n_2/2}]^{\varepsilon}\right)
\end{align*}
for some $C \in \Q^\times$.
The above equality was proved independently by Harder \cite[Theorem 8.7]{HR2020} and Weselmann \cite[\S\,9.6.12]{HR2020} with $\gamma(0,\pi_1 \times \pi_2^\vee,\psi_\infty)$ replaced by $(2\pi\sqrt{-1})^{-n_1n_2/2}$.
Finally, an immediate computation, together with the condition that $e(\pi_1) = e(\pi_2)$, shows that 
\[
\gamma(0,\pi_1 \times \pi_2^\vee,\psi_\infty) \in (2\pi\sqrt{-1})^{-n_1n_2/2}\cdot \Q^\times.
\]
Indeed, this is (\ref{E:well-defined pf 4}) and please refer to the computation therein.
This completes the proof.
\end{proof}

\section{Eisenstein cohomology classes}\label{S:Eisenstein}
%\cite[\S\,II.1.7 and \S\,IV.1.11]{MW1995}
\subsection{Eisenstein series}\label{SS:Eisenstein}

Let $P = P_{(n_1,...,n_k)}$ be a standard parabolic subgroup of $\GL_n$ of type $(n_1,...,n_k)$.
Let $\itPi_1 \otimes \cdots\otimes \itPi_k$ be a cuspidal automorphic representation of $M_P(\A)\cong \GL_{n_1}(\A)\times\cdots\times\GL_{n_k}(\A)$.
%Let $\itPi$ be an isobaric automorphic representation of $\GL_n(\A)$. We have
%\[
%\itPi = \itPi_1 \boxplus \cdots \boxplus \itPi_k
%\]
%for some cuspidal automorphic representation $\itPi_i$ of $\GL_{n_i}(\A)$ for $1 \leq i \leq k $. 
%Following \cite{LLS2021}, we say $\itPi$ is $tamely$ $isobaric$ if there exists $s \in \R$ such that $\itPi_i \otimes |\mbox{ }|_\A^s$ is unitary for $1 \leq i \leq k$. In this case, $\itPi$ is fully induced (cf.\,\cite{Bernstein1984} and \cite{Vogan1986}).
%We recall the realization of $\itPi$ in the space $\mathcal{A}(\GL_n)$ of automorphic forms on $\GL_n(\A)$ using cuspidal Eisenstein series as follows: Let $P=P_{(n_1,\cdots,n_k)}$ be the standard parabolic subgroup of $\GL_n$ of type $(n_1,\cdots,n_k)$.
Let 
\[
I_P^{\GL_n}(\otimes_{i=1}^k\itPi_i)
\]
be the space of smooth and ${\rm O}_n(\R)$-finite functions $h : U_P(\A)M_P(\Q)\backslash\GL_n(\A) \rightarrow \C$ such that for all $k \in {\rm O}_n(\R)\times\GL_n(\widehat{\Z})$, the function
\[
M_P(\A)\longrightarrow \C,\quad m\longmapsto \delta_P(m)^{-1} \cdot h(mk)
\]
is a cusp form in $\itPi_1 \otimes \cdots\otimes \itPi_k$,
where $\delta_P$ is the square-root of the modulus character of $P(\A)$. %, and $V_{\itPi_i}$ is the space of $\itPi_i$ realized in the space of cusp forms on $\GL_{n_i}(\A)$ for $1 \leq i \leq k$.
It is clear that $I_P^{\GL_n}(\otimes_{i=1}^k\itPi_i)$ can be identified with %the space of ${\rm O}_n(\R)$-finite functions in 
the induced representation ${\rm Ind}_{P(\A)}^{\GL_n(\A)}(\otimes_{i=1}^k\itPi_i)$. %, which we denote by ${\rm Ind}_{P(\A)}^{\GL_n(\A)}(\otimes_{i=1}^k\itPi_i)_{{\rm O}_n(\R)}$.
For $h \in I_P^{\GL_n}(\otimes_{i=1}^k\itPi_i)$ and $\underline{s} = (s_1,...,s_k) \in \C^k$, let $h_{\underline{s}} \in I_P^{\GL_n}(\otimes_{i=1}^k(\itPi_i\otimes|\mbox{ }|_\A^{s_i}))$  defined by
\begin{align}\label{E:flat section}
h_{\underline{s}}(g) = \prod_{i=1}^k|\det m_i|_\A^{s_i}\cdot h(g)
\end{align}
for $g=umk \in U_P(\A)M_P(\A)({\rm O}_n(\R)\times\GL_n(\widehat{\Z}))$ with $m = {\rm diag}(m_1,...,m_k)$.
We define the Eisenstein series 
\[
E(g,h_{\underline{s}}) = \sum_{\gamma \in P(\Q)\backslash\GL_n(\Q)}h_{\underline{s}}(\gamma g),\quad g\in\GL_n(\A)
\]
which converges absolutely when ${\rm Re}(s_i-s_j)$ is sufficiently large for $1 \leq i <j\leq k$. 
By the result of Langlands \cite{Langlands1976}, it admits meromorphic continuation to $\underline{s} \in \C^k$ (cf.\,\cite[IV.1.8]{MW1995}). %Moreover, since $\itPi$ is tamely isobaric, the Eisenstein series is holomorphic at $\lambda=0$ (cf.\,\cite[\S\,IV.1.11]{MW1995}). 
%We then realize $\itPi$ in $\mathcal{A}(\GL_n)$ as the image of the $((\frak{g}_{n},{\rm O}_{n}(\R))\times \GL_n(\A_f))$-equivariant homomorphism 
For $w \in (W_n^P)^{-1}$ with $wM_Pw^{-1}=M_{P'}$ for some standard parabolic subgroup $P'$, there exists a unique permutation $\tau_w$ on $\{1,...,k\}$ such that
\begin{align}\label{E:permuation}
w\cdot{\rm diag}(m_1,...,m_k)\cdot w^{-1} = {\rm diag}(m_{\tau_w(1)},...,m_{\tau_w(k)}),\quad {\rm diag}(m_1,...,m_k) \in M_P.
\end{align}
One can easily verify that the association $w \mapsto \tau_w$ is injective.
Recall we have the intertwining operator%an $((\frak{g}_{n},{\rm O}_{n}(\R))\times \GL_n(\A_f))$-equivariant isomorphism
\[
M_w^{P'}
 : I_P^{\GL_n}(\otimes_{i=1}^k (\itPi_i\otimes |\mbox{ }|_\A^{s_i})) \longrightarrow I_{P'}^{\GL_n}(\otimes_{i=1}^k (\itPi_{\tau_w(i)}\otimes|\mbox{ }|_\A^{s_{\tau_w(i)}}))
\]
defined by the meromorphic continuation of the intertwining integral
\begin{align}\label{E:local intertwining}
\int_{(U_{P'}(\A)\cap wU_P(\A)w^{-1})\backslash U_{P'}(\A)}h_{\underline{s}}(w^{-1}ug)\,du^{\rm Tam}
\end{align}
which converges absolutely when ${\rm Re}(s_i-s_j)$ is sufficiently large for all $1 \leq i < j \leq k$ such that $\tau_w^{-1}(i)>\tau_w^{-1}(j)$ (cf.\,\cite[II.1.6]{MW1995}). 
For an automorphic form $\varphi$ on $\GL_n(\A)$, recall the constant term of $\varphi$ along $P$ is defined by
\[
\varphi_P(g) = \int_{U_P(\Q)\backslash U_P(\A)}\varphi(g)\,du^{\rm Tam},\quad g \in \GL_n(\A).
\] 
Here $du^{\rm Tam}$ is the Tamagawa measure on $U_P(\A)$.
%For $h \in I_P^{\GL_n}(\otimes_{i=1}^k\itPi_i)$ and $w \in W_n$ such that $wM_Pw^{-1}=M_P$, let 
%\[
%M_w^P(h) : U_P(\A)M_P(\Q)\backslash \GL_n(\A) \longrightarrow \C
%\]
%be the function defined so that $M_w^P(h)(g)$ is the evaluation at $\lambda=0$ of the meromorphic function 
%Note that the meromorphic function is holomorphic at $\lambda=0$, since $\itPi$ is tamely isobaric (cf.\,\cite[\S\,IV.1.11]{MW1995}). 
Then the constant term of the Eisenstein series $E(h_{\underline{s}})$ along $P$ is given by (cf.\,\cite[II.1.7]{MW1995})
\begin{align}\label{E:constant term}
E_P(h_{\underline{s}}) = \sum_{ w \in W_n^P\cap (W_n^P)^{-1}\atop wM_Pw^{-1}=M_P} M_w^P(h_{\underline{s}}).
\end{align}
In the following proposition, we give sufficient conditions under which the Eisenstein series $E(h_{\underline{s}})$ is holomorphic at $\underline{s} =0$.

\begin{lemma}\label{P:Eisenstein}
For $1\leq i \leq k$, let $e(\itPi_i)$ be the unique real number such that $\itPi_i \otimes |\mbox{ }|_\A^{-e(\itPi_i)}$ is unitary.
Assume the following conditions are satisfied for all $1 \leq i<j \leq k$:
\begin{itemize}
\item[(1)] $\itPi_i \otimes |\mbox{ }|_\A^{-e(\itPi_i)}$ and $\itPi_j\otimes |\mbox{ }|_\A^{-e(\itPi_j)}$ are inequivalent.
\item[(2)] $e(\itPi_i) \geq e(\itPi_j)$.
%\item[(2)] Either $e(\itPi_i) = e(\itPi_j)$ or $e(\itPi_i)-e(\itPi_j)\geq 1$ holds.
\end{itemize}
Then the Eisenstein series $E(h_{\underline{s}})$ is holomorphic at $\underline{s}=0$ for all $h \in I_P^{\GL_n}(\otimes_{i=1}^k\itPi_i )$.
\end{lemma}

\begin{proof}
%By the formula for the constant terms of Eisenstein series (cf.\,\cite[\S\,II.1.7]{MW1995}), 
It suffices to show that the intertwining operators are holomorphic at $\underline{s}=0$ (cf.\,Remark in \cite[IV.4.4]{MW1995}).
Let $w \in (W_n^P)^{-1}$ such that $wM_Pw^{-1}=M_{P'}$ for some standard parabolic subgroup $P'$.
Fix isomorphisms 
\[
I_P^{\GL_n}(\otimes_{i=1}^k\itPi_i ) \cong \bigotimes_v{\rm Ind}_{P(\Q_v)}^{\GL_n(\Q_v)}(\otimes_{i=1}^k\itPi_{i,v}),\quad I_{P'}^{\GL_n}(\otimes_{i=1}^k \itPi_{\tau_w(i)}) \cong \bigotimes_v{\rm Ind}_{P'(\Q_v)}^{\GL_n(\Q_v)}(\otimes_{i=1}^k \itPi_{\tau_w(i),v}),
\]
where the restricted tensor product on the right-hand sides are defined with respect to spherical sections $h_p^\circ$ and $\tilde{h}_p^{\circ}$ normalized to take the value $1$ on the identity for all primes $p$ at which $\itPi_p$ is unramified.
For each place $v$ of $\Q$, let $M_{w,v}^{P'}$ be the local intertwining operator defined analogous to (\ref{E:local intertwining}).
Let $h = \otimes_v h_v$ be a pure tensor.
There is a finite set $S$ of places of $\Q$ including the archimedean place such that $\itPi_v$ is unramified and $h_v = h_v^\circ$ for $v \notin S$.
By the unramified computation in \cite[\S\,4.3]{Shahidi2010} and Theorem 4.2.2 in \textit{loc.\ cit.}, we have
\begin{align}\label{E:Eisenstein pf 1}
M_w^{P'}(h_{\underline{s}}) = \prod_{1 \leq i < j \leq k \atop \tau_w^{-1}(i)>\tau_w^{-1}(j)}\frac{L(s_i-s_j, \itPi_{i}\times \itPi_j^\vee)}{L(s_i-s_j+1, \itPi_{i}\times \itPi_j^\vee)}\cdot \left(\bigotimes_{v \notin S}\tilde{h}_{v,\underline{s}}^\circ\right) \otimes \left(\bigotimes_{v \in S}N_{w,v}^{P'}(h_{v,\underline{s}})\right),
\end{align}
where $N_{w,v}^{P'}$ is the normalized intertwining operator defined by
\[
N_{w,v}^{P'} =  \prod_{1 \leq i < j \leq k \atop \tau_w^{-1}(i)>\tau_w^{-1}(j)}\frac{L(s_i-s_j+1, \itPi_{i,v}\times \itPi_{j,v}^\vee)}{L(s_i-s_j, \itPi_{i,v}\times \itPi_{j,v}^\vee)}\cdot M_{w,v}^{P'}.
\]
By conditions (1) and (2), and the results of Jacquet--Shalika \cite[Theorem 5.3]{JS1981} and Shahidi \cite[Theorem 5.2]{Shahidi1981}, the ratio of $L$-functions in (\ref{E:Eisenstein pf 1}) is holomorphic at $\underline{s}=0$.
By the result of Moeglin--Waldspurger \cite[Proposition I.10]{MW1989}, condition (2) implies that $\otimes_{v \in S}N_{w,v}^{P'}$ is holomorphic at $\underline{s}=0$.
Therefore, $M_w^{P'}(h_{\underline{s}})$ is holomorphic at $\underline{s}=0$.
\end{proof}

%\begin{rmk}
%If $\itPi_i$ is essentially tempered at all places for $1 \leq i \leq k$, then we can replace condition (2) by the weaker condition that $e(\itPi_1)\geq \cdots\geq e(\itPi_k)$.
%\end{rmk}

\subsection{Whittaker functionals}\label{SS:Whittaker}

We keep the notation of \S\,\ref{SS:Eisenstein}. 
We denote by ${}^a{\rm Ind}_{P(\A)}^{\GL_n(\A)}$ (resp.\,${}^a{\rm Ind}_{P(\Q_v)}^{\GL_n(\Q_v)}$) the algebraic induction for $\GL_n(\A)$ (resp.\,$\GL_n(\Q_v)$) defined without the normalizing factor $\delta_P$. For $1\leq i \leq k$, let $\delta_i = \delta_P \vert_{\GL_{n_i}(\A)}$ and we write $\itPi_i\delta_i:=\itPi_i\otimes\delta_i$. Then we have
\[
{\rm Ind}_{P(\A)}^{\GL_n(\A)}(\otimes_{i=1}^k\itPi_i) = {}^a{\rm Ind}_{P(\A)}^{\GL_n(\A)}(\otimes_{i=1}^k\itPi_i\delta_i).
\]
For $f \in {\rm Ind}_{P(\A)}^{\GL_n(\A)}\left(\otimes_{i=1}^k\mathcal{W}(\itPi_i)\right)$ and $\underline{s} = (s_1,...,s_k) \in \C^k$, let $f_{\underline{s}} \in {\rm Ind}_{P(\A)}^{\GL_n(\A)}\left(\otimes_{i=1}^k\mathcal{W}(\itPi_i\otimes|\mbox{ }|_\A^{s_i})\right)$ be defined as in (\ref{E:flat section}), and ${\mathbb W}(f_{\underline{s}})$ be the Whittaker function of the induced representation ${\rm Ind}_{P(\A)}^{\GL_n(\A)}(\otimes_{i=1}^k(\itPi_i\otimes|\mbox{ }|_\A^{s_i}))$ with respect to $\psi_{U_n}$ defined by %defined by the Jacquet integral (cf.\,\cite{Jacquet1967})
\begin{align}\label{E:Jacquet integral global}
{\mathbb W}(g,f_{\underline{s}}) = \int_{w_P\overline{U}_P(\A)w_P^{-1}}f_{\underline{s}}(w_P^{-1}ug)({\bf 1}_{n_1},...,{\bf 1}_{n_k})\overline{\psi_{U_n}(u)}\,du^{\rm Tam}, \quad g \in \GL_n(\A).
\end{align}
Here $\overline{U}_P$ is the transpose of $U_P$, $w_P \in W_n$ is the Weyl element represented by
\[
\bp  & & {\bf 1}_{n_k}\\&\iddots&\\ {\bf 1}_{n_1} &  & \ep,
\]
and $du^{\rm Tam}$ is the Tamagawa measure. 
Note that the integral converges absolutely when ${\rm Re}(s_i-s_j)$ is sufficiently large for $1 \leq i <j\leq k$ and admits holomorphic continuation to $\underline{s} \in \C^k$ (cf.\,\cite[\S\,3]{Shahidi2010}). 
For $h \in I_P^{\GL_n}(\otimes_{i=1}^k\itPi_i)$, let $f(h) \in {\rm Ind}_{P(\A)}^{\GL_n(\A)}\left(\otimes_{i=1}^k\mathcal{W}(\itPi_i)\right)$ defined by composing with the global Whittaker functional (\ref{E:Whittaker functional}) with $U_n$ replaced by $U_n\cap M_P$, that is, for $g \in \GL_n(\A)$, $f(h)(g) \in \otimes_{i=1}^k\mathcal{W}(\itPi_i)$ is given by
\[
f(h)(g): m\longmapsto \delta_P(m)^{-1}\int_{(U_n\cap M_P)(\Q)\backslash (U_n\cap M_P)(\A)} h(umg)\overline{\psi_{U_n}(u)}\,du^{\rm Tam},\quad m \in M_P(\A).
\]
We then have the following equality between Whittaker functions as holomorphic functions in $\underline{s} \in \C^k$ outside the singularities of $E(h_{\underline{s}})$ (cf.\,the proof of \cite[Theorem 7.1.2]{Shahidi2010}):
\begin{align}\label{E:Fourier coeff.}
W(E(h_{\underline{s}})) = {\mathbb W}(f(h_{\underline{s}})).
\end{align}
Here $W(E(h_{\underline{s}}))$ is defined as in (\ref{E:Whittaker functional}) with $\varphi$ therein replaced by $E(h_{\underline{s}})$.

Let $v$ be a place of $\Q$. For $f_v \in {\rm Ind}_{P(\Q_v)}^{\GL_n(\Q_v)}\left(\otimes_{i=1}^k\mathcal{W}(\itPi_{i,v})\right)$, let ${\mathbb W}(f_{v,\underline{s}})$ be the Whittaker function %of ${\rm Ind}_{P(\Q_v)}^{\GL_n(\Q_v)}(\otimes_{i=1}^k(\itPi_{i,v}\otimes|\mbox{ }|_v^{s_i}))$ 
on $\GL_n(\Q_v)$ defined in a similar way as in (\ref{E:Jacquet integral global}), with the Haar measure $du_v$ on $w_P\overline{U}_P(\A)w_P^{-1}$ chosen so that 
\[
\begin{cases}
\mbox{${\rm vol}(w_P\overline{U}_P(\Z_p)w_P^{-1},du_p)=1$} & \mbox{ if $v=p$ is finite},\\
\mbox{$du_\infty$ is the product measure of the Lebesgue measures on $\R$} & \mbox{ if $v = \infty$}.
\end{cases}
\]
Similarly, ${\mathbb W}(f_{v,\underline{s}})$ admits holomorphic continuation to $\underline{s} \in \C^k$.
By \cite[Theorem 15.4.1]{Wallach1992} and \cite[Corollaries 3.4.9 and 3.6.11]{Shahidi2010}, we then have a non-zero equivariant homomorphism
\begin{align}\label{E:Jacquet integral}
{\mathbb W} : {\rm Ind}_{P(\Q_v)}^{\GL_n(\Q_v)}\left(\otimes_{i=1}^k\mathcal{W}(\itPi_{i,v})\right) \longrightarrow \mathcal{W}\left({\rm Ind}_{P(\Q_v)}^{\GL_n(\Q_v)}(\otimes_{i=1}^k\itPi_{i,v})\right),\quad {\mathbb W}(f_v) = {\mathbb W}(f_{v,\underline{s}})\vert_{\underline{s}=0}.
\end{align}
Assume $v=p$ is finite. For $\sigma \in {\rm Aut}(\C)$, we have the $\sigma$-linear $M_P(\Q_p)$-equivariant isomorphism $\otimes_{i=1}^kt_\sigma : \otimes_{i=1}^k\mathcal{W}((\itPi_i \delta_i)_p) \rightarrow \otimes_{i=1}^k\mathcal{W}({}^\sigma\!(\itPi_i\delta_i)_p)$ defined as in (\ref{E:sigma-linear Whittaker}), where ${}^\sigma\!(\itPi_i\delta_i)_p$ is the $\sigma$-conjugate of $(\itPi_i\delta_i)_p$ and $u_{\sigma}$ therein is replaced by its $p$-component $u_{\sigma,p}\in\Z_p^\times$.
This defines a $\sigma$-linear $\GL_n(\Q_p)$-equivariant isomorphism 
\[
{}^a{\rm Ind}_{ P(\Q_p)}^{\GL_n(\Q_p)}(\otimes_{i=1}^kt_\sigma) : {}^a{\rm Ind}_{ P(\Q_p)}^{\GL_n(\Q_p)}\left(\otimes_{i=1}^k\mathcal{W}((\itPi_i\delta_i)_p)\right) \longrightarrow {}^a{\rm Ind}_{ P(\A_f)}^{\GL_n(\Q_p)}\left(\otimes_{i=1}^k\mathcal{W}({}^\sigma\!(\itPi_i \delta_i)_p)\right).
\]
We have the following lemma on the Galois-equivariance property of the local Whittaker functional, which will be used in the proof of Theorem \ref{T:Raghuram}.

\begin{lemma}\label{L:Galois equiv. zeta}
Let $p$ be a prime number.
For $\sigma \in {\rm Aut}(\C)$, we have the following equality of $\sigma$-linear homomorphisms from ${}^a{\rm Ind}_{P(\Q_p)}^{\GL_n(\Q_p)}\left(\otimes_{i=1}^k\mathcal{W}((\itPi_{i}\delta_{i})_p)\right)$ to $\mathcal{W}\left({}^a{\rm Ind}_{P(\Q_p)}^{\GL_n(\Q_p)}\left(\otimes_{i=1}^k{}^\sigma\!(\itPi_{i}\delta_{i})_p\right)\right)$:
\[
t_\sigma \circ {\mathbb W} = \left(\prod_{1 \leq i < j \leq k}{}^\sigma\!\omega_{(\itPi_{j}\delta_{j})_p}(u_{\sigma,p})^{-n_i}\right)\cdot {\mathbb W}\circ {}^a{\rm Ind}_{P(\Q_p)}^{\GL_n(\Q_p)}\left(\otimes_{i=1}^kt_\sigma\right).
\]
%Here $t_\sigma$ is the $\sigma$-linear isomorphism defined as in (\ref{E:sigma-linear Whittaker}), and
Here $u_{\sigma,p} \in \Z_p^\times$ is the unique element such that $\sigma(\psi_p(x)) = \psi_p(u_{\sigma,p}x)$ for all $x \in \Q_p$.
\end{lemma}

\begin{proof}
We drop the subscript $p$ for brevity. Since the Whittaker integrals ${\mathbb W}$ commute with parabolic induction in stages, we may assume $k=2$.
For $N \geq 1$, let $d_{\sigma,N} = {\rm diag}(u_\sigma^{-N+1},u_\sigma^{-N+2},...,1) \in \GL_{N}(\Z_p)$. 
Let $f \in {}^a{\rm Ind}_{P(\Q_p)}^{\GL_n(\Q_p)}\left(\mathcal{W}(\itPi_1\delta_1)\otimes \mathcal{W}(\itPi_2\delta_2)\right)$ and $g \in \GL_n(\Q_p)$.
By \cite[Theorem 3.4.7]{Shahidi2010}, the Whittaker integral ${\mathbb W}(d_{\sigma,n}g,f)$ is a stable integral, that is, there exists a sufficiently small integer $m$ such that
\[
{\mathbb W}(d_{\sigma,n}g,f) = \int_{{\rm M}_{n_2,n_1}(p^{m}\Z_p)} f\left( \bp 0 & {\bf 1}_{n_1}\\{\bf 1}_{n_2} & 0\ep \bp {\bf 1}_{n_2} & x \\ 0 & {\bf 1}_{n_1}\ep d_{\sigma,n}g \right)({\bf 1}_{n_1}, {\bf 1}_{n_2})\overline{\psi(x_{n_2,1})}\,dx.
\]
Since ${\rm M}_{n_2,n_1}(p^{m}\Z_p)$ is compact, the integration commutes with $\sigma$-conjugation. 
Therefore, we have
\begin{align*}
&(t_\sigma\circ{\mathbb W})(g,f)\\
& = \int_{{\rm M}_{n_2,n_1}(p^{m}\Z_p)} \sigma \left( f\left( \bp 0 & {\bf 1}_{n_1}\\{\bf 1}_{n_2} & 0\ep \bp {\bf 1}_{n_2} & x \\ 0 & {\bf 1}_{n_1}\ep \bp u_\sigma^{-n_1}\cdot d_{\sigma,n_2} & 0 \\ 0 & d_{\sigma,n_1}\ep g \right)({\bf 1}_{n_1}, {\bf 1}_{n_2})\right) \overline{\psi(u_\sigma x_{n_2,1})}\,dx\\
& = \int_{{\rm M}_{n_2,n_1}(p^{m}\Z_p)} \sigma \left( f\left( \bp  d_{\sigma,n_1} & 0 \\ 0 & u_\sigma^{-n_1}\cdot d_{\sigma,n_2} \ep\bp 0 & {\bf 1}_{n_1}\\{\bf 1}_{n_2} & 0\ep \bp {\bf 1}_{n_2} & u_\sigma^{n_1}\cdot d_{\sigma,n_2}^{-1}x d_{\sigma,n_1} \\ 0 & {\bf 1}_{n_1}\ep g \right)({\bf 1}_{n_1}, {\bf 1}_{n_2})\right) \\
&\quad\quad\quad\quad\quad\quad\quad\quad\quad\quad\quad\quad\quad\quad\quad\quad\quad\quad\quad\quad\quad\quad\quad\quad\quad\quad\quad\quad\quad\quad\quad\quad\quad\quad\quad\quad\quad\overline{\psi(u_\sigma x_{n_2,1})}\,dx\\
& = {}^\sigma\!\omega_{\itPi_2\delta_2}(u_\sigma)^{-n_1}\cdot \int_{{\rm M}_{n_2,n_1}(p^{m}\Z_p)} \sigma \left( f\left( \bp 0 & {\bf 1}_{n_1}\\{\bf 1}_{n_2} & 0\ep \bp {\bf 1}_{n_2} &  y \\ 0 & {\bf 1}_{n_1}\ep g \right)(d_{\sigma,n_1}, d_{\sigma,n_2})\right)\overline{\psi(y_{n_2,1})}\,dy\\
& = {}^\sigma\!\omega_{\itPi_2\delta_2}(u_\sigma)^{-n_1}\cdot {\mathbb W}\left(g,{}^a{\rm Ind}_{P(\Q_p)}^{\GL_n(\Q_p)}\left(t_\sigma\otimes t_\sigma\right)(f)\right).
\end{align*}
Here in the third equality we make a change of variable from $u_\sigma^{n_1}\cdot d_{\sigma,n_2}^{-1}x d_{\sigma,n_1}$ to $y$.
This completes the proof.
\end{proof}

\begin{rmk}\label{R:Gauss sum}
The occurrence of the product of Gauss sums in Definition \ref{D:Betti-Whittaker} is a consequence of Lemma \ref{L:Galois equiv. zeta}.
\end{rmk}

\subsection{Algebraicity of Eisenstein cohomology classes}

Let 
\[
\itPi = \itPi_1 \boxplus \cdots \boxplus \itPi_k
\]
be a regular algebraic automorphic representation of $\GL_n(\A)$.
We assume the cuspidal summands $\itPi_1,...,\itPi_k$ are arranged so that conditions in Definition \ref{D:Betti-Whittaker} are satisfied.
%Assume further the following condition is also satisfied:
%\begin{align}\label{E:exponent condition}
%\mbox{For $1 \leq i < j \leq k$, we have either $e(\itPi_i) = e(\itPi_j)$ or $e(\itPi_i)-e(\itPi_j)\geq 1$.}
%\end{align}
Let $P=P_{(n_1,...,n_k)}$.
By Lemma \ref{P:Eisenstein}, we have an $((\frak{g}_n,{\rm O}_n(\R))\times\GL_n(\A_f))$-equivariant homomorphism
\begin{align}\label{E:Eisenstein realization}
E:I_P^{\GL_n}(\otimes_{i=1}^k\itPi_i ) \longrightarrow \mathcal{A}(\GL_n),\quad h \longmapsto E(h):=E(h_{\underline{s}})\vert_{\underline{s}=0}.
\end{align}
%For $1 \leq i \leq k$, we have $(\itPi_i\delta_i)_\infty \in \Omega_{\mu^{(i)}}$ for some $\mu^{(i)} \in X^+(T_{n_i})$.
By Lemma \ref{L:combinatorial}, we have $\itPi_\infty \in \Omega_{\mu}$ for some $\mu \in X^+(T_{n})$. % uniquely determined by $\mu^{(1)},\cdots,\mu^{(k)}$ as described in the lemma.
Let $\varepsilon \in \{\pm1\}$ if $n$ is even, and $\varepsilon = \varepsilon(\itPi_\infty)$ if $n$ is odd.
%Note that $\varepsilon(\itPi_\infty) = \varepsilon(\itPi_{k,\infty})$ if $n$ (hence $n_k$) is odd.
%In (\ref{E:Betti-Whittaker}), we have defined a $\GL_{n_i}(\A_f)$-equivariant homomorphism $\Phi_{\itPi_i}^\varepsilon$ from $\mathcal{W}(\itPi_i)$ to $H^{b_{n_i}}(\mathcal{S}_{n_i},\mathcal{M}_{\mu^{(i)},\C})$ for each $i$. 
%In the following, we define analogous equivariant homomorphism $\Phi_{\itPi}^\varepsilon$ from $\mathcal{W}\left({\rm Ind}_{P(\A)}^{\GL_n(\A)}(\otimes_{i=1}^k\itPi_i)\right)$ to $H^{b_n}(\mathcal{S}_{n},\mathcal{M}_{\mu,\C})$.
Let 
\[
E^{\bullet} : H^{\bullet}(\frak{g}_n,K_n^\circ;I_P^{\GL_n}(\otimes_{i=1}^k\itPi_i)\otimes M_{\mu,\C}) \longrightarrow H^{\bullet}(\mathcal{S}_{n},\mathcal{M}_{\mu,\C})
\]
be the $(\pi_0(\GL_n(\R))\times\GL_n(\A_f))$-equivariant homomorphism induced from $E$ in (\ref{E:Eisenstein realization}), which is called the Eisenstein map associated to $\itPi$.
The cohomology classes in the image of $E^\bullet$ are called the Eisenstein cohomology classes associated to $\itPi$.
The equivariant isomorphism (cf.\,(\ref{E:Whittaker isomorphism}))
\begin{align*}
{\rm Ind}_{P(\A)}^{\GL_n(\A)}(\otimes_{i=1}^k \Upsilon_{\itPi_i} ):{\rm Ind}_{P(\A)}^{\GL_n(\A)}\left(\otimes_{i=1}^k\mathcal{W}(\itPi_i)\right) &\longrightarrow  {\rm Ind}_{P(\A)}^{\GL_n(\A)}(\otimes_{i=1}^k\itPi_i)
%\mathbb{W} : {\rm Ind}_{P(\A)}^{\GL_n(\A)}(\otimes_{i=1}^k\mathcal{W}(\itPi_i)) &\longrightarrow \mathcal{W}\left({\rm Ind}_{P(\A)}^{\GL_n(\A)}(\otimes_{i=1}^k\itPi_i)\right)
\end{align*}
induces an $(\pi_0(\GL_n(\R))\times\GL_n(\A_f))$-equivariant isomorphism
\begin{align*}
&{\rm Ind}_{P(\A)}^{\GL_n(\A)}(\otimes_{i=1}^k \Upsilon_{\itPi_i})^{\bullet}:
H^{\bullet}(\frak{g}_n,K_n^\circ;{\rm Ind}_{P(\A)}^{\GL_n(\A)}(\otimes_{i=1}^k\mathcal{W}(\itPi_i))\otimes M_{\mu,\C}  ) \\
&\quad\quad\quad\quad\quad\quad\quad\quad\quad\quad\quad\quad\quad\quad\quad\quad\quad\quad\quad\quad\longrightarrow  H^{\bullet}(\frak{g}_n,K_n^\circ;I_P^{\GL_n}(\otimes_{i=1}^k\itPi_i)\otimes M_{\mu,\C}).
\end{align*}
%and
%\begin{align*} 
%&\mathbb{W}^{b_n} : H^{b_n}(\frak{g}_n,K_n^\circ;{\rm Ind}_{P(\A)}^{\GL_n(\A)}(\otimes_{i=1}^k\mathcal{W}(\itPi_i))\otimes M_{\mu,\C}  )\\
%&\quad\quad\quad\quad\quad\quad\quad\quad\quad\quad\quad\quad\longrightarrow H^{b_n}\left(\frak{g}_n,K_n^\circ;\mathcal{W}\left({\rm Ind}_{P(\A)}^{\GL_n(\A)}(\otimes_{i=1}^k\itPi_i)\right)\otimes M_{\mu,\C}\right).
%\end{align*}
For $1 \leq i \leq k$, we have fixed a generator $[(\itPi_i\delta_i)_\infty]^\varepsilon$ for the relative Lie algebra cohomology in bottom degree fixed in Definition \ref{D:generator}. 
Let $I_{\otimes_{i=1}^k\mathcal{W}(\itPi_{i,\infty})}^\varepsilon$ be the isomorphism defined in (\ref{E:Delorme}) with $V_{\itPi_{i,\infty}} = \mathcal{W}(\itPi_{i,\infty})$. Then we have a generator
\[
\left(I_{\otimes_{i=1}^k\mathcal{W}(\itPi_{i,\infty})}^\varepsilon\right)^{-1}\left(\otimes_{i=1}^k[(\itPi_i\delta_i)_\infty]^\varepsilon\right) \in H^{b_n}(\frak{g}_n,K_n^\circ;{\rm Ind}_{P(\R)}^{\GL_n(\R)}(\otimes_{i=1}^k\mathcal{W}(\itPi_{i,\infty}))\otimes M_{\mu,\C})[\varepsilon].
\]
Let 
\[
\Phi_\itPi^\varepsilon : {}^a{\rm Ind}_{P(\A_f)}^{\GL_n(\A_f)}\left(\otimes_{i=1}^k\mathcal{W}((\itPi_{i}\delta_i)_f)\right) \longrightarrow H^{b_n}(\mathcal{S}_n,\mathcal{M}_{\mu,\C})
\]
be the $\GL_n(\A_f)$-equivariant homomorphism defined by
\begin{align}\label{E:Betti-Whittaker Eisenstein}
\Phi_\itPi^\varepsilon = E^{b_n}\circ{\rm Ind}_{P(\A)}^{\GL_n(\A)}(\otimes_{i=1}^k \Upsilon_{\itPi_i})^{b_n}\circ\left(\left(I_{\otimes_{i=1}^k\mathcal{W}(\itPi_{i,\infty})}^\varepsilon\right)^{-1}\left(\otimes_{i=1}^k[(\itPi_i\delta_i)_\infty]^\varepsilon\right)\otimes\,\cdot\,\right).
\end{align}
%Let $\itPi = \itPi_1 \boxplus\cdots\boxplus\itPi_k$ be a cohomological tamely isobaric automorphic representation of $\GL_n(\A)$, where $\itPi_i$ is a cuspidal automorphic representation of $\GL_{n_i}(\A)$ for $1 \leq i \leq k$. Let $P$ be the standard parabolic subgroup of $\GL_n$ of type $(n_1,\cdots,n_k)$.
%Then we have
%\[
%\itPi \cong {}^a{\rm Ind}_{P(\A)}^{\GL_n(\A)}(\itPi_1\rho_1\otimes \cdots \otimes \itPi_k\rho_k).
%\]
%Recall here the subscript refers to the ${\rm O}_n(\R)$-finite part, and $\itPi_i\rho_i$ is defined as in \S\,\ref{SS:BW periods} for $1 \leq i \leq k$.
For $\sigma \in {\rm Aut}(\C)$, we have the $\sigma$-linear $\GL_n(\A_f)$-equivariant isomorphism 
\[
{}^a{\rm Ind}_{ P(\A_f)}^{\GL_n(\A_f)}(\otimes_{i=1}^kt_\sigma) : {}^a{\rm Ind}_{ P(\A_f)}^{\GL_n(\A_f)}\left(\otimes_{i=1}^k\mathcal{W}((\itPi_i\delta_i)_f)\right) \longrightarrow {}^a{\rm Ind}_{ P(\A_f)}^{\GL_n(\A_f)}\left(\otimes_{i=1}^k\mathcal{W}({}^\sigma\!(\itPi_i \delta_i)_f)\right),
\]
where $t_\sigma: \mathcal{W}((\itPi_i\delta_i)_f) \rightarrow \mathcal{W}({}^\sigma\!(\itPi_i\delta_i)_f)$ is defined in (\ref{E:sigma-linear Whittaker}) for each $1 \leq i \leq k$.
Following is the main result of this section, where we prove that the algebraicity of the Eisenstein cohomology classes in bottom degree can be expressed in terms of the Betti--Whittaker periods of the cuspidal summands.
%The idea of the proof is to give a cohomological interpretation of the Langlands--Shahidi method for the Rankin--Selberg $L$-functions.
The idea of the proof goes back to Harder \cite[\S\,II]{Harder1987} for $\GL_2$ over a number field (see also \cite[\S\,4.5]{Mahnkopf2005}, \cite[\S\,5]{GH2016}, \cite[Chapter 6]{HR2020}, \cite[\S\,2]{GL2021}).
We consider the restriction to the boundary component associated to $P$ of the Borel--Serre compactification, and construct a splitting homomorphism of the restriction map in the bottom degree by using the constant term of Eisenstein series along $P$.

\begin{thm}\label{T:period relation}
Let $\itPi = \itPi_1\boxplus\cdots\boxplus\itPi_k$ be a regular algebraic automorphic representation of $\GL_n(\A)$ satisfying the conditions in Definition \ref{D:Betti-Whittaker}. % and (\ref{E:exponent condition}).
Let $\varepsilon \in \{\pm1\}$ if $n$ is even, and $\varepsilon = \varepsilon(\itPi_\infty)$ if $n$ is odd.
Then we have
\begin{align*}
\sigma^{b_n}\circ\left( \frac{\Phi_\itPi^\varepsilon}{\prod_{i=1}^k p(\itPi_i\delta_i,\varepsilon)} \right) = \left(\frac{\Phi_{{}^\sigma\!\itPi}^\varepsilon}{\prod_{i=1}^k p({}^\sigma\!(\itPi_i\delta_i),\varepsilon)}\right)\circ {}^a{\rm Ind}_{ P(\A_f)}^{\GL_n(\A_f)}(\otimes_{i=1}^kt_\sigma),\quad \sigma \in {\rm Aut}(\C).
%& = \frac{p({}^\sigma\!\itPi,\varepsilon)}{\prod_{1 \leq i < j \leq k}G({}^\sigma\!\omega_{\itPi_j\rho_j})^{n_i}L^{(\infty)}(1,{}^\sigma\!(\itPi_i\rho_i)\rho_i^{-1} \times {}^\sigma\!(\itPi_j^\vee\rho_j^{-1})\rho_j)\cdot \prod_{i=1}^{k} p({}^\sigma\!(\itPi_i\rho_i),\varepsilon)}.
\end{align*} 
%Here $p(\itPi,\varepsilon) \in \C^\times$ is defined by
%\begin{align}
%p(\itPi,\varepsilon) = \prod_{1 \leq i < j \leq k}G(\omega_{\itPi_j\rho_j})^{n_i}L^{(\infty)}(1,\itPi_i \times \itPi_j^\vee)\cdot \prod_{i=1}^k p(\itPi_i\rho_i,\varepsilon).
%\end{align}
\end{thm}

\begin{proof}

Recall the topologial space $\mathcal{S}_n$ in (\ref{E:orbifold}). 
Let $\overline{\mathcal{S}}_n$ be the Borel--Serre compactification of $\mathcal{S}_n$ (cf.\,\cite{BS1973} and \cite{Rohlfs1996}) and $\partial_P\mathcal{S}_n$ %$\partial_P\mathcal{S}_n = P(\Q)\backslash \GL_n(\A)/K_n^\circ$
 be the boundary component of $\overline{\mathcal{S}}_n$ associated to a standard parabolic subgroup $P$.
Let $\mu \in X^+(T_n)$ and $\mathcal{M}_\mu$ be the sheaf of $\Q$-vector spaces on $\mathcal{S}_n$ associated to $M_\mu$ (cf.\,\S\,\ref{SS:rational sheaf}).
The closed embedding $\partial_P\mathcal{S}_n \hookrightarrow \overline{\mathcal{S}}_n$ defines the $(\pi_0(\GL_n(\R))\times \GL_n(\A_f))$-equivariant restriction map
\begin{align*}
{r}_P:H^{b_n}(\mathcal{S}_n,\mathcal{M}_{\mu}) \longrightarrow H^{b_n}(\partial_P\mathcal{S}_n,\mathcal{M}_{\mu}).
\end{align*}
We denote by $r_{P,\C}$ the corresponding map at the transcendental level (cf.\,(\ref{E:transcendental})). Then we have
\begin{align}\label{E:period relation proof 1}
\sigma^{b_n}\circ r_{P,\C} = r_{P,\C}\circ \sigma^{b_n}, \quad  \sigma\in{\rm Aut}(\C).
\end{align}
%By \cite[Satz 1.10]{Schwermer1983}, the map $\frak{r}_{P,\C}$ factors through itself composing with the equivariant map obtained from the constant term map
%\[
%C^\infty(\GL_n(\Q)\backslash\GL_n(\A))\longrightarrow C^\infty(U_P(\A)M_P(\A)\backslash\GL_n(\A)),\quad \varphi\longmapsto \varphi_P.
%\]
For the boundary cohomology, we have a canonical $(\pi_0(\GL_n(\R))\times \GL_n(\A_f))$-equivariant isomorphism (cf.\,\cite[Proposition 4.3]{HR2020}) %between $H^{b_n}(\partial_P\mathcal{S}_n,\mathcal{M}_{\mu})$ and 
%\begin{align*}
%\bigoplus_{w \in W_n^P}{}^a{\rm Ind}_{ P(\A_f)}^{\GL_n(\A_f)}\left[ H^{b_n-\ell(w)}(\mathcal{S}_n^{M_P},\otimes_{i=1}^k \mathcal{M}_{\mu_w,i})^{\pi_0(K_{n}^{M_P})}\right],
%\end{align*}
\begin{align*}
I_{\partial_P} : H^{b_n}(\partial_P\mathcal{S}_n,\mathcal{M}_{\mu}) \longrightarrow  \bigoplus_{w \in W_n^P \atop q_1+\cdots+q_k = b_n-\ell(w)}{}^a{\rm Ind}_{ P(\A_f)}^{\GL_n(\A_f)}\left[ \left(\bigotimes_{i=1}^k H^{q_i}(\mathcal{S}_{n_i}, \mathcal{M}_{\mu_{w,i}})\right)^{\pi_0(K_{n}^{M_P})}\right],
\end{align*}
where %$\mathcal{S}_n^{M_P} = M_P(\Q)\backslash M_P(\A)/K_n^{M_P,\circ}$ and 
$(\mu_{w,1},...,\mu_{w,k}) = w(\mu+\rho_n)-\rho_n$ for $w \in W_n^P$.
%Recall $K_n^{M_P} =K_n^\circ \cap M_P(\R)$.
%$K_n^{M_P} = K_n^\circ \cap M_P(\R)$ and $K_n^{M_P,\circ}$ is its connected component.
%The natural surjection $\mathcal{S}_{n_1}\times \cdots \times \mathcal{S}_{n_k} \twoheadrightarrow \mathcal{S}_n^{M_P}$ induces canonical $\pi_0(M_P(\R)) \times M_P(\A_f)$-equivariant homomorphism
%\[
%H^{b_n-\ell(w)}(\mathcal{S}_n^{M_P},\otimes_{i=1}^k \mathcal{M}_{\mu_w,i}) \longrightarrow \bigoplus_{q_1+\cdots+ q_k = b_n-\ell(w)}\left(\bigotimes_{i=1}^k H^{q_i}(\mathcal{S}_{n_i},\mathcal{M}_{\mu_{w,i}})\right)
%\]
%for $w \in W_n^P$. We thus obtain a $\pi_0(\GL_n(\R))\times \GL_n(\A_f)$-equivariant homomorphism
%\begin{align*}
%I_{\partial_P} : H^{b_n}(\partial_P\mathcal{S}_n,\mathcal{M}_{\mu}) \longrightarrow  \bigoplus_{w \in W_n^P\atop q_1+\cdots+ q_k = b_n-\ell(w)}{}^a{\rm Ind}_{ P(\A_f)}^{\GL_n(\A_f)}\left[ \left( \bigotimes_{i=1}^k H^{q_i}(\mathcal{S}_{n_i},\mathcal{M}_{\mu_{w,i}})\right)^{\pi_0(K_{n}^{M_P})}\right].
%\end{align*}}
We denote by $I_{\partial_P,\C}$ the corresponding map at the transcendental level, and it is clear that
\begin{align}\label{E:period relation proof 2}
{}^a{\rm Ind}_{P(\A_f)}^{\GL_n(\A_f)}(\otimes_{i=1}^k\sigma^{b_{n_i}})\circ I_{\partial_P,\C} = I_{\partial_P,\C}\circ \sigma^{b_n}, \quad  \sigma\in{\rm Aut}(\C).
\end{align}

We now assume $P=P_{(n_1,...,n_k)}$ and $\mu$ is the dominant weight so that $\itPi_\infty \in \Omega_\mu$. For $1\leq i \leq k$, we have $(\itPi_i\delta_i)_\infty \in \Omega_{\mu^{(i)}}$ for some $\mu^{(i)} \in X^+(T_{n_i})$. 
Let $w_0 \in W_n^P$ be the unique Weyl element in Lemma \ref{L:combinatorial} such that $(\mu_{w_0,1},...,\mu_{w_0,k}) = (\mu^{(1)},...,\mu^{(k)})$. 
%We denote by $H^{b_n}(\mathcal{S}_n,\mathcal{M}_{\mu,\C})[\varepsilon\times\itPi_f]$ the image of $H^{b_n}(\frak{g}_n,K_n^\circ;\itPi \otimes M_{\mu,\C})[\varepsilon]$ under $\Psi_\itPi$ in (\ref{E:Psi}).
We define the $\GL_n(\A_f)$-equivariant homomorphism
\[
r_\itPi^\varepsilon : H^{b_n}(\mathcal{S}_n,\mathcal{M}_{\mu,\C})\longrightarrow {}^a{\rm Ind}_{P(\A_f)}^{\GL_n(\A_f)} \left( \bigotimes_{i=1}^k H^{b_{n_i}}(\mathcal{S}_{n_i},\mathcal{M}_{\mu^{(i)},\C})[\varepsilon\times(\itPi_i\delta_i)_f]\right)%{}^a{\rm Ind}_{ P(\A_f)}^{\GL_n(\A_f)}\left[ H^{b_n-\ell(w_0)}(\mathcal{S}_n^{M_P},\otimes_{i=1}^k \mathcal{M}_{\mu_w,i})^{\pi_0(K_{n}^{M_P})}\right]
\]
%\[
%\frak{r}_\itPi^\varepsilon : H^{b_n}(\mathcal{S}_n,\mathcal{M}_{\mu,\C})[\varepsilon\times\itPi_f] \longrightarrow \bigoplus_{q_1+\cdots+q_k=b_n-\ell(w_0)}{}^a{\rm Ind}_{P(\A_f)}^{\GL_n(\A_f)}\left[ \left( \bigotimes_{i=1}^k H^{q_i}(\mathcal{S}_{n_i},\mathcal{M}_{\mu_{i}})\right)^{\pi_0(K_{n}^{M_P})}\right]
%\]
by projecting the image of $I_{\partial_P}\circ r_P$ to the summand indexed by $w_0$, and then following it by the surjective homomorphism
\begin{align*}
&\bigoplus_{q_1+\cdots+q_k = b_{n_1}+\cdots+b_{n_k}}{}^a{\rm Ind}_{ P(\A_f)}^{\GL_n(\A_f)}\left[ \left(\bigotimes_{i=1}^k H^{q_i}(\mathcal{S}_{n_i}, \mathcal{M}_{\mu^{(i)}})\right)^{\pi_0(K_{n}^{M_P})}\right] \\
&\quad\quad\quad\quad\quad\quad\quad\quad\quad\quad\quad\quad\longrightarrow{}^a{\rm Ind}_{P(\A_f)}^{\GL_n(\A_f)} \left( \bigotimes_{i=1}^k H^{b_{n_i}}(\mathcal{S}_{n_i},\mathcal{M}_{\mu^{(i)},\C})[\varepsilon\times(\itPi_i\delta_i)_f]\right). 
\end{align*}
Let $\omega$ be a class in $H^{b_n}(\frak{g}_n,K_n^\circ;I_P^{\GL_n}(\otimes_{i=1}^k\itPi_i)\otimes M_{\mu,\C})$ represented by
\begin{align*}
\sum_{\alpha,\beta}h_{\alpha,\beta}\otimes X_\alpha^*\otimes {\bf v}_\beta
\end{align*}
for some $h_{\alpha,\beta}\in I_P^{\GL_n}(\otimes_{i=1}^k\itPi_i)$, $X_\alpha^* \in \extp^{b_n}(\frak{g}_{n,\C}/\frak{k}_{n,\C})^*$, and ${\bf v}_\beta \in M_{\mu,\C}$.
By \cite[Satz 1.10]{Schwermer1983} %(see also \cite[\S\,9.2.1]{Harderbook})
 and the formula for constant term in (\ref{E:constant term}), the class $r_{P,\C}\circ E^{b_n}(\omega)$ is then represented by
\begin{align}\label{E:period relation proof 3}
\sum_{ w \in W_n^P\cap (W_n^P)^{-1}\atop wM_Pw^{-1}=M_P}\sum_{\alpha,\beta}M_w^P(h_{\alpha,\beta})\otimes X_{\alpha}^*\otimes {\bf v}_\beta.
\end{align}
For $w \in W_n^P\cap (W_n^P)^{-1}$ with $wM_Pw^{-1}=M_P$, by Lemma \ref{L:combinatorial} there exists a unique $w' \in W_n^P$ such that $\ell(w') = b_n - \sum_{i=1}^kb_{n_i}$ and $(\itPi_{\tau_w(i)}\delta_i)_\infty \in \Omega_{\mu_{w',i}}$ for $1 \leq i \leq k$. %is cohomological with coefficients in $M_{\mu_{w',i},\C}$ for $1 \leq i \leq k$. 
Here $\tau_w$ is the permutation defined in (\ref{E:permuation}). 
It is clear that the association $w \mapsto w'$ is injective and we have $1 \mapsto w_0$ by definition.
If $w \mapsto w'$, then we have the $(\pi_0(\GL_n(\R))\times \GL_n(\A_f))$-equivariant isomorphism constructed in Proposition \ref{L:Delorme}
\begin{align*}
I_{\otimes_{i=1}^k \itPi_{\tau_w(i)}} & :  H^{b_n}(\frak{g}_n,K_n^\circ;I_P^{\GL_n}(\otimes_{i=1}^k \itPi_{\tau_w(i)}) \otimes M_{\mu,\C}) \\
&\quad\quad\quad\quad\quad\quad\quad\longrightarrow {}^a{\rm Ind}_{P(\A_f)}^{\GL_n(\A_f)}\left[ \left( \bigotimes_{i=1}^k H^{b_{n_i}}(\frak{g}_{n_i},K_{n_i}^\circ;\itPi_{\tau_w(i)}\delta_i \otimes M_{\mu_{w',i},\C})\right)^{\pi_0(K_{n}^{M_P})}\right],
\end{align*}
where we take the automorphic realizations as the underlying spaces.
%by \cite[III, Theorem 3.3]{BW2000}.
Recall $I_{\otimes_{i=1}^k \itPi_{\tau_w(i)}}^\varepsilon$ is the isomorphism obtained by restricting to the $\varepsilon$-isotypic components on both sides.
By (\ref{E:period relation proof 3}) and the constructions of $I_{\partial_P,\C}$ and $I_{\otimes_{i=1}^k \itPi_{\tau_w(i)}}$, we see that the image of $I_{\partial_P,\C}\circ r_{P,\C}\circ E^{b_n}(\omega)$ in the summand indexed by $w'$ is equal to
\[
{}^a{\rm Ind}_{P(\A_f)}^{\GL_n(\A_f)}(\otimes_{i=1}^k \iota^{b_{n_i}})\circ I_{\otimes_{i=1}^k \itPi_{\tau_w(i)}}^\varepsilon\left(\sum_{\alpha,\beta}M_w^P(h_{\alpha,\beta})\otimes X_{\alpha}^*\otimes {\bf v}_\beta\right).
\]
On the other hand, if $w' \in W_n^P$ is not in the image of the above association, then $I_{\partial_P,\C}\circ r_{P,\C}\circ E^{b_n}(\omega)$ has zero contribution to the summand indexed by $w'$.
For our purpose here, we can consider the projection of $I_{\partial_P,\C}\circ r_{P,\C}$ to an arbitrary summand indexed by some $w' \in W_n^P$ which lies in the image of the association $w\mapsto w'$. % whenever the intertwining operator $M_w^P$ dose not have a zero at $\underline{s} = 0$.
%This condition is obviously satisfied when $w'=w_0$.
In conclusion, we have the following commutative diagram of $\GL_n(\A_f)$-equivariant homomorphisms:
\begin{equation}\label{E:diagram 0}
\begingroup
\small
\begin{tikzcd}[row sep=normal, column sep=large]
H^{b_n}(\frak{g}_n,K_n^\circ;I_P^{\GL_n}(\otimes_{i=1}^k \itPi_{i}) \otimes M_{\mu,\C})[\varepsilon] \arrow[rr, "E^{b_n}"] \arrow[dd, "I_{\otimes_{i=1}^k \itPi_{i}}^\varepsilon"]& &H^{b_n}(\mathcal{S}_n,\mathcal{M}_{\mu,\C})\arrow[dd, "r_\itPi^\varepsilon"]\\
& &\\
{}^a{\rm Ind} \left(\bigotimes_{i=1}^k H^{b_{n_i}}(\frak{g}_{n_i},K_{n_i}^\circ;\itPi_{i}\delta_i \otimes M_{\mu^{(i)},\C})[\varepsilon]\right) \arrow[rr, "{}^a{\rm Ind}\left(\otimes_{i=1}^k \iota^{b_{n_i}}\right)"] & &{}^a{\rm Ind}\left( \bigotimes_{i=1}^k H^{b_{n_i}}(\mathcal{S}_{n_i},\mathcal{M}_{\mu^{(i)},\C})[\varepsilon\times(\itPi_i\delta_i)_f]\right).
\end{tikzcd}
\small
\endgroup
\end{equation}
Here in the diagram ${}^a{\rm Ind}$ is the abbreviation for ${}^a{\rm Ind}_{P(\A_f)}^{\GL_n(\A_f)}$. In particular, $r_{\itPi}^\varepsilon\circ E^{b_n}$ is injective.
%In particular, $\frak{r}_\itPi^\varepsilon \circ E^*$ is injective.

To prove the assertion of the theorem, let $\sigma \in {\rm Aut}(\C)$ and consider the following diagram (non-commutative):
\begin{equation}\label{E:diagram}
\begingroup
\tiny
\begin{tikzcd}[row sep=normal, column sep=small]
& & &
H^{b_n}(\mathcal{S}_n,\mathcal{M}_{\mu,\C})  \arrow[dddddddl, "\sigma^{b_n}"] \arrow[dd, "r_\itPi^\varepsilon"]
\\
&&&
\\
& 
{}^a{\rm Ind}\left(\otimes_{i=1}^k\mathcal{W}((\itPi_i\delta_i)_f)\right)\arrow[uurr, "\Phi_\itPi^\varepsilon"] \arrow[dddddddl, "{}^a{\rm Ind}\left(\bigotimes_{i=1}^k t_\sigma\right)"] \arrow[rr, "{}^a{\rm Ind}\left(\bigotimes_{i=1}^k \Phi_{\itPi_i\delta_i}^\varepsilon\right)"']
& &
{}^a{\rm Ind} \left( \bigotimes_{i=1}^k H^{b_{n_i}}(\mathcal{S}_{n_i},\mathcal{M}_{\mu^{(i)},\C})[\varepsilon\times(\itPi_i\delta_i)_f]\right)\arrow[dddddddl, "{}^a{\rm Ind}\left(\bigotimes_{i=1}^k\sigma^{b_{n_i}}\right)"] 
\\ 
&&&
\\
&&&
\\
&&&
\\
&&&
\\
& &
H^{b_n}(\mathcal{S}_n,\mathcal{M}_{\mu,\C}) \arrow["r_{{}^\sigma\!\itPi}^\varepsilon"', dd]
&
\\
&&&
\\
{}^a{\rm Ind}\left(\otimes_{i=1}^k\mathcal{W}({}^\sigma\!(\itPi_i\delta_i)_f)\right)\arrow[uurr, "\Phi_{{}^\sigma\!\itPi}^\varepsilon"] \arrow[rr, "\,{}^a{\rm Ind}\left(\bigotimes_{i=1}^k\Phi_{{}^\sigma\!(\itPi_i\delta_i)}^\varepsilon\right)"']
& &
{}^a{\rm Ind} \left( \bigotimes_{i=1}^k H^{b_{n_i}}(\mathcal{S}_{n_i},\mathcal{M}_{\mu^{(i)},\C})[\varepsilon\times{}^\sigma\!(\itPi_i\delta_i)_f]\right)
&
\end{tikzcd}
\tiny
\endgroup
\end{equation}
%It is clear that the assertion in the theorem holds if and only if there exists $C \in \Q^\times$ such that the up face of the diagram is commutative after we multiply $\Phi_\itPi^\varepsilon$ and $\Phi_{{}^\sigma\!\itPi}^\varepsilon$ by $C\cdot p(\itPi,\varepsilon)^{-1}$ and $C\cdot p({}^\sigma\!\itPi,\varepsilon)^{-1}$, respectively.
%We prove this later assertion by consider the remaining five faces of diagram (\ref{E:diagram}) as follows:
%\\
%\\
The right face of diagram (\ref{E:diagram}) is commutative by (\ref{E:period relation proof 1}) and (\ref{E:period relation proof 2}).
By definition and normalization of the Betti--Whittaker periods of $\itPi_i\delta_i$ and ${}^\sigma\!(\itPi_i\delta_i)$ for $1 \leq i\leq k$ in Lemma \ref{L:BW period}, the down face of diagram (\ref{E:diagram}) is commutative after we multiply ${}^a{\rm Ind}\left(\otimes_{i=1}^k\left(\Phi_{\itPi_i\delta_i}^\varepsilon\right)\right)$ and ${}^a{\rm Ind}\left(\otimes_{i=1}^k\left(\Phi_{{}^\sigma\!(\itPi_i\delta_i)}^\varepsilon\right)\right)$ by $\prod_{i=1}^k p(\itPi_i\delta_i,\varepsilon)^{-1}$ and $\prod_{i=1}^{k} p({}^\sigma\!(\itPi_i\delta_i),\varepsilon)^{-1}$, respectively.
For the back face of diagram (\ref{E:diagram}), by the definition of $\Phi_\itPi^\varepsilon$ and $\Phi_{\itPi_i\delta_i}^\varepsilon$, it can be separate into three parts: The isomorphism
\[
(\otimes_{i=1}^k [(\itPi_i\delta_i)_\infty]^\varepsilon) \otimes \,\,\cdot\,\,: {}^a{\rm Ind}(\otimes_{i=1}^k\mathcal{W}((\itPi_i\delta_i)_f)) \longrightarrow {}^a{\rm Ind} \left(\bigotimes_{i=1}^k H^{b_{n_i}}(\frak{g}_{n_i},K_{n_i}^\circ;\mathcal{W}(\itPi_{i}\delta_i) \otimes M_{\mu^{(i)},\C})[\varepsilon]\right),
\]
the commutative diagram in (\ref{E:diagram 0}), and
\[
\begingroup
\small
\begin{tikzcd}[row sep=normal, column sep=normal]
H^{b_n}(\frak{g}_n,K_n^\circ;{\rm Ind}(\otimes_{i=1}^k \mathcal{W}(\itPi_i))\otimes M_{\mu,\C})[\varepsilon]\arrow[rrr, "{}^a{\rm Ind}\left(\otimes_{i=1}^k\Upsilon_{\itPi_i\delta_i}\right)^{b_n}"]\arrow[d, "I_{\otimes_{i=1}^k \mathcal{W}(\itPi_{i,\infty})}^\varepsilon"]&&&H^{b_n}(\frak{g}_n,K_n^\circ;I_P^{\GL_n}(\otimes_{i=1}^k \itPi_{i}) \otimes M_{\mu,\C})[\varepsilon]\arrow[d, "I_{\otimes_{i=1}^k \itPi_{i}}^\varepsilon"]\\
{}^a{\rm Ind} \left(\bigotimes_{i=1}^k H^{b_{n_i}}(\frak{g}_{n_i},K_{n_i}^\circ;\mathcal{W}(\itPi_{i}\delta_i) \otimes M_{\mu^{(i)},\C})[\varepsilon]\right) \arrow[rrr, "{}^a{\rm Ind}\left(\otimes_{i=1}^k\Upsilon_{\itPi_i\delta_i}^{b_{n_i}}\right)"] &&&{}^a{\rm Ind} \left(\bigotimes_{i=1}^k H^{b_{n_i}}(\frak{g}_{n_i},K_{n_i}^\circ;\itPi_{i}\delta_i \otimes M_{\mu^{(i)},\C})[\varepsilon]\right). 
\end{tikzcd}
\small
\endgroup
\]
The above diagram is commutative by the construction of $I_{\otimes_{i=1}^k \mathcal{W}(\itPi_{i,\infty})}^\varepsilon$ and $I_{\otimes_{i=1}^k \itPi_{i}}^\varepsilon$.
Therefore, the back face of diagram (\ref{E:diagram}) is commutative.
Similarly, the front face of diagram (\ref{E:diagram}) is also commutative.
Finally, as in the proof of \cite[Proposition 1.6]{Grobner2018b}, %the minimality of $b_n$ and conditions (1) and (2) in Definition \ref{D:Betti-Whittaker} imply that 
the image of $E^{b_n}$ is equal to a summand in the decomposition \cite[Theorem 2.3]{FS1998} for $H^{b_n}(\mathcal{S}_n,\mathcal{M}_\mu)$.
In the notation of \textit{loc.\ cit.}, the summand is indexed by the associate classes $\{P\}$ and $\{\varphi_P\}$, where $\varphi_P$ consisting of $\otimes_{i=1}^k\itPi_{\tau_w(i)}$ for $w \in W_n^P\cap (W_n^P)^{-1}$ with $wM_Pw^{-1}=M_P$.
It then follows from \cite[Theorem 4.3]{FS1998} that the images of $\sigma^{b_n}\circ\Phi_\itPi^\varepsilon$ and $\Phi_{{}^\sigma\!\itPi}^\varepsilon$ are equal.
Therefore, we conclude from the above discussion and the injectivity of $r_{{}^\sigma\!\itPi}^\varepsilon\circ E^{b_n}$ that the up face of diagram (\ref{E:diagram}) is commutative after multiply $\Phi_\itPi^\varepsilon$ and $\Phi_{{}^\sigma\!\itPi}^\varepsilon$ by $\prod_{i=1}^k p(\itPi_i\delta_i,\varepsilon)^{-1}$ and $\prod_{i=1}^{k} p({}^\sigma\!(\itPi_i\delta_i),\varepsilon)^{-1}$, respectively.
This completes the proof.
\end{proof}

\section{Proof of main result}\label{S:proof main}

The aim of this section is to prove our main result Theorem \ref{T:Main}. 
One of the key ingredients is an algebraicity result for $\GL_n \times \GL_{n-1}$ proved in Theorem \ref{T:Raghuram} below.
In \S\,\ref{SS:DC}, we show that Conjecture \ref{C:main 1} is compatible with Deligne's conjecture on critical motivic $L$-values.

\subsection{Rankin--Selberg $L$-functions}\label{SS:RS}

Let $\itSigma$ and $\itPi$ be algebraic automorphic representations of $\GL_n(\A)$ and $\GL_{n'}(\A)$, respectively. %Denote by $L(s,\itSigma \times \itPi)$ the Rankin--Selberg $L$-function for $\itSigma \times \itPi$.
%Recall we have the functional equation
%\begin{align}\label{E:fe}
%L(1-s,\itSigma^\vee \times \itPi^\vee) = \varepsilon(s,\itSigma\times\itPi)\cdot L(s,\itSigma \times \itPi),
%\end{align}
%where $\varepsilon(s,\itSigma\times\itPi) = \prod_v \varepsilon(s,\itSigma_v\times\itPi_v,\psi_v)$ is the product of the local $\varepsilon$-factor for $\itSigma_v\times\itPi_v$ with respect to $\psi_v$.
A critical point for the Rankin--Selberg $L$-function $L(s,\itSigma \times \itPi)$ is a half-integer $m_0 \in \Z+\tfrac{n+n'}{2}$ such that the local factors $L(s,\itSigma_\infty \times \itPi_\infty)$ and $L(1-s,\itSigma_\infty^\vee \times \itPi_\infty^\vee)$ are holomorphic at $s=m_0$.
For instance, if $\itSigma_\infty \in \Omega(n)$ and $\itPi_\infty \in \Omega(n')$ are essentially tempered with infinity types $(\underline{\kappa};\,{\sf w})$ and $(\underline{\ell};\,{\sf u})$, respectively, and $n$ is even, then the set of critical points is given by (cf.\,\cite[Proposition 7.7]{HR2020})
\begin{align}\label{E:critical range}
\left \{ m_0 \in \Z+\tfrac{n'}{2} \,\left\vert\, \tfrac{2-{\sf w}-{\sf u}-d(\underline{\kappa},\underline{\ell})}{2} \leq m_0 \leq \tfrac{-{\sf w}-{\sf u}+d(\underline{\kappa},\underline{\ell})}{2} \right\}\right.,
\end{align}
where 
\begin{align}\label{E:distance}
d(\underline{\kappa},\underline{\ell}) = \begin{cases}
\min\{|\kappa_i-\ell_j|\} & \mbox{ if $n'$ is even},\\
\min\{|\kappa_i-\ell_j|,|\kappa_i-1|\} & \mbox{ if $n'$ is odd}.
\end{cases}
\end{align}
In this case, $L(s,\itSigma \times \itPi)$ admits critical points if and only if $d(\underline{\kappa},\underline{\ell}) \geq 1$. Also, it admits a central critical point if and only if $d(\underline{\kappa},\underline{\ell}) \geq 1$ and ${\sf w}+{\sf u} \equiv n+n'+1\,({\rm mod}\,2)$.
When $\itSigma_\infty \in \Omega(n)$ and $\itPi_\infty \in \Omega(n')$ are essentially tempered and $L(s,\itSigma \times \itPi)$ admits critical points, the cuspidal summands of $\itSigma$ and $\itPi^\vee$ are non-isomorphic up to central twist, hence $L(s,\itSigma \times \itPi)$ is entire. %(cf.\,\cite[Propositions 3.3 and 3.6]{JS1981b}).
Moreover, it is non-vanishing at non-central critical points (cf.\,\cite[Theorem 1.3]{JS1976}, \cite[Theorem 5.3]{JS1981}, \cite[Theorem 5.2]{Shahidi1981}).

\subsection{Algebraicity of critical values for $\GL_n \times \GL_{n-1}$}\label{SS:Raghuram}

Let $\itSigma$ and $\itPi$ be regular algebraic automorphic representations of $\GL_n(\A)$ and $\GL_{n-1}(\A)$, respectively. 
In this section, we generalize the results of Raghuram \cite{Raghuram2009}, \cite{Raghuram2016} on the algebraicity of critical values of $L(s,\itSigma \times \itPi)$, which is a key ingredient in the proof of our main theorem in \S\,\ref{SS:main}. In \textit{loc.\ cit.}, Raghuram considered the case when both $\itSigma$ and $\itPi$ are cuspidal.
We allow $\itPi$ to be non-cuspidal and consider Eisenstein series on $\GL_{n-1}$ associated to it (see also %\cite{Mahnkopf2005}, 
\cite{GH2016}, \cite{GL2021}, \cite{LLS2022} when $\itPi_\infty$ is essentially tempered).
%We then prove a result on the existence of regular algebraic cuspidal automorphic representations, which will also be used in the proof.
%We begin with the following result of Raghuram \cite{Raghuram2009}, where the algebraicity is expressed in terms of Betti--Whittaker periods. 
%We refer to the introduction of $loc.$ $cit.$ for a brief survey on algebraicity results for $\GL_n \times \GL_{n-1}$ in the literature.  
Assume $\itSigma_\infty \in \Omega_\mu$ and $\itPi_\infty \in \Omega_\lambda$ for some $\mu \in X^+(T_n)$ and $\lambda \in X^+(T_{n-1})$. 
Let $(\underline{\kappa};\,\underline{\sf w})$ and $(\underline{\ell};\,\underline{\sf u})$ be the infinity types of $\itSigma_\infty$ and $\itPi_\infty$ respectively defined in (\ref{E:parameter}).
We say $(\itSigma_\infty,\itPi_\infty)$ is balanced/interlaced if
\begin{align}\label{E:balanced}
\begin{cases}
\kappa_1>\ell_1>\kappa_2>\ell_2>\cdots>\ell_{r-1}>\kappa_r & \mbox{ if $n=2r$},\\
\kappa_1>\ell_1>\kappa_2>\ell_2>\cdots>\ell_{r-1}>\kappa_r>\ell_r & \mbox{ if $n=2r+1$}.
\end{cases}
\end{align}
When $(\itSigma_\infty,\itPi_\infty)$ is balanced, it is easy to verify that $m+\tfrac{1}{2} \in \Z+\tfrac{1}{2}$ is critical for $L(s,\itSigma \times \itPi)$ if and only if 
\[
\mu_i \geq -\lambda_{n-i}+m \geq \mu_{i+1},\quad 1 \leq i \leq n-1.
\]
By the branching law from $\GL_{n}$ to $\GL_{n-1}$, the above inequality holds if and only if 
\[
{\rm Hom}_{\GL_{n-1}(\Q)}(M_{\mu-m}\otimes M_{\lambda},\Q) \neq 0.
\]
In this case, the space of invariant functionals is one-dimensional and any non-zero functional $\phi$ induces a pairing 
\[
Z(\cdot,\cdot): H^{b_n}(\frak{g}_n,K_n^\circ;\mathcal{W}(\itSigma_\infty\otimes|\mbox{ }|^m) \otimes M_{\mu-m,\C}) \times H^{b_{n-1}}(\frak{g}_{n-1},K_{n-1}^\circ;\mathcal{W}(\itPi_\infty) \otimes M_{\lambda,\C}) \longrightarrow \C
\]
as follows:
For $W \in \mathcal{W}(\itSigma_\infty\otimes|\mbox{ }|^m)$ and $W' \in \mathcal{W}(\itPi_\infty)$, let $Z(s,W,W')$ be the local zeta integral defined by
\begin{align}\label{E:local zeta}
Z(s,W,W')= \int_{U_{n-1}(\R)\backslash \GL_{n-1}(\R)}W\left(\bp g & 0 \\ 0 & 1 \ep\right)W'({\rm diag}(-1,1,...,(-1)^{n-1})g)|\det g|^s\,dg.
\end{align}
Here $dg$ is the quotient measure of the Haar measure on $\GL_{n-1}(\R)$ defined by the invariant measure on $\GL_{n-1}(\R) / {\rm SO}_{n-1}(\R)$ corresponding to $\extp_{1 \leq i \leq j \leq n-1}e_{ij}^* \in \extp^{n(n-1)/2} (\frak{g}_{n-1}/\frak{so}_{n-1})^*$ and the Haar measure on ${\rm SO}_{n-1}(\R)$ with ${\rm vol}({\rm SO}_{n-1}(\R))=1$, by the Haar measure on $U_{n-1}(\R)$ given by the product of Lebesgue measures, where $e_{ij}$ is the $n-1$ by $n-1$ matrix with $1$ in the $(i,j)$-entry and zeros otherwise and $\{e_{ij}^*\,\vert\,1 \leq i,j\leq n-1\}$ is the corresponding dual basis ordered lexicographically.
The integral converges absolutely for ${\rm Re}(s)$ sufficiently large and admits meromorphic continuation to $s \in \C$. Since $L(s,\itSigma_\infty \times \itPi_\infty)$ is holomorphic at $s=m+\tfrac{1}{2}$, the local zeta integral is holomorphic at $s=0$ by the result of Jacuqet \cite[Theorem 2.1-(ii)]{Jacquet2009}. 
The functionals $\phi$ and $Z(s,\cdot,\cdot)\vert_{s=0}$, together with the diagonal embedding
\begin{align}\label{E:diagonal embedding}
\GL_{n-1}(\R) \longrightarrow \GL_{n}(\R)\times\GL_{n-1}(\R),\quad g \longmapsto \left(  \bp g & 0 \\ 0 & 1 \ep, g\right),
\end{align}
then induce a homomorphism
\[
H^\bullet\left(\frak{g}_n \times \frak{g}_{n-1},K_n^\circ \times K_{n-1}^\circ;\,(\mathcal{W}(\itSigma_\infty\otimes|\mbox{ }|^m)\otimes \mathcal{W}(\itPi_\infty))\otimes (M_{\mu-m,\C}\otimes M_{\lambda,\C})\right)\longrightarrow H^\bullet \left(\frak{g}_{n-1},{\rm SO}_{n-1}(\R);\,\C\right).
\]
This defines the pairing $Z(\cdot,\cdot)$ when $\bullet = b_n+b_{n-1} = \tfrac{n(n-1)}{2}$ by the K\"{u}nneth formula and the minimality of $b_n$ and $b_{n-1}$. Here we identify $H^{n(n-1)/2} (\frak{g}_{n-1},{\rm SO}_{n-1}(\R);\C) = \extp^{n(n-1)/2} (\frak{g}_{n-1}/\frak{so}_{n-1})^*$ with $\C$ by mapping $\extp_{1 \leq i \leq j \leq n-1}e_{ij}^*$ to $1$.
We recall the following non-vanishing result for the pairing $Z(\cdot,\cdot)$. 
When $\itSigma_\infty$ and $\itPi_\infty$ are essentially tempered, this was first proved by Binyong Sun in \cite{Sun2017}. More generally, the non-vanishing in the general case follows from the archimedean period relations established by Li--Liu--Sun \cite{LLS2022}, and has recently been proved in full generality by Jin--Liu--Sun in \cite{JLS2025a}.
Alternatively, following the arguments of Sun in \cite{Sun2017}, the pairing $Z(\cdot,\cdot)$ is already non-vanishing at the level of cochain complexes when restricted to minimal $K$-types. It therefore suffices to show that these minimal $K$-types contribute non-trivially to the corresponding cohomology groups. We give a proof of this fact in Proposition \ref{L:cohomology}, as it may be of independent interest.

\begin{thm}\label{T:nonvanishing}
Assume $(\itSigma_\infty,\itPi_\infty)$ is balanced. Let $m+\tfrac{1}{2}$ be a critical point for $L(s,\itSigma \times \itPi)$. Then we have
\[
Z([\itSigma_\infty \otimes |\mbox{ }|^m]^{\varepsilon_m},[\itPi_\infty]^{\varepsilon_m'}) \neq 0.
\]
Here $\varepsilon_m=\varepsilon(\itSigma_\infty\otimes|\mbox{ }|^m)$ if $n$ is odd, $\varepsilon_m'=\varepsilon(\itPi_\infty)$ if $n$ is even, and $\varepsilon_m\varepsilon_m'=(-1)^{n}$.
\end{thm}

Following is the main result of this section on the algebraicity of the critical values of $L(s,\itSigma \times \itPi)$ in terms of the Betti--Whittaker periods of $\itSigma$ and $\itPi$.

\begin{thm}%[Raghuram]
\label{T:Raghuram}
Let $\itSigma$ and $\itPi$ be regular algebraic automorphic representations of $\GL_{n}(\A)$ and $\GL_{n-1}(\A)$, respectively.
Assume $\itSigma$ is cuspidal, $\itPi$ satisfies the conditions in Definition \ref{D:Betti-Whittaker}, %and (\ref{E:exponent condition}), 
and $(\itSigma_\infty,\itPi_\infty)$ is balanced.
For any critical point $m+\tfrac{1}{2} \in \Z+\tfrac{1}{2}$ for $L(s,\itSigma \times \itPi)$, we have
\begin{align*}
\sigma \left(\frac{L^S(m+\tfrac{1}{2},\itSigma \times \itPi)\cdot Z([\itSigma_\infty \otimes |\mbox{ }|^m]^{\varepsilon_m},[\itPi_\infty]^{\varepsilon_m'})}{G(\omega_\itPi)\cdot p(\itSigma\otimes|\mbox{ }|_\A^m,\varepsilon_m)\cdot p(\itPi,\varepsilon_m')} \right) = \frac{L^S(m+\tfrac{1}{2},{}^\sigma\!\itSigma \times {}^\sigma\!\itPi)\cdot Z([\itSigma_\infty \otimes |\mbox{ }|^m]^{\varepsilon_m},[\itPi_\infty]^{\varepsilon_m'})}{G(\omega_{{}^\sigma\!\itPi})\cdot p({}^\sigma\!\itSigma\otimes|\mbox{ }|_\A^m,\varepsilon_m)\cdot p({}^\sigma\!\itPi,\varepsilon_m')}
\end{align*}
for all $\sigma \in {\rm Aut}(\C)$. Here $S$ is a sufficiently large finite set of places of $\Q$ containing the archimedean place, $\varepsilon_m=\varepsilon(\itSigma_\infty\otimes|\mbox{ }|^m)$ if $n$ is odd, $\varepsilon_m'=\varepsilon(\itPi_\infty)$ if $n$ is even, and $\varepsilon_m\varepsilon_m'=(-1)^{n}$.
\end{thm}

\begin{proof}
We identify $\GL_{n-1}$ with a subgroup of $\GL_{n}$ by the map $g \mapsto \bp g & 0 \\ 0 & 1\ep $.
Consider the topological space
\[
\widetilde{\mathcal{S}}_{n-1} = \GL_{n-1}(\Q)\backslash \GL_{n-1}(\A) / {\rm SO}_{n-1}(\R).
\]
Then the inclusion $\GL_{n-1} \subset \GL_n$ induces a proper map $j:\widetilde{S}_{n-1} \rightarrow \mathcal{S}_n$. We also have natural surjection $\pi:\widetilde{S}_{n-1} \rightarrow \mathcal{S}_{n-1}$.  
These maps then induce homomorphisms between sheaf cohomology groups:
\[
H_c^\bullet(\mathcal{S}_n,\mathcal{M}_{\mu-m}) \longrightarrow H_c^\bullet(\widetilde{\mathcal{S}}_{n-1},j^*\mathcal{M}_{\mu-m}),\quad H^\bullet(\mathcal{S}_{n-1},\mathcal{M}_{\lambda}) \longrightarrow H^\bullet(\widetilde{\mathcal{S}}_{n-1},\pi^*\mathcal{M}_{\lambda}).
\]
Here $\mu$ and $\lambda$ are chosen so that $\itSigma_\infty \in \Omega_\mu$ and $\itPi_\infty \in \Omega_\lambda$.  
It is easy to verify that $j^*\mathcal{M}_{\mu-m}$ and $\pi^*\mathcal{M}_\lambda$ are isomorphic to the locally constant sheaves on $\widetilde{\mathcal{S}}_{n-1}$ associated to $M_{\mu-m}\vert_{\GL_{n-1}(\Q)}$ and $M_\lambda$ respectively.
As we mentioned above, the conditions $(\itSigma_\infty,\itPi_\infty)$ is balanced and $m+\tfrac{1}{2}$ is critical imply that $M_\lambda^\vee$ appears in $M_{\mu-m}\vert_{\GL_{n-1}(\Q)}$ with multiplicity one.
Hence $\pi^*\mathcal{M}_\mu^\vee$ also appears in $j^*\mathcal{M}_{\mu-m}$ with multiplicity one.
%the constant sheaf $\underline{\Q}$ on $\widetilde{\mathcal{S}}_{n-1}$.
% appears in $j^*\mathcal{M}_{\mu-m} \otimes \mathcal{M}_{\lambda}$ with multiplicity one.
Since 
$
b_n+b_{n-1}  = {\rm dim}_\R\,(\frak{g}_{n-1}/\frak{so}_{n-1})$,
composing with the Poincar\'{e} duality map, we then obtain a homomorphism 
\[
\<\cdot,\cdot\> : H_c^{b_n}(\mathcal{S}_n,\mathcal{M}_{\mu-m})\times H^{b_{n-1}}(\mathcal{S}_{n-1},\mathcal{M}_{\lambda}) \longrightarrow \Q
\]
which depends implicitly on a choice of non-zero $\phi \in {\rm Hom}_{\GL_{n-1}(\Q)}(M_{\mu-m}\otimes M_{\lambda},\Q)$.
We denote by $\<\cdot,\cdot\>_\C$ the corresponding map at the transcendental level. Then it is clear that
\begin{align}\label{E:Raghuram pf 1}
\sigma(\<\cdot,\cdot\>_\C) = \<\cdot,\cdot\>_\C\circ (\sigma^{b_n}\times\sigma^{b_{n-1}}),\quad \sigma \in {\rm Aut}(\C).
\end{align}
Let $\omega$ and $\omega'$ be classes in $H_c^{b_n}(\mathcal{S}_n,\mathcal{M}_{\mu-m})$ and $H^{b_{n-1}}(\mathcal{S}_{n-1},\mathcal{M}_{\lambda})$ represented respectively by
\begin{align*}
\sum_{\alpha,\beta}\varphi_{\alpha,\beta}\otimes X_\alpha^*\otimes {\bf v}_\beta,\quad \sum_{\alpha',\beta'}\varphi'_{\alpha',\beta'}\otimes Y_{\alpha'}^*\otimes {\bf v}'_{\beta'}
\end{align*}
for some {$\varphi_{\alpha,\beta} \in \mathcal{A}_0(\GL_n(\A)) $}, $\varphi'_{\alpha',\beta'} \in \mathcal{A}(\GL_{n-1}(\A)) $, $X_\alpha^* \in \extp^{b_n}(\frak{g}_{n,\C}/\frak{k}_{n,\C})^*$, $Y_{\alpha'}^* \in \extp^{b_{n-1}}(\frak{g}_{n-1,\C}/\frak{k}_{n-1,\C})^*$, ${\bf v}_\beta \in M_{\mu-m,\C}$, ${\bf v}'_{\beta'} \in M_{\lambda,\C}$.
We have
\begin{align}\label{E:Raghuram pf 2}
\begin{split}
\<\omega,\omega'\>_\C &= \sum_{\alpha,\beta,\alpha',\beta'}\left(\iota(X_\alpha^*), {\rm pr}(Y_{\alpha'}^*)\right)\cdot {\phi_\C({\bf v}_\beta\otimes{\bf v}'_{\beta'})}\cdot\int_{\GL_{n-1}(\Q)\backslash\GL_{n-1}(\A)}\varphi_{\alpha,\beta}(g)\varphi'_{\alpha',\beta'}(g)\,dg^{\rm Tam}.
\end{split}
\end{align}
Here $\iota : \extp^{b_n}(\frak{g}_{n,\C}/\frak{k}_{n,\C})^* \rightarrow \extp^{b_n}(\frak{g}_{n-1,\C}/\frak{so}_{n-1,\C})^*$ is induced from the inclusion $\GL_{n-1}(\R) \subset \GL_n(\R)$, ${\rm pr}: \extp^{b_{n-1}}(\frak{g}_{n-1,\C}/\frak{k}_{n-1,\C})^* \rightarrow \extp^{b_{n-1}}(\frak{g}_{n-1,\C}/\frak{so}_{n-1,\C})^*$ is the natural projection, and $\left(\iota(X_\alpha^*), {\rm pr}(Y_{\alpha'}^*)\right)$ is the complex number uniquely determined by
\[
\iota(X_\alpha^*)\wedge {\rm pr}(Y_{\alpha'}^*) =  \left(\iota(X_\alpha^*), {\rm pr}(Y_{\alpha'}^*)\right)\cdot \extp_{1 \leq i \leq j \leq n-1}e_{ij}^*.
\]
On the other hand, for {$\varphi \in \mathcal{A}_0(\GL_n(\A)) $}, $\varphi' \in \mathcal{A}(\GL_{n-1}(\A))$, an unfolding argument as in the proof of \cite[Proposition 3.3]{JS1981b} shows that 
\begin{align}\label{E:Raghuram pf 3}
\begin{split}
&\int_{\GL_{n-1}(\Q)\backslash\GL_{n-1}(\A)}\varphi(g)\varphi'(g)|\det g|_\A^s\,dg^{\rm Tam} \\
&= \int_{U_{n-1}(\A)\backslash \GL_{n-1}(\A)}W(g,\varphi)W({\rm diag}(-1,1,...,(-1)^{n-1})g,\varphi')|\det g|_\A^s\,dg^{\rm Tam}.
\end{split}
\end{align}
Write $\itPi = \itPi_1\boxplus \cdots \boxplus \itPi_k$ for some cuspidal automorphic representation $\itPi_i$ of $\GL_{n_i}(\A)$ arranged so that conditions in Definition \ref{D:Betti-Whittaker} are satisfied. % and (\ref{E:exponent condition}) are satisfied. Put $P = P_{(n_1,\cdots,n_k)}$. 
Let $S$ be any finite set of places of $\Q$ containing the archimedean place such that $\itSigma_p$ and $\itPi_p$ are unramified for $p \notin S$. For $p \notin S$, let 
\[
W_{\itSigma_p \otimes |\mbox{ }|_p^m}^\circ \in \mathcal{W}(\itSigma_p\otimes |\mbox{ }|
_p^m),\quad W_{(\itPi_i\delta_i)_p}^\circ \in \mathcal{W}((\itPi_i\delta_i)_p)
\]
be the $\GL_{n}(\Z_p)$-invariant and $\GL_{n-1}(\Z_p)$-invariant Whittaker functions respectively normalized so that $W_{\itSigma_p \otimes |\mbox{ }|_p^m}^\circ({\bf 1}_n) = W_{(\itPi_i\delta_i)_p}^\circ({\bf 1}_{n-1}) =1$ for $1 \leq i \leq k$. 
Let $f_p^\circ \in {}^a{\rm Ind}_{P(\Q_p)}^{\GL_{n-1}(\Q_p)}\left(\otimes_{i=1}^k \mathcal{W}((\itPi_i\delta_i)_p)\right)$ be the $\GL_{n-1}(\Z_p)$-invariant section normalized so that $f_p^\circ({\bf 1}_{n-1}) = \otimes_{i=1}^k W_{(\itPi_i\delta_i)_p}^\circ$.
By the unramified computation in \cite[Proposition 7.1.4]{Shahidi2010}, we have
\begin{align}\label{E:Raghuram pf 4}
{\mathbb W}({\bf 1}_{n-1},f_{p,\underline{s}}^\circ) = \prod_{1\leq i <  j \leq k}L(s_i-s_j+1,\itPi_{i,p} \times \itPi_{j,p}^\vee)^{-1}.
\end{align}
Assume
\[
\varphi = \Upsilon_{\itSigma \otimes |\mbox{ }|
_\A^m}(W),\quad  \varphi' = E\circ {\rm Ind}_{P(\A)}^{\GL_{n-1}(\A)}(\otimes_{i=1}^k\Upsilon_{\itPi_i})(f)
\]
for some
$
W \in \mathcal{W}(\itSigma\otimes |\mbox{ }|
_\A^m)$ and $f \in {\rm Ind}_{P(\A)}^{\GL_{n-1}(\A)}\left(\otimes_{i=1}^k\mathcal{W}(\itPi_i)\right)$.
By the definition of $\Upsilon$ and (\ref{E:Fourier coeff.}), we have $W(\varphi) = W$ and $W(\varphi') = \mathbb{W}(f)$.
If we assume further that 
\[
W = \prod_{v \in S}W_{v} \cdot \prod_{p\notin S}W_{\itSigma_p \otimes |\mbox{ }|_p^m}^\circ,\quad f = \left(\bigotimes_{v \in S}f_v\right) \otimes \left( \bigotimes_{p \notin S}f_p^\circ\right),
\]
then the right-hand side of (\ref{E:Raghuram pf 3}) is {equal to} (by \cite[\S\,1.4]{JS1981b} and (\ref{E:Raghuram pf 4}))
\[
\frac{L^S(s+m+\tfrac{1}{2},\itSigma\times\itPi)}{\prod_{1 \leq i < j \leq k}L^S(1,\itPi_i \times \itPi_j^\vee)}\cdot \prod_{v \in S} Z(s,W_v,\mathbb{W}(f_v)),
\]
where $Z(s,W_v,\mathbb{W}(f_v))$ is the local zeta integral defined similarly as in (\ref{E:local zeta}).
By Lemma \ref{L:GL_n generator}, there exists $C \in \Q^\times$ such that
\[
{\mathbb W}^{b_{n-1}}\circ \left(I^{\varepsilon_m'}_{\otimes_{i=1}^k\mathcal{W}(\itPi_{i,\infty})}\right)^{-1}\left( \otimes_{i=1}^k[(\itPi_i\delta_i)_\infty]^{\varepsilon_m'}\right) = C\cdot [\itPi_\infty]^{\varepsilon_m'}.
\]
We thus conclude from (\ref{E:Raghuram pf 2}) and the definition of $\Phi_\itPi^{\varepsilon_m'}$ and $Z([\itSigma_\infty \otimes |\mbox{ }|^m]^{\varepsilon_m},[\itPi_\infty]^{\varepsilon_m'})$ that 
\begin{align}\label{E:Raghuram pf 5}
\begin{split}
\left\<\Phi_{\itSigma \otimes |\mbox{ }|_\A^m}^{\varepsilon_m}(W),\,\Phi_{\itPi}^{\varepsilon_m'}(f)\right\>_\C 
&= C\cdot \frac{L^S(m+\tfrac{1}{2},\itSigma\times\itPi)}{\prod_{1 \leq i < j \leq k}L^S(1,\itPi_i \times \itPi_j^\vee)}\\
&\times Z([\itSigma_\infty \otimes |\mbox{ }|^m]^{\varepsilon_m},[\itPi_\infty]^{\varepsilon_m'})\cdot\prod_{p \in S\smallsetminus \{\infty\}} \left.Z(s,W_p,\mathbb{W}(f_p))\right\vert_{s=0}
\end{split}
\end{align}
for all
\[
W = \prod_{p \in S\smallsetminus\{\infty\}}W_{p} \cdot \prod_{p\notin S}W_{\itSigma_p \otimes |\mbox{ }|_p^m}^\circ,\quad f = \left(\bigotimes_{p \in S\smallsetminus\{\infty\}}f_p \right) \otimes \left(\bigotimes_{p \notin S}f_p^\circ\right).
\]
Let $\sigma \in {\rm Aut}(\C)$. We write
\[
{}^\sigma\!f = {}^a{\rm Ind}_{ P(\A_f)}^{\GL_{n-1}(\A_f)}(\otimes_{i=1}^kt_\sigma)(f),\quad {}^\sigma\!f_p = {}^a{\rm Ind}_{ P(\A_f)}^{\GL_{n-1}(\A_f)}(\otimes_{i=1}^kt_\sigma)(f_p).
\]
Similarly as above, we have
\begin{align}\label{E:Raghuram pf 6}
\begin{split}
\left\<\Phi_{{}^\sigma\!\itSigma \otimes |\mbox{ }|_\A^m}^{\varepsilon_m}(t_\sigma W),\,\Phi_{{}^\sigma\!\itPi}^{\varepsilon_m'}({}^\sigma\!f)\right\>_\C
& = C\cdot \frac{L^S(m+\tfrac{1}{2},{}^\sigma\!\itSigma\times{}^\sigma\!\itPi)}{\prod_{1 \leq i < j \leq k}L^S(1,{}^\sigma\!(\itPi_i\delta_i)\delta_i^{-1} \times {}^\sigma\!(\itPi_j\delta_j)^\vee\delta_j)}\\
&\times Z([\itSigma_\infty \otimes |\mbox{ }|^m]^{\varepsilon_m},[\itPi_\infty]^{\varepsilon_m'})\cdot\prod_{p \in S\smallsetminus \{\infty\}} \left.Z\left(s,t_\sigma W_p,\mathbb{W}({}^\sigma\!f_p)\right)\right\vert_{s=0}.
\end{split}
\end{align}
For the left-hand sides of (\ref{E:Raghuram pf 5}) and (\ref{E:Raghuram pf 6}), by Lemma \ref{L:BW period}, Theorem \ref{T:period relation}, and (\ref{E:Raghuram pf 1}), we have
\[
\sigma \left(\frac{\left\<\Phi_{\itSigma \otimes |\mbox{ }|_\A^m}^{\varepsilon_m}(W),\,\Phi_{\itPi}^{\varepsilon_m'}(f)\right\>_\C}{p(\itSigma\otimes|\mbox{ }|_\A^m,\varepsilon_m)\cdot \prod_{i=1}^k p(\itPi_i\delta_i,\varepsilon_m')}\right) = \frac{\left\<\Phi_{{}^\sigma\!\itSigma \otimes |\mbox{ }|_\A^m}^{\varepsilon_m}(t_\sigma W),\,\Phi_{{}^\sigma\!\itPi}^{\varepsilon_m'}({}^\sigma\!f)\right\>_\C}{p({}^\sigma\!\itSigma\otimes|\mbox{ }|_\A^m,\varepsilon_m)\cdot \prod_{i=1}^k p({}^\sigma\!(\itPi_i\delta_i),\varepsilon_m')}.
\]
For the right-hand sides of (\ref{E:Raghuram pf 5}) and (\ref{E:Raghuram pf 6}), by \cite[Lemma 3.1]{Chen2023}, (\ref{E:Galois Gauss sum}), and Lemma \ref{L:Galois equiv. zeta}, we have
\begin{align*}
&\sigma \left( G(\omega_\itPi)\prod_{1 \leq i < j \leq k}G(\omega_{\itPi_j\delta_j})^{n_i}\cdot\prod_{p \in S\smallsetminus \{\infty\}} \left.Z(s,W_p, \mathbb{W}(f_p))\right\vert_{s=0}\right)\\
& = G(\omega_{{}^\sigma\!\itPi})\prod_{1 \leq i < j \leq k}G({}^\sigma\!\omega_{\itPi_j\delta_j})^{n_i}\cdot\prod_{p \in S\smallsetminus \{\infty\}} \left.Z(s,t_\sigma W_p,\mathbb{W}({}^\sigma\!f_p))\right\vert_{s=0}.
\end{align*}
Note that for all prime numbers $p$, the local factor $L(s,\itPi_{i,p}\times\itPi_{j,p}^\vee)$ is holomorphic at $s=1$ by condition (3) in Definition \ref{D:Betti-Whittaker}, and we have (cf.\,\cite[Proposition 3.17]{Raghuram2009} and \cite[Lemma 3.1]{Chen2023})
\[
\sigma (L(1,\itPi_{i,p}\times\itPi_{j,p}^\vee)) = L(1,{}^\sigma\!(\itPi_i\delta_i)_p\delta_{i,p}^{-1} \times {}^\sigma\!(\itPi_j\delta_j)_p^\vee\delta_{j,p}).
\]
The theorem then follows if we choose $W_p$ and $f_p$ such that $\left.Z(s,W_p, \mathbb{W}(f_p))\right\vert_{s=0}$ is non-zero for $p \in S\smallsetminus \{\infty\}$ (see the proof of \cite[Theorem 2.7-(ii)]{JPSS1983}).
This completes the proof.
%The above Galois-equivariance of local zeta integrals can be proved by arguments as in the proof of Lemma \ref{L:Galois equiv. zeta}.
\end{proof}

\subsection{Main Theorem}\label{SS:main}

\begin{thm}\label{T:Main}
Conjecture \ref{C:main 1} holds if the following conditions are satisfied:
\begin{itemize}
\item[(1)] $nn'$ is even.
\item[(2)] $\itSigma_\infty \in \Omega(n)$ and $\itPi_\infty \in \Omega(n')$ are essentially tempered.
\item[(3)] Let $(\underline{\kappa};\,{\sf w})$ and $(\underline{\ell};\,{\sf u})$ be the infinity types of $\itSigma_\infty$ and $\itPi_\infty$. We have 
\begin{align}\label{E:regularity 1}
\min\{\kappa_i,\ell_j\} \geq
\begin{cases}
 3 & \mbox{ if $n$ and $n'$ are even},\\
 5 & \mbox{ if $n$ or $n'$ is odd},
\end{cases}
\end{align}
and
\begin{align}\label{E:regularity 2}
\mbox{$\underline{\kappa}\sqcup\underline{\ell}$ is}
\begin{cases}
\mbox{$4$-regular} & \mbox{ if ${\sf w}$ or ${\sf u}$ is odd, and ${\sf w}+{\sf u} \equiv n+n'\,({\rm mod}\,2)$},\\
\mbox{$5$-regular} & \mbox{ if ${\sf w}+{\sf u} \equiv n+n'+1\,({\rm mod}\,2)$},\\
\mbox{$6$-regular} & \mbox{ if ${\sf w}$ and ${\sf u}$ are even, and ${\sf w}+{\sf u} \equiv n+n'\,({\rm mod}\,2)$}.
\end{cases}
\end{align}
Here $\underline{\kappa}\sqcup\underline{\ell}$ is the $(\lfloor\tfrac{n}{2}\rfloor+\lfloor\tfrac{n'}{2}\rfloor)$-tuple in descending order with entries from $\underline{\kappa}$ and $\underline{\ell}$.
\end{itemize}
\end{thm}

%We have the following remark on the conditions in Theorem \ref{T:main 1}.
%\begin{rmk}
%Condition (1) in Theorem \ref{T:main 1} is equivalent to the existence of integers $t_1,t_2$ such that the isobaric sum $(\itSigma \otimes |\mbox{ }|_\A^{t_1/2}) \boxplus (\itPi^\vee \otimes |\mbox{ }|_\A^{t_2/2})$ is cohomological and tamely isobaric.
%Condition (2) in Theorem \ref{T:main 1} is necessary for the existence of a cohomological cuspidal automorphic representation of $\GL_{n+n'+1}(\A)$ such that Raghuram's result \cite{Raghuram2009} is applicable to the corresponding Rankin--Selberg $L$-function for $\GL_{n+n'+1} \times \GL_{n+n'}$.
%Conditions (1)-(5) in Theorem \ref{T:main 2} are slightly artificial in nature. 
%It guarantees the existence of cohomological tamely isobaric automorphic representation $\itPi'$ of $\GL_{n-n'-1}(\A)$ and integers $t_1,t_2$ such that $\itPsi = (\itPi \otimes |\mbox{ }|_\A^{t_1/2}) \boxplus (\itPi' \otimes |\mbox{ }|_\A^{t_2/2})$ is cohomological and tamely isobaric and Raghuram's result \cite{Raghuram2009} is applicable to $L(s,\itSigma \times \itPsi)$.
%If we assume $\itSigma$ and $\itSigma'$ are cuspidal, then Conjecture \ref{C:main 1} also holds under conditions (1)-(5) in Theorem \ref{T:main 2}.
%The main difficulty to lift these conditions is to generalize Raghuram's result to cohomological tamely isobaric automorphic representations. 
%We refer to \cite[\S\,6.6]{GH2016} for a discussion.
%\end{rmk}

\begin{proof}

%Let $\itSigma, \itSigma'$ (resp.\,$\itPi,\itPi'$) be algebraic automorphic representations of $\GL_n(\A)$ (resp.\,$\GL_{n'}(\A)$) such that %the archimedean components are essentially tempered with
%\[
%\itSigma_\infty = \itSigma_\infty',\quad \itPi_\infty = \itPi_\infty'.
%\]
%Let $(\underline{\kappa};\,{\sf w})$ and $(\underline{\ell};\,{\sf u})$ be the infinity types of $\itSigma$ and $\itPi$, respectively. Consider the following assertion:
%\\
%\\
%For all critical points $m_0 \in \Z+\tfrac{n+n'}{2}$ for $L(s,\itSigma\times\itPi)$ with $L(m_0,\itSigma \times \itPi')\cdot L(m_0,\itSigma' \times \itPi) \neq 0$, we have
%\begin{align}\label{E:main 1 proof 1}
%\frac{L(m_0,\itSigma \times \itPi)\cdot L(m_0,\itSigma' \times \itPi')}{L(m_0,\itSigma \times \itPi')\cdot L(m_0,\itSigma' \times \itPi)} \sim 1.
%\end{align}
%\\
%\\
%In this section, we prove assertion (\ref{E:main 1 proof 1}) under conditions (1) and (2) in Theorem \ref{T:main 1}.
By condition (1), we may assume $n$ is even. %and, after twisting $\itSigma$ and $\itPi$ by integral powers of $|\mbox{ }|_\A$, we have
%\[
%\begin{cases}
%{\sf w}+{\sf u}=0 & \mbox{ if $n'$ is even},\\
%{\sf w}+{\sf u}=1 & \mbox{ if $n'$ is odd}.
%\end{cases}
%\]
Put $\delta \in \{0,1\}$ with $\delta \equiv n'\,({\rm mod}\,2)$.
Let $m+\tfrac{\delta}{2}$ be a non-central and right-half critical point for $L(s,\itSigma \times \itPi)$ defined by
\begin{align}\label{E:main proof 0}
m+\tfrac{\delta}{2} = -\tfrac{({\sf w}+{\sf u})}{2} + \begin{cases}
1 & \mbox{ if ${\sf w}+{\sf u}\equiv n'\,({\rm mod}\,2)$},\\
\tfrac{5}{2} & \mbox{ if ${\sf w}+{\sf u}\equiv n'+1\,({\rm mod}\,2)$}.
\end{cases}
\end{align}
Consider a regular algebraic automorphic representation of $\GL_{n+n'}(\A)$ defined by 
\[
\tau(\itSigma,\itPi) = (\itSigma\otimes|\mbox{ }|_\A^{-\delta/2}) \boxplus (\itPi^\vee\otimes |\mbox{ }|_\A^{-m+1-\delta}).
\]
By condition (2) and our choice of $m$, it is clear that $\tau(\itSigma,\itPi)$ satisfies the conditions in Definition \ref{D:Betti-Whittaker}. % and (\ref{E:exponent condition}).
By conditions (2) and (3), there exists an essentially tempered $\pi \in \Omega(n+n'+1)$ with infinity type $(\underline{\gamma};\,\delta)$ such that $(\pi,\tau(\itSigma,\itPi)_\infty)$ is balanced and $d(\underline{\gamma},\underline{\kappa}) \geq 2$, $d(\underline{\gamma},\underline{\ell}) \geq 2$ (cf.\,(\ref{E:distance}) for notation).
Indeed, write $\underline{\kappa} \sqcup\underline{\ell} = (\kappa_1',...,\kappa_{r}')$ with $r = {\lfloor \tfrac{n+n'}{2}\rfloor}$. If ${\sf w}$, ${\sf u}$, $n'$ are even, then we define $\underline{\gamma} = (\gamma_1,...,\gamma_r)$ by $\gamma_i = \kappa_i'+3$. If ${\sf w}$, ${\sf u}$ are odd and $n'$ is even, then we define $\underline{\gamma} = (\gamma_1,...,\gamma_r)$ by $\gamma_i = \kappa_i'+2$. If ${\sf u}$ is even and ${\sf w}$, $n'$ are odd, then we define
$\underline{\gamma} = (\gamma_1,...,\gamma_{r+1})$ by $\gamma_i = \kappa_i'+2$ for $1 \leq i \leq r$ and $\gamma_{r+1} = \kappa'_r-2$.
In the remaining cases when ${\sf w}+{\sf u} \equiv n+n'+1\,({\rm mod}\,2)$, we can define $\underline{\gamma}$ in a similar way.
%Note that $\underline{\gamma}$ is $4$-regular and $\min\{\gamma_i\} \geq 5$ if $n'$ is even.
%Note that it is cohomological and tamely isobaric by condition (1), and the archimedean component $\tau(\itSigma,\itPi)_\infty$ has infinity type $(\underline{\kappa} \sqcup\underline{\ell},\,{\sf w}-\delta)$. We write
%\[
%\underline{\kappa} \sqcup\underline{\ell} = (\kappa_1',\cdots,\kappa_{r'}'), \quad r' = {\lfloor \tfrac{n+n'}{2}\rfloor}.
%\]
%Let $r = {\lfloor \tfrac{n+n'+1}{2}\rfloor}$ and define $\underline{\lambda} = (\lambda_1,\cdots,\lambda_r) \in \Z^r$ by
%\[
%\begin{cases}
%\lambda_i = \kappa_i'+2 & \mbox{ if ${\sf w}$ or ${\sf u}$ is odd},\\
%\lambda_i = \kappa_i'+3 & \mbox{ otherwise},\\
%\end{cases}
%\]
%for $1 \leq i \leq r'$ and $\lambda_r = \kappa_{r'}'-2$ if $n'$ is odd. %Note that $\lambda_r \geq 5$ when $n'$ is even by our assumption on $\kappa_{r'}'$ in condition (2). 
%Also the parity condition (1) implies that $(\underline{\lambda};\,\delta)$ is the infinity type of an irreducible admissible, cohomological, essentially tempered $(\frak{g}_{n+n'+1},{\rm O}_{n+n'+1}(\R))$-module.
Let $\itPsi$ be a regular algebraic cuspidal automorphic representation of $\GL_{n+n'+1}(\A)$ such that $\itPsi_\infty = \pi$. 
For the existence of $\itPsi$, we refer to \cite[Lemma 4.5]{Chen2023} and the reference therein (resp.\,\cite[Theorem A.1]{RW2004}) when $n'$ is even (resp.\,odd).
Instead, it also follows from the recent result of Darshan and Raghuram \cite{DR2024}.
%It is clear that the Rankin--Selberg $L$-function $L(s,\itPsi \times \tau(\itSigma,\itPi))$ admits exactly $2$ non-central critical points (cf.\,(\ref{E:critical range})), and $(\itPsi_\infty,\tau(\itSigma,\itPi)_\infty)$ is balanced by the definition of $\underline{\lambda}$.
%Fix a non-central critical point $m+\tfrac{1}{2} \in \Z+\tfrac{1}{2}$ for $L(s,\itPsi \times \tau(\itSigma,\itPi))$.
Consider the Rankin--Selberg $L$-function
\[
L(s,\itPsi \times \tau(\itSigma,\itPi)) = L(s-\tfrac{\delta}{2}, \itPsi \times \itSigma)\cdot L(s-m+1-\delta,\itPsi\times\itPi^\vee).
\]
By (\ref{E:critical range}) and condition (2), $M+\tfrac{1}{2}  \in \Z+\tfrac{1}{2}$ is a critical point for $L(s,\itPsi \times \tau(\itSigma,\itPi))$ if and only if 
\[
\tfrac{1-{\sf w}-d(\underline{\gamma},\underline{\kappa})}{2} \leq M \leq \tfrac{-1-{\sf w}+d(\underline{\gamma},\underline{\kappa})}{2},\quad m+\tfrac{\delta}{2}+\tfrac{-1+{\sf u}-d(\underline{\gamma},\underline{\ell})}{2} \leq M \leq 
m+\tfrac{\delta}{2}+\tfrac{-3+{\sf u}+d(\underline{\gamma},\underline{\ell})}{2}.
\]
By our choice of $m$ and the condition that $d(\underline{\gamma},\underline{\kappa}) \geq 2$, $d(\underline{\gamma},\underline{\ell}) \geq 2$, we see that the set of critical points is non-empty.
Moreover, when ${\sf w}+{\sf u} \equiv n+n'\,({\rm mod}\,2)$, $\tau(\itSigma,\itPi)_\infty$ is essentially tempered and $L(s,\itPsi \times \tau(\itSigma,\itPi))$ admits at least two critical points.
When ${\sf w}+{\sf u} \equiv n+n'+1\,({\rm mod}\,2)$, the parity of $d(\underline{\gamma},\underline{\kappa})$ and $d(\underline{\gamma},\underline{\ell})$ are different. 
Therefore, in any case, there exists an integer $M$ such that $M+\tfrac{1-\delta}{2}$ and $M+\tfrac{3}{2}-m-\delta$ are non-central critical points for $L(s,\itPsi \times \itSigma)$ and $L(s,\itPsi \times \itPi^\vee)$ respectively. In particular, $L(M+\tfrac{1}{2},\itPsi \times \tau(\itSigma,\itPi))$ is non-zero.
By Theorem \ref{T:Raghuram}, we have
\begin{align*}
&L^S(M+\tfrac{1}{2},\itPsi \times \tau(\itSigma,\itPi)) \cdot Z([\itPsi_\infty\otimes|\mbox{ }|^M]^{\varepsilon_M}, [\tau(\itSigma,\itPi)_\infty]^{\varepsilon_M'})\\
&\sim G(\omega_\itSigma\omega_\itPi^{-1})\cdot p(\itPsi\otimes|\mbox{ }|_\A^M,\varepsilon_M)\cdot p(\tau(\itSigma,\itPi),\varepsilon_M'),
\end{align*}
where $S$ is a sufficiently large finite set of places of $\Q$ containing the archimedean place, 
$\varepsilon_M=\varepsilon(\itPsi_\infty\otimes|\mbox{ }|^M)$ if $n'$ is even, $\varepsilon_M'=\varepsilon(\tau(\itSigma,\itPi)_\infty)$ if $n'$ is odd, and $\varepsilon_M\varepsilon_M'=(-1)^{n'+1}$.
By the period relation in Lemma \ref{L:period relation}, we have
\[
p(\tau(\itSigma,\itPi),\varepsilon_M') \sim G(\omega_\itPi)^{-n}\cdot L^{(\infty)}(m+\tfrac{\delta}{2}, \itSigma \times \itPi)\cdot p(\itSigma \otimes |\mbox{ }|_\A^{(n'-\delta)/2},\varepsilon_M')\cdot p(\itPi^\vee\otimes|\mbox{ }|_\A^{-m+1-\delta-n/2},\varepsilon_M').
\]
%Note that
%\[
%L(s,\itPsi \times \tau(\itSigma,\itPi)) = L(s-\tfrac{\delta}{2},\itPsi \times \itSigma)\cdot L(s,\itPsi \times \itPi^\vee).
%\]
Therefore, we conclude that
\begin{align}\label{E:main proof 1}
\begin{split}
&L^{(\infty)}(m+\tfrac{\delta}{2}, \itSigma \times \itPi) \\
&\sim  \frac{L^S(M+\tfrac{1-\delta}{2},\itPsi \times \itSigma )\cdot L^S(M+\tfrac{3}{2}-m-\delta,\itPsi \times \itPi^\vee)\cdot Z([\itPsi_\infty\otimes|\mbox{ }|^M]^{\varepsilon_M}, [\tau(\itSigma,\itPi)_\infty]^{\varepsilon_M'})}{G(\omega_\itSigma)\cdot G(\omega_{\itPi})^{-n-1}\cdot p(\itPsi\otimes|\mbox{ }|_\A^M,\varepsilon_M)\cdot p(\itSigma \otimes |\mbox{ }|_\A^{(n'-\delta)/2},\varepsilon_M')\cdot p(\itPi^\vee\otimes|\mbox{ }|_\A^{-m+1-\delta-n/2},\varepsilon_M')}. 
\end{split}
\end{align}
%Since $L^{(\infty)}(m+\tfrac{1}{2},\itPsi \times \tau(\itSigma,\itPi))$ is non-zero, it then follows immediately from the assumptions $\itSigma_\infty = \itSigma_\infty'$ and $\itPi_\infty = \itPi_\infty'$ that assertion (\ref{E:main 1 proof 1}) holds for the critical point $m_0 = 1-\tfrac{\delta}{2}$.
%Note that condition (1) implies that all critical points for $L(s,\itSigma \times \itPi)$ are non-central (cf.\,(\ref{E:critical range})).
%In particular, all critical values are non-zero. 
Note that the right-hand side is non-zero by Theorem \ref{T:nonvanishing}.
We denote by 
\[
\mathcal{R}(\itSigma \times \itPi) = \mathcal{R}_{m,\itPsi,M,S}(m,\itSigma \times \itPi)
\]
the ratio on the right-hand side of (\ref{E:main proof 1}) depending on the choices of $m,\itPsi,M,S$.
Let $\itSigma'$ and $\itPi'$ be other algebraic automorphic representations of $\GL_n(\A)$ and $\GL_{n'}(\A)$ respectively such that $\itSigma_\infty = \itSigma_\infty'$ and $\itPi_\infty = \itPi_\infty'$.
We enlarge $S$ so that it contains the primes at which $\itSigma'$ or $\itPi'$ are ramified.
Then it is clear that (\ref{E:main proof 1}) also holds for the Rankin--Selberg $L$-functions $L(s,\itSigma' \times \itPi')$, $L(s,\itSigma \times \itPi')$, $L(s,\itSigma' \times \itPi)$ with the same choices of $m,\itPsi,M,S$.
We define the quantities $\mathcal{R}(\itSigma' \times \itPi')$, $\mathcal{R}(\itSigma \times \itPi')$, $\mathcal{R}(\itSigma' \times \itPi)$ with respect to the same choices of $m,\itPsi,M,S$.
Then it is clear that
\[
\frac{\mathcal{R}(\itSigma \times \itPi)\cdot \mathcal{R}(\itSigma' \times \itPi')}{\mathcal{R}(\itSigma \times \itPi')\cdot \mathcal{R}(\itSigma' \times \itPi)}=1.
\]
Therefore, it follows from (\ref{E:main proof 1}) that Conjecture \ref{C:main 1} holds for the specific right-half critical point $m+\tfrac{\delta}{2}$ in (\ref{E:main proof 0}).
By the result of Harder and Raghuram \cite[Theorem 7.21]{HR2020} and decomposing the four $L$-functions into products of Rankin--Selberg $L$-functions of regular algebraic cuspidal automorphic representations, we then deduce that Conjecture \ref{C:main 1} holds for all right-half critical points and the possible central critical point, that is, for critical points greater than $-\tfrac{({\sf w}+{\sf u})}{2}$.
Replacing $\itSigma,\itSigma',\itPi,\itPi'$ by their contragredients, Conjecture \ref{C:main 1} also holds for the ratio
\[
\frac{L(s,\itSigma^\vee \times \itPi^\vee)\cdot L(s,(\itSigma')^\vee \times (\itPi')^\vee)}{L(s,\itSigma^\vee \times (\itPi')^\vee)\cdot L(s,(\itSigma')^\vee \times \itPi^\vee)}
\]
at the right-half critical points for $L(s,\itSigma^\vee \times \itPi^\vee)$ which correspond to the left-half critical points for $L(s,\itSigma \times \itPi)$ by the global functional equation.
Therefore, we are reduced to show that
\[
 \frac{L(m_0,\itSigma \times \itPi)\cdot L(m_0,\itSigma' \times \itPi')}{L(m_0,\itSigma \times \itPi')\cdot L(m_0,\itSigma' \times \itPi)} \sim \frac{L(1-m_0,\itSigma^\vee \times \itPi^\vee)\cdot L(1-m_0,(\itSigma')^\vee \times (\itPi')^\vee)}{L(1-m_0,\itSigma^\vee \times (\itPi')^\vee)\cdot L(1-m_0,(\itSigma')^\vee \times \itPi^\vee)}
\]
for all non-central critical points $m_0$ for $L(s,\itSigma \times \itPi)$. %, where $\varepsilon(s,\itSigma_f \times \itPi_f,\psi_f) = \prod_p\varepsilon(s,\itSigma_p \times \itPi_p,\psi_p)$ is the product of local $\varepsilon$-factor over finite places.
This follows from an automorphic analogue of the motivic period relation \cite[(5.1.7)]{Deligne1979} observed by Deligne which relates Deligne's periods $c^\pm(M)$ with $c^\mp(M^\vee)$ for a pure motive $M$ over $\Q$. 
We have proved the automorphic analogue in \cite[Theorem 3.2]{Chen2023} based on the global functional equation and the Galois-equivariance property of the local $\varepsilon$-factors for the Rankin--Selberg convolutions.
%But this follows directly from \cite[Theorem 3.2]{Chen2023} by decomposing the $L$-functions into product of Rankin--Selberg $L$-functions of regular algebraic cuspidal automorphic representations.
This completes the proof.
\end{proof}

\begin{rmk}\label{R:regularity}
If the regularity conditions (\ref{E:regularity 1}) or (\ref{E:regularity 2}) dose not hold, then either (i) there does not exist $\itPsi$ such that $(\itPsi_\infty,\tau(\itSigma,\itPi)_\infty)$ is balanced, or (ii) such $\itPsi$ do exist, but $L(s, \itPsi \times \itSigma)$ and $L(s, \itPsi \times \itPi^\vee)$ admit only central critical points, or $L(s, \itPsi \times \tau(\itSigma, \itPi))$ does not admit any critical points.
In case (ii), the arguments in the proof remain valid, provided that the central critical values are non-vanishing.
For instance, subject to the non-vanishing assumption, Conjecture A (Conjecture \ref{C:Deligne Sym}) holds for $\kappa \geq 3$.
\end{rmk}

\subsection{Compatibility with Deligne's conjecture}\label{SS:DC}

In this section, we verify in Proposition \ref{L:Yoshida} that Conjecture \ref{C:main 1} is compatible with Deligne's conjecture \cite{Deligne1979}, based on Clozel's conjecture \cite{Clozel1990} and Yoshida's computation of motivic periods \cite{Yoshida2001}.
Let $M$ be a motive over $\Q$ with coefficients in a number field $\E$. Associated to $M$, we have the Betti realization $H_B(M)$, the de Rham realization $H_{dR}(M)$, and the $\frak{l}$-adic realization $H_{\frak{l}}(M)$ for each finite place $\frak{l}$ of $\E$.
The dual motive $M^\vee$ of $M$ is also over $\Q$ with coefficients in $\E$ whose realizations are dual to that of $M$.
The Betti realization $H_B(M)$ is a finite dimensional $\E$-vector space equipped with an action of archimedean Frobenius $F_\infty$ and a Hodge decomposition:
\[
H_B(M) \otimes_\Q\C = \bigoplus_{p,q \in \Z}H_B^{p,q}(M),
\]
where $H_B^{p,q}(M)$ is a free $\E\otimes_\Q\C$-module of rank $h^{p,q}_M$ such that $F_\infty(H_B^{p,q}(M)) = H_B^{q,p}(M)$.
Let $H_B^\pm(M)$ be the $\pm$-eigenspace of $H_B(M)$ under the action of $F_\infty$, and write $d^\pm(M) = {\rm dim}_\E\,H_B^\pm(M)$.
We say $M$ is \textit{pure} of weight $w$ if $H_B^{p,q}(M)$ is non-zero only when $p+q=w$. We say $M$ is \textit{regular} if the non-zero Hodge numbers are all equal to $1$.
When $M$ is pure of weight $w$, the archimedean local factor of $M$ is defined by
\[
L_\infty(M,s) = \prod_{p<\tfrac{w}{2}}\Gamma_\C(s-p)^{h_M^{p,w-p}}\cdot\begin{cases}
1 & \mbox{ if $w$ is odd},\\
\Gamma_\R(s-\tfrac{w}{2})^{h^{w/2,w/2,+}_M}\cdot\Gamma_\R(s-\tfrac{w}{2}+1)^{h^{w/2,w/2,-}_M} & \mbox{ if $w$ is even}.
\end{cases}
\]
Here $h^{w/2,w/2,+}_M$ and $h^{w/2,w/2,-}_M$ are the rank of the eigenspaces of $H_B^{w/2,w/2}(M)$ under the action of $F_\infty$ with eigenvalues $(-1)^{w/2}$ and $(-1)^{w/2+1}$, respectively, and
\[
\Gamma_\R(s) = \pi^{-s/2}\Gamma(\tfrac{s}{2}),\quad \Gamma_\C(s) = 2(2\pi)^{-s}\Gamma(s).
\] 
A critical point for $M$ is an integer $m \in \Z$ such that $L_\infty(M,s)$ and $L_\infty(M^\vee,1-s)$ are holomorphic at $s=m$. Note that $h_M^{p,q} = h_{M^\vee}^{-q,-p}$ and $h_M^{w/2,w/2,\pm} = h_{M^\vee}^{-w/2,-w/2,\pm}$ if $w$ is even. Therefore, $M$ admits critical points if and only if $h_M^{w/2,w/2,+}=0$ or $h_M^{w/2,w/2,-}=0$ when $w$ is even.
The de Rham realization $H_{dR}(M)$ is a finite dimensional $\E$-vector space equipped with a decreasing Hodge filtration $\{F_{dR}^i(M)\}_{i \in \Z}$. 
We have the comparison isomorphism
\[
I_\infty : H_B(M) \otimes_\Q\C \longrightarrow H_{dR}(M) \otimes_\Q\C
\]
such that $I_\infty\circ(F_\infty\otimes c) = ({\rm id}\otimes c)\circ I_\infty$ and $I_\infty\left(\oplus_{p \geq i}H_B^{p,q}(M)\right) = F_{dR}^i(M) \otimes_\Q\C$ for all $i \in \Z$.
Assume $M$ is pure and admits critical points. Then there are subspaces $F^\pm(M)$ of $H_{dR}(M)$ from the Hodge filtration such that ${\rm dim}_\E\,F^\pm(M) = d^\pm(M)$. In this case, $I_\infty$ induces isomorphisms
\[
I_\infty^\pm : H^\pm_B(M) \otimes_\Q\C \hookrightarrow H_B(M)\otimes_\Q\C \longrightarrow H_{dR}(M)\otimes_\Q\C\twoheadrightarrow H_{dR}(M) / F^\mp(M) \otimes_\Q\C.
\]
With respect to fixed $\E$-rational bases on both sides, the Deligne's periods $c^\pm(M)$ are defined by
\[
c^\pm(M) := \det(I_\infty^\pm) \in (\E\otimes_\Q\C)^\times.
\]
The $\frak{l}$-adic realizations $\{H_{\frak{l}}(M)\}_{\frak{l}}$ form a strictly compatible system of Galois representations of $\Gal(\overline{\Q}/\Q)$ on $\E_{\frak l}$-vector spaces $H_{\frak l}(M)$. Under assumptions \cite[\S\,1.2.1]{Deligne1979} on rationality and independency of local factors, we can define the $\E\otimes_\Q\C$-valued $L$-function $L(M,s)$ of $M$. 
Note that
\[
{\rm dim}_\E\,H_B(M) = {\rm dim}_\E\,H_{dR}(M) = {\rm dim}_{\E_{\frak l}}\,H_{\frak l}(M).  
\]
This common dimension is denoted by $d(M)$, called the rank of $M$.
In \cite[Conjecture 2.8]{Deligne1979}, Deligne proposed the following:
\begin{conj}[Deligne]\label{C:Deligne}
Assume $M$ is pure and admits critical points. For a critical point $m\in\Z$ for $L(M,s)$, we have
\begin{align}\label{E:Deligne}
\frac{L(M,m)}{(2\pi\sqrt{-1})^{d^{(-1)^m}(M)m}\cdot c^{(-1)^m}(M)} \in \E.
\end{align}
\end{conj}

{%An automorphic representation $\itPi \cong \bigotimes_v \itPi_v$ of $\GL_n(\A)$ is called $algebraic$ if it is isobaric and the infinitesimal character of $\itPi_\infty$ belongs to $(\Z+\tfrac{n-1}{2})^n$.
Let $\itPi \cong \otimes_v \itPi_v$ be an algebraic automorphic representation of $\GL_n(\A)$ (cf.\,\cite[Definition 1.8]{Clozel1990}).
For $\sigma \in {\rm Aut}(\C)$, let ${}^\sigma\!\itPi: = \itPi_\infty\otimes {}^\sigma\!\itPi_f$ be the $\sigma$-conjugate of $\itPi$.
The rationality field $\Q(\itPi)$ of $\itPi$ is define to be the fixed field of $\left\{ \sigma \in {\rm Aut}(\C)\,\vert\,{}^\sigma\!\itPi \cong \itPi\right\}$.
%Clozel conjectured in \cite{Clozel1990} that if $\itPi$ is algebraic, then ${}^\sigma\!\itPi$ is automorphic for all $\sigma \in {\rm Aut}(\C)$ and $\Q(\itPi)$ is a number field. The assertion was proved in $loc.$ $cit.$ when $\itPi$ is regular algebraic and cuspidal, based on the study of cuspidal cohomology of the locally symmetric spaces associated to $\GL_n$. Moreover, 
We have the following conjecture proposed by Clozel in \cite[Conjecture 4.5]{Clozel1990}.}

\begin{conj}[Clozel]\label{C:Clozel}
Let $\itPi$ be an algebraic cuspidal automorphic representation of $\GL_n(\A)$. There exists an absolutely irreducible, pure motive $M_\itPi$ over $\Q$ of rank $n$ with coefficients in a number field $\E$ containing $\Q(\itPi)$ satisfying the following properties:
\begin{itemize}
\item[(1)] Let $\phi_{\itPi_\infty}$ be the Langlands parameter associated to $\itPi_\infty$ by the local Langlands correspondence. Write
\[
\phi_{\itPi_\infty} = (h_1\cdot \phi_{\kappa_1}\oplus\cdots\oplus h_r\cdot\phi_{\kappa_r}\oplus h_+\cdot{\bf 1}\oplus h_-\cdot {\rm sgn}) \otimes |\mbox{ }|^{{\sf w}/2}
\]
for some $\kappa_1,\cdots,\kappa_r \geq 2$, ${\sf w}\in\Z$, and some multiplicities $h_1,...,h_r,h_\pm$,
where $\phi_\kappa$ corresponds to the discrete series representation of $\GL_2(\R)$ with minimal weight $\kappa\geq 2$.
Then $M_\itPi$ is pure of weight $-{\sf w}-n+1$ and the Hodge numbers are given by
\begin{align*}
h_{M_\itPi}^{(-\kappa_i-{\sf w}-n+2)/2,(\kappa_i-{\sf w}-n)/2} &= h_{M_\itPi}^{(\kappa_i-{\sf w}-n)/2, (-\kappa_i-{\sf w}-n+2)/2} = h_i,\quad 1 \leq i \leq r,\\
h_{M_\itPi}^{(-{\sf w}-n+1)/2,(-{\sf w}-n+1)/2,\pm} &= h_\pm.
\end{align*}
\item[(2)] Let $p\neq l$ be prime numbers and $\frak{l}$ be a finite place of $\E$ lying over $l$. Then we have
\[
W\!D\left((H_{\frak{l}}(M_\itPi) \otimes_{\E_\frak{l}} \overline{\Q}_l)\vert_{\Gal(\overline{\Q}_p/\Q_p)}\right)^{F-ss} \cong \iota_l^{-1}{\rm rec}_{\Q_p}\left(\itPi_p \otimes |\mbox{ }|_p^{(n-1)/2}\right).
\]
Here $W\!D(\cdot)$ denotes the associated Weil--Deligne representation, $F-ss$ denotes the Frobenius semisimplification, ${\rm  rec}_{\Q_p}$ denotes the local Langlands correspondence, and $\iota_l : \overline{\Q}_l \rightarrow \C$ is an isomorphism.
\end{itemize}
In particular, we have $L_\infty(M_\itPi,s) = L(s+\tfrac{n-1}{2},\itPi_\infty)$ and 
\[
L(M_\itPi,s) = \left( L^{(\infty)}(s+\tfrac{n-1}{2},{}^\sigma\!\itPi)\right)_{\sigma : \E\rightarrow \C}.
\]
\end{conj}

%\begin{rmk}\label{R:GL_2}
%When $n=2$ and $\itPi$ is regular algebraic, the conjectural motive $M_\itPi$ was constructed by Deligne \cite{Deligne1971} and Scholl \cite{Scholl1990}.
%Under the canonical isomorphism $\E\otimes_\Q\C \cong \prod_{\sigma : \E\rightarrow\C}\C$, we write
%\[
%c^\pm(M_{\itPi}) = \left( c^\pm({}^\sigma\!\itPi)\right)_{\sigma : \E\rightarrow\C}.
%\]
%In this case, Conjecture \ref{C:Deligne} is known (cf.\,\cite[\S\,7]{Deligne1979}) and we have the following period relation between the Deligne's periods and the Betti--Whittaker periods of $\itPi$:
%\begin{align}\label{E:period GL_2}
%\sigma \left(\frac{c^\pm(\itPi)}{(2\pi\sqrt{-1})^{(\kappa+{\sf w})/2}\cdot p(\itPi,\pm)}\right) = \frac{c^\pm({}^\sigma\!\itPi)}{(2\pi\sqrt{-1})^{(\kappa+{\sf w})/2}\cdot p({}^\sigma\!\itPi,\pm)},\quad \sigma \in {\rm Aut}(\C).
%\end{align}
%Here $(\kappa;\,{\sf w})$ is the infinity type of $\itPi$.
%Here we need to use the existence of non-zero twisted central critical values when there is only one critical value which is zero (cf.\,\cite{Rohrlich1989}).
%\end{rmk}

Let $\itPi$ be an algebraic automorphic representation of $\GL_n(\A)$ with cuspidal summands $\itPi_1,...,\itPi_k$. Then $\itPi_i \otimes |\mbox{ }|_\A^{(n-n_i)/2}$ is algebraic and cuspidal for $1 \leq i \leq k$ by definition, where $n_i$ is the degree of $\itPi_i$. %It is clear that $\itPi_\infty$ is essentially tempered if and only if $\itPi$ is tamely isobaric. 
Let $M_\itPi$ be the conjectural motive over $\Q$ of rank $n$ with coefficients in some number field containing $\prod_{i=1}^k\Q(\itPi_i \otimes|\mbox{ }|_\A^{(n-n_i)/2})$ defined by
\[
M_\itPi = M_{\itPi_1 \otimes|\mbox{ }|_\A^{(n-n_1)/2}} \oplus \cdots \oplus M_{\itPi_k \otimes|\mbox{ }|_\A^{(n-n_k)/2}}.
\]
Then $M_\itPi$ is pure if and only if $\itPi_\infty$ is essentially tempered. 

\begin{prop}\label{L:Yoshida}
If Conjectures \ref{C:Deligne} and \ref{C:Clozel} hold, then Conjecture \ref{C:main 1} holds.
\end{prop}

\begin{proof}
The assertion is a direct consequence of a result of Yoshida in \cite{Yoshida2001}. We use freely the terminology therein. 
Let $M$ and $N$ be pure motives over $\Q$ with coefficients in a number field $\E$. Denote by $X_M \in {\rm M}_{d(M),d(M)}(\E\otimes_\Q\C)$ and $X_N \in {\rm M}_{d(N),d(N)}(\E\otimes_\Q\C)$ the period matrices of $M$ and $N$, respectively (well-defined up to $\E^\times$).
In \cite[Proposition 12]{Yoshida2001}, Yoshida explicitly computed the Deligne's periods $c^\pm(M\otimes N)$ of the tensor motive $M
\otimes N$, and show that there exist polynomial functions $f^\pm$ and $g^\pm$ on ${\rm M}_{d(M),d(M)}$ and ${\rm M}_{d(N),d(N)}$ respectively whose types are determined by the Hodge numbers of $M$ and $N$ (including the ranks of $\pm$-eigenspaces of the middle Hodge types), such that
\[
c^\pm(M\otimes N ) \in f^\pm(X_M)g^\pm(X_N)\cdot \E^\times.
\]
Therefore, if $M'$ and $N'$ are other pure motives over $\Q$ with coefficients in $\E$ such that the Hodge numbers of $M'$ and $N'$ are equal to that of $M$ and $N$, respectively, then we have
\begin{align}\label{E:Yoshida} 
\frac{c^\pm(M\otimes N)\cdot c^\pm(M'\otimes N')}{c^\pm(M\otimes N')\cdot c^\pm(M'\otimes N)} \in \frac{f^\pm(X_M)g^\pm(X_N) \cdot f^\pm(X_{M'})g^\pm(X_{N'})}{f^\pm(X_M)g^\pm(X_{N'})\cdot f^\pm(X_{M'})g^\pm(X_N)}\cdot \E^\times = \E^\times.
\end{align}
Let $\itSigma, \itSigma'$ (resp.\,$\itPi,\itPi'$) be algebraic automorphic representations of $\GL_n(\A)$ (resp.\,$\GL_{n'}(\A)$) such that %the archimedean components are essentially tempered with
\[
\itSigma_\infty = \itSigma_\infty',\quad \itPi_\infty = \itPi_\infty'.
\]
When $\itSigma_\infty$ and $\itPi_\infty$ are essentially tempered, we see that Conjecture \ref{C:main 1} follows from Conjecture \ref{C:Deligne} and (\ref{E:Yoshida}) by taking
\[
M = M_\itSigma,\quad M' = M_{\itSigma'},\quad N = M_\itPi,\quad N' = M_{\itPi'}.
\]
The general case follows from this special case by decomposing $\itSigma,\itSigma',\itPi,\itPi'$ into isobaric sums of algebraic automorphic representations with essentially tempered archimedean components according to the exponents of the summands.
\end{proof}

\section{Applications to Deligne's conjecture}\label{S:applications}

This section is devoted to applications of the cross-ratio formula Theorem \ref{T:Main} to Conjecture \ref{C:Deligne}. 
In general, Deligne's conjecture is widely open. 
For motives associated to algebraic Hecke characters, the conjecture was known when the base field is totally real (cf.\,\cite[Proposition 3.1]{Shimura1978}) or a CM-field (notably by Blasius \cite{Blasius1986}), and is settled by Kufner in \cite{Kufner2024} for arbitrary base field based on the constructions of King and Sprang in \cite{KS2024} (an alternative proof is also given by Jin, Liu, and Sun in \cite{JLS2025} by following the strategy of Harder and Schappacher in \cite{HS1985}).
%It is known for motives of rank $1$ associated to algebraic Hecke characters over totally real or CM fields (cf.\,\cite{Blasius1986}), 
The conjecture is also known for motives associated to regular algebraic cuspidal automorphic representations of $\GL_2(\A)$ (cf.\,\cite[\S\,7]{Deligne1979}, \cite{Scholl1990}). In \cite{GHL2021}, \cite{GHL2025}, Grobner, Harris, and Lin proved the conjecture for tensor products of two automorphic motives associated to regular algebraic conjugate self-dual cuspidal automorphic representations over CM-fields under various assumptions (see also \cite{GL2021}). A key ingredient in their proof is an automorphic analogue of Deligne's conjecture \cite[Theorem 1]{GHL2025}, established by Grobner--Harris--Lin--Raghuram, using a generalization of Theorem \ref{T:Main} to CM-fields.
For some other known cases in the literature, we have results for Rankin--Selberg $L$-functions for $\GL_n \times \GL_2$ in the unbalanced case, tensor product $L$-functions for $\GL_2$, and symmetric power $L$-functions for $\GL_2$ recalled in \S\,\ref{SS:Deligne for GL_n x GL_2}, \S\,\ref{SS:Blasius}, and \S\,\ref{SS:Deligne Sym} respectively.
In these known cases, the conjecture is proven based on integral representations of the $L$-function in question, and a cohomological interpretation of the integral representation (or algebraicity of the Eisenstein series involved in the construction).
As applications of our main result Theorem \ref{T:Main}, we provide another approach to prove Deligne's conjecture for these motives in higher rank cases where integral representations of the $L$-functions are not yet available.
The main results of this section are Theorems \ref{T:RS for GL_n x GL_2}, \ref{T:tensor product for GL_2}, and \ref{T:Sym odd}. 
We prove Theorems \ref{T:tensor product for GL_2} and \ref{T:Sym odd} by mathematical induction and the strategy is as follows:
\begin{itemize}
\item Show that the automorphic/motivic $L$-function in question is equal to (a factor of) a Rankin--Selberg $L$-function $L(s,\itSigma\times\itPi)$ for some regular algebraic automorphic representations $\itSigma$ and $\itPi$ with essentially tempered archimedean components.
\item Choose auxiliary algebraic automorphic representations $\itSigma'$ and $\itPi'$ with $\itSigma_\infty = \itSigma'_\infty$ and $\itPi_\infty = \itPi'_\infty$, such that the Rankin--Selberg $L$-functions $L(s,\itSigma \times \itPi')$, $L(s,\itSigma' \times \itPi)$, and $L(s,\itSigma' \times \itPi')$ decompose into products of automorphic/motivic $L$-functions of the relevant type with smaller degree.
\item By Theorem \ref{T:Main} and the induction hypothesis, we are reduced to the cases for $\GL_2 \times \GL_1$, $\GL_2 \times \GL_2$, and $\GL_2 \times \GL_2 \times \GL_2$ which are already known thanks to \cite{Shimura1976}, \cite{Shimura1977}, and \cite{GH1993}.
%In the proof of Theorem \ref{T:RS for GL_n x GL_2}, we also need the results \cite{Harris1997}, \cite{Harris2021}, and \cite{GL2016} for $\GL_n \times \GL_1$ over CM-fields.
\end{itemize}
The strategy for Theorem \ref{T:RS for GL_n x GL_2} is similar, except in the second step we show that Deligne's conjecture holds for the auxiliary Rankin--Selberg $L$-function $L(s,\itSigma \times \itPi')$ by reducing the problem from $\GL_n \times \GL_2$ over $\Q$ to $\GL_n \times \GL_1$ over imaginary quadratic fields.

\subsection{Automorphic tensor products}\label{SS:5.1}

For irreducible automorphic representations $\itPi_1$ and $\itPi_2$ of $\GL_n(\A)$ and $\GL_{n'}(\A)$, we denote by $\itPi_1 \boxtimes \itPi_2$ the automorphic tensor product of $\itPi_1$ and $\itPi_2$. It is an irreducible admissible $((\frak{g}_{nn'},{\rm O}_{nn'}(\R))\times \GL_{nn'}(\A_f))$-module defined by the restricted tensor product
\[
\itPi_1 \boxtimes \itPi_2 = \bigotimes_v \itPi_{1,v}\boxtimes \itPi_{2,v}.
\]
Here $\itPi_{1,v}\boxtimes \itPi_{2,v}$ is the irreducible admissible representation of $\GL_{nn'}(\Q_v)$ (resp.\,$(\frak{g}_{nn'},{\rm O}_{nn'}(\R))$-module) if $v$ is finite (resp.\,archimedean) defined by the local Langlands correspondence with respect to the tensor representation
\[
\GL_n(\C) \times \GL_{n'}(\C) \longrightarrow \GL_{nn'}(\C).
\]
As special cases of the Langlands functoriality conjecture, automorphic tensor products are expected to be automorphic.
Thanks to the works of Ramakrishnan \cite{Rama2000}, Dieulefait \cite{Dieulefait2020}, and Newton--Thorne \cite{NT2021}, \cite{NT2021b}, all the automorphic tensor products that appear in \S\,\ref{S:applications} will be automorphic with essentially tempered archimedean components. Moreover, they will be locally tempered everywhere up to a twist, hence must be isobaric (cf.\,\cite[Lemme 1.5]{Clozel1990}).

%\subsection{Regular algebraic cuspidal automorphic representations of $\GL_2(\A)$ of CM-type}\label{SS:5.2}
\subsection{Rankin--Selberg $L$-functions for $\GL_n \times \GL_2$}\label{SS:Deligne for GL_n x GL_2}

Let $\itSigma$ and $\itPi$ be regular algebraic automorphic representations of $\GL_n(\A)$ and $\GL_2(\A)$ with essentially tempered archimedean components.
Let $(\underline{\kappa};\,{\sf w})$ and $(\kappa';\,{\sf w}')$ be the infinity types of $\itSigma$ and $\itPi$ respectively. 
Note that $\itPi$ must be cuspidal by the regularity and essentially temperedness on $\itPi_\infty$.
%Let $f_{\itPi}$ be the normalized newform of $\itPi$ in the sense of Casselman \cite{Casselman1973}, and $\Vert f_{\itPi} \Vert$ be its Petersson norm defined by
%\begin{align}\label{E:Petersson norm}
%\Vert f_{\itPi} \Vert = \int_{\A^\times\GL_2(\Q)\backslash\GL_2(\A)}|f_{\itPi}(g)|^2 |\det g|_\A^{-{\sf w}}\,dg^{\rm Tam},
%\end{align}
%where $dg^{\rm Tam}$ is the Tamagawa measure on $\A^\times\backslash \GL_2(\A)$.
We denote by $M_{\itPi}$ the motive associated to $\itPi$ constructed by Deligne \cite{Deligne1971} and Scholl \cite{Scholl1990} satisfying the properties in Conjecture \ref{C:Clozel}. Under the canonical isomorphism $\Q(\itPi)\otimes_\Q\C \cong \prod_{\sigma : \Q(\itPi)\rightarrow\C}\C$, we write
\[
c^\pm(M_{\itPi}) = \left( c^\pm({}^\sigma\!\itPi)\right)_{\sigma : \Q(\itPi)\rightarrow\C}.
\]
In \cite[Proposition 3.1]{Yoshida1994}, Yoshida explicitly computed the (conjectural) Deligne's periods of the tensor product motive $M_\itSigma \otimes M_{\itPi}$ under the unbalanced condition that $\kappa' > \kappa_1 = \max\{\kappa_i\}$.
%Similar computation can be carried out when $n$ is odd (cf.\,\cite{Bhagwat2014}).
The resulting formula can be expressed in terms of the Gauss sum $G(\omega_\itSigma)$ and the motivic periods $c^\pm(\itPi)$, which is irrelevant to the existence of $M_\itSigma$.
We have the following refinement of Conjecture \ref{C:Deligne} for $M_\itSigma \otimes M_{\itPi}$:

\begin{conj}\label{C:Deligne for GL_n x GL_2}
Assume $\kappa' > \kappa_1$.
Put $\delta \in \{0,1\}$ with $\delta \equiv n \,({\rm mod}\,2)$.
For a critical point $m+\tfrac{\delta}{2} \in \Z+\tfrac{n}{2}$ and $\sigma \in {\rm Aut}(\C)$, we have
\[
\sigma\left(\frac{L^{(\infty)}(m+\tfrac{\delta}{2},\itSigma \times \itPi)}{(2\pi\sqrt{-1})^{nm}\cdot q^{(-1)^m}(\itSigma \times \itPi)} \right) = \frac{L^{(\infty)}(m+\tfrac{\delta}{2},{}^\sigma\!\itSigma \times {}^\sigma\!\itPi)}{(2\pi\sqrt{-1})^{nm}\cdot q^{(-1)^m}({}^\sigma\!\itSigma \times {}^\sigma\!\itPi)}.
\]
%Here
%\begin{align*}
%q^\pm(\itPi \times \itPi') &= (2\pi\sqrt{-1})^{r\delta+n(\kappa'+{\sf w}+{\sf w}')/2}\cdot (\sqrt{-1})^{r{\sf w}'}\cdot G(\omega_\itPi\omega_{\itPi'}^r)\cdot \Vert f_{\itPi'}\Vert^r\\
%&\times  \begin{cases}
%1 & \mbox{ if $n=2r$},\\
% p(\itPi',\pm(-1)^{r+{\sf w}/2}\omega_{\itPi,\infty}(-1)) & \mbox{ if $n=2r+1$}.
%\end{cases}
%\end{align*}
Here
\begin{align*}
q^\pm(\itSigma \times \itPi) &= (2\pi\sqrt{-1})^{\lfloor\frac{n}{2}\rfloor(\delta-1)+n{\sf w}/2}\cdot  G(\omega_\itSigma)\cdot (c^+(\itPi)\cdot c^-(\itPi))^{\lfloor\frac{n}{2}\rfloor}\cdot  \begin{cases}
1 & \mbox{ if $n$ is even},\\
 c^{\pm\varepsilon}(\itPi) & \mbox{ if $n$ is odd},
\end{cases}
\end{align*}
with $\varepsilon = (-1)^{\lfloor\frac{n}{2}\rfloor+{\sf w}/2}\omega_{\itSigma,\infty}(-1)$ if $n$ is odd.
\end{conj}

In \cite{FM2014} and \cite{FM2016}, Furusawa and Morimoto proved algebraicity results for ${\rm SO}(V) \times \GL_2$, which imply Conjecture \ref{C:Deligne for GL_n x GL_2} under the following global and regularity assumptions:
\begin{itemize}
\item  A twist of $\itSigma$ by some integral power of $|\mbox{ }|_\A^{1/2}$ is self-dual and descends to an automorphic representation of ${\rm SO}(V)(\A)$ for some anisotropic quadratic space $V$ over $\Q$.
In particular, $n$ must be even and $\itSigma$ is essentially self-dual.
\item $\kappa'\geq \kappa_1+{\rm dim}_\Q\,V+2$.
\end{itemize}

The following lemma is on the multiplicative property of the period $q^\pm(\itSigma \times \itPi)$.

\begin{lemma}\label{L:mult. period relation}
Assume $\kappa' > \kappa_1$.
Put $\delta \in \{0,1\}$ with $\delta \equiv n \,({\rm mod}\,2)$.
If $\itSigma = \itSigma_1 \boxplus \itSigma_2$ for some isobaric automorphic representations $\itSigma_i$ of $\GL_{n_i}(\A)$ such that $n_1$ is even, then we have
\[
q^\pm(\itSigma \times \itPi) \sim q^\pm(\itSigma_1 \otimes |\mbox{ }|_\A^{\delta/2} \times \itPi) \cdot q^\pm(\itSigma_2 \times \itPi).
\]
\end{lemma}

\begin{proof}
The regularity implies that $\itSigma_i \otimes |\mbox{ }|_\A^{(n-n_i)/2}$ is regular algebraic for $i=1,2$. In particular, $\itSigma_1 \otimes |\mbox{ }|_\A^{\delta/2}$ and $\itSigma_2$ are regular algebraic since $n_1$ is even. 
It is clear the unbalanced condition also holds for $\itSigma_1 \otimes |\mbox{ }|_\A^{\delta/2} \times \itPi$ and $\itSigma_2 \times \itPi$. 
By the essentially temperedness of $\itSigma_\infty$, the infinity types of $\itSigma_1 \otimes |\mbox{ }|_\A^{\delta/2}$ and $\itSigma_2$ are of the form $(\underline{\kappa}^{(1)};\,{\sf w}+\delta)$ and $(\underline{\kappa}^{(2)};\,{\sf w})$ respectively for some $\underline{\kappa}^{(1)}$ and $\underline{\kappa}^{(2)}$. 
Note that $\omega_{\itSigma,\infty}(-1)=(-1)^{n_1/2}\omega_{\itSigma_2,\infty}(-1)$ if $n$ is odd, since $\omega_{\itSigma}=\omega_{\itSigma_1}\omega_{\itSigma_2}$ and ${\sf w}$ is even in this case (cf.\,\S\,\ref{SS:cohomological}).
Also by (\ref{E:Galois Gauss sum}) we have
\[
G(\omega_{\itSigma}) \sim G\left(\omega_{\itSigma_1 \otimes |\mbox{ }|_\A^{\delta/2} }\right)G(\omega_{\itSigma_2}).
\]
The period relation then follows immediately. This completes the proof.
\end{proof}

By Lemma \ref{L:mult. period relation}, to prove Conjecture \ref{C:Deligne for GL_n x GL_2}, it suffices to prove the case when $\itSigma$ is cuspidal.
In Theorem \ref{T:RS for GL_n x GL_2} below, we prove the conjecture under some conditions. Following the strategy outlined in the beginning of \S\,\ref{S:applications}, in the second step we choose some $\itPi'$ of CM-type. 
The conjecture for the auxiliary Rankin--Selberg $L$-function $L(s,\itSigma \times \itPi')$ can be deduced from Theorem \ref{T:Harris GL} and Lemma \ref{L:CM periods} below by the method of base change and automorphic induction.
We recall these results in the following paragraphs.

Let $\K$ be an imaginary quadratic extension of $\Q$.
Let $c\in\Gal(\K/\Q)$ be the non-trivial automorphism. We identify $\K_\infty$ with $\C$ by the natural inclusion $\K\subset\C$.
Let $\chi$ be a Hecke character of $\A_\K^\times$.
Assume $\chi$ is algebraic, that is, there exist $a,b\in\Z$ such that  
\[
\chi_\infty(z) = z^a\overline{z}^b.
\]
%We call $(\kappa;\,{\sf w})$ the infinity type of $\chi$.
In this case, let $\K(\chi)$ be the number field generated over $\K$ by the values of $\chi$ on the finite ideles $\A_{\K,f}^\times$. As in \cite[Appendix]{HK1991} and \cite[\S\,1]{Harris1993}, we denote by 
\[
p(\widecheck{\chi},\tau) \in \C^\times/\K(\chi)^\times
\]
the CM-period of $\widecheck{\chi}:=(\chi^c)^{-1}$ defined with respect to the CM-type $\{\tau\}\subset\Gal(\K/\Q)$. Moreover, we normalize the family of periods $\{p({}^\sigma\!\widecheck{\chi},\sigma\circ\tau)\}_{\sigma \in {\rm Aut}(\C)}$ in a compatible way as in \cite[Lemma 1.3 and Remark 1.3.1]{Harris1993}.
Then, when $a \neq b$, we have
\begin{align}\label{E:Harris GL}
\sigma \left ( \frac{L^{(\infty)}(m,\chi)}{(2\pi\sqrt{-1})^m\cdot p(\widecheck{\chi},\tau_\chi)}\right) =  \frac{L^{(\infty)}(m,{}^\sigma\!\chi)}{(2\pi\sqrt{-1})^m\cdot p({}^\sigma\!\widecheck{\chi},\tau_\chi)}
\end{align}
for all critical points $m \in \Z$ for $L(s,\chi)$ and $\sigma \in {\rm Aut}(\C/\K)$.
Here $\tau_\chi = {\rm id}$ (resp.\,$\tau_\chi = c$) if $a<b$ (resp.\,$a>b$).
More generally, we have the following result on the algebraicity of twisted standard $L$-functions for $\GL_n(\A_\K)$ due to Harris and Guerberoff--Lin.
We only state the result in the unbalanced case and follow the statement in \cite[Theorem 5.13]{GHL2021}.

\begin{thm}[Harris \cite{Harris1997}, \cite{Harris2021}, Guerberoff--Lin \cite{GL2016}]\label{T:Harris GL}
Let $\itPsi$ be a regular algebraic conjugate self-dual cuspidal automorphic representation of $\GL_n(\A_\K)$ with central character $\omega_\itPsi$. Let $\{z^{a_i}\overline{z}^{-a_i}\}_{1 \leq i \leq n}$ be the infinity type of $\itPsi$.
We assume the following assumption is satisfied:
\begin{itemize}
\item There exists a $n$-dimensional non-degenerate definite Hermitian space $V$ over $\K$ such that $\itPsi$ descends to an automorphic representation of $U(V)(\A)$.
\end{itemize}
Let $\chi$ be an algebraic Hecke character of $\A_\K^\times$ with infinity type $z^a\overline{z}^b$ such that $2a_i+a-b > 0$ for all $1 \leq i \leq n$.
For a critical point $m_0 \in \Z+\tfrac{n-1}{2}$ and $\sigma \in {\rm Aut}(\C/\K)$, we have
\begin{align*}
&\sigma \left( \frac{L^{(\infty)}(m_0,\itPsi\otimes\chi)}{(2\pi\sqrt{-1})^{nm_0}\cdot p(\widecheck{\omega}_\itPsi,c)\cdot p(\widecheck{\chi},c)^n}\right) = \frac{L^{(\infty)}(m_0,{}^\sigma\!\itPsi\otimes{}^\sigma\!\chi)}{(2\pi\sqrt{-1})^{nm_0}\cdot p(\widecheck{\omega}_{{}^\sigma\!\itPsi},c)\cdot p({}^\sigma\!\widecheck{\chi},c)^{n}}.
\end{align*}
Here $L(s,\itPsi\otimes\chi)$ is the standard $L$-function of $\itPsi\otimes\chi$. %and $\widecheck{\omega}_\itPsi = ({\omega}_\itPsi^c)^{-1}$, $\widecheck{\chi} = (\chi^c)^{-1}$.
\end{thm}

\begin{rmk}
We have replace the period $P^{(0)}(\itPsi,{\rm id})$ in \cite[Theorem 5.13]{GHL2021} by the CM-period $p(\widecheck{\omega}_\itPsi,c)$. This period relation is stated in \cite[Theorem 2.6]{GHL2021} and proved by Lin in \cite[Lemma 7.6.1]{Lin2015}.
\end{rmk}

\begin{rmk}\label{R:Harris GL}
Let $V$ be a definite Hermitian space over $\K$. By Arthur's multiplicity formula for $U(V)(\A)$ announced by Kaletha--M\'inguez--Shin--White in \cite[Theorem* 1.7.1]{KMSW2014}, the obstruction to descending $\itPsi$ to $U(V)(\A)$ is purely local.
In particular, if either $n$ is odd or $n \equiv 0 \,({\rm mod}\,4)$, then the assumption always holds.
When $n \equiv 2 \,({\rm mod}\,4)$, the assumption holds if we assume further that there exists a prime $p$ such that $\itPsi_p$ is a discrete series representation and $p$ is {\color{black}non-split} in $\K$.
%Let $0 \leq i \leq n$. If $n \equiv 2i \,({\rm mod}\,4)$, then the assertion holds for all algebraic Hecke characters $\chi$ with ${}^\sharp I(\itPsi,\chi)=i$. 
%If $n\equiv 2i+2 \,({\rm mod}\,4)$, then we need to assume the existence of prime $p$ such that $\itPsi_p$ is a discrete series representation and $p$ splits in $\K$.
\end{rmk}

For an algebraic Hecke character $\chi$ of $\A_\K^\times$ with infinity type $z^a\overline{z}^b$, we denote by $I_\K(\chi|\mbox{ }|_{\A_\K}^{1/2})$ the automorphic induction of $\chi|\mbox{ }|_{\A_\K}^{1/2}$ to $\GL_2(\A)$. 
Recall that the central character of $I_\K(\chi|\mbox{ }|_{\A_\K}^{1/2})$ is equal to $\omega_{\K/\Q}(\chi\vert_{\A^\times})|\mbox{ }|_\A$, where $\omega_{\K/\Q}$ is the quadratic Hecke character of $\A^\times$ associated to $\K/\Q$ by class field theory.
If $a \neq b$, then $I_\K(\chi|\mbox{ }|_{\A_\K}^{1/2})$ is regular algebraic and cuspidal with infinity type 
\[
(|a-b|+1;\,a+b+1).
\]
%In this case, we write $f_\chi = f_{I_\K(\chi|\mbox{ }|_{\A_\K}^{1/2})}$ for the normalized newform of $I_\K(\chi|\mbox{ }|_{\A_\K}^{1/2})$.
%In the following lemma, we establish some period relations for cohomological cuspidal automorphic representations of $\GL_2(\A)$ of CM-type. 
In the following lemma, we have a period relation for regular algebraic cuspidal automorphic representations of $\GL_2(\A)$ of CM-type.

\begin{lemma}\label{L:CM periods}
Let $\chi$ be an algebraic Hecke characters of $\A_\K^\times$ with infinity types $z^a\overline{z}^b$ such that $a \neq b$.
For $\sigma \in {\rm Aut}(\C/\K)$, we have
\[
\sigma \left( \frac{c^\pm(I_\K(\chi|\mbox{ }|_{\A_\K}^{1/2}))}{(2\pi\sqrt{-1})\cdot p(\widecheck{\chi},\tau_\chi)}\right) = \frac{c^\pm(I_\K({}^\sigma\!\chi|\mbox{ }|_{\A_\K}^{1/2}))}{(2\pi\sqrt{-1})\cdot p({}^\sigma\!\widecheck{\chi},\tau_\chi)}.
\]
Here $\tau_\chi = {\rm id}$ (resp.\,$\tau_\chi = c$) if $a<b$ (resp.\,$a>b$).
\end{lemma}

%\begin{lemma}
%Let $\chi$ and $\eta$ be algebraic Hecke characters of $\A_\K^\times$ with infinity types $z^a\overline{z}^b$ and $z^{c}\overline{z}^{d}$, respectively. Assume $a\neq b$ %and $c\neq d$.
%The following period relations hold:
%\begin{itemize}
%\item[(1)] For $\sigma \in {\rm Aut}(\C/\K)$, we have
%\[
%\sigma \left( \frac{p(I_\K(\chi|\mbox{ }|_{\A_\K}^{1/2}),\pm)}{(2\pi\sqrt{-1})^{-\max\{a,b\}}\cdot p(\widecheck{\chi},\tau_\chi)}\right) = \frac{p(I_\K({}^\sigma\!\chi|\mbox{ }|_{\A_\K}^{1/2}),\pm)}{(2\pi\sqrt{-1})^{-\max\{a,b\}}\cdot p({}^\sigma\!\widecheck{\chi},\tau_\chi)}.
%\]
%Here $\tau_\chi = {\rm id}$ (resp.\,$\tau_\chi = c$) if $a<b$ (resp.\,$a>b$).
%\item[(2)]
%Assume
%\[
%a-b+2 \geq c-d \geq 2.
%\]
%For $n,m\in\Z$ and $\sigma \in {\rm Aut}(\C)$, we have
%\begin{align*}
%\sigma\left(\frac{\Vert f_{\chi\eta}\Vert^n\cdot \Vert f_{\chi\eta^c}\Vert^m}{\pi^{n-m}\cdot \Vert f_{\chi}\Vert^{n+m}\cdot \Vert f_{\eta}\Vert^{n-m}}\right) = \frac{\Vert f_{{}^\sigma\!\chi{}^\sigma\!\eta}\Vert^n\cdot \Vert f_{{}^\sigma\!\chi{}^\sigma\!\eta^c}\Vert^m}{\pi^{n-m}\cdot \Vert f_{{}^\sigma\!\chi}\Vert^{n+m}\cdot \Vert f_{{}^\sigma\!\eta}\Vert^{n-m}}.
%\end{align*}
%\end{itemize}
%\end{lemma}
\begin{proof}
Note that for any finite order Hecke character $\xi$ of $\A^\times$, we have
\[
L(s,\chi\cdot(\xi\circ{\rm N}_{\K/\Q})) = L(s-\tfrac{1}{2},I_\K(\chi|\mbox{ }|_{\A_\K}^{1/2})\otimes\xi).
\]
On the other hand, by the properties of CM-periods in \cite[Proposition 1.4]{Harris1993} and \cite[(1.10.10)]{Harris1997}, we have
\[
\sigma \left( \frac{p(\widecheck{\chi}\cdot\widecheck{(\xi\circ{\rm N}_{\K/\Q})},\tau)}{G(\xi)\cdot p(\widecheck{\chi},\tau)}\right) = \frac{p({}^\sigma\!\widecheck{\chi}\cdot\widecheck{({}^\sigma\!\xi\circ{\rm N}_{\K/\Q})},\sigma\circ\tau)}{G({}^\sigma\!\xi)\cdot p({}^\sigma\!\widecheck{\chi},\sigma\circ\tau)}
\]
for all $\sigma \in {\rm Aut}(\C)$ and $\tau \in \Gal(\K/\Q)$.
The assertion then follows immediately from (\ref{E:Deligne}) for $M=M_{I_\K(\chi|\mbox{ }|_{\A_\K}^{1/2})}$ and (\ref{E:Harris GL}) by choosing $\xi$ so that all the critical values of $L(s,\chi\cdot(\xi\circ{\rm N}_{\K/\Q}))$ are non-zero (cf.\,\cite{Rohrlich1989}).
\end{proof}

%\begin{rmk}
%By using (1) and the period relations in \cite[Proposition 1.4]{Harris1993} and \cite[(1.10.10)]{Harris1997}, one can prove that (2) holds for $\sigma \in {\rm Aut}(\C/\K)$ assuming only $a-b>c-d>0$. 
%\end{rmk}

We prove Conjecture \ref{C:Deligne for GL_n x GL_2} for essentially self-dual $\itSigma$ under some regularity conditions on the infinity types, and a local condition when $n \equiv 2\,({\rm mod}\,4)$. More precisely: 
\begin{thm}\label{T:RS for GL_n x GL_2}
Conjecture \ref{C:Deligne for GL_n x GL_2} holds %for the Rankin--Selberg $L$-function $L(s,\itPi\times\itPi')$ 
if the following conditions are satisfied:
\begin{itemize}
\item[(1)] $\itSigma$ is essentially self-dual.
\item[(2)] If $n\equiv 2\,({\rm mod}\,4)$, then $\itSigma_p$ is an essentially discrete series representation for some prime $p$.
%\item[(3)] ${\sf w}+{\sf w}' \equiv n\,({\rm mod}\,2)$.
\item[(3)] $\min\{\kappa_i\} \geq 3$ if $n$ is even, $\min\{\kappa_i\} \geq 5$ if $n$ is odd, and 
\begin{align*}
\mbox{$\underline{\kappa}\sqcup\{\kappa'\}$ is}
\begin{cases}
\mbox{$4$-regular} & \mbox{ if ${\sf w}$ or ${\sf w}'$ is odd, and ${\sf w}+{\sf w}' \equiv n\,({\rm mod}\,2)$},\\
\mbox{$5$-regular} & \mbox{ if ${\sf w}+{\sf w}' \equiv n+1\,({\rm mod}\,2)$},\\
\mbox{$6$-regular} & \mbox{ if ${\sf w}$ and ${\sf w}'$ are even.}%, and ${\sf w}+{\sf w}' \equiv n\,({\rm mod}\,2)$}.
\end{cases}
\end{align*}
\end{itemize}
\end{thm}

\begin{proof}
As we have mentioned, by Lemma \ref{L:mult. period relation}, we may assume $\itSigma$ is cuspidal.
%Let $\itSigma$ and $\itPi$ be cohomological cuspidal automorphic representations of $\GL_n(\A)$ and $\GL_2(\A)$ with infinity types $(\underline{\kappa};\,{\sf w})$ and $(\kappa';\,{\sf w}')$, respectively. We assume $\kappa'>\kappa_1$ and the conditions (1)-(4) in Theorem \ref{T:RS for GL_n x GL_2} are satisfied.
%In this section, we prove that Conjecture \ref{C:Deligne for GL_n x GL_2} holds for the Rankin--Selberg $L$-function $L(s,\itSigma \times \itPi)$.
Put $\delta \in \{0,1\}$ with $\delta \equiv n \,({\rm mod}\,2)$ and $r = \lfloor\tfrac{n}{2}\rfloor$.
%We rewrite $q^\pm(\itSigma \times \itPi)$ in terms of the Betti--Whittaker periods of $\itPi$ (cf.\,(\ref{E:period GL_2})):
%\begin{align*}
%q^\pm(\itSigma \times \itPi) &= (2\pi\sqrt{-1})^{r\delta+n({\sf w}+\kappa'+{\sf w}')/2}\cdot (\sqrt{-1})^{r{\sf w}'}\cdot G(\omega_\itSigma)\cdot G(\omega_{\itPi})^r\cdot \Vert f_{\itPi}\Vert^r\\
%&\times  \begin{cases}
%1 & \mbox{ if $n$ is even},\\
% p(\itPi,\pm(-1)^{r+{\sf w}/2}\omega_{\itSigma,\infty}(-1)) & \mbox{ if $n$ is odd}.
%\end{cases}
%\end{align*}
We choose auxiliary regular algebraic automorphic representations $\itSigma'$ and $\itPi'$ of $\GL_n(\A)$ and $\GL_2(\A)$ respectively as follows:
\begin{itemize}
\item For $1 \leq i \leq r$, let $\itPi_i$ be a regular algebraic cuspidal automorphic representation of $\GL_2(\A)$ with infinity type $(\kappa_i;\,{\sf w}+\delta)$. If $n$ is odd, let $\xi$ be a finite order Hecke character of $\A^\times$ such that $\xi_\infty(-1) = (-1)^r\omega_{\itSigma,\infty}(-1)$. We define $\itSigma'$ by
\[
\itSigma' = (\itPi_1\otimes|\mbox{ }|_\A^{-\delta/2}) \boxplus \cdots \boxplus (\itPi_r\otimes|\mbox{ }|_\A^{-\delta/2}) \boxplus \begin{cases}
\emptyset & \mbox{ if $n$ is even},\\
\xi|\mbox{ }|_\A^{{\sf w}/2} & \mbox{ if $n$ is odd}.
\end{cases}
\]
\item Let $\K$ be an imaginary quadratic extension of $\Q$ such that if $n \equiv 2 \,({\rm mod}\,4)$, then there exists a prime number $p$ which is inert in $\K$ and $\itSigma_p$ is an essentially discrete series representation.
We then let $\itPi' = I_\K(\chi'|\mbox{ }|_{\A_\K}^{1/2})$ for some algebraic Hecke character $\chi'$ of $\A_\K^\times$ with infinity type $z^{(\kappa'+{\sf w}'-2)/2}\overline{z}^{(-\kappa'+{\sf w}')/2}$.
\end{itemize}
It is clear that $\itSigma'_\infty$ and $\itPi_\infty'$ are essentially tempered and $\itSigma_\infty = \itSigma_\infty'$ and $\itPi_\infty = \itPi'_\infty$. 
By the result of Harder and Raghuram \cite[Theorem 7.21]{HR2020}, we only need to consider non-central critical points.
Let $m+\tfrac{\delta}{2} \in \Z+\tfrac{n}{2}$ be a non-central critical point for $L(s,\itSigma\times\itPi)$.
Conjecture \ref{C:Deligne for GL_n x GL_2} is known when $n=1$ (cf.\,\cite[\S\,7]{Deligne1979}).
For $n=2$, the conjecture follows from the result of Shimura \cite[Theorem 3]{Shimura1976}.
Therefore, by Lemma \ref{L:mult. period relation}, Conjecture \ref{C:Deligne for GL_n x GL_2} holds for $\itSigma'\times\itPi'$ and $\itSigma'\times\itPi$, that is,
\begin{align}\label{E:RS proof 1}
\begin{split}
L^{(\infty)}(m+\tfrac{\delta}{2},\itSigma'\times\itPi')&\sim (2\pi\sqrt{-1})^{nm}\cdot q^{(-1)^m}(\itSigma'\times\itPi'),\\
L^{(\infty)}(m+\tfrac{\delta}{2},\itSigma'\times\itPi)&\sim (2\pi\sqrt{-1})^{nm}\cdot q^{(-1)^m}(\itSigma'\times\itPi). \end{split}
\end{align}
%Here we have used the following period relation due to Shimura \cite[Theorem 1-(iv)]{Shimura1977}:
%\begin{align}\label{E:RS proof 6}
%\sigma\left(\frac{c^+(\itPi)\cdot c^-(\itPi)}{(2\pi\sqrt{-1})^{1+\kappa+{\sf w}}\cdot(\sqrt{-1})^{{\sf w}}\cdot G(\omega_\itPi)\cdot\Vert f_\itPi \Vert}\right) = \frac{c^+({}^\sigma\!\itPi)\cdot c^-({}^\sigma\!\itPi)}{(2\pi\sqrt{-1})^{1+\kappa+{\sf w}}\cdot(\sqrt{-1})^{{\sf w}}\cdot G({}^\sigma\!\omega_\itPi)\cdot\Vert f_{{}^\sigma\!\itPi} \Vert}.
%\end{align}
By condition (1), we have $\itSigma^\vee = \itSigma \otimes \chi$ for some algebraic Hecke character $\chi$ of $\A^\times$. Fix an algebraic Hecke character $\eta$ of $\A_\K^\times$ such that $\eta_\infty(z) = z^{-{\sf w}}$ and $\eta\vert_{\A^\times} = \chi$ (cf.\,\cite[Lemma 4.1.4]{CHT2008}). Let ${\rm BC}_{\K/\Q}(\itSigma)$ be the base change of $\itSigma$ to $\GL_n(\A_\K)$. Put $\itPsi = {\rm BC}_{\K/\Q}(\itSigma)\otimes\eta$. Then $\itPsi$ is regular algebraic and conjugate self-dual. Moreover, the infinity type of $\itPsi$ is given by $\{z^{a_i}\overline{z}^{-a_i}\}_{1\leq i \leq n}$ with
\[
a_i = \tfrac{\kappa_i-1-{\sf w}}{2},\quad a_i+a_{n-i+1} = -{\sf w}
\]
for $1 \leq i \leq r$ and $a_{r+1} =-\tfrac{\sf w}{2}$ if $n$ is odd.
By the adjointness between automorphic induction and base change, we have
\begin{align*}
L(s,\itSigma \times \itPi') &= L(s+\tfrac{1}{2}, {\rm BC}_{\K/\Q}(\itSigma) \otimes \chi')\\
& = L(s+\tfrac{1}{2},\itPsi \otimes \eta^{-1}\chi').
\end{align*}
We assume further that $\K$ is chosen so that ${\rm BC}_{\K/\Q}(\itSigma)$ is cuspidal.
By condition (2) and the condition that $p$ is inert in $\K$, the assumption in Theorem \ref{T:Harris GL} holds for ${\rm BC}_{\K/\Q}(\itSigma)$ (cf.\,Remark \ref{R:Harris GL}). 
It then follows from the assumption $\kappa'>\kappa_1$ and Theorem \ref{T:Harris GL} that
\[
L^{(\infty)}(m+\tfrac{\delta}{2},\itSigma \times \itPi') \sim_\K (2\pi\sqrt{-1})^{n(m+(\delta+1)/2)}\cdot p(\widecheck{\omega}_\itPsi,c)\cdot p(\widecheck{\eta}^{-1}\widecheck{\chi}',c)^n.
\]
Here $\sim_\K$ means the ratio is equivariant under ${\rm Aut}(\C/\K)$ (recall the notation in \S\,\ref{SS:notation}, here we replace ${\rm Aut}(\C)$ by ${\rm Aut}(\C/\K)$).
By Lemma \ref{L:CM periods}, we have
\begin{align*}
p(\check{\chi}',c)^n &\sim_\K (2\pi\sqrt{-1})^{-n}\cdot c^\pm(\itPi')^n.
\end{align*}
%\begin{align*}
%p(\check{\chi}',c)^n &\sim_\K (2\pi\sqrt{-1})^{r+n(\kappa'+{\sf w}'-2)/2}\cdot (\sqrt{-1})^{r{\sf w}'}\cdot G(\omega_{\itPi'})^{r}\cdot \Vert f_{\itPi'} \Vert^r\\
%& \times \begin{cases}
%1 & \mbox{ if $n$ is even},\\
%(2\pi\sqrt{-1})^{-(\kappa'+{\sf w}')/2}\cdot c^\pm(\itPi') & \mbox{ if $n$ is odd}.
%\end{cases}
%\end{align*}
By \cite[Proposition 1.4-(b)]{Harris1993}, for algebraic Hecke characters $\chi_1,\chi_2$ of $\A_\K^\times$, we have
\[
p(\chi_1\chi_2,\tau) \sim_\K p(\chi_1,\tau)\cdot p(\chi_2,\tau),\quad \tau \in \Gal(\K/\Q).
\]
In particular, since $\omega_\itPsi = (\omega_\itSigma\circ {\rm N}_{\K/\Q})\cdot \eta^n$, we have
\begin{align*}
p(\widecheck{\omega}_\itPsi,c)\cdot p(\widecheck{\eta}^{-1},c)^n &\sim_\K p(\widecheck{(\omega_\itSigma\circ {\rm N}_{\K/\Q})},c)\\
&\sim_\K (2\pi\sqrt{-1})^{n{\sf w}/2}\cdot G(\omega_\itSigma).
\end{align*}
Here we have used \cite[(1.10.9) and (1.10.10)]{Harris1997}.
We thus conclude that
\begin{align}\label{E:RS proof 2}
L^{(\infty)}(m+\tfrac{\delta}{2},\itSigma\times\itPi')&\sim _\K(2\pi\sqrt{-1})^{nm}\cdot q^{(-1)^m}(\itSigma\times\itPi').
\end{align}
Note that condition (3) is the regularity assumption on $\itSigma_\infty$ and $\itPi_\infty$ in Theorem \ref{T:Main}. Therefore, by (\ref{E:RS proof 1}), (\ref{E:RS proof 2}), and Theorem \ref{T:Main}, we have
\begin{align}\label{E:RS proof 3}
\begin{split}
L^{(\infty)}(m+\tfrac{\delta}{2},\itSigma \times \itPi) &\sim_\K (2\pi\sqrt{-1})^{nm}\cdot \frac{q^{(-1)^m}(\itSigma\times\itPi')\cdot q^{(-1)^m}(\itSigma'\times\itPi)}{q^{(-1)^m}(\itSigma'\times\itPi')}\\
& = (2\pi\sqrt{-1})^{nm}\cdot q^{(-1)^m}(\itSigma\times\itPi).
\end{split}
\end{align}
We take another imaginary quadratic extension $\K'$ of $\Q$ satisfying the same assumptions as above such that $\K\cdot \Q(\itSigma,\itPi) \cap \K'\cdot \Q(\itSigma,\itPi) = \Q(\itSigma,\itPi)$. Then (\ref{E:RS proof 3}) also holds for $\K'$. In particular, we see that
\begin{align}\label{E:RS proof 4}
\frac{L^{(\infty)}(m+\tfrac{\delta}{2},\itSigma \times \itPi)}{(2\pi\sqrt{-1})^{nm}\cdot q^{(-1)^m}(\itSigma\times\itPi)} \in \Q(\itSigma,\itPi).
\end{align}
Assume $\K \cap \Q(\itSigma,\itPi) = \Q$.
Let $\sigma \in {\rm Aut}(\C)$. Write $\sigma\vert_{\K\Q(\itSigma,\itPi)^{\rm Gal}} = \sigma_1\circ \sigma_2$ for some $\sigma_1 \in \Gal(\K\Q(\itSigma,\itPi)^{\rm Gal}/\K)$ and $\sigma_2 \in \Gal(\K\Q(\itSigma,\itPi)^{\rm Gal}/\Q(\itSigma,\itPi))$, where $\Q(\itSigma,\itPi)^{\rm Gal}$ is the Galois closure of $\Q(\itSigma,\itPi)$ over $\Q$.
Then we have
\begin{align}\label{E:RS proof 5}
\begin{split}
\sigma \left( \frac{L^{(\infty)}(m+\tfrac{\delta}{2},\itSigma \times \itPi)}{(2\pi\sqrt{-1})^{nm}\cdot q^{(-1)^m}(\itSigma\times\itPi)}\right) &= \sigma_1\left( \frac{L^{(\infty)}(m+\tfrac{\delta}{2},\itSigma \times \itPi)}{(2\pi\sqrt{-1})^{nm}\cdot q^{(-1)^m}(\itSigma\times\itPi)}\right)\\
& = \frac{L^{(\infty)}(m+\tfrac{\delta}{2},{}^{\sigma_1}\!\itSigma \times {}^{\sigma_1}\!\itPi)}{(2\pi\sqrt{-1})^{nm}\cdot q^{(-1)^m}({}^{\sigma_1}\!\itSigma\times{}^{\sigma_1}\!\itPi)}\\
& = \frac{L^{(\infty)}(m+\tfrac{\delta}{2},{}^{\sigma}\!\itSigma \times {}^{\sigma}\!\itPi)}{(2\pi\sqrt{-1})^{nm}\cdot q^{(-1)^m}({}^{\sigma}\!\itSigma\times{}^{\sigma}\!\itPi)}.
\end{split}
\end{align}
Here the first and second equalities follow from (\ref{E:RS proof 4}) and (\ref{E:RS proof 3}), respectively.
This completes the proof.
\end{proof}
%\begin{rmk}
%If we assume the validity of Conjecture \ref{C:main 1}, then we can prove Conjecture \ref{C:Deligne for GL_n x GL_2} under assumptions (1) and (2), and the regularity condition that $\underline{\kappa}\sqcup \kappa'$ is $2$-regular.
%\end{rmk}

%\subsection{Proof of Theorem \ref{T:tensor product for GL_2}}
\subsection{Tensor product $L$-functions for $\GL_2$}\label{SS:Blasius}

Let $\itPi_1,...,\itPi_n$ be regular algebraic cuspidal automorphic representations of $\GL_2(\A)$.
For $1 \leq i \leq n$, let $(\kappa_i;\,{\sf w}_i)$ be the infinity type of $\itPi_i$.
Let $f_{\itPi_i}$ be the normalized newform of $\itPi_i$, and $\< f_{\itPi_i}, f_{\itPi_i}\>$ be its Petersson norm defined by
\[
\< f_{\itPi_i}, f_{\itPi_i}\> = \int_{\A^\times\GL_2(\Q)\backslash\GL_2(\A)}|f_{\itPi_i}(g)|^2 |\det g|_\A^{-{\sf w}_i}\,dg^{\rm Tam},
\]
where $dg^{\rm Tam}$ is the Tamagawa measure on $\A^\times\backslash \GL_2(\A)$.
We consider the automorphic tensor product $\itPi_1 \boxtimes \cdots \boxtimes \itPi_n$ of $\itPi_1,...,\itPi_n$ defined in \S\,\ref{SS:5.1}.
%It is an irreducible admissible $((\frak{g}_{2^n},{\rm O}_{2^n}(\R))\times\GL_{2^n}(\A_f))$-module defined by the restricted tensor product 
%\[
%\itPi_1 \boxtimes \cdots \boxtimes \itPi_n = \bigotimes_v \itPi_{1,v}\boxtimes \cdots \boxtimes\itPi_{n,v}.
%\]
%Here, for each place $v$ of $\Q$, $\itPi_{1,v}\boxtimes \cdots \boxtimes\itPi_{n,v}$ is the irreducible admissible representation of $\GL_{2^n}(\Q_v)$ (resp.\,$(\frak{g}_{2^n},{\rm O}_{2^n}(\R))$-module) if $v$ is finite (resp.\,archimedean) defined by the local Langlands correspondence with respect to the tensor representation
%\[
%\GL_2(\C) \times \cdots \times \GL_2(\C) \longrightarrow \GL((\C^2)^{\otimes n}).
%\]
%As a special instance of Langlands functoriality conjecture, we expect $\itPi_1 \boxtimes \cdots \boxtimes \itPi_n$ to be automorphic. 
For $n=2$, the functoriality of the automorphic tensor product was proved by Ramakrishnan \cite{Rama2000}. For $n \geq 3$, Dieulefait \cite{Dieulefait2020} proved the functoriality under the assumptions that $\itPi_i$ has level $1$ for all $1 \leq i \leq n$ and the infinitesimal character of $\itPi_{1,\infty}\boxtimes \cdots \boxtimes \itPi_{n,\infty}$ is regular.
%Consider the tensor product $L$-function defined by the Euler product
%\[
%L(s,\itPi_1 \times\cdots\times\itPi_n) = \prod_v L(s,\itPi_{1,v}\boxtimes \cdots \boxtimes \itPi_{n,v}),
%\]
%which converges absolutely for ${\rm Re}(s)>1-\tfrac{1}{2}\sum_{i=1}^n{\sf w}_i$.
%A critical point for $L(s,\itPi_1 \times\cdots\times\itPi_n)$ is a half-integer in $\Z+\tfrac{n}{2}$ which is not a pole of the archimedean local factors $L(s,\itPi_{1,\infty} \boxtimes\cdots\boxtimes\itPi_{n,\infty})$ and $L(1-s,\itPi_{1,\infty}^\vee \boxtimes\cdots\boxtimes\itPi_{n,\infty}^\vee)$.
Note that the tensor product $L$-function $L(s,\itPi_1 \times\cdots\times\itPi_n)$ admits critical points if and only if
\begin{align}\label{E:critical tensor}
\sum_{i=1}^n\varepsilon_i(\kappa_i-1) \neq 0
\end{align}
for all $\underline{\varepsilon} = (\varepsilon_1,...,\varepsilon_n) \in \{\pm1\}^n$.
Moreover, we have the following lemma on the archimedean component.
\begin{lemma}\label{L:archimedean tensor}
Assume (\ref{E:critical tensor}) holds and $|\sum_{i=1}^n(\varepsilon_i-\varepsilon_i')(\kappa_i-1)| \neq 0$ for all $(\varepsilon_1,...,\varepsilon_n)\neq(\varepsilon_1',...,\varepsilon_n')$ in $\{\pm1\}^n$.
Then $(\itPi_{1,\infty}\boxtimes \cdots \boxtimes \itPi_{n,\infty})\otimes|\mbox{ }|^{(1+\delta)/2} \in \Omega(2^n)$ and is essentially tempered with infinity type
\begin{align*}%\label{E:infinity type tensor}
\left( \left(\left|\sum_{i=1}^{n}\varepsilon_i(\kappa_i-1)\right|+1\right)_{\underline{\varepsilon} \in \{\pm1\}^n/\{\pm1\}};\,\sum_{i=1}^{n}{\sf w}_i+1+\delta\right),
\end{align*}
where $\delta \in \{0,1\}$ with $\delta \equiv n \,({\rm mod}\,2)$. 
\end{lemma}

\begin{proof}
For $\kappa \geq 2$, let $\phi_\kappa : W_\R \rightarrow \GL_2(\C)$ be the irreducible representation of the Weil group $W_\R$ of $\R$ associated to the discrete series representation $D_\kappa$ by the local Langlands correspondence. For $\kappa > \ell \geq 2$, we have
\[
\phi_\kappa \otimes \phi_\ell \cong \phi_{\kappa+\ell-1} \oplus \phi_{\kappa-\ell+1}.
\]
Note that the Langlands parameter of $\itPi_{1,\infty}\boxtimes \cdots \boxtimes \itPi_{n,\infty}$ is the tensor representation
\[
(\phi_{\kappa_1} \otimes \cdots \otimes \phi_{\kappa_n})\otimes |\mbox{ }|^{\sum_{i=1}^n{\sf w}_i/2}.
\]
The assertion then follows by mathematical induction on $n$.
\end{proof}

In \cite{Blasius1987}, when $n \geq 2$, Blasius explicitly computed Deligne's period of the tensor product motive $M_{\itPi_1}\otimes \cdots \otimes M_{\itPi_n}$ in terms of products of Deligne's periods of $M_{\itPi_1},...,M_{\itPi_n}$ which are related to the Petersson norms $\< f_{\itPi_1}, f_{\itPi_1}\>,...,\< f_{\itPi_n}, f_{\itPi_n}\>$ by the period relation due to Shimura \cite[Theorem 1-(iv)]{Shimura1977}. 
More precisely, we have the following refinement of Conjecture \ref{C:Deligne} for $M_{\itPi_1}\otimes \cdots \otimes M_{\itPi_n}$:

\begin{conj}[Blasius]\label{C:Blasius}
Assume $n \geq 2$. % and (\ref{E:critical tensor}) holds.
Put $\delta \in \{0,1\}$ with $\delta \equiv n \,({\rm mod}\,2)$.
The tensor product $L$-function admits meromorphic continuation to $s\in\C$ and is entire if it admits critical points. For a critical point $m+\tfrac{\delta}{2} \in \Z+\tfrac{n}{2}$ and $\sigma \in {\rm Aut}(\C)$, we have
\begin{align*}
\sigma \left( \frac{L^{(\infty)}(m+\tfrac{\delta}{2},\itPi_1\times\cdots\times\itPi_n)}{(2\pi\sqrt{-1})^{2^{n-1}m}\cdot q(\itPi_1\times\cdots\times\itPi_n)}\right) = \frac{L^{(\infty)}(m+\tfrac{\delta}{2},{}^\sigma\!\itPi_1\times\cdots\times{}^\sigma\!\itPi_n)}{(2\pi\sqrt{-1})^{2^{n-1}m}\cdot q({}^\sigma\!\itPi_1\times\cdots\times{}^\sigma\!\itPi_n)}.
\end{align*}
Here
\[
q(\itPi_1\times\cdots\times\itPi_n) = (2\pi\sqrt{-1})^{2^{n-2}(\delta+\sum_{i=1}^n{\sf w}_i)}\cdot \prod_{i=1}^n G(\omega_{\itPi_i})^{2^{n-2}}\cdot (\pi^{\kappa_i} \cdot \<f_{\itPi_i},f_{\itPi_i}\>)^{2^{n-2}-t_i}
\]
with $t_i$ equal to the cardinality of the set
\begin{align}\label{E:cardinality}
\left\{(\varepsilon_1,...,\varepsilon_n)\in \{\pm1\}^n\,\left\vert\, 2(\kappa_i-1)<\sum_{j=1}^n \varepsilon_j(\kappa_j-1),\quad \varepsilon_i=1\right\}\right. ,\quad 1 \leq i \leq n.
\end{align}
%for $1 \leq i \leq n$.
\end{conj}

For $n=2$, the conjecture was proved by Shimura in \cite[Theorem 3]{Shimura1976}. When $n=3$, there are two cases according to $\kappa_1+\kappa_2+\kappa_3 \geq 2\max\{\kappa_1,\kappa_2,\kappa_3\}$ or $2\max\{\kappa_1,\kappa_2,\kappa_3\} > \kappa_1+\kappa_2+\kappa_3$. The former (resp.\,later) case is called the balanced (resp.\,unbalanced) case.
In the balanced case, the conjecture was proved by Garrett and Harris in \cite{GH1993} under some conditions which were lifted by the author in \cite{Chen2021d}.
We refer to the introduction of \cite{Chen2021d} for a survey on known results in this case.
In the unbalanced case, the conjecture is partially proved. In \cite{HK1991}, Harris and Kudla prove the case for central critical values when $\omega_{\itPi_1}\omega_{\itPi_2}\omega_{\itPi_3}=1$. In \cite{FM2014} and \cite{FM2016}, as an application to their main result, Furusawa and Morimoto proved the conjecture under some local and global assumptions.
By Theorem \ref{T:RS for GL_n x GL_2} applied to $\itSigma = (\itPi_2 \boxtimes \itPi_3)\otimes|\mbox{ }|_\A^{-1/2}$ and $\itPi=\itPi_1$, in the notation therein we have
\[
\underline{\kappa} = (\kappa_2+\kappa_3-1,|\kappa_2-\kappa_3|+1),\quad \kappa'=\kappa_1.
\]
Thus the conjecture for $n=3$ in the unbalanced case holds under the following condition:
\begin{align}\label{E:unbalanced}
\begin{cases}
\mbox{$\kappa_1 \geq \kappa_2+\kappa_3+3 \geq |\kappa_2-\kappa_3|+9 \geq 11$} & \mbox{ if ${\sf w}_1$ is odd and ${\sf w}_2+{\sf w}_3$ is even},\\
\mbox{$\kappa_1 \geq \kappa_2+\kappa_3+4 \geq |\kappa_2-\kappa_3|+11 \geq 13$} & \mbox{ if ${\sf w}_1+{\sf w}_2+{\sf w}_3$ is even},\\
\mbox{$\kappa_1 \geq \kappa_2+\kappa_3+5 \geq |\kappa_2-\kappa_3|+13 \geq 15$} & \mbox{ if ${\sf w}_1$ is even and ${\sf w}_2+{\sf w}_3$ is odd}.
\end{cases}
\end{align}
For arbitrary $n \geq 2$, we prove the conjecture under a regularity condition on the infinity types, together with an assumption on functoriality if $n\geq 5$. More precisely: 

\begin{thm}\label{T:tensor product for GL_2}
 %and (\ref{E:critical tensor}) holds.
Conjecture \ref{C:Blasius} holds %for the tensor product $L$-function $L(s,\itPi_1 \times \cdots\times\itPi_n)$ 
if the following conditions are satisfied:
\begin{itemize}
%\item[(1)] ${\sf w}_1+\cdots+{\sf w}_n \equiv n\,({\rm mod}\,2)$.
\item[(1)] $|\sum_{i=1}^n\varepsilon_i(\kappa_i-1)| \geq 6$ and 
$|\sum_{i=1}^n(\varepsilon_i-\varepsilon_i')(\kappa_i-1)| \geq 6$ for all $(\varepsilon_1,...,\varepsilon_n)\neq(\varepsilon_1',...,\varepsilon_n')$ in $\{\pm1\}^n$. %with $(\varepsilon_1,\cdots,\varepsilon_n)\neq \pm(\varepsilon_1',\cdots,\varepsilon_n')$.
\item[(2)] $\itPi_1\boxtimes \cdots\boxtimes \itPi_k$ is automorphic for all $1 \leq k \leq n-2$.
\end{itemize}
In particular, condition (2) holds if either $n \leq 4$ or $\itPi_1,...,\itPi_{n-2}$ have level $1$.
\end{thm}

\begin{proof}
%Let $\itPi_1,\cdots,\itPi_n$ be cohomological cuspidal automorphic representations of $\GL_2(\A)$ with $n\geq 2$. In this section, we show that Conjecture \ref{C:Blasius} holds for the tensor product $L$-function $L(s,\itPi_1 \times \cdots \times \itPi_n)$ when conditions (1)-(3) in Theorem \ref{T:tensor product for GL_2} are satisfied. 
First note that condition (2) implies that $L(s,\itPi_1 \times \cdots \times \itPi_n)$ is entire. Indeed, by (\ref{E:critical tensor}), the cuspidal summands of $\itPi_1 \boxtimes \cdots \boxtimes \itPi_{n-2}$ and $\itPi_{n-1}^\vee\boxtimes\itPi_n^\vee$ must be non-isomorphic. 
We will prove the conjecture by induction on $n$. 
It holds unconditionally when $n=2$ or $n=3$ in the balanced case, and holds for $n=3$ under (\ref{E:unbalanced}) in the unbalanced case (thus under condition (1) here). 
Now let $n \geq 4$ and assume the conjecture holds for $n-1$.
Put $\delta \in \{0,1\}$ with $\delta \equiv n \,({\rm mod}\,2)$. 
We choose some auxiliary automorphic representations $\itPi(\underline{\varepsilon})$, $\itPi_{n-1}'$, $\itPi_{n-1}''$ as follows:
\begin{itemize}
\item For $\underline{\varepsilon} \in \{\pm1\}^{n-2}/\{\pm1\}$, let $\itPi(\underline{\varepsilon})$ be a regular algebraic cuspidal automorphic representation of $\GL_2(\A)$ with infinity type 
\[
\left(\left|\sum_{i=1}^{n-2}\varepsilon_i(\kappa_i-1)\right|+1;\,\sum_{i=1}^{n-2}{\sf w}_i+1+\delta\right).
\]
\item Let $\itPi_{n-1}'$ and $\itPi_{n-1}''$ be regular algebraic cuspidal automorphic representations of $\GL_2(\A)$ with infinity types
\[
(\kappa_{n-1}+\kappa_n-1;\,{\sf w}_{n-1}+{\sf w}_n-1), \quad (|\kappa_{n-1}-\kappa_n|+1;\,{\sf w}_{n-1}+{\sf w}_n-1)
\]
respectively. 
\end{itemize}
Let $\itSigma$ and $\itSigma'$ be regular algebraic automorphic representations of $\GL_{2^{n-2}}(\A)$ defined by
\begin{align*}
\itSigma = (\itPi_1\boxtimes\cdots\boxtimes\itPi_{n-2})\otimes |\mbox{ }|_\A^{(1+\delta)/2},\quad
\itSigma'  = \bigboxplus_{\underline{\varepsilon}\in \{\pm1\}^{n-2}/\{\pm1\}}\itPi(\underline{\varepsilon}).
\end{align*}
Let $\itPi$ and $\itPi'$ be regular algebraic automorphic representations of $\GL_4(\A)$ defined by
\begin{align*}
\itPi = (\itPi_{n-1}\boxtimes\itPi_n)\otimes |\mbox{ }|_\A^{-1/2},\quad
\itPi'  = \itPi_{n-1}' \boxplus \itPi_{n-1}''.
\end{align*}
Then $\itSigma_\infty$ and $\itPi_\infty$ are essentially tempered with $\itSigma_\infty = \itSigma_\infty'$ and $\itPi_\infty = \itPi'_\infty$ by Lemma \ref{L:archimedean tensor}. 
By \cite[Theorem 7.21]{HR2020}, we only need to consider non-central critical points.
By (\ref{E:distance}) and Lemma \ref{L:archimedean tensor}, condition (1) implies that $L(s,\itPi_1 \times \cdots \times \itPi_n)$ admits non-central critical points.
Let $m+\tfrac{\delta}{2} \in \Z+\tfrac{n}{2}$ be a non-central critical point.
It is clear that we have the following equality and factorization:
\begin{align*}
L(s,\itSigma \times \itPi) & = L(s+\tfrac{\delta}{2},\itPi_1 \times \cdots \times \itPi_n),\\
L(s,\itSigma \times \itPi') & = L(s+\tfrac{1+\delta}{2},\itPi_1 \times \cdots \times \itPi_{n-2}\times\itPi_{n-1}') \cdot L(s+\tfrac{1+\delta}{2},\itPi_1 \times \cdots \times \itPi_{n-2}\times\itPi_{n-1}'').
\end{align*}
Recall the sequence $(t_1,...,t_n) \in \Z^n$ defined in (\ref{E:cardinality}) with respect to $(\kappa_1,...,\kappa_n)$. Similarly, let $(t_1',...,t_{n-1}')$ and $(t_1'',...,t_{n-1}'')$ in $\Z^{n-1}$ defined with respect to $(\kappa_1,...,\kappa_{n-2},\kappa_{n-1}+\kappa_n-1)$ and $(\kappa_1,...,\kappa_{n-2},|\kappa_{n-1}-\kappa_n|+1)$, respectively.
{We have assumed that the tuple $(\itPi_1,\ldots,\itPi_n)$ satisfies conditions (1) and (2).
It follows immediately that conditions (1) and (2) are also satisfied by the tuples $(\itPi_1,\ldots,\itPi_{n-2},\itPi_{n-1}')$ and $(\itPi_1,\ldots,\itPi_{n-2},\itPi_{n-1}'')$.}
Therefore, by the induction hypothesis, we have 
\begin{align}\label{E:tensor proof 1}
\begin{split}
&L^{(\infty)}(m,\itSigma\times\itPi')\\
&\sim(2\pi\sqrt{-1})^{2^{n-1}m+2^{n-2}(\delta+\sum_{i=1}^n{\sf w}_i)}\cdot \prod_{i=1}^{n-2} G(\omega_{\itPi_i})^{2^{n-2}}\cdot (\pi^{\kappa_i}\cdot \<f_{\itPi_i},f_{\itPi_i}\>)^{2^{n-2}-t_i'-t_i''}\\
&\times G(\omega_{\itPi_{n-1}'})^{2^{n-3}}\cdot(\pi^{\kappa_{n-1}+\kappa_n-1}\cdot \<f_{\itPi_{n-1}'},f_{\itPi_{n-1}'}\>)^{2^{n-3}-t_{n-1}'}\\
&\times G(\omega_{\itPi_{n-1}''})^{2^{n-3}}\cdot  (\pi^{|\kappa_{n-1}-\kappa_n|+1}\cdot \<f_{\itPi_{n-1}''},f_{\itPi_{n-1}''}\>)^{2^{n-3}-t_{n-1}''}.
\end{split}
\end{align}
On the other hand, we have the following factorizations:
\begin{align*}
L(s,\itSigma' \times \itPi') &= \prod_{\underline{\varepsilon}\in \{\pm1\}^{n-2}/\{\pm1\}} L(s,\itPi(\underline{\varepsilon}) \times \itPi_{n-1}')\cdot L(s,\itPi(\underline{\varepsilon}) \times \itPi_{n-1}''),\\
L(s,\itSigma' \times \itPi) &= \prod_{\underline{\varepsilon}\in \{\pm1\}^{n-2}/\{\pm1\}}L(s-\tfrac{1}{2},\itPi(\underline{\varepsilon}) \times \itPi_{n-1} \times \itPi_n).
\end{align*}
For each $\underline{\varepsilon}$, it is clear that condition (1) holds for the triple product $L$-function $L(s,\itPi(\underline{\varepsilon}) \times \itPi_{n-1} \times \itPi_n)$. Therefore, Conjecture \ref{C:Blasius} holds for these $L$-functions when $n=2,3$ and we have
\begin{align}\label{E:tensor proof 2}
\begin{split}
&\frac{L^{(\infty)}(m,\itSigma' \times \itPi')}{L^{(\infty)}(m,\itSigma'\times \itPi)}\\
&\sim 
\frac{G(\omega_{\itPi_{n-1}'})^{2^{n-3}}\cdot G(\omega_{\itPi_{n-1}''})^{2^{n-3}}}{G(\omega_{\itPi_{n-1}})^{2^{n-2}}\cdot G(\omega_{\itPi_n})^{2^{n-2}}}\\
&\times\prod_{ \underline{\varepsilon}\in \{\pm1\}^{n-2}/\{\pm1\}} \frac{(\pi^{\kappa_{n-1}+\kappa_n-1}\cdot \< f_{\itPi_{n-1}'},f_{\itPi_{n-1}'} \>)^{1-t_{n-1}'(\underline{\varepsilon})}\cdot (\pi^{|\kappa_{n-1}-\kappa_n|+1}\cdot \< f_{\itPi_{n-1}''},f_{\itPi_{n-1}''} \>)^{1-t_{n-1}''(\underline{\varepsilon})}}{(\pi^{\kappa_{n-1}}\cdot \< f_{\itPi_{n-1}},f_{\itPi_{n-1}} \>)^{2-t_{n-1}(\underline{\varepsilon})}\cdot (\pi^{\kappa_{n}}\cdot \< f_{\itPi_{n}},f_{\itPi_{n}} \>)^{2-t_n(\underline{\varepsilon})}}.
\end{split}
\end{align}
Here
\begin{align*}
t_{n-1}(\underline{\varepsilon})  &= {}^\sharp\!\left\{ (\varepsilon',\varepsilon_n)\in\{\pm1\}^2 \, \left\vert \, \kappa_{n-1}-1 < \varepsilon'\left|\sum_{i=1}^{n-2}\varepsilon_i(\kappa_i-1)\right| + \varepsilon_n(\kappa_n-1)\right\}\right.,\\
t_n(\underline{\varepsilon})  &= {}^\sharp\!\left\{ (\varepsilon',\varepsilon_{n-1})\in\{\pm1\}^2 \, \left\vert \, \kappa_{n}-1 < \varepsilon'\left|\sum_{i=1}^{n-2}\varepsilon_i(\kappa_i-1)\right| + \varepsilon_{n-1}(\kappa_{n-1}-1)\right\}\right.,\\
t_{n-1}'(\underline{\varepsilon}) & = {}^\sharp\!\left\{ \varepsilon' \in\{\pm1\}\,\left\vert\, \kappa_{n-1}+\kappa_n-2 < \varepsilon'\left|\sum_{i=1}^{n-2}\varepsilon_i(\kappa_i-1)\right| \right\}\right.,\\
t_{n-1}''(\underline{\varepsilon}) & = {}^\sharp\!\left\{ \varepsilon' \in\{\pm1\}\,\left\vert\, |\kappa_{n-1}-\kappa_{n-1}| < \varepsilon'\left|\sum_{i=1}^{n-2}\varepsilon_i(\kappa_i-1)\right| \right\}\right..
\end{align*}
Finally, it is clear from definition that
\begin{align*}
\sum_{\underline{\varepsilon}\in \{\pm1\}^{n-2}/\{\pm1\}}t_{n-1}(\underline{\varepsilon}) &= t_{n-1},\quad \sum_{\underline{\varepsilon}\in \{\pm1\}^{n-2}/\{\pm1\}}t_{n}(\underline{\varepsilon}) = t_{n},\\
\sum_{\underline{\varepsilon}\in \{\pm1\}^{n-2}/\{\pm1\}}t_{n-1}'(\underline{\varepsilon}) &= t_{n-1}',\quad \sum_{\underline{\varepsilon}\in \{\pm1\}^{n-2}/\{\pm1\}}t_{n-1}'(\underline{\varepsilon}) = t_{n-1}''.
\end{align*}
{Also, it is straightforward to verify that the regularity condition (1) implies the regularity conditions \eqref{E:regularity 1} and \eqref{E:regularity 2} on $\itSigma_\infty$ and $\itPi_\infty$, by Lemma \ref{L:archimedean tensor}.}
We thus conclude from (\ref{E:tensor proof 1}), (\ref{E:tensor proof 2}), and Theorem \ref{T:Main} that Conjecture \ref{C:Blasius} holds for $L^{(\infty)}(m+\tfrac{\delta}{2},\itPi_1\times\cdots\times\itPi_n)$.
This completes the proof.
\end{proof}

%\begin{rmk}
%If we assume the validity of Conjecture \ref{C:main 1}, then we can prove that Conjecture \ref{C:Blasius} holds for all non-central critical points under assumption (3).
% and the regularity condition $|\sum_{i=1}^n(\varepsilon_i-\varepsilon_i')(\kappa_i-1)| \geq 2$ for all $(\varepsilon_1,\cdots,\varepsilon_n)$ and $(\varepsilon_1',\cdots,\varepsilon_n')$ in $\{\pm1\}^n$ with $(\varepsilon_1,\cdots,\varepsilon_n) \neq \pm(\varepsilon_1',\cdots,\varepsilon_n')$. 
%\end{rmk}

%\subsection{Proof of Theorem \ref{T:Sym odd}}
\subsection{Symmetric power $L$-functions for $\GL_2$}\label{SS:Deligne Sym}

Let $\itPi$ be a regular algebraic cuspidal automorphic representation of $\GL_2(\A)$ with infinity type $(\kappa;\,{\sf w})$.
For $n \geq 1$, let ${\rm Sym}^n\itPi$ be the functorial lift of $\itPi$ to $\GL_{n+1}(\A)$ with respect to the symmetric $n$-th power representation of $\GL_2(\C)$.
The existence of the lifts was proved recently by Newton and Thorne in \cite{NT2021} and \cite{NT2021b} (we refer to \textit{loc.\ cit.} for a survey of known cases for $n\leq 8$).
Note that ${\rm Sym}^n\itPi$ is regular algebraic with essentially tempered archimedean component (cf.\,\cite[Theorem 5.3]{Raghuram2009}) and infinity type
\begin{align}\label{E:infinity type sym}
\left(\left((n-2i)(\kappa-1)+1\right)_{0 \leq i\leq {\lfloor \frac{n-1}{2}\rfloor}};\,n{\sf w}\right).
\end{align}
Moreover, it is cuspidal for all $n$ if $\itPi$ is not of CM-type. In \cite[Proposition 7.7]{Deligne1979}, Deligne explicitly computed the motivic periods of ${\rm Sym}^nM_\itPi$ in terms of $\delta(M_\itPi)$ and $c^\pm(M_\itPi)$.
In \cite{Blasius1997}, Blasius investigated the behavior of Deligne's periods upon twisting by Artin motives.
More precisely, for a finite order Hecke character $\chi$ of $\A^\times$ and a pure motive $M$ over $\Q$ with coefficients in a number field $\mathbb{E}$,
by \cite[M.17, Corollary 4]{Blasius1997} we have
\[
c^\pm(M\otimes{\rm Art}_\chi) \in \left(G({}^\sigma\!\chi)^{d^{\pm\chi_\infty(-1)}(M)}\right)_{\sigma:\E\rightarrow\C} \cdot c^{\pm\chi_\infty(-1)}(M)\cdot \E^\times.
\]
In particular, we have the following refinement of Conjecture \ref{C:Deligne} for the motive ${\rm Sym}^nM_\itPi \otimes {\rm Art}_\chi$: %where ${\rm Art}_\chi$ denotes the Artin motive of a finite order Hecke character $\chi$:

\begin{conj}\label{C:Deligne Sym}
Put $\delta \in \{0,1\}$ with $\delta \equiv n \,({\rm mod}\,2)$ and $r = \lfloor\tfrac{n}{2}\rfloor$.
Let $\chi$ be a finite order Hecke character of $\A^\times$. For a critical point $m+\tfrac{\delta}{2} \in \Z+\tfrac{n}{2}$ for the twisted symmetric $n$-th power $L$-function $L(s,{\rm Sym}^n\itPi\otimes\chi)$ and $\sigma \in {\rm Aut}(\C)$, we have
\begin{align*}
&\sigma \left( \frac{L^{(\infty)}(m+\tfrac{\delta}{2},{\rm Sym}^n\itPi\otimes\chi)}{(2\pi\sqrt{-1})^{d_n^\pm m}\cdot G(\chi)^{d_n^\pm}\cdot q^{\pm}({\rm Sym}^{n}\itPi)} \right) = \frac{L^{(\infty)}(m+\tfrac{\delta}{2},{\rm Sym}^n{}^\sigma\!\itPi\otimes{}^\sigma\!\chi)}{(2\pi\sqrt{-1})^{d_n^\pm m}\cdot G({}^\sigma\!\chi)^{d_n^\pm}\cdot q^{\pm}({\rm Sym}^{n}{}^\sigma\!\itPi)}.
\end{align*}
Here $\pm = (-1)^{m+r}\chi_\infty(-1)$, $d_n^+ = r+1,\,d_n^- = r$ (resp.\,$d_n^\pm = r+1$) if $n$ is even (resp.\,odd), and
\begin{align*}
&q^\pm({\rm Sym}^n\itPi) \\
&= 
\begin{cases}
(2\pi\sqrt{-1})^{r(r\pm1){\sf w}/2-r(r+1)/2}\cdot G(\omega_\itPi)^{r(r\pm1)/2}\cdot (c^+(\itPi)\cdot c^-(\itPi))^{r(r+1)/2} & \mbox{ if $n$ is even},\\
(2\pi\sqrt{-1})^{r(r+1)({\sf w}-1)/2}\cdot G(\omega_\itPi)^{r(r+1)/2}\cdot c^{\pm}(\itPi)^{(r+1)(r+2)/2}\cdot c^{\mp}(\itPi)^{r(r+1)/2} & \mbox{ if $n$ is odd}.
\end{cases}
\end{align*}
\end{conj}

To compare the statements in Conjectures A and \ref{C:Deligne Sym}, note that if $f$ is the normalized elliptic newform associated to $\itPi$, then we have
\begin{align*}
\delta(M_f) &= ((2\pi\sqrt{-1})^{1-\kappa}\cdot G({}^\sigma\!\omega_{\itPi}))_{\sigma : \Q(\itPi) \rightarrow \C},\\
c^\pm(M_f) &\in (2\pi\sqrt{-1})^{-(\kappa+{\sf w})/2}\cdot c^{\pm(-1)^{(\kappa+{\sf w})/2}}(M_\itPi)\cdot \Q(\itPi)^\times,\\
L(s,{\rm Sym}^n(f)) &= L^{(\infty)}(s-\tfrac{(\kappa+{\sf w}-1)n}{2},{\rm Sym}^n\itPi).
\end{align*}
We refer to the paragraph below Conjecture A for a brief survey on the results in the literature toward Conjecture \ref{C:Deligne Sym} when $n \leq 6$.
For arbitrary $n$, we prove the following:

\begin{thm}\label{T:Sym odd}
Conjecture \ref{C:Deligne Sym} holds when $\kappa \geq 5$.
\end{thm}

%\begin{rmk}\label{R:assumptions}
%If we assume the validity of Conjecture \ref{C:main 1}, then we can prove Conjecture \ref{C:Deligne Sym} for all $\kappa \geq 4$.
%\end{rmk}

\begin{rmk}
As a consequence of the theorem, we obtain period relations between the Betti--Whittaker periods of the symmetric power liftings of $\itPi$ and the motivic periods of $\itPi$. 
We refer to \cite[Conjecture 3.4]{Chen2021g} for the precise statement.
\end{rmk}

\begin{rmk}
When $n$ is odd, an automorphic variant of Conjecture \ref{C:Deligne Sym} is given by Grobner--Raghuram \cite[Proposition 8.1.4]{GR2014}, where the algebraicity is expressed in terms of the Betti--Shalika periods of ${\rm Sym}^n\itPi$.
\end{rmk}

We divide the proof of Theorem \ref{T:Sym odd} into two parts. In \S\,\ref{SS:Deligne Sym odd}, we first establish the conjecture for the case where $n$ is odd. Building on this result, along with some preliminaries in \S\,\ref{SS:Deligne Sym even}, we then prove the conjecture for the case where $n$ is even in \S\,\ref{SS:Deligne Sym even 2}.

%We divide the proof of Theorem \ref{T:Sym odd} into two parts. In \S\,\ref{SS:Deligne Sym odd}, we first prove the conjecture when $n$ is odd. Based on this case, we prove in 
%\S\,\ref{SS:Deligne Sym even} and \S\,\ref{SS:Deligne Sym even 2} for the case when $n$ is even.

%Let $\itPi$ be a cohomological cuspidal automorphic representation of $\GL_2(\A)$ with infinity type $(\kappa;\,{\sf w})$.
%For $n \geq 1$, let ${\rm Sym}^n \itPi$ be the symmetric $n$-th power functorial lift of $\itPi$ established by Newton and Thorne \cite{NT2021}, \cite{NT2021b}. Then ${\rm Sym}^n \itPi$ is a cohomological tamely isobaric automorphic representation of $\GL_{n+1}(\A)$. Moreover, its infinity type is given by
%\[
%\left(\left((n-2i)(\kappa-1)+1\right)_{0 \leq i\leq {\lfloor \tfrac{n-1}{2}\rfloor}};\,n{\sf w}\right).
%\]
%In this section, we show that Conjecture \ref{C:Deligne Sym} holds for ${\rm Sym}^n\itPi$ when $n=2r+1$ is odd and $\kappa\geq 5$ is odd. 
\subsubsection{Symmetric odd power $L$-functions for $\GL_2$}\label{SS:Deligne Sym odd}

As in the proof of Theorem \ref{T:tensor product for GL_2}, we prove Theorem \ref{T:Sym odd} when $n=2r+1$ by induction on $r$ using Theorem \ref{T:Main}.
In the induction step, we will encounter critical values of Rankin--Selberg $L$-functions for ${\rm Sym}^{N}\itPi \times \itPi'$ with $N=r,r\pm1,r-2$ and some regular algebraic cuspidal automorphic representation $\itPi'$ of $\GL_2(\A)$.
Conjecture \ref{C:Deligne} for the motive ${\rm Sym}^{N}M_\itPi \otimes M_{\itPi'}$ can also be proved by follow the strategy outlined in the beginning of \S\,\ref {S:applications}.
The precise statement is as follows.

\begin{thm}\label{T:Sym x GL_2}
Let $\itPi'$ be a regular algebraic cuspidal automorphic representation of $\GL_2(\A)$ with infinity type $(\kappa';\,{\sf w}')$. 
Let $N \geq 1$ be an integer. 
Put $\delta \in \{0,1\}$ with $\delta\equiv N+1\,({\rm mod}\,2)$.
Assume 
\[
\kappa'>(N-2)(\kappa-1)+1
\] 
and 
\begin{align}\label{E:regularity 3}
\begin{split}
&\min\{|N(\kappa-1)+1-\kappa'|,|(N-2)(\kappa-1)+1-\kappa'|, 2(\kappa-1)\}\\
& \geq \begin{cases} 4 & \mbox{ if ${\sf w}$ and ${\sf w}'$ are odd},\\
5 & \mbox{ if $N{\sf w}+{\sf w}' \equiv N\,({\rm mod}\,2)$},\\
6 & \mbox{ if ${\sf w}$ is even and ${\sf w}' \equiv N+1\,({\rm mod}\,2)$}.
\end{cases}
\end{split}
\end{align}
Then, for a critical point $m+\tfrac{\delta}{2} \in \Z+\tfrac{N+1}{2}$ for $L(s,{\rm Sym}^N\itPi\times\itPi')$ and $\sigma \in {\rm Aut}(\C)$, we have
\[
\sigma \left( \frac{L^{(\infty)}(m+\tfrac{\delta}{2},{\rm Sym}^N\itPi\times\itPi')}{(2\pi\sqrt{-1})^{(N+1)m}\cdot q^{(-1)^m}({\rm Sym}^N\itPi\times\itPi')}\right) = \frac{L^{(\infty)}(m+\tfrac{\delta}{2},{\rm Sym}^N{}^\sigma\!\itPi\times{}^\sigma\!\itPi')}{(2\pi\sqrt{-1})^{(N+1)m}\cdot q^{(-1)^m}({\rm Sym}^N{}^\sigma\!\itPi\times{}^\sigma\!\itPi')}.
\]
Here $q^\pm({\rm Sym}^N\itPi\times\itPi')$ are defined as follows:
\begin{itemize}
\item[(i)] If $\kappa' > N(\kappa-1)+1$, then
\begin{align*}
q^\pm({\rm Sym}^N\itPi\times\itPi') &= (2\pi\sqrt{-1})^{\lfloor\frac{N+1}{2}\rfloor(\delta-1)+N(N+1){\sf w}/2}\cdot G(\omega_\itPi)^{N(N+1)/2}\\
&\times (c^+(\itPi')\cdot c^-(\itPi'))^{\lfloor\frac{N+1}{2}\rfloor} \cdot \begin{cases}
c^{\pm(-1)^{N/2}}(\itPi') & \mbox{ if $N$ is even},\\
1 & \mbox{ if $N$ is odd}.
\end{cases}
\end{align*}
\item[(ii)] If $N(\kappa-1)+1>  \kappa' > (N-2)(\kappa-1)+1$, then
\begin{align*}
q^\pm({\rm Sym}^N\itPi\times\itPi')
& = (2\pi\sqrt{-1})^{\lfloor \frac{N}{2} \rfloor\delta-N-\lfloor \frac{N-1}{2} \rfloor+N(N-1){\sf w}/2+{\sf w}'}\cdot G(\omega_\itPi)^{N(N-1)/2}\cdot G(\omega_{\itPi'})\\
&\times (c^+(\itPi)\cdot c^-(\itPi))^N \cdot (c^+(\itPi')\cdot c^-(\itPi'))^{\lfloor \frac{N-1}{2} \rfloor}\cdot
\begin{cases}
c^{\pm(-1)^{N/2}}(\itPi') & \mbox{ if $N$ is even},\\
1 & \mbox{ if $N$ is odd}.
\end{cases}
\end{align*}
%\begin{align*}
%q^\pm({\rm Sym}^n\itPi\times\itPi')
%& = (2\pi\sqrt{-1})^{\lfloor \tfrac{n}{2} \rfloor\delta+n\kappa+n(n+1){\sf w}/2+((n-1)\kappa'+(n+1){\sf w}')/2}\cdot (\sqrt{-1})^{n{\sf w}+\lfloor \tfrac{n-1}{2} \rfloor{\sf w}'}\\
%&\times G(\omega_\itPi)^{n(n+1)/2}\cdot G(\omega_{\itPi'})^{\lfloor \tfrac{n+1}{2} \rfloor}\cdot\Vert f_\itPi\Vert^n\cdot\Vert f_{\itPi'}\Vert^{\lfloor \tfrac{n-1}{2} \rfloor}\\
%&\times
%\begin{cases}
%p(\itPi',\pm(-1)^{n/2}) & \mbox{ if $n$ is even},\\
%1 & \mbox{ if $n$ is odd}.
%\end{cases}
%\end{align*}
\end{itemize}
\end{thm}

\begin{proof}
We prove the assertion by induction on $N$. 
%Therefore, $\itPi$ will be fixed and $\itPi'$ varies with $n$.
When $N=1$, Cases (i) and (ii) are special cases of Conjecture \ref{C:Blasius} and were proved by Shimura \cite[Theorem 3]{Shimura1976}. 
When $N=2$, Case (i) is a special case of Theorem \ref{T:RS for GL_n x GL_2} and Case (ii) is a direct consequence of the result of Garrett and Harris \cite{GH1993} on triple product $L$-functions for $\GL_2$ in the balanced case (cf.\,\cite[Theorem 6]{Morimoto2021}). Note that the assumptions in \cite[Theorem 6]{Morimoto2021} can be further weakened by our result \cite[Theorem 1.2]{Chen2021d} and we only need to assume 
\[
2\kappa-1 \geq \kappa'+2 \geq 5
\]
in Case (ii) when $N=2$. Indeed, we have
\[
L(s,\itPi \times \itPi \times \itPi') = L(s, {\rm Sym}^2\itPi \times \itPi')\cdot L(s,\itPi'\otimes \omega_\itPi).
\]
The above condition is equivalent to saying that $L(s,\itPi \times \itPi \times \itPi')$ admits non-central critical points by (\ref{E:distance}). 
Therefore, Case (ii) when $N=2$ then follows from Deligne's conjecture for $L(s,\itPi \times \itPi \times \itPi')$ and $L(s,\itPi'\otimes \omega_\itPi)$.
In the induction step, the arguments for Cases (i) and (ii) are similar, we give details for Case (ii) and leave it to the readers to verify Case (i).
Note that in the induction step, the choices of $\itSigma'$ and $\itPsi'$ below are the same for both Case (i) and Case (ii).
We also stress that in the induction step, we only use the induction hypothesis and Theorem \ref{T:Main} in Case (i). Whereas in Case (ii), we also need the result for Case (i).
Let $N \geq 3$ and assume the assertion for Case (ii) holds for all $N-1 \geq N' \geq 1$.
Let $\itPi'$ and $\delta$ be as in the assertion satisfying the regularity assumptions.
We choose auxiliary automorphic representations $\itPi_j$, $\itPi_1'$, $\itPi_2'$, and $\chi$ as follows:
\begin{itemize}
\item For $0 \leq j \leq {\lfloor \tfrac{N-3}{2} \rfloor}$, let $\itPi_j$ be a regular algebraic cuspidal automorphic representation of $\GL_2(\A)$ with infinity type
\[
((N-1-2j)(\kappa-1)+1;\,(N-1){\sf w}+1-\delta).
\]
\item If $N$ is even and $j=\tfrac{N-2}{2}$, let 
\[
\itPi_{j} = \itPi \otimes |\mbox{ }|_\A^{(N-2){\sf w}/2}.
\]
\item Let $\itPi_1'$ and $\itPi_2'$ be regular algebraic automorphic representations of $\GL_2(\A)$ with infinity types
\[
(\kappa+\kappa'-1;\,{\sf w}+{\sf w}'-1),\quad (\kappa'-\kappa+1;\,{\sf w}+{\sf w}'-1)
\]
respectively.
\item If $N$ is odd, let $\chi$ be a finite order Hecke character of $\A^\times$ such that $\chi_\infty(-1) = (-1)^{(N-1)(1+{\sf w})/2}$.
\end{itemize}
Let $\itSigma$ and $\itSigma'$ be regular algebraic automorphic representations of $\GL_N(\A)$ defined by
\begin{align*}
\itSigma = {\rm Sym}^{N-1}\itPi,\quad \itSigma'  = \begin{cases}
\bigboxplus_{j=0}^{{\lfloor \frac{N-2}{2} \rfloor}} \itPi_j & \mbox{ if $N$ is even},\\
\bigboxplus_{j=0}^{{\lfloor \frac{N-2}{2} \rfloor}} (\itPi_j\otimes|\mbox{ }|_\A^{-1/2})\boxplus \chi|\mbox{ }|_\A^{(N-1){\sf w}/2} & \mbox{ if $N$ is odd}.
\end{cases}
\end{align*}
Let $\itPsi$ and $\itPsi'$ be regular algebraic automorphic representations of $\GL_4(\A)$ defined by
\begin{align*}
\itPsi = (\itPi\boxtimes\itPi')\otimes|\mbox{ }|_\A^{-1/2},\quad\itPsi'  = \itPi_1'\boxplus\itPi_2'.
\end{align*}
The infinity types of $\itSigma_\infty$ and $\itPsi_\infty$ are given in (\ref{E:infinity type sym}) and Lemma \ref{L:archimedean tensor} respectively.
Also by our condition on $\chi_\infty(-1)$, we have $\varepsilon(\itSigma_\infty) = \varepsilon(\itSigma_\infty')$ if $N$ is odd (cf.\,(\ref{E:signature})).
Thus $\itSigma_\infty$ and $\itPsi_\infty$ are essentially tempered with $\itSigma_\infty = \itSigma'_\infty$ and $\itPsi_\infty = \itPsi_\infty'$.
By \cite[Theorem 7.21]{HR2020}, we only need to consider non-central critical points.
By (\ref{E:distance}), our regularity assumptions imply that $L(s,{\rm Sym}^N\itPi \times \itPi')$ admits non-central critical points. 
Let $m+\tfrac{\delta}{2} \in \Z+\tfrac{N+1}{2}$ be a non-central critical point. Note that  
\[
({\rm Sym}^{N-1} \itPi) \boxtimes \itPi = {\rm Sym}^N\itPi \boxplus ({\rm Sym}^{N-2}\itPi \otimes \omega_\itPi).
\]
In particular, we have the factorization of $L$-function:
\[
L(s,\itSigma \times \itPsi) = L(s-\tfrac{1}{2},{\rm Sym}^N\itPi \times \itPi')\cdot L(s-\tfrac{1}{2},{\rm Sym}^{N-2}\itPi \times \itPi'\otimes \omega_\itPi).
\]
Also it is clear that
\[
L(s,\itSigma \times \itPsi') = L(s,{\rm Sym}^{N-1}\itPi \times \itPi_1')\cdot L(s,{\rm Sym}^{N-1}\itPi \times \itPi_2').
\]
Applying the result for Case (i) to $L(m+\tfrac{\delta}{2},{\rm Sym}^{N-2}\itPi \times \itPi'\otimes \omega_\itPi)$ and $L(m+\tfrac{\delta+1}{2},{\rm Sym}^{N-1}\itPi \times \itPi_1')$, we obtain
\begin{align}\label{E:Sym x GL_2 proof 1}
\begin{split}
&\frac{L^{(\infty)}(m+\tfrac{\delta+1}{2},\itSigma \times \itPsi)}{L^{(\infty)}(m+\tfrac{\delta+1}{2},\itSigma \times \itPsi') }\\ &\sim %\frac{L^{(\infty)}(m+\tfrac{\delta}{2},{\rm Sym}^n \itPi \times \itPi')}{L^{(\infty)}(m+\tfrac{1+\delta}{2},{\rm Sym}^{n-1} \itPi \times \itPi_2')}\\
%&\times (2\pi\sqrt{-1})^{-m-{\lfloor \tfrac{n}{2} \rfloor}+n-\delta-(n(\kappa+{\sf w})+\kappa'+{\sf w}')/2}\cdot (\sqrt{-1})^{{\lfloor \tfrac{n}{2} \rfloor}(1+{\sf w})+(n-1){\sf w}'}\\
%&\times G(\omega_{\itPi'})^{{\lfloor \tfrac{n-1}{2} \rfloor}}\cdot G(\omega_{\itPi_1'})^{-{\lfloor \tfrac{n}{2} \rfloor}}\cdot \Vert f_{\itPi'} \Vert^{{\lfloor \tfrac{n-1}{2} \rfloor}}\cdot \Vert f_{\itPi_1'} \Vert^{-{\lfloor \tfrac{n}{2} \rfloor}}\\
%&\times \begin{cases}
%p(\itPi',(-1)^{m+(n-2)/2}) & \mbox{ if $n$ is even},\\
%p(\itPi_1',(-1)^{m+(n-1)/2})^{-1} & \mbox{ if $n$ is odd}.
%\end{cases}
\frac{L^{(\infty)}(m+\tfrac{\delta}{2},{\rm Sym}^N \itPi \times \itPi')}{L^{(\infty)}(m+\tfrac{\delta+1}{2},{\rm Sym}^{N-1} \itPi \times \itPi_2')}\\
&\times (2\pi\sqrt{-1})^{-m-{\lfloor \frac{N}{2} \rfloor}}\cdot \frac{(c^+(\itPi')\cdot c^-(\itPi'))^{\lfloor \frac{N-1}{2} \rfloor}}{(c^+(\itPi_1')\cdot c^-(\itPi_1'))^{\lfloor \frac{N}{2} \rfloor}} \cdot \begin{cases}
c^{(-1)^{m-1+N/2}}(\itPi') & \mbox{ if $N$ is even},\\
c^{(-1)^{m+(N-1)/2}}(\itPi_1')^{-1} & \mbox{ if $N$ is odd}.
\end{cases}
\end{split}
\end{align}
Here we have used the period relation (cf.\,\cite[(3.6) and (3.9)]{Yoshida1994})
\[
c^\pm(\itPi' \otimes \omega_\itPi)\sim (2\pi\sqrt{-1})^{\sf w}\cdot G(\omega_{\itPi})\cdot c^{\pm(-1)^{{\sf w}}}(\itPi').
\]
On the other hand, we have the following factorizations:
\begin{align*}
L(s,\itSigma' \times \itPsi) &= \prod_{j=0}^{{\lfloor \frac{N-2}{2} \rfloor}}L(s+\tfrac{\delta-2}{2},\itPi_j \times \itPi \times \itPi')\cdot \begin{cases}
1 & \mbox{ if $N$ is even},\\
L(s+\tfrac{(N-1){\sf w}-1}{2},\itPi \times \itPi'\otimes\chi) & \mbox{ if $N$ is odd}.
\end{cases}\\
L(s,\itSigma' \times \itPsi') &= \prod_{j=0}^{{\lfloor \frac{N-2}{2} \rfloor}}L(s+\tfrac{\delta-1}{2},\itPi_j \times (\itPi_1' \boxplus \itPi_2'))\cdot \begin{cases}
1 & \mbox{ if $N$ is even},\\
L(s+\tfrac{(N-1){\sf w}}{2},(\itPi_1' \boxplus \itPi_2')\otimes\chi) & \mbox{ if $N$ is odd}.
\end{cases}
\end{align*}
For $0 \leq j \leq {\lfloor \tfrac{N-2}{2} \rfloor}$, Conjecture \ref{C:Blasius} holds for the triple product $L$-function $L(s,\itPi_j \times \itPi \times \itPi')$.
Indeed, when $N$ is even and $j=\tfrac{N-2}{2}$, we have
\[
L(s,\itPi_j \times \itPi \times \itPi') = L(s+\tfrac{(N-2){\sf w}}{2},{\rm Sym}^2\itPi \times \itPi')\cdot L(s+\tfrac{(N-2){\sf w}}{2},\itPi' \otimes \omega_\itPi).
\]
Therefore, Conjecture \ref{C:Blasius} holds by Case (i) for $L(s,{\rm Sym}^2\itPi \times \itPi')$ and Conjecture \ref{C:Deligne} for $M=M_{\itPi' \otimes \omega_\itPi}$.
When $1 \leq j \leq {\lfloor \tfrac{N-3}{2} \rfloor}$, we are in the unbalanced cases and the conditions in (\ref{E:unbalanced}) are satisfied by the regularity assumptions here, thus Conjecture \ref{C:Blasius} also holds.
When $j=0$, we are in the balanced case and Conjecture \ref{C:Blasius} holds by \cite{GH1993} and \cite{Chen2021d}.
Thus we have
\begin{align*}
&L^{(\infty)}(m+\tfrac{\delta+1}{2},\itSigma' \times \itPsi) \\
%&\sim (2\pi\sqrt{-1})^{4(m+\delta-1/2)+2(N{\sf w}+{\sf w}'+1-\delta)}\cdot G(\omega_{\itPi_0}\omega_{\itPi}\omega_{\itPi'})^2\\
%&\times (\pi^{(N-1)(\kappa-1)+1}\cdot \<f_{\itPi_0},f_{\itPi_0}\>)\cdot (\pi^{\kappa}\cdot \<f_{\itPi},f_{\itPi}\>)\cdot(\pi^{\kappa'}\cdot \<f_{\itPi'},f_{\itPi'}\>)\\
%&\times \prod_{j=1}^{\lfloor \frac{N-2}{2}\rfloor}(2\pi\sqrt{-1})^{4(m+\delta-1/2)+2(N{\sf w}+{\sf w}'+1-\delta)}\cdot G(\omega_{\itPi_j}\omega_{\itPi}\omega_{\itPi'})^2\cdot (\pi^{\kappa'}\cdot \<f_{\itPi'},f_{\itPi'}\>)^2\\
%&\times \begin{cases}
%1 & \mbox{ if $N$ is even},\\
%(2\pi\sqrt{-1})^{2(m+(N-1){\sf w}/2)+{\sf w}+{\sf w}'}\cdot G(\omega_{\itPi}\omega_{\itPi'}\chi^2)\cdot (\pi^{\kappa'}\cdot \<f_{\itPi'},f_{\itPi'}\>) & \mbox{ if $N$ is odd},
%\end{cases}\\
&\sim (2\pi\sqrt{-1})^{(4(m+\delta-1/2)+2(N{\sf w}-\delta))(1+\lfloor \frac{N-2}{2} \rfloor)-{\sf w}+{\sf w}'}\cdot
G(\omega_{\itPi})^{N-1}\cdot G(\omega_{\itPi'})
\cdot \prod_{j=0}^{\lfloor \frac{N-2}{2} \rfloor}G(\omega_{\itPi_j})^2\\
&\times(\pi^{(N-1)(\kappa-1)+1}\cdot \<f_{\itPi_0},f_{\itPi_0}\>)\cdot(c^+(\itPi)\cdot c^-(\itPi))\cdot(c^+(\itPi')\cdot c^-(\itPi'))^{N-1}\\
&\times \begin{cases}
1 & \mbox{ if $N$ is even},\\
(2\pi\sqrt{-1})^{2(m+(N-1){\sf w}/2)+{\sf w}-1}\cdot G(\chi^2) & \mbox{ if $N$ is odd}.
\end{cases}
\end{align*}
Here we have used the period relations (cf.\,\cite[Theorem 1-(iv)]{Shimura1977} and \cite[(3.6)]{Yoshida1994})
\begin{align}\label{E:period relation Shimura}
\begin{split}
c^+(\itPi)\cdot c^-(\itPi) &\sim (2\pi\sqrt{-1})^{{\sf w}+1}\cdot G(\omega_\itPi)\cdot (\pi^\kappa\cdot\<f_\itPi,f_\itPi\>),\\
c^\pm(\itPi \otimes \chi) &\sim G(\chi)\cdot c^{\pm\chi_\infty(-1)}(\itPi),
\end{split}
\end{align}
and its analogue for $\itPi'$.
By Conjecture \ref{C:Blasius} for $n=2$ (cf.\,\cite[Theorem 3]{Shimura1976}) and Conjecture \ref{C:Deligne} for $M=M_{\itPi_1'\otimes\chi}$ and $M=M_{\itPi_2'\otimes\chi}$, we have
\begin{align*}
&L^{(\infty)}(m+\tfrac{\delta+1}{2},\itSigma' \times \itPsi') \\
%&\sim (2\pi\sqrt{-1})^{4(m+\delta)+2(N{\sf w}+{\sf w}'-\delta)}\cdot G(\omega_{\itPi_0}^2\omega_{\itPi_1'}\omega_{\itPi_2'})\cdot(\pi^{(N-1)(\kappa-1)+1}\cdot\<f_{\itPi_0},f_{\itPi_0}\>)\cdot (\pi^{\kappa+\kappa'-1}\cdot\<f_{\itPi_1'},f_{\itPi_1'}\>)\\
%&\prod_{j=1}^{\lfloor \frac{N-2}{2}\rfloor} (2\pi\sqrt{-1})^{4(m+\delta)+2(N{\sf w}+{\sf w}'-\delta)}\cdot G(\omega_{\itPi_j}^2\omega_{\itPi_1'}\omega_{\itPi_2'})\cdot (\pi^{\kappa+\kappa'-1}\cdot\<f_{\itPi_1'},f_{\itPi_1'}\>)\cdot(\pi^{\kappa'-\kappa+1}\cdot\<f_{\itPi_2'},f_{\itPi_2'}\>)\\
%&\times \begin{cases}
%1 & \mbox{ if $N$ is even},\\
%(2\pi\sqrt{-1})^{2(m+(N-1){\sf w}/2)}\cdot G(\chi^2)\cdot c^{(-1)^{m+(N-1){\sf w}/2}}(\itPi_1')\cdot c^{m+(N-1){\sf w}/2}(\itPi_2') & \mbox{ if $N$ is odd},
%\end{cases}\\
&\sim (2\pi\sqrt{-1})^{(4(m+\delta)+2((N-1){\sf w}-\delta))(1+\lfloor \frac{N-2}{2}\rfloor)+{\sf w}+{\sf w}'}\cdot G(\omega_{\itPi_2'})\cdot\prod_{j=0}^{\lfloor \frac{N-2}{2}\rfloor}G(\omega_{\itPi_j})^2\\
&\times (\pi^{(N-1)(\kappa-1)+1}\cdot\<f_{\itPi_0},f_{\itPi_0}\>)\cdot (c^+(\itPi_1')\cdot c^-(\itPi_1'))^{\lfloor \frac{N}{2}\rfloor} \cdot (c^+(\itPi_2')\cdot c^-(\itPi_2'))^{\lfloor \frac{N-2}{2}\rfloor} \\
&\times \begin{cases}
1 & \mbox{ if $N$ is even},\\
(2\pi\sqrt{-1})^{2(m+(N-1){\sf w}/2)}\cdot G(\chi^2)\cdot c^{(-1)^{m+(N-1)/2}}(\itPi_1')\cdot c^{(-1)^{m+(N-1)/2}}(\itPi_2') & \mbox{ if $N$ is odd}.
\end{cases}
\end{align*}
Here we also used the period relation (\ref{E:period relation Shimura}) for $\itPi_1'$ and $\itPi_2'$.
Hence we conclude that
\begin{align}\label{E:Sym x GL_2 proof 2}
\begin{split}
&\frac{L^{(\infty)}(m+\tfrac{\delta+1}{2},\itSigma' \times \itPsi')}{L^{(\infty)}(m+\tfrac{\delta+1}{2},\itSigma' \times \itPsi)}\\
&\,\sim (2\pi\sqrt{-1})^{N+(-N+2){\sf w}}\cdot G(\omega_{\itPi})^{-N+1}\cdot G(\omega_{\itPi'})^{-1}\cdot G(\omega_{\itPi_2'})\\
&\times \frac{(c^+(\itPi_1')\cdot c^-(\itPi_1'))^{\lfloor \frac{N}{2} \rfloor}\cdot(c^+(\itPi_2')\cdot c^-(\itPi_2'))^{\lfloor \frac{N-2}{2} \rfloor}}{(c^+(\itPi)\cdot c^-(\itPi))\cdot (c^+(\itPi')\cdot c^-(\itPi'))^{N-1}}
\cdot\begin{cases}
1 & \mbox{ if $N$ is even},\\
c^{(-1)^{m+(N-1)/2}}(\itPi_1')\cdot c^{(-1)^{m+(N-1)/2}}(\itPi_2') & \mbox{ if $N$ is odd}.
\end{cases}
\end{split}
\end{align}
{Note that the three cases of the regularity assumptions considered in (\ref{E:regularity 3}) correspond, case by case and in the same order, to those in (\ref{E:regularity 2}) imposed on $\itSigma_\infty$ and $\itPsi_\infty$.}
Therefore, by Theorem \ref{T:Main} we have
\begin{align}\label{E:Sym x GL_2 proof 3}
\frac{L(m+\tfrac{\delta+1}{2},\itSigma \times \itPsi)\cdot L(m+\tfrac{\delta+1}{2},\itSigma' \times \itPsi')}{L(m+\tfrac{\delta+1}{2},\itSigma \times \itPsi')\cdot L(m+\tfrac{\delta+1}{2},\itSigma' \times \itPsi)} \sim1.
\end{align}
In (\ref{E:Sym x GL_2 proof 1}), when $N$ is even, we rewrite $c^{(-1)^{m-1+N/2}}(\itPi')$ as $(c^+(\itPi')\cdot c^-(\itPi'))\cdot c^{(-1)^{m+N/2}}(\itPi')^{-1}$.
We then deduce from (\ref{E:Sym x GL_2 proof 1}), (\ref{E:Sym x GL_2 proof 2}), and (\ref{E:Sym x GL_2 proof 3}) that
\begin{align*}
&\frac{L^{(\infty)}(m+\tfrac{\delta}{2},{\rm Sym}^N \itPi \times \itPi')}{L^{(\infty)}(m+\tfrac{\delta+1}{2},{\rm Sym}^{N-1} \itPi \times \itPi_2')}\\
&\sim(2\pi\sqrt{-1})^{m+{\lfloor \frac{N}{2} \rfloor}-N+(N-2){\sf w}}\cdot G(\omega_{\itPi})^{N-1}\cdot G(\omega_{\itPi'})\cdot G(\omega_{\itPi_2'})^{-1}\\
&\times \frac{(c^+(\itPi)\cdot c^-(\itPi))\cdot (c^+(\itPi')\cdot c^-(\itPi'))^{\lfloor \frac{N-1}{2} \rfloor}}{(c^+(\itPi_2')\cdot c^-(\itPi_2'))^{\lfloor \frac{N-2}{2} \rfloor}}\cdot\begin{cases}
c^{(-1)^{m+N/2}}(\itPi') & \mbox{ if $N$ is even},\\
c^{(-1)^{m+(N-1)/2}}(\itPi_2')^{-1} & \mbox{ if $N$ is odd}.
\end{cases}
\end{align*}
Note that we are in Case (ii) for the $L$-function in the denominator of the above ratio.
It is easy to verify that the right-hand side is equal to
\[
\frac{(2\pi\sqrt{-1})^{(N+1)m}\cdot q^{(-1)^m}({\rm Sym}^N\itPi \times \itPi')}{(2\pi\sqrt{-1})^{N(m+\delta)}\cdot q^{(-1)^{m+\delta}}({\rm Sym}^{N-1}\itPi \times \itPi_2')}.
\]
{Finally, the regularity assumptions on the tuple
$({\rm Sym}^N \itPi,\itPi')$ immediately imply the corresponding assumptions for $({\rm Sym}^{N-1}\itPi,\itPi_2')$.}
By the induction hypothesis, we have
\[
L^{(\infty)}(m+\tfrac{\delta+1}{2},{\rm Sym}^{N-1} \itPi \times \itPi_2')\sim(2\pi\sqrt{-1})^{N(m+\delta)}\cdot q^{(-1)^{m+\delta}}({\rm Sym}^{N-1}\itPi \times \itPi_2').
\]
We thus conclude that 
\[
L^{(\infty)}(m+\tfrac{\delta}{2},{\rm Sym}^N \itPi \times \itPi')\sim(2\pi\sqrt{-1})^{(N+1)m}\cdot q^{(-1)^m}({\rm Sym}^N\itPi \times \itPi').
\]
This completes the proof.
\end{proof}

\begin{rmk}
Since ${\rm Sym}^N\itPi$ is essentially self-dual, the result for Case (i) is a special case of Conjecture \ref{C:Deligne for GL_n x GL_2}.
Contrary to condition (2) in Theorem \ref{T:RS for GL_n x GL_2}, here we do not need to assume local condition on finite place.
\end{rmk}

Now we begin the proof of Theorem \ref{T:Sym odd} when $n$ is odd. 
%We rewrite $q^\pm({\rm Sym}^{2r+1}\itPi)$ in terms of the Betti--Whittaker periods of $\itPi$ (cf.\,(\ref{E:period GL_2})):
%\begin{align*}
%q^\pm({\rm Sym}^{2r+1}\itPi) &= (2\pi\sqrt{-1})^{(r+1)((r+1)\kappa+(2r+1){\sf w}-r)/2}\cdot G(\omega_\itPi)^{r(r+1)/2}\\
%&\times p(\itPi,\pm(-1)^r)^{(r+1)(r+2)/2}\cdot p(\itPi,\mp(-1)^r)^{r(r+1)/2}.
%\end{align*}
We prove the assertion by induction on $n=2r+1$. 
When $r=0$, Conjecture \ref{C:Deligne Sym} is known, as explained by Deligne in \cite[\S\,7]{Deligne1979}.
When $r=1$, Conjecture \ref{C:Deligne Sym} was proved in \cite[Theorem 6.2]{GH1993} and \cite[Theorem 1.6]{Chen2021d} under the assumption that $\kappa \geq 3$. Let $r \geq 2$ and assume Conjecture \ref{C:Deligne Sym} holds for all $r-1 \geq r' \geq 1$. Put $\delta \in \{0,1\}$ with $\delta \equiv r \,({\rm mod}\,2)$.
We choose auxiliary automorphic representations $\itPi'$, $\itPi''$, and $\eta$ as follows:
\begin{itemize}
\item Let $\itPi'$ and $\itPi''$ be regular algebraic cuspidal automorphic representations of $\GL_2(\A)$ with infinity types
\begin{align*}
(\kappa';\,{\sf w}') = ((r+1)(\kappa-1)+1;\,(r+1){\sf w}+\delta),\quad
(\kappa'';\,{\sf w}'') = (r(\kappa-1)+1;\,r{\sf w}+1-\delta).
\end{align*}
\item Let $\eta$ be a finite order Hecke character of $\A^\times$ such that $\eta_\infty(-1) = (-1)^{1+{\sf w}}$.
\end{itemize}
Let $\itSigma$ and $\itSigma'$ be regular algebraic automorphic representations of $\GL_{r+2}(\A)$ defined by
\[
\itSigma = {\rm Sym}^{r+1}\itPi,\quad \itSigma' = ({\rm Sym}^{r-1}\itPi\otimes \eta|\mbox{ }|_\A^{{\sf w}})\boxplus(\itPi' \otimes |\mbox{ }|_\A^{-\delta/2}).
\]
Let $\itPsi$ and $\itPsi'$ be regular algebraic automorphic representations of $\GL_{r+1}(\A)$ defined by
\[
\itPsi = {\rm Sym}^r\itPi \otimes \chi,\quad \itPsi' = ({\rm Sym}^{r-2}\itPi\otimes\chi\eta|\mbox{ }|_\A^{{\sf w}})\boxplus(\itPi'' \otimes |\mbox{ }|_\A^{(\delta-1)/2}).
\]
Recall the infinity types of the symmetric powers of $\itPi$ in (\ref{E:infinity type sym}). Also the condition on $\eta_\infty(-1)$ implies that $\varepsilon(\itSigma_\infty) = \varepsilon(\itSigma_\infty')$ (resp.\,$\varepsilon(\itPsi_\infty) = \varepsilon(\itPsi_\infty')$) if $r$ is odd (resp.\,even) by (\ref{E:signature}).
Thus $\itSigma_\infty$ and $\itPsi_\infty$ are essentially tempered with $\itSigma_\infty = \itSigma'_\infty$ and $\itPsi_\infty = \itPsi'_\infty$. 
By \cite[Theorem 7.21]{HR2020}, we only need to consider non-central critical points.
Indeed, $\kappa \geq 5$ implies that the symmetric odd power $L$-function admits at least two consecutive non-central right-half critical points.
Therefore, if we can show that Conjecture \ref{C:Deligne Sym} holds for $L(s,{\rm Sym}^{2r+1}\itPi\otimes\chi)$ at all non-central critical points, then the relative periods of ${\rm Sym}^{2r+1}\itPi$ defined in \cite[Definition 5.3]{HR2020} are  given by the ratios between $q^+({\rm Sym}^{2r+1}\itPi)$ and $q^-({\rm Sym}^{2r+1}\itPi)$. 
Let $m+\tfrac{1}{2}\in\Z+\tfrac{1}{2}$ be a non-central critical point for $L(s,{\rm Sym}^{2r+1}\itPi\otimes \chi)$.
Note that
\[
({\rm Sym}^{r+1}\itPi \otimes \chi) \boxtimes {\rm Sym}^r\itPi = \bigboxplus_{j=0}^{r}{\rm Sym}^{2(r-j)+1}\itPi \otimes \omega_\itPi^j\chi.
\]
In particular, we have
\begin{align*}
L(s,\itSigma \times \itPsi) &= \prod_{j=0}^r L(s,{\rm Sym}^{2(r-j)+1}\itPi \otimes \omega_\itPi^j\chi)\\
&=L(s,{\rm Sym}^{2r+1}\itPi \otimes \chi)\cdot\prod_{j=0}^{r-1} L(s,{\rm Sym}^{2(r-j)-1}\itPi \otimes \omega_\itPi^{j+1}\chi).
\end{align*}
Similarly, we have
\[
L(s,\itSigma' \times \itPsi) = \prod_{j=0}^{r-1}L(s+{\sf w},{\rm Sym}^{2(r-j)-1}\itPi \otimes \omega_\itPi^j\chi\eta)\cdot L(s-\tfrac{\delta}{2},{\rm Sym}^r\itPi \times \itPi'\otimes\chi).
\]
{Note that the assumption $\kappa \geq 5$ implies that the regularity conditions in Theorem~\ref{T:Sym x GL_2} are satisfied for the tuple $({\rm Sym}^r \itPi, \itPi' \otimes \chi)$.
More precisely, if $\kappa$ is odd (resp.\ even), then we are in the first (resp.\,second) case of \eqref{E:regularity 3}, and the corresponding condition is satisfied if and only if $\kappa \geq 5$ (resp.\ $\kappa \geq 6$).}
By Theorem \ref{T:Sym x GL_2}, we have
\[
L^{(\infty)}(m+\tfrac{1-\delta}{2},{\rm Sym}^r\itPi \times \itPi'\otimes\chi) \sim(2\pi\sqrt{-1})^{(r+1)m}\cdot q^{(-1)^m}({\rm Sym}^r\itPi \times \itPi' \otimes \chi).
\]
Put 
$
\varepsilon = (-1)^{m+r}\chi_\infty(-1)$.
By the induction hypothesis, we have
\[
\prod_{j=0}^{r-1}\frac{L^{(\infty)}(m+\tfrac{1}{2},{\rm Sym}^{2(r-j)-1}\itPi\otimes\omega_{\itPi}^{j+1}\chi)}{L^{(\infty)}(m+{\sf w}+\tfrac{1}{2},{\rm Sym}^{2(r-j)-1}\itPi\otimes\omega_{\itPi}^j\chi\eta)}\sim G(\omega_\itPi\eta^{-1})^{r(r+1)/2}\cdot\prod_{j=0}^{r-1}\frac{q^{(-1)^{j-1}\varepsilon}({\rm Sym}^{2(r-j)-1}\itPi)}{q^{(-1)^j\varepsilon}({\rm Sym}^{2(r-j)-1}\itPi)}.
\]
Therefore, we have
\begin{align}\label{E:Sym odd proof 1}
\begin{split}
\frac{L(m+\tfrac{1}{2},\itSigma \times \itPsi)}{L(m+\tfrac{1}{2},\itSigma' \times \itPsi)}&\sim\frac{L^{(\infty)}(m+\tfrac{1}{2},{\rm Sym}^{2r+1}\itPi\otimes\chi)}{(2\pi\sqrt{-1})^{(r+1)m}\cdot q^{(-1)^m}({\rm Sym}^r\itPi \times \itPi' \otimes \chi)}\\
&\times G(\omega_\itPi\eta^{-1})^{r(r+1)/2}\cdot\prod_{j=0}^{r-1}\frac{q^{(-1)^{j-1}\varepsilon}({\rm Sym}^{2(r-j)-1}\itPi)}{q^{(-1)^j\varepsilon}({\rm Sym}^{2(r-j)-1}\itPi)}.
\end{split}
\end{align}
On the other hand, we also have
\begin{align*}
{\rm Sym}^{r+1}\itPi \boxtimes ({\rm Sym}^{r-2}\itPi\otimes \chi\eta) = \bigboxplus_{j=0}^{r-2}{\rm Sym}^{2(r-j)-1}\itPi \otimes \omega_\itPi^j\chi\eta.
\end{align*}
Hence we have the following factorizations:
\begin{align*}
L(s,\itSigma \times \itPsi') & = \prod_{j=0}^{r-2}L(s+{\sf w},{\rm Sym}^{2(r-j)-1}\itPi\otimes\omega_\itPi^j\chi\eta)\cdot L(s+\tfrac{\delta-1}{2},{\rm Sym}^{r+1}\itPi \times \itPi''),\\
L(s,\itSigma' \times \itPsi') & = \prod_{j=0}^{r-2}L(s+2{\sf w},{\rm Sym}^{2(r-j)-3}\itPi\otimes\omega_\itPi^j\chi\eta^2)\cdot L(s-\tfrac{1}{2},\itPi'\times \itPi'')\\
&\times L(s+{\sf w}+\tfrac{\delta-1}{2},{\rm Sym}^{r-1}\itPi\times \itPi''\otimes\eta)\cdot L(s+{\sf w}-\tfrac{\delta}{2},{\rm Sym}^{r-2}\itPi\times \itPi'\otimes\chi\eta).
\end{align*}
Similarly, the assumption $\kappa \geq 5$ implies that all the Rankin--Selberg $L$-functions appearing in the factorizations of $L(s,\itSigma \times \itPsi')$ or $L(s,\itSigma' \times \itPsi')$, except for the symmetric odd power $L$-functions, satisfy the regularity conditions in Theorem \ref{T:Sym x GL_2}.
Therefore, by Theorem \ref{T:Sym x GL_2} we have
\begin{align*}
&\frac{L^{(\infty)}(m,\itPi'\times \itPi'')\cdot L^{(\infty)}(m+{\sf w}+\tfrac{\delta}{2},{\rm Sym}^{r-1}\itPi\times \itPi''\otimes \eta)\cdot L^{(\infty)}(m+{\sf w}+\tfrac{1-\delta}{2},{\rm Sym}^{r-2}\itPi\times \itPi'\otimes\chi\eta)}{L^{(\infty)}(m+\tfrac{\delta}{2},{\rm Sym}^{r+1}\itPi \times \itPi'')}\\
& \sim (2\pi\sqrt{-1})^{(r-1)m+(2r-1){\sf w}}\\
&\times \frac{q^{(-1)^m}(\itPi'\times \itPi'')\cdot q^{(-1)^{m+{\sf w}}}({\rm Sym}^{r-1}\itPi\times \itPi''\otimes \eta)\cdot q^{(-1)^{m+{\sf w}}}({\rm Sym}^{r-2}\itPi\times \itPi'\otimes\chi\eta)}{q^{(-1)^m}({\rm Sym}^{r+1}\itPi \times \itPi'')}.
\end{align*}
By the induction hypothesis, we have
\begin{align*}
&\prod_{j=0}^{r-2}\frac{L^{(\infty)}(m+2{\sf w}+\tfrac{1}{2},{\rm Sym}^{2(r-j)-3}\itPi\otimes\omega_\itPi^j\chi\eta^2)}{L^{(\infty)}(m+{\sf w}+\tfrac{1}{2},{\rm Sym}^{2(r-j)-1}\itPi\otimes\omega_\itPi^j\chi\eta)}\\
&= \frac{L^{(\infty)}(m+2{\sf w}+\tfrac{1}{2},\itPi\otimes\omega_\itPi^{r-2}\chi\eta^2)}{L^{(\infty)}(m+{\sf w}+\tfrac{1}{2},{\rm Sym}^{2r-1}\itPi\otimes\chi\eta)}\cdot \prod_{j=0}^{r-3}\frac{L^{(\infty)}(m+2{\sf w}+\tfrac{1}{2},{\rm Sym}^{2(r-j)-3}\itPi\otimes\omega_\itPi^j\chi\eta^2)}{L^{(\infty)}(m+{\sf w}+\tfrac{1}{2},{\rm Sym}^{2(r-j)-3}\itPi\otimes\omega_\itPi^{j+1}\chi\eta)}\\
& \sim (2\pi\sqrt{-1})^{(1-r)m}\cdot G(\omega_\itPi\eta^{-1})^{-(r-1)(r-2)/2}\cdot G(\chi)^{1-r}\cdot \frac{c^{(-1)^r\varepsilon}(\itPi)}{q^{\varepsilon}({\rm Sym}^{2r-1}\itPi)}\cdot \prod_{j=0}^{r-3} \frac{q^{(-1)^{j}\varepsilon}({\rm Sym}^{2(r-j)-3}\itPi)}{q^{(-1)^{j-1}\varepsilon}({\rm Sym}^{2(r-j)-3}\itPi)}.
\end{align*}
Therefore, we have
\begin{align}\label{E:Sym odd proof 2}
\begin{split}
\frac{L(m+\tfrac{1}{2},\itSigma' \times \itPsi')}{L(m+\tfrac{1}{2},\itSigma \times \itPsi')}&\sim(2\pi\sqrt{-1})^{(2r-1){\sf w}}\cdot G(\omega_\itPi\eta^{-1})^{-(r-1)(r-2)/2}\cdot G(\chi)^{1-r}\\
&\times \frac{q^{(-1)^m}(\itPi'\times \itPi'')\cdot q^{(-1)^{m+{\sf w}}}({\rm Sym}^{r-1}\itPi\times \itPi''\otimes\eta)\cdot q^{(-1)^{m+{\sf w}}}({\rm Sym}^{r-2}\itPi\times \itPi'\otimes\chi\eta)}{q^{(-1)^m}({\rm Sym}^{r+1}\itPi \times \itPi'')}\\
&\times \frac{c^{(-1)^r\varepsilon}(\itPi)}{q^{\varepsilon}({\rm Sym}^{2r-1}\itPi)}\cdot \prod_{j=0}^{r-3} \frac{q^{(-1)^{j}\varepsilon}({\rm Sym}^{2(r-j)-3}\itPi)}{q^{(-1)^{j-1}\varepsilon}({\rm Sym}^{2(r-j)-3}\itPi)}.
\end{split}
\end{align}
Since $\kappa \geq 5$, the conditions in Theorem \ref{T:Main} on the infinity types of $\itSigma$ and $\itPsi$ are satisfied.
{We refer to (\ref{E:infinity type sym}) for the infinity types of $\itSigma$ and $\itPsi$.}
By Theorem \ref{T:Main}, we have
\begin{align}\label{E:Sym odd proof 3}
\frac{L(m+\tfrac{1}{2},\itSigma \times \itPsi)\cdot L(m+\tfrac{1}{2},\itSigma' \times \itPsi')}{L(m+\tfrac{1}{2},\itSigma' \times \itPsi)\cdot L(m+\tfrac{1}{2},\itSigma \times \itPsi')} \sim1.
\end{align}
We thus deduce from (\ref{E:Sym odd proof 1})-(\ref{E:Sym odd proof 3}) that
\begin{align*}
&L^{(\infty)}(m+\tfrac{1}{2},{\rm Sym}^{2r+1}\itPi\otimes\chi)\\
&\sim (2\pi\sqrt{-1})^{(r+1)m-(2r-1){\sf w}}\cdot G(\omega_\itPi\eta^{-1})^{-2r+1}\cdot G(\chi)^{r-1}\\
&\times\frac{q^{(-1)^m}({\rm Sym}^{r+1}\itPi \times \itPi'')\cdot q^{(-1)^m}({\rm Sym}^r\itPi \times \itPi'\otimes\chi)}{q^{(-1)^m}(\itPi'\times \itPi'')\cdot q^{(-1)^{m+{\sf w}}}({\rm Sym}^{r-1}\itPi\times \itPi''\otimes\eta)\cdot q^{(-1)^{m+{\sf w}}}({\rm Sym}^{r-2}\itPi\times \itPi'\otimes\chi\eta)}\\
&\times \frac{q^{\varepsilon}({\rm Sym}^{2r-1}\itPi)}{c^{(-1)^r\varepsilon}(\itPi)}\cdot \prod_{j=0}^{r-3} \frac{q^{(-1)^{j-1}\varepsilon}({\rm Sym}^{2(r-j)-3}\itPi)}{q^{(-1)^{j}\varepsilon}({\rm Sym}^{2(r-j)-3}\itPi)}\cdot \prod_{j=0}^{r-1}\frac{q^{(-1)^{j}\varepsilon}({\rm Sym}^{2(r-j)-1}\itPi)}{q^{(-1)^{j-1}\varepsilon}({\rm Sym}^{2(r-j)-1}\itPi)}.
\end{align*}
Note that ${\rm Sym}^r\itPi \times \itPi'\otimes\chi$, ${\rm Sym}^{r-1}\itPi\times \itPi''\otimes\eta$, ${\rm Sym}^{r-2}\itPi\times \itPi'\otimes\chi\eta$ are in Case (i) of Theorem \ref{T:Sym x GL_2}, and ${\rm Sym}^{r+1}\itPi \times \itPi''$, $\itPi' \times \itPi''$ are in Case (ii).
By a direct computation and the period relation (\ref{E:period relation Shimura}) for $c^\pm(\itPi'\otimes\chi)$, $c^\pm(\itPi''\otimes\eta)$, and $c^\pm(\itPi' \otimes\chi\eta)$, we have
\begin{align*}
&\frac{q^{\pm}({\rm Sym}^{r+1}\itPi \times \itPi'')\cdot q^{\pm}({\rm Sym}^r\itPi \times \itPi'\otimes\chi)}{q^{\pm}(\itPi'\times \itPi'')\cdot q^{\pm(-1)^{{\sf w}}}({\rm Sym}^{r-1}\itPi\times \itPi''\otimes\eta)\cdot q^{\pm(-1)^{{\sf w}}}({\rm Sym}^{r-2}\itPi\times \itPi'\otimes\chi\eta)}\\
&\sim(2\pi\sqrt{-1})^{-r+(3r-1){\sf w}}\cdot G(\omega_\itPi)^{3r-1}\cdot G(\eta)^{-2r+1} \cdot G(\chi)^2\cdot (c^+(\itPi)\cdot c^-(\itPi))^{r+1}.
\end{align*}
Note that 
\[
\frac{q^\pm({\rm Sym}^{2j+1}\itPi)}{q^\mp({\rm Sym}^{2j+1}\itPi)} = \left(\frac{c^\pm(\itPi)}{c^\mp(\itPi)}\right)^{j+1}.
\]
Hence we have
\begin{align*}
\prod_{j=0}^{r-3} \frac{q^{(-1)^{j-1}\varepsilon}({\rm Sym}^{2(r-j)-3}\itPi)}{q^{(-1)^{j}\varepsilon}({\rm Sym}^{2(r-j)-3}\itPi)}\cdot \prod_{j=0}^{r-1}\frac{q^{(-1)^{j}\varepsilon}({\rm Sym}^{2(r-j)-1}\itPi)}{q^{(-1)^{j-1}\varepsilon}({\rm Sym}^{2(r-j)-1}\itPi)}%& = \prod_{j=1}^{r-1}\frac{q^{-\varepsilon}({\rm Sym}^{2j-1}\itPi)\cdot q^{-\varepsilon}({\rm Sym}^{2j+1}\itPi)}{q^{\varepsilon}({\rm Sym}^{2j-1}\itPi)\cdot q^{\varepsilon}({\rm Sym}^{2j+1}\itPi)}\\
& = \begin{cases}
\displaystyle{\frac{c^\varepsilon(\itPi)}{c^{-\varepsilon}(\itPi)}}& \mbox{ if $r$ is even},\\
1 & \mbox{ if $r$ is odd}.
\end{cases}
\end{align*}
We conclude that
\begin{align*}
L^{(\infty)}(m+\tfrac{1}{2},{\rm Sym}^{2r+1}\itPi\otimes\chi) &\sim (2\pi\sqrt{-1})^{(r+1)m+r({\sf w}-1)}\cdot G(\omega_\itPi)^r\cdot G(\chi)^{r+1}\\
&\times q^{\varepsilon}({\rm Sym}^{2r-1}\itPi)\cdot c^\varepsilon(\itPi)^{r+1}\cdot c^{-\varepsilon}(\itPi)^r\\
&= (2\pi\sqrt{-1})^{(r+1)m}\cdot G(\chi)^{r+1}\cdot q^\varepsilon({\rm Sym}^{2r+1}\itPi).
\end{align*}
This completes the proof of Theorem \ref{T:Sym odd} when $n$ is odd.
\qed

\subsubsection{Symmetric even power $L$-functions for $\GL_2$: preliminaries}\label{SS:Deligne Sym even}

%In this section, we introduce some preliminaries to the proof of Theorem \ref{T:Sym odd} in \S\,\ref {SS:Deligne Sym even 2} for symmetric even power $L$-functions of modular forms. 
%A key observation due to Morimoto \cite{Morimoto2021} is that symmetric even power $L$-functions are related to Asai $L$-functions which we shall now recall.
In this section, we introduce some preliminaries for the proof of Theorem \ref{T:Sym odd} in \S\,\ref{SS:Deligne Sym even 2}, concerning symmetric even power $L$-functions of modular forms. 
A key observation due to Morimoto \cite{Morimoto2021} is that symmetric even power $L$-functions are related to Asai $L$-functions, which we now recall.
%A key observation, due to Morimoto \cite{Morimoto2021}, is that symmetric even power $L$-functions are closely related to Asai $L$-functions, which we now briefly recall.
Let $\E/\F$ be a quadratic extension of number fields and $\pi$ be a cuspidal automorphic representation of $({\rm Res}_{\E/\F}\GL_{n,\F})(\A_\F) = \GL_n(\A_\E)$.
For each place $v$ of $\F$, let $\pi_v$ be the local component of $\pi$ at $v$, which is a representation of $({\rm Res}_{\E_v/\F_v}\GL_{n,\F_v})(\F_v) = \GL_n(\E_v)$.
Let 
\[
\phi_{\pi_v} : W\!D_{\F_v} \longrightarrow (\GL_n(\C) \times \GL_n(\C))\rtimes \Z/2\Z
\] 
be the Langlands parameter of $\pi_v$, where the action of $\Z/2\Z$ is trivial if $v$ splits in $\E$ and is given by permutation of components otherwise.
Let ${\rm As}(\phi_{\pi_v}) : W\!D_{\F_v} \rightarrow \GL_{n^2}(\C)$ be the associated Asai representation. Recall that when $v$ splits in $\E$, ${\rm As}(\phi_{\pi_v})$ is given by composing $\phi_{\pi_v}$ with the tensor representation $\GL_n(\C) \times \GL_n(\C) \rightarrow \GL_{n^2}(\C)$. If $v$ is non-split in $\E$, then $\phi_{\pi_v}$ can be canonically identified with a Langlands parameter $W\!D_{\E_v} \rightarrow \GL_n(\C)$ for $\GL_n(\E_v)$ (cf.\,\cite[\S\,5]{Borel1979}). In this case the Asai representation is the multiplicative induction of $\phi_{\pi_v}$ from $W\!D_{\E_v}$ to $W\!D_{\F_v}$ recalled in (\ref{E:Asai rep}) below. We denote by ${\rm As}^+(\phi_{\pi_v}) = {\rm As}(\phi_{\pi_v})$ and ${\rm As}^-(\phi_{\pi_v}) = {\rm As}(\phi_{\pi_v})\otimes \omega_{\E_v/\F_v}$.
For a Hecke character $\chi$ of $\A_\F^\times$, we then have the twisted ($\pm$) Asai $L$-function of $\pi$ defined by the Euler product
\[
L(s,\pi,{\rm As}^\pm\otimes\chi) := \prod_v L(s,{\rm As}^\pm(\phi_{\pi_v})\otimes\chi_v), \quad {\rm Re}(s) \gg0.
\]
The Asai $L$-function admits meromorphic continuation to the whole complex plane and satisfies a functional equation.
In the special case when $\pi = {\rm BC}_{\E/\F}({\rm Sym}^{n-1}\tau)$ is the quadratic base change of the symmetric $(n-1)$-th power functorial lift of a cuspidal automorphic representation $\tau$ of $\GL_2(\A_\F)$, we have
\begin{align}\label{E:Asai and symmetric}
L(s,{\rm BC}_{\E/\F}({\rm Sym}^{n-1}\tau), {\rm As}^{(-1)^{n}}\otimes \omega_\tau^{-n+1}) = \prod_{j=0}^{n-1} L(s,{\rm Sym}^{2j}\tau \otimes \omega_\tau^{-j}\omega_{\E/\F}^{j+1}).
\end{align}
The equality follows from a comparison of the local Langlands parameters that define the local factors on both sides. 
When $v$ splits in $\E$, this follows directly from the Clebsch–Gordan formula for the tensor representation ${\rm Sym}^{n-1}(\C^2) \otimes {\rm Sym}^{n-1}(\C^2)$ of $\GL_2(\C)$.
For non-split places, we include a proof for the reader's convenience.

\begin{lemma}\label{L:Asai rep}
Let $G$ be a group with a subgroup $H$ of index $2$. Fix a set of representatives $\{1,c\}$ for $G/H$. For a representation $(\rho,V)$ of $H$, let ${\rm As}(\rho)$ be the associated Asai representation of $G$ on $V \otimes V$ defined by
\begin{align}\label{E:Asai rep}
\begin{split}
{\rm As}(\rho)(h)(v_1\otimes v_2) &:= \rho(h)v_1 \otimes \rho(c^{-1}hc)v_2, \quad h \in H\\
{\rm As}(\rho)(c)(v_1\otimes v_2) &:= \rho(c^2)v_2 \otimes v_1.
\end{split}
\end{align}
Let $\phi$ be a $2$-dimensional complex representation of $G$. Assume $\phi(c)$ is semisimple. Then we have
\[
{\rm As}\left(({\rm Sym}^n\phi)\vert_H\right) \cong \bigoplus_{i=0}^n {\rm Sym}^{2n-2i}\phi \otimes \det(\phi)^i\omega_{G/H}^i
\]
for all $n \geq 1$, where $\omega_{G/H}$ is the quadratic character of $G$ with kernel $H$.
\end{lemma}

\begin{proof}
We write ${\rm S}^i = {\rm Sym}^i$. Let $V$ be the representation space of $\phi$. After conjugation by the linear isomorphism 
\[
{\rm S}^nV \otimes {\rm S}^nV \longrightarrow {\rm S}^nV \otimes {\rm S}^nV,\quad v_1\otimes v_2 \longmapsto v_1 \otimes ({\rm S}^n\phi)(c)v_2,
\]
the Asai representation is given by
\begin{align*}
{\rm As}\left(({\rm S}^n\phi)\vert_H\right)(h)(v_1\otimes v_2) &= ({\rm S}^n\phi)(h)v_1 \otimes ({\rm S}^n\phi)(h)v_2, \quad h \in H\\
{\rm As}\left(({\rm S}^n\phi)\vert_H\right)(\rho)(c)(v_1\otimes v_2) &= ({\rm S}^n\phi)(c)v_2 \otimes ({\rm S}^n\phi)(c)v_1.
\end{align*}
Therefore, by the the Clebsch–Gordan formula for $\GL_2(\C)$, we have
\begin{align*}
{\rm As}\left(({\rm S}^n\phi)\vert_H\right)\vert_H = ({\rm S}^n\phi \otimes {\rm S}^n\phi)\vert_H
\cong \bigoplus_{i=0}^n ({\rm S}^{2n-2i}\phi \otimes \det(\phi)^i)\vert_H.
\end{align*}
For $1 \leq i \leq n$, let $V_i$ be the subspace of ${\rm S}^nV \otimes {\rm S}^nV$ that realizing the representation $({\rm S}^{2n-2i}\phi \otimes \det(\phi)^i)\vert_H$.
It remains to verify the action of $c$ on $V_i$ under the Asai representation.
Let $e_1$ and $e_2$ be linearly independent eigenvectors of $\phi(c)$ with eigenvalues $\alpha$ and $\beta$. We identify $V$ with $\C^2$ via $z_1e_1 + z_2 e_2 \mapsto (z_1,z_2)$. 
Consider the weight vectors of $V_i$ under the action of the diagonal torus of $\GL_2(\C)$.
By the explicit formula of the Clebsch–Gordan coefficients for ${\rm S}^n(\C^2) \otimes {\rm S}^n(\C^2)$, a vector in $V_i$ of weight $(2n-i-j,i+j)$ is a linear combination of tensor products of the form
\[
v_1 \otimes v_2 +(\pm1)^i v_2\otimes v_1
\]
for some weight vectors $v_1$ and $v_2$ in ${\rm S}^nV$ of weight $(n-k_1,k_1)$ and $(n-k_2,k_2)$ respectively such that $k_1+k_2 = i+j$.
For tensor product of the above form, we have
\[
{\rm As}\left(({\rm S}^n\phi)\vert_H\right)(c) (v_1 \otimes v_2 +(\pm1)^i v_2\otimes v_1) = (\pm1)^i\alpha^{2n-i-j}\beta^{i+j}\cdot (v_1 \otimes v_2 +(\pm1)^i v_2\otimes v_1).
\]
Therefore, the action of ${\rm As}\left(({\rm S}^n\phi)\vert_H\right)(c)$ on weight vectors of $V_i$ coincides with that of $({\rm Sym}^{2n-2i}\phi \otimes \det(\phi)^i\omega_{G/H}^i)(c)$. This completes the proof.
\end{proof}

Now we sketch the idea of the proof of Theorem \ref{T:Sym odd} for symmetric even power $L$-functions. Let $\itPi$ be a regular algebraic cuspidal automorphic representation of $\GL_2(\A)$ with infinity type $(\kappa;\,{\sf w})$. For an imaginary quadratic extension $\K/\Q$, consider the base change liftings $\itSigma = {\rm BC}_{\K/\Q}({\rm Sym}^r\itPi)$ and $\itSigma' = {\rm BC}_{\K/\Q}({\rm Sym}^{r-1}\itPi)$. By (\ref{E:Asai and symmetric}), we have
\[
L(s,\itSigma, {\rm As}^{(-1)^{r+1}}\otimes \omega_\itPi^{-r}) = L(s,{\rm Sym}^{2r}\itPi \otimes \omega_\itPi^{-r}\omega_{\K/\Q}^{r+1})\cdot L(s,\itSigma', {\rm As}^{(-1)^{r}}\otimes \omega_\itPi^{-r+1})
\] 
and $L(s,\itSigma', {\rm As}^{(-1)^{r1}}\otimes \omega_\itPi^{-r+1})$ is a product of twisted symmetric even power $L$-functions of $\itPi$ of smaller degrees.
By the result of Grobner and Lin recalled in Theorem \ref{T:GL} below, for $\kappa \geq 4$, the algebraicity of the critical values of the Rankin--Selberg $L$-function $L(s,\itSigma \times \itSigma')$ for $\GL_{r+1} (\A_\K) \times \GL_r(\A_\K)$ can be expressed in terms of the product of Asai $L$-values
\[
L(1,\itSigma, {\rm As}^{(-1)^{r+1}}\otimes \omega_\itPi^{-r})\cdot L(1,\itSigma', {\rm As}^{(-1)^{r}}\otimes \omega_\itPi^{-r+1}).
\] 
On the other hand, $L(s,\itSigma \times \itSigma')$ is a product of twisted symmetric odd power $L$-functions of $\itPi$ (cf.\,(\ref{E:RS = symmetric odd}) below). Therefore, if $\kappa \geq 5$, then the algebraicity of critical values of $L(s,\itSigma \times \itSigma')$ can be expressed in terms of products of Deligne's periods of $\itPi$.
By induction on $r$, one then verify that Deligne's conjecture holds for the right-half critical value $L(1,{\rm Sym}^{2r}\itPi \otimes \omega_\itPi^{-r}\omega_{\K/\Q}^{r+1})$ if $\kappa \geq 5$.
For an arbitrary finite order Hecke character $\chi$ of $\A^\times$ and critical value $L(m,{\rm Sym}^{2r}\itPi \otimes \chi)$, we need an algebraicity result for $\GL_{2r+1}(\A) \times \GL_1(\A)$ which can be used to study the algebraicity of the ratio
\[
\frac{L(m,{\rm Sym}^{2r}\itPi \otimes \chi)}{L(1,{\rm Sym}^{2r}\itPi \otimes \omega_\itPi^{-r}\omega_{\K/\Q}^{r+1})}.
\]
One possible choice is the result for twisted standard $L$-functions for $\GSp_{2r}$ over totally real number fields, as recalled in Theorem \ref{T:Liu} below, which is proved based on the doubling method. Nonetheless, the result is not applicable to the critical value $L(1,{\rm Sym}^{2r}\itPi \otimes \omega_\itPi^{-r}\omega_{\K/\Q}^{r+1})$ as the base field is $\Q$ in this case. 
To address this issue, we consider the base change liftings of symmetric powers of \( \itPi \) to a totally real number field \( \F \neq \Q \) and the CM-field \( \E = \K\F \). By applying similar arguments as above, with \( \itPi \) and \( \K/\Q \) replaced by \( {\rm BC}_{\F/\Q}\itPi \) and \( \E/\F \), and subsequently varying \( \K \) and \( \F \), we can deduce Deligne's conjecture for symmetric even power \( L \)-functions associated with \( \itPi \). We refer to \S\,\ref{SS:Deligne Sym even 2} for the details.

%In the induction step, we consider base change of the symmetric power lifts to CM or totally real fields.
%The proof is based on Theorem \ref{T:Sym odd} for symmetric odd powers, together with some algebraicity results of Liu \cite{Liu2019b} and Grobner--Lin \cite{GL2021} which we shall now recall.
Let $\F$ be a totally real number field with $[\F:\Q]=d$.
We have the following result due to Liu \cite{Liu2019b} on the algebraicity of critical values of the twisted standard $L$-functions for $\GSp_{2n}(\A_\F)$.

\begin{thm}[Liu]\label{T:Liu}
Let $\itPsi$ be a cohomological cuspidal automorphic representation of $\GSp_{2n}(\A_\F)$ such that $\itPsi_v$ is a holomorphic discrete series representation for each $v \mid \infty$. 
There exists a famility of non-zero complex numbers $(\Omega({}^\sigma\!\itPsi))_{\sigma \in {\rm Aut}(\C)}$ satisfying the following property: Let $\chi$ be a finite order Hecke character of $\A_\F^\times$ with parallel signature and $m \in \Z_{\geq 1}$ be a critical point of the twisted standard $L$-function $L(s,\itPsi , \,{\rm std}\otimes \chi)$ such that $m\neq1$ if $\F=\Q$ and $\chi^2=1$. 
For $\sigma \in {\rm Aut}(\C)$, we have
\begin{align*}
&\sigma\left(  \frac{L^S(m,\itPsi,{\rm std}\otimes\chi)}{(2\pi\sqrt{-1})^{(n+1)dm}\cdot G(\chi\vert_{\A^\times})^{n+1}\cdot \Omega(\itPsi)}\right) = \frac{L^S(m,{}^\sigma\!\itPsi,{\rm std}\otimes{}^\sigma\!\chi)}{(2\pi\sqrt{-1})^{(n+1)dm}\cdot  G({}^\sigma\!\chi\vert_{\A^\times})^{n+1}\cdot \Omega({}^\sigma\!\itPsi)}.
\end{align*}
%if $m \geq 1$; and
%\begin{align*}
%&\sigma\left(  \frac{L^S(m,\itPsi,{\rm std}\otimes\chi)}{|D_\F|^{n/2}\cdot(2\pi\sqrt{-1})^{ndm}\cdot (\sqrt{-1})^{d{\sf v}}\cdot G(\chi)^{n}\cdot p(\itPsi)}\right)\\
%& = \frac{L^S(m,{}^\sigma\!\itPsi,{\rm std}\otimes{}^\sigma\!\chi)}{|D_\F|^{n/2}\cdot(2\pi\sqrt{-1})^{ndm}\cdot (\sqrt{-1})^{d{\sf v}}\cdot  G({}^\sigma\!\chi)^{n}\cdot p({}^\sigma\!\itPsi)}
%\end{align*}
%if $m\leq 0$.
Here $S$ is a sufficiently large set of places containing the archimedean places.
Moreover, there exists a number field $\E(\itPsi)$ containing the rationality field of $\itPsi$ such that $\Omega({}^\sigma\!\itPsi) = \Omega(\itPsi)$ if $\sigma \in {\rm Aut}(\C/\E(\itPsi))$.
\end{thm}

\begin{rmk}
The result in \cite[Corollary 0.0.2]{Liu2019b} is stated for \( \F = \Q \), with \( \GSp_{2n} \) and \( {\rm Aut}(\C) \) replaced by \( \Sp_{2n} \) and \( {\rm Aut}(\C/\Q(\zeta_N)) \) for some cyclotomic field \( \Q(\zeta_N) \). The generalization to arbitrary \( \F \) is straightforward (cf.\,\cite[Theorem 28.8]{Shimura2000}). 
By considering the Hilbert--Siegel modular varieties associated with \( \GSp_{2n,\F} \) instead of the connected Shimura varieties associated with \( \Sp_{2n,\F} \), we can establish the full Galois-equivariance property over \( {\rm Aut}(\C) \) by following the standard arguments as in \cite[\S\,3]{GL2016} for the standard \( L \)-functions of unitary groups. One of the key ingredients in \textit{loc.\ cit.} is the Galois-equivariance property of holomorphic degenerate Eisenstein series, given in \cite[Proposition 3.4.1]{GL2016}. We refer to \cite[Proposition 2.6]{Chen2021d} for the analogous result for Hilbert--Siegel Eisenstein series.
Liu's breakthrough lies in the explicit computation of the archimedean doubling local zeta integrals.

%The result in \cite[Corollary 0.0.2]{Liu2019b} is stated for $\F=\Q$ with $\GSp_{2n}$ and ${\rm Aut}(\C)$ replaced by $\Sp_{2n}$ and ${\rm Aut}(\C/\Q(\zeta_N))$ for some cyclotomic field $\Q(\zeta_N)$. The generalization to arbitrary $\F$ is straightforward (cf.\,\cite[Theorem 28.8]{Shimura2000}). By considering the Hilbert--Siegel modular varieties associated to $\GSp_{2n,\F}$ instead of the connected Shimura varieties associated to $\Sp_{2n,\F}$, we can prove the full Galois-equivariance property over ${\rm Aut}(\C)$ by following the standard arguments as in \cite[\S\,3]{GL2016} for the standard $L$-functions of unitary groups. 
%One of the key ingredients in $loc.$ $cit.$ is the Galois-equivariance property of holomorphic degenerate Eisenstein series in \cite[Proposition 3.4.1]{GL2016}. We refer to \cite[Proposition 2.6]{Chen2021d} for the analogous result for Hilbert--Siegel Eisenstein series.
%The breakthrough of Liu is on the explicit computation of the archimedean doubling local zeta integrals. 
\end{rmk}

\begin{rmk}
Other results in the literature on the algebraicity of twisted standard $L$-functions were mainly for scalar-valued Hilbert--Siegel cusp forms (cf.\,\cite{Harris1981}, \cite{Sturm1981}, \cite{Mizumoto1991}, \cite{Shimura2000}, \cite{BS2000}). For our purpose here, we need to consider vector-valued Hilbert--Siegel cusp forms. We refer to \cite{Kozima2000} for certain vector-valued cases and \cite{PSS2020} for $n=2$ and ${\rm sgn}(\chi)=-1$. We also refer to the recent result \cite{HPSS2021} under some parity conditions on $m$ and $\chi$.
\end{rmk}

Let $\E/\F$ be a totally imaginary quadratic extension with non-trivial automorphism $c \in \Gal(\E/\F)$.
Let $\omega_{\E/\F}$ be the quadratic Hecke character of $\A_\F^\times$ associated to $\E/\F$ by the class field theory.
For each archimedean place $v$ of $\F$, we fix an embedding $\iota_v : \E \rightarrow \C$ extending $v$.
Let $\itSigma$ be a regular algebraic cuspidal automorphic representation of $\GL_n(\A_\E)$. For each archimedean place $v$ of $\F$, the embedding $\iota_v$ determines an isomorphism
\begin{align}\label{E:complex infinity type}
\itSigma_v \cong {\rm Ind}_{B_n(\C)}^{\GL_n(\C)}(\chi_{\kappa_{1,v}} \otimes \cdots \otimes \chi_{\kappa_{n,v}}) \otimes |\mbox{ }|_\C^{{\sf w}/2}
\end{align}
for some integers $\kappa_{1,v}>\cdots>\kappa_{n,v}$ and ${\sf w}$ such that
\[
\kappa_{1,v}\equiv \cdots \equiv \kappa_{n,v} \equiv n+1+{\sf w}\,({\rm mod}\,2).
\]
Here $B_n$ is the standard Borel subgroup of $\GL_n$ and $\chi_\kappa(z) = (z/\overline{z})^{\kappa/2}$ for $\kappa \in \Z$. We write $\underline{\kappa}_v = (\kappa_{1,v},...,\kappa_{n,v})$ and call $(\underline{\kappa}_v;\,{\sf w})$ the infinity type of $\itSigma_v$.
Note that ${\sf w}$ is independent of the choice of $v$.
Assume further that $\itSigma$ is essentially conjugate self-dual, that is, there exists a Hecke character $\chi$ of $\A_\F^\times$ with parallel signature ${\rm sgn}(\chi)$ such that
\[
\itSigma^c = \itSigma^\vee \otimes \chi\circ {\rm N}_{\E/\F}.
\]
Note that $\chi$ is necessary algebraic. After replacing $\chi$ by $\omega_{\E/\F}\cdot\chi$ if necessary, we may assume that ${\rm sgn}(\chi) = (-1)^{\sf w}$, that is, $\chi_v(-1) = (-1)^{\sf w}$ for all $v \mid \infty$.
In this case, the twisted Asai $L$-function 
\[
L(s,\itSigma,{\rm As}^{(-1)^n}\otimes\chi^{-1})
\]
is holomorphic and non-vanishing at $s=1$ (cf.\,\cite[Corollary 2.5.9]{Mok2015}).
We have the following result of Grobner and Lin \cite[Theorem C]{GL2021} on the algebraicity of the critical values of Rankin--Selberg $L$-functions for $\GL_n(\A_\E) \times \GL_{n-1}(\A_\E)$ in terms of critical values of twisted Asai $L$-functions.

\begin{thm}[Grobner--Lin]\label{T:GL}
Let $\itSigma$ and $\itSigma'$ be regular algebraic essentially conjugate self-dual cuspidal automorphic representations of $\GL_n(\A_\E)$ and $\GL_{n-1}(\A_\E)$, respectively. 
For each archimedean place $v$ of $\F$, let $(\underline{\kappa}_v;\,{\sf w})$ and $(\underline{\kappa}_v';\,{\sf w}')$ be the infinity types of $\itSigma_v$ and $\itSigma_v'$, respectively. 
Let $\chi$ and $\chi'$ be algebraic Hecke characters of $\A_\F^\times$ with parallel signatures $(-1)^{\sf w}$ and $(-1)^{{\sf w}'}$ respectively such that
\[
\itSigma^c = \itSigma^\vee \otimes \chi\circ {\rm N}_{\E/\F},\quad (\itSigma')^c = (\itSigma')^\vee \otimes \chi'\circ {\rm N}_{\E/\F}.
\]
Assume for each archimedean place $v$ of $\F$, the following conditions are satisfied:
\begin{itemize}
\item[(1)] $\underline{\kappa}_v$ and $\underline{\kappa}_v'$ are $6$-regular.
\item[(2)] $(\itSigma_v,\itSigma'_v)$ is balanced, that is, 
\[
\kappa_{1,v} > -\kappa'_{n-1,v} > \kappa_{2,v} > -\kappa_{n-2,v}'>\cdots>-\kappa_{1,v}'>\kappa_{n,v}.
\]
\end{itemize}
Then, for a critical point $m+\tfrac{1}{2} \in \Z+\tfrac{1}{2}$ for $L(s, \itSigma \times \itSigma')$ and $\sigma \in {\rm Aut}(\C/\E^{\Gal})$, we have
\begin{align*}
\sigma\left( \frac{L^{(\infty)}(m+\tfrac{1}{2}, \itSigma \times \itSigma')}{(2\pi\sqrt{-1})^{n(n-1)dm}\cdot q(\itSigma \times \itSigma')} \right) = \frac{L^{(\infty)}(m+\tfrac{1}{2}, {}^\sigma\!\itSigma \times {}^\sigma\!\itSigma')}{(2\pi\sqrt{-1})^{n(n-1)dm}\cdot q({}^\sigma\!\itSigma \times {}^\sigma\!\itSigma')}.
\end{align*}
Here $\E^{\Gal}$ is the Galois closure of $\E$ over $\Q$ and
\begin{align*}
q(\itSigma \times \itSigma') &= (2\pi\sqrt{-1})^{n(n-1)d({\sf w}+{\sf w}')/2-n(n+1)d/2}\cdot G(\chi\vert_{\A^\times})^{n(n-1)/2}\cdot G((\chi')\vert_{\A^\times})^{(n-1)(n-2)/2}\cdot G(\omega_{\itSigma'}\vert_{\A^\times})\\
&\times L^{(\infty)}(1,\itSigma,{\rm As}^{(-1)^n}\otimes\chi^{-1})\cdot L^{(\infty)}(1,\itSigma',{\rm As}^{(-1)^{n-1}}\otimes(\chi')^{-1}).
\end{align*}
\end{thm}

\begin{rmk}
In \cite{GL2021}, the result is proved when $\itSigma$ and $\itSigma'$ are conjugate self-dual. The general case can be proved in a similar way. For the (conjectural) ${\rm Aut}(\C)$-equivariant version, one needs to modify the period $q(\itSigma \times \itSigma')$ by multiplying by certain integral powers of $G(\omega_{\E/\F})$ and by the square-root of the absolute discriminant of $\F/\Q$. We refer to \cite[Corollary 7.9]{Chen2021} for the case $\GL_3 \times \GL_2$ based on explicit computations.
\end{rmk}

\subsubsection{\color{black}Symmetric even power $L$-functions for $\GL_2$: proof}\label{SS:Deligne Sym even 2}

Let \( \itPi \) be a regular algebraic cuspidal automorphic representation of \( \GL_2(\A) \) with infinity type \( (\kappa; {\sf w}) \). 
If \( \itPi \) is of CM-type, then Conjecture \ref{C:Deligne Sym} is known to hold (cf.\,\cite[\S\,4]{RS2007c}). 
For the remainder of the discussion, it is assumed that \( \itPi \) is non-CM and that \( \kappa \geq 5 \). 
Note that the symmetric power lifts of \( \itPi \) are all cuspidal since $\itPi$ is non-CM.

The proof of Theorem \ref{T:Sym odd} for ${\rm Sym}^n\itPi$ with \( n = 2r \) is now presented. 
The assertion is established by induction on \( r \). 
It should be noted that Conjecture \ref{C:Deligne Sym} remains meaningful even when \( n = 0 \), in which case it reduces to Deligne's conjecture for finite order Hecke characters of \( \A^\times \), a result that is well-known.
Let \( r \geq 1 \), and assume that Conjecture \ref{C:Deligne Sym} holds for all \( r' \) such that \( 0 \leq r' \leq r - 1 \). 
The proof begins with the following lemma, which follows as a consequence of Theorem \ref{T:Liu} and the existence of weak automorphic descent from \( \GL_{2r+1}(\A_\F) \) to \( \GSp_{2r}(\A_\F) \) for the symmetric \( 2r \)-th representation of \( \GL_2(\C) \).

%Let $\itPi$ be a regular algebraic cuspidal automorphic representation of $\GL_2(\A)$ with infinity type $(\kappa;\,{\sf w})$. 
%If $\itPi$ is of CM-type, then Conjecture \ref{C:Deligne Sym} holds (cf.\,\cite[\S\,4]{RS2007c}).
%We assume $\itPi$ is non-CM and $\kappa \geq 5$. 
%Then the symmetric power lifts of $\itPi$ are all cuspidal.
%Now we begin the proof of Theorem \ref{T:Sym odd} when $n=2r$. 
%We prove the assertion by induction on $r$. 
%Note that Conjecture \ref{C:Deligne Sym} makes sense even when $n=0$. In this case, it is just Deligne's conjecture for finite order Hecke characters of $\A^\times$, which is well-known.
%Ler $r \geq 1$ and assume Conjecture \ref{C:Deligne Sym} holds for $r-1 \geq r' \geq 0$.
%We begin with the following lemma which is a consequence of Theorem \ref{T:Liu} and the existence of weak automorphic descent 
%from $\GL_{2r+1}(\A_\F)$ to $\GSp_{2r}(\A_\F)$ for the symmetric $2r$-th representation of $\GL_2(\C)$.

\begin{lemma}\label{P:auxiliary}
Assume $\F/\Q$ is cyclic of prime degree and the base change ${\rm BC}_{\F/\Q}({\rm Sym}^{2r}\itPi)$ to $\GL_{2r+1}(\A_\F)$ is cuspidal.
Let $\chi_1,\chi_2$ be finite order Hecke characters of $\A^\times$ with $\chi_{1,\infty}(-1) = \chi_{2,\infty}(-1)$.
For a critical point $m \in \Z_{\geq 1}-r{\sf w}$ for $L(s,{\rm Sym}^{2r}\itPi \otimes \chi_1)$ and $\sigma \in {\rm Aut}(\C)$, we have
\begin{align*}
\prod_{i=1}^d\sigma \left( \frac{L^{(\infty)}(m,{\rm Sym}^{2r}\itPi \otimes\chi_1\omega_{\F/\Q}^i)\cdot G(\chi_1)^{-r-1}}{L^{(\infty)}(m,{\rm Sym}^{2r}\itPi \otimes\chi_2\omega_{\F/\Q}^i)\cdot G(\chi_2)^{-r-1}}\right) = \prod_{i=1}^d\frac{L^{(\infty)}(m,{\rm Sym}^{2r}{}^\sigma\!\itPi \otimes{}^\sigma\!\chi_1\omega_{\F/\Q}^i)\cdot G({}^\sigma\!\chi_1)^{-r-1}}{L^{(\infty)}(m,{\rm Sym}^{2r}{}^\sigma\!\itPi \otimes{}^\sigma\!\chi_2\omega_{\F/\Q}^i)\cdot G({}^\sigma\!\chi_2)^{-r-1}}.
\end{align*}
Here $\omega_{\F/\Q}$ is a non-trivial Hecke character of $\A^\times$ trivial on ${\rm N}_{\F/\Q}(\A_\F^\times)$. 
\end{lemma}

\begin{proof}
The existence of base change lifting was proved by Arthur and Clozel \cite{AC1989}.
Since $\itPi^\vee = \itPi \otimes \omega_\itPi^{-1}$, we see that ${\rm Sym}^{2r}\itPi \otimes \omega_{\itPi}^{-r}$ is self-dual with trivial central character.
Therefore, ${\rm BC}_{\F/\Q}({\rm Sym}^{2r}\itPi)\otimes \omega_{\itPi}^{-r}\circ{\rm N}_{\F/\Q}$ is also self-dual with trivial central character.
Moreover, we have the factorization of the symmetric square $L$-function of ${\rm BC}_{\F/\Q}({\rm Sym}^{2r}\itPi)\otimes \omega_{\itPi}^{-r}\circ{\rm N}_{\F/\Q}$:
\begin{align*}
L(s,{\rm BC}_{\F/\Q}({\rm Sym}^{2r}\itPi)\otimes \omega_{\itPi}^{-r}\circ{\rm N}_{\F/\Q},{\rm Sym}^2 ) = \prod_{i=0}^r L(s,{\rm BC}_{\F/\Q}({\rm Sym}^{4i}\itPi)\otimes \omega_{\itPi}^{-2i}\circ{\rm N}_{\F/\Q}).
\end{align*}
Since the twisted symmetric power $L$-functions are holomorphic and non-zero at $s=1$ by \cite[Corollary 7.1.5]{BLGG2011}, we see that $L(s,{\rm BC}_{\F/\Q}({\rm Sym}^{2r}\itPi)\otimes \omega_{\itPi}^{-r}\circ{\rm N}_{\F/\Q},{\rm Sym}^2 )$ has a pole at $s=1$ contributed from the factor $i=0$.
By Arthur's multiplicity formula \cite[Theorem 1.5.2]{Arthur2013} for $\Sp_{2r}(\A_\F)$, the self-dual cuspidal automorphic representation ${\rm BC}_{\F/\Q}({\rm Sym}^{2r}\itPi)\otimes \omega_{\itPi}^{-r}\circ{\rm N}_{\F/\Q}$ descends, with respect to the standard representation ${\rm SO}_{2r+1}(\C) \rightarrow \GL_{2r+1}(\C)$, to a discrete automorphic representation $\itPsi'$ of $\Sp_{2r}(\A_\F)$.
Since ${\rm BC}_{\F/\Q}({\rm Sym}^{2r}\itPi)\otimes \omega_{\itPi}^{-r}\circ{\rm N}_{\F/\Q}$ is cohomological, the
descents at archimedean places are discrete series representations (cf.\,\cite[Theorem 5.6-(1)]{RS2018}).
Note that the condition ${\rm BC}_{\F/\Q}({\rm Sym}^{2r}\itPi)$ is cuspidal is equivalent to saying that the global Arthur parameter of $\itPsi'$ is stable.
Therefore, by the multiplicity formula again, we can choose $\itPsi'$ freely from its global Arthur packet, that is, there is no global obstruction.
In particular, we may assume that $\itPsi_v'$ is a holomorphic discrete series representation for all $v \mid \infty$.
Also it follows from the temperedness of $\itPsi_v'$ for $v \mid \infty$ and \cite[Theorem 4.3]{Wallach1984} that $\itPsi'$ must be cuspidal.
By the results of Patrikis \cite[Corollary 3.1.6 and Proposition 3.1.14]{Patrikis2019}, there exists a cohomological cuspidal automorphic representation $\itPsi$ of $\GSp_{2r}(\A_\F)$ satisfying the following conditions:
\begin{itemize}
\item For all $v \mid \infty$, $\itPsi_v$ is a holomorphic discrete series representation such that $\itPsi_v \vert_{\Sp_{2r}(\F_v)}$ contains $\itPsi_v'$.
\item $\itPsi'$ is a weak functorial lift of $\itPsi$ with respect to the $L$-homomorphism ${\rm GSpin}_{2r+1}(\C) \rightarrow {\rm SO}_{2r+1}(\C)$. 
\end{itemize}
In particular, for any finite order Hecke character $\chi$ of $\A^\times$, the set of critical points for $L(s,\itPsi, {\rm std}\otimes \chi\circ{\rm N}_{\F/\Q})$ and $L(s,{\rm Sym}^{2r}\itPi \otimes \omega_{\itPi}^{-r}\chi)$ are equal and we have
\begin{align*}
L^S(s,\itPsi, {\rm std}\otimes \chi\circ{\rm N}_{\F/\Q}) 
& = L^S(s,{\rm BC}_{\F/\Q}({\rm Sym}^{2r}\itPi)\otimes (\omega_{\itPi}^{-r}\chi)\circ{\rm N}_{\F/\Q})\\
& = \prod_{i=1}^d L^S(s,{\rm Sym}^{2r}\itPi \otimes \omega_{\itPi}^{-r}\chi\omega_{\F/\Q}^i).
\end{align*}
for some sufficiently large finite set $S$ of places containing $\infty$.
The assertion then follows immediately from Theorem \ref{T:Liu} with $\chi = \chi_1\omega_\itPi^r|\mbox{ }|_\A^{-r{\sf w}}$ and $\chi = \chi_2\omega_\itPi^r|\mbox{ }|_\A^{-r{\sf w}}$. Note that $\F \neq \Q$ by our assumption.
This completes the proof.
\end{proof}

Assume $\F$ is chosen so that the assumptions in Lemma \ref{P:auxiliary} are satisfied with odd prime degree $d=[\F:\Q]$, and let $\K$ be an imaginary quadratic extension of $\Q$ such that the base change ${\rm BC}_{\K\F/\Q}({\rm Sym}^{r}\itPi)$ and ${\rm BC}_{\K\F/\Q}({\rm Sym}^{r-1}\itPi)$ are cuspidal.
We apply Theorem \ref{T:GL} to $\E=\E^{\Gal} = \K\F$ and the regular algebraic essentially conjugate self-dual cuspidal automorphic representations $\itSigma$ and $\itSigma'$ of $\GL_{r+1}(\A_{\K\F})$ and $\GL_r(\A_{\K\F})$, respectively, defined by
\[
\itSigma = {\rm BC}_{\K\F/\Q}({\rm Sym}^r \itPi),\quad \itSigma' = {\rm BC}_{\K\F/\Q}({\rm Sym}^{r-1} \itPi).
\]
More precisely, we have $\itSigma^c = \itSigma$, $(\itSigma')^c = \itSigma'$, and 
\[
\itSigma = \itSigma^\vee \otimes \omega_\itPi^r\circ {\rm N}_{\K\F/\Q},\quad \itSigma' = (\itSigma')^\vee \otimes \omega_\itPi^{r-1}\circ {\rm N}_{\K\F/\Q},
\]
and the infinity types (\ref{E:complex infinity type}) of $\itSigma_v$ and $\itSigma_v'$ for any archimedean place $v$ of $\F$ are given by
\[
\left( ((r-2i)(\kappa-1))_{0 \leq i \leq r};\,r{\sf w}\right),\quad \left( ((r-2i-1)(\kappa-1))_{0 \leq i \leq r-1};\,(r-1){\sf w}\right).
\]
Thus the assumptions (1) and (2) in Theorem \ref{T:GL} are satisfied since $\kappa \geq 4$.
By the adjointness property of the base change liftings for $\K\F/\F$ and $\F/\Q$ in \cite[Proposition 3.1-(3) and (5)]{PR1999}, for any cuspidal automorphic representation $\pi_i$ of $\GL_{n_i}(\A)$ for $i=1,2$, we have
\begin{align*}%\label{E:adjointness}
%\begin{split}
L(s,{\rm BC}_{\K\F/\Q}(\pi_1)\times{\rm BC}_{\K\F/\Q}(\pi_2)) & = L(s, {\rm BC}_{\F/\Q}(\pi_1) \times {\rm BC}_{\F/\Q}(\pi_2)) \cdot L(s, {\rm BC}_{\F/\Q}(\pi_1) \times {\rm BC}_{\F/\Q}(\pi_2)\otimes \omega_{\K\F/\F})\\
& = \prod_{i=1}^d L(s,\pi_1 \times \pi_2 \otimes \omega_{\F/\Q}^i)\cdot L(s,\pi_1 \times \pi_2 \otimes \omega_{\K/\Q}\omega_{\F/\Q}^i).
%\end{split}
\end{align*}
Here in the first equality we use ${\rm BC}_{\K\F/\Q} = {\rm BC}_{\K\F/\F}\circ{\rm BC}_{\F/\Q}$, and in the second one we use $\omega_{\K/\Q}\circ{\rm N}_{\F/\Q} = \omega_{\K\F/\F}$ which follows easily from the assumption that $d$ is odd.
Note that 
\[
{\rm Sym}^{r}\itPi \boxtimes {\rm Sym}^{r-1}\itPi = \bigboxplus_{j=0}^{r}{\rm Sym}^{2(r-j)-1}\itPi \otimes \omega_\itPi^j.
\]
Therefore, let $\pi_1 = {\rm Sym}^{r}\itPi$ and $\pi_2 = {\rm Sym}^{r-1}\itPi$%in (\ref{E:adjointness})
, we then have
\begin{align}\label{E:RS = symmetric odd}
L(s,\itSigma \times \itSigma') %&= \prod_{i=1}^dL(s,{\rm Sym}^r \itPi \times {\rm Sym}^{r-1} \itPi\otimes \omega_{\F/\Q}^i)\cdot L(s,{\rm Sym}^r \itPi \times {\rm Sym}^{r-1} \itPi\otimes \omega_{\K/\Q}\omega_{\F/\Q}^i)\\
& = \prod_{i=1}^d\prod_{j=0}^{r-1} L(s,{\rm Sym}^{2(r-j)-1}\itPi \otimes \omega_\itPi^j\omega_{\F/\Q}^i)\cdot L(s,{\rm Sym}^{2(r-j)-1}\itPi \otimes \omega_\itPi^j\omega_{\K/\Q}\omega_{\F/\Q}^i).
\end{align}
As we have proved in \S\,\ref{SS:Deligne Sym odd}, Conjecture \ref{C:Deligne Sym} holds for symmetric odd power $L$-functions of $\itPi$ when $\kappa \geq 5$. 
Therefore, for any critical point $m+\tfrac{1}{2} \in \Z+\tfrac{1}{2}$ for $L(s,\itSigma \times \itSigma')$, we have
\begin{align*}
L^{(\infty)}(m+\tfrac{1}{2},\itSigma \times \itSigma') &\sim 
\left((2\pi\sqrt{-1})^{r(r+1)m+r(r-1)(r+1){\sf w}/3}\cdot G(\omega_\itPi)^{r(r-1)(r+1)/3}\cdot G(\omega_{\K/\Q})^{r(r+1)/2} \right)^d\\
&\times \left ( \prod_{j=0}^{r-1}q^+({\rm Sym}^{2(r-j)-1}\itPi)\cdot q^-({\rm Sym}^{2(r-j)-1}\itPi) \right)^d.
\end{align*}
%For any cuspidal automorphic representation $\pi$ of $\GL_2(\A_\F)$ and $N \geq 1$, by Lemma \ref{L:Asai rep} we have
%\begin{align*}%\label{E:Asai and symmetric}
%L(s,{\rm BC}_{\K\F/\F}({\rm Sym}^N\pi), {\rm As}^{(-1)^{N+1}}\otimes \omega_\pi^{-N}) = \prod_{j=0}^N L(s,{\rm Sym}^{2j}\pi \otimes \omega_\pi^{-j}\omega_{\K\F/\F}^{j+1}).
%\end{align*}
By (\ref{E:Asai and symmetric}) with $\tau = {\rm BC}_{\F/\Q}(\itPi)$ and $n=r+1$ or $n=r$, together with the adjointness property of the base change liftings for $\F/\Q$ in \cite[Proposition 3.1-(3) and (5)]{PR1999} and $\omega_{\K/\Q}\circ{\rm N}_{\F/\Q} = \omega_{\K\F/\F}$, we obtain
\begin{align*}
L(s,\itSigma,{\rm As}^{(-1)^{r+1}}\otimes \omega_\itPi^{-r}\circ{\rm N}_{\F/\Q}) & =L(s,\itSigma',{\rm As}^{(-1)^{r}}\otimes \omega_\itPi^{-r+1}\circ{\rm N}_{\F/\Q})\\
&\times \prod_{i=1}^d L(s,{\rm Sym}^{2r}\itPi \otimes \omega_\itPi^{-r}\omega_{\K/\Q}^{r+1}\omega_{\F/\Q}^i),\\
L(s,\itSigma',{\rm As}^{(-1)^{r}}\otimes \omega_\itPi^{-r+1}\circ{\rm N}_{\F/\Q}) & = \prod_{i=1}^d\prod_{j=0}^{r-1} L(s,{\rm Sym}^{2j}\itPi \otimes \omega_\itPi^{-j}\omega_{\K/\Q}^{j+1}\omega_{\F/\Q}^i).
\end{align*}
By the induction hypothesis, Conjecture \ref{C:Deligne Sym} holds for ${\rm Sym}^{2j}\itPi$ with $0 \leq j \leq r-1$. 
Since $\omega_{\C/\R}$ is odd, $s=1$ is always a right-half critical point for the twisted symmetric even power $L$-functions of $\itPi$ appearing in the above decompositions. 
We thus have $\pm = +$ in Conjecture \ref{C:Deligne Sym}, and therefore
\begin{align*}
L^{(\infty)}(1,\itSigma',{\rm As}^{(-1)^{r}}\otimes \omega_\itPi^{-r+1}\circ{\rm N}_{\F/\Q})^2& \sim \left( (2\pi\sqrt{-1})^{r(r+1)/2-r(r-1)(r+1){\sf w}/3}\cdot G(\omega_\itPi)^{-r(r-1)(r+1)/3}\right)^{2d}\\
&\times \left( \prod_{j=0}^{r-1}q^+({\rm Sym}^{2j}\itPi)\right)^{2d}.
\end{align*}
Here we have use the obvious facts that $G(\omega_{\K/\Q})^2 \sim 1$ and $G(\omega_{\F/\Q})^d \sim 1$.
Combining with Theorem \ref{T:GL} for $L^{(\infty)}(m+\tfrac{1}{2},\itSigma \times \itSigma')$ (take $m+\tfrac{1}{2}$ to be non-central), we deduce that
\begin{align*}
&\prod_{i=1}^d L^{(\infty)}(1,{\rm Sym}^{2r}\itPi \otimes \omega_\itPi^{-r}\omega_{\K/\Q}^{r+1}\omega_{\F/\Q}^i)\\ &\sim_{\K\F}
\left( (2\pi\sqrt{-1})^{-r(r+1){\sf w}/2-(r+1)(r-2)/2}\cdot G(\omega_\itPi)^{-r(r+1)/2}\cdot G(\omega_{\K/\Q})^{r(r+1)/2} \right)^d\\
&\times \left( \prod_{j=0}^{r-1}\frac{q^+({\rm Sym}^{2j+1}\itPi)\cdot q^-({\rm Sym}^{2j+1}\itPi)}{q^+({\rm Sym}^{2j}\itPi)^2} \right)^d.
\end{align*}
Here $\sim_{\K\F}$ means the ratio is equivariant under ${\rm Aut}(\C/\K\F)$ (recall the notation in \S\,\ref{SS:notation}, here we replace ${\rm Aut}(\C)$ by ${\rm Aut}(\C/\K\F)$).
By the formula for the periods $q^\pm$ of symmetric powers, the right-hand side of the above relation is equal to 
\[
\left( (2\pi\sqrt{-1})^{(r+1)(1-r{\sf w})}\cdot G(\omega_\itPi^{-r})^{r+1}\cdot q^+({\rm Sym}^{2r}\itPi)\cdot G(\omega_{\K/\Q})^{r(r+1)/2} \right)^d.
\]
On the other hand, for any finite order Hecke character $\chi$ of $\A^\times$ such that $\chi_\infty(-1) = (-1)^{r+1}$, by Lemma \ref{P:auxiliary} we have
\[
\prod_{i=1}^d L^{(\infty)}(1,{\rm Sym}^{2r}\itPi \otimes \omega_\itPi^{-r}\omega_{\K/\Q}^{r+1}\omega_{\F/\Q}^i) \sim \prod_{i=1}^d L^{(\infty)}(1,{\rm Sym}^{2r}\itPi \otimes \omega_\itPi^{-r}\chi\omega_{\F/\Q}^i) \cdot G(\chi\omega_{\K/\Q}^{r+1})^{-r-1}. 
\]
Note that $G(\omega_{\K/\Q})\sim_\K 1$ and $G(\omega_{\F/\Q})\sim_\F 1$. We thus conclude that
\begin{align}\label{E:Sym even proof 1}
\begin{split}
\prod_{i=1}^d \frac{L^{(\infty)}(1,{\rm Sym}^{2r}\itPi \otimes \omega_\itPi^{-r}\chi\omega_{\F/\Q}^i)}{(2\pi\sqrt{-1})^{(r+1)(1-r{\sf w})}\cdot G(\omega_\itPi^{-r}\chi\omega_{\F/\Q}^i)^{r+1}\cdot q^+({\rm Sym}^{2r}\itPi)} \sim_{\K\F} 1.
\end{split}
\end{align}
As we have explained in the proof of Lemma \ref{P:auxiliary} (with base field $\Q$ this time), there exists a cohomological cuspidal automorphic representation $\itPsi$ of $\GSp_{2r}(\A)$ such that $\itPsi_\infty$ is a holomorphic discrete series representation and satisfies the following conditions:
For any Hecke character $\chi$ of $\A^\times$, we have
\begin{itemize}
\item $L(s,\itPsi_\infty,{\rm std}\otimes\chi_\infty) = L(s,{\rm Sym}^{2r}\itPi_\infty \otimes \omega_{\itPi,\infty}^{-r}\chi_\infty)$.
\item $L^S(s,\itPsi,{\rm std}\otimes\chi) = L^S(s,{\rm Sym}^{2r}\itPi \otimes \omega_{\itPi}^{-r}\chi)$ for some sufficiently large set $S$ of places containing $\infty$.
\end{itemize}
We fix one such $\itPsi$, and let $(\Omega({}^\sigma\!\itPsi))_{\sigma \in {\rm Aut}(\C)}$ and $\E(\itPsi)$ be as in the statement of Theorem \ref{T:Liu}.
We enlarge $\E(\itPsi)$ so that it also contains $\Q(\itPi)$.
Assume $\chi$ is chosen so that $\chi_\infty(-1)=(-1)^{r+1}$ and $\chi^2 \notin \<\omega_{\F/\Q}\>$. By Theorem \ref{T:Liu} and (\ref{E:Sym even proof 1}), we have
\begin{align}\label{E:Sym even proof 2}
\left( \frac{\Omega(\itPsi)}{(2\pi\sqrt{-1})^{-r(r+1){\sf w}}\cdot q^+({\rm Sym}^{2r}\itPi)} \right)^d \sim_{\K\F} 1.
\end{align}
Choose another imaginary quadratic extension \(\K'/\Q\) and a totally real extension \(\F'/\Q\) with \([\F':\Q] = [\F:\Q]\), such that \(\K\F\E(\itPsi) \cap \K'\F'\E(\itPsi) = \E(\itPsi)\) and the base change liftings of \({\rm Sym}^r\itPi\) or \({\rm Sym}^{r-1}\itPi\) to \(\F'\) or \(\K'\F'\) are all cuspidal. Under these assumptions, equation \((\ref{E:Sym even proof 2})\) also holds for \({\rm Aut}(\C/\K'\F')\), allowing us to deduce that its left-hand side actually belongs to \(\E(\itPsi)\). 
Assume further that \(\K\F \cap \E(\itPsi) = \Q\). Using the same arguments as in \((\ref{E:RS proof 5})\), it follows that the equivalence \((\ref{E:Sym even proof 2})\) actually holds under \({\rm Aut}(\C)\) rather than just \({\rm Aut}(\C/\K\F)\). By choosing \(d = 3, 5\) (or any two distinct odd prime numbers), we then obtain the period relation:
\[
\frac{\Omega(\itPsi)}{(2\pi\sqrt{-1})^{-r(r+1){\sf w}} \cdot q^+({\rm Sym}^{2r}\itPi)} \sim 1.
\]
From Theorem \ref{T:Liu} again, it then follows that Conjecture \ref{C:Deligne Sym} holds for \(L^{(\infty)}(m,{\rm Sym}^{2r}\itPi \otimes \omega_\itPi^{-r}\chi)\), for any finite order Hecke character \(\chi\) of \(\A^\times\) and any critical point \(m \in \Z_{\geq 1}\) for \(L(s,{\rm Sym}^{2r}\itPi \otimes \omega_\itPi^{-r}\chi)\) such that \(m \neq 1\) if \(\chi_\infty(-1) = (-1)^{r+1}\) and \(\chi^2 = 1\). 
Now consider a quadratic Hecke character \(\chi\) of \(\A^\times\) satisfying \(\chi_\infty(-1) = (-1)^{r+1}\). To demonstrate that Conjecture \ref{C:Deligne Sym} also holds for \(L^{(\infty)}(1,{\rm Sym}^{2r}\itPi \otimes \omega_\itPi^{-r}\chi)\), we consider equation \((\ref{E:Sym even proof 1})\). 
Since \(\omega_{\F/\Q}\) has odd prime order \(d\), Conjecture \ref{C:Deligne Sym} holds for the critical values \(L^{(\infty)}(1,{\rm Sym}^{2r}\itPi \otimes \omega_\itPi^{-r}\chi\omega_{\F/\Q}^i)\) with \(1 \leq i \leq d-1\). It then follows from \((\ref{E:Sym even proof 1})\) that:
\[
\frac{L^{(\infty)}(1,{\rm Sym}^{2r}\itPi \otimes \omega_\itPi^{-r}\chi)}{(2\pi\sqrt{-1})^{(r+1)(1-r{\sf w})} \cdot G(\omega_\itPi^{-r}\chi)^{r+1} \cdot q^+({\rm Sym}^{2r}\itPi)} \sim_{\K\F} 1.
\]
By varying $\K$ and $\F$ as above, we see that Conjecture \ref{C:Deligne Sym} holds for $L^{(\infty)}(1,{\rm Sym}^{2r}\itPi \otimes \omega_\itPi^{-r}\chi)$.
As for the left-half critical points (critical points $m \leq 0$ for $L(s,{\rm Sym}^{2r}\itPi \otimes \omega_\itPi^{-r}\chi)$), it follows from the global functional equation for $\GL_{2r+1}(\A)$ applied to the cuspidal automorphic representation ${\rm Sym}^{2r}\itPi$.
Indeed, by \cite[Theorem 3.2]{Chen2023} (see also \cite[Proposition A.4]{Chen2021}), we have
\begin{align*}
&L^{(\infty)}(m,{\rm Sym}^{2r}\itPi\otimes \omega_\itPi^{-r}\chi) \\
&\sim \gamma(m,{\rm Sym}^{2r}\itPi_\infty \otimes \omega_{\itPi,\infty}^{-r}\chi_\infty,\psi_\infty)\cdot G(\chi)^{2r+1}\cdot L^{(\infty)}(1-m,{\rm Sym}^{2r}\itPi\otimes \omega_\itPi^{-r}\chi^{-1})
\end{align*}
for all left-half critical points $m$ for $L(s,{\rm Sym}^{2r}\itPi \otimes \omega_\itPi^{-r}\chi)$.
By a direct computation using the formulas of archimedean local factors (cf.\,\cite[(3.6) and (3.7)]{Knapp1994}), we have
\[
\gamma(m,{\rm Sym}^{2r}\itPi_\infty \otimes \omega_{\itPi,\infty}^{-r}\chi_\infty,\psi_\infty) \in (2\pi\sqrt{-1})^{(2r+1)m-r-1} \cdot \Q^\times.
\]
This completes the proof of Theorem \ref{T:Sym odd}.
\qed

\appendix
\section{Minimal ${\rm O}_n(\R)$-types and cohomological contributions}\label{S:appendix}

In this appendix, we prove that the minimal $K$-types contribute non-trivially to the relevant relative Lie algebra cohomology in the bottom degree. As explained before Theorem~\ref{T:nonvanishing}, this reduces the non-vanishing of the pairing $Z(\cdot,\cdot)$ to a statement at the level of cochain complexes. We include the proof of Proposition~\ref{L:cohomology} as it may be of independent interest.

Let $\frak{g}=\frak{g}_n$ be the complexified Lie algebra of $\GL_n(\R)$ and identify it with the algebra of $n$ by $n$ complex matrices. 
Let $\frak{t}\subset\frak{g}$, $\frak{u}\subset \frak{g}$, and $\frak{so}_n\subset \frak{g}$ be the complexified Lie algebras of $T_n(\R)$, $U_n(\R)$, and ${\rm SO}_n(\R)$ respectively. 
Define $k_\infty \in {\rm SU}_n(\R)$ by
\[
k_\infty = \begin{cases}
\displaystyle{{\rm diag}\left(\frac{1}{\sqrt{2}}\bp 1 & \sqrt{-1} \\  \sqrt{-1}  & 1 \ep,...,\frac{1}{\sqrt{2}}\bp 1 & \sqrt{-1} \\  \sqrt{-1}  & 1 \ep\right)} & \mbox{ if $n$ is even},\\
\displaystyle{{\rm diag}\left(\frac{1}{\sqrt{2}}\bp 1 & \sqrt{-1} \\  \sqrt{-1}  & 1 \ep,...,\frac{1}{\sqrt{2}}\bp 1 & \sqrt{-1} \\  \sqrt{-1}  & 1 \ep,1\right)} & \mbox{ if $n$ is odd}.
\end{cases}
\]
For $1 \leq i,j \leq n$, let $e_{ij} \in \frak{g}$ be the matrix with $1$ in the $(i,j)$-entry and zeros otherwise, and define $X_{ij} \in \frak{g}$ by
\[
X_{ij} = k_\infty \cdot e_{ij}\cdot k_\infty^{-1}.
\]
Let $\frak{l} = k_\infty\cdot \frak{t}\cdot k_\infty^{-1}$ and $\frak{n} \subset \frak{g}$ be the nilpotent subalgebra defined by
\[
\frak{n} = \bigoplus_{i=1}^{\lfloor\frac{n}{2}\rfloor} \C\cdot X_{2i-1,2i} \bigoplus_{{1 \leq i \leq \lfloor\frac{n}{2}\rfloor} \atop{ 2i<j\leq n}}\left( \C\cdot X_{j,2i}\oplus \C\cdot {}^tX_{j,2i}\right).
\]
Then 
\[
\frak{b} := \frak{l}\ltimes\frak{n}
\]
is a $\theta$-stable Borel subalgebra of $\frak{g}$, where $\theta$ is the Cartan involution given by $\theta(x) = -{}^tx$.
Let 
\[
\frak{b}^c = \frak{b}\cap \frak{so}_n,\quad \frak{l}^c = \frak{l}\cap \frak{so}_n,\quad \frak{n}^c = \frak{n}\cap \frak{so}_n.
\]
Then $\frak{b}^c = \frak{l}^c\ltimes\frak{n}^c$ is a Borel subalgebra of $\frak{so}_n$. We identify $(\frak{l}^c)^*$ with $\C^{\lfloor\frac{n}{2}\rfloor}$ by
\[
\sum_{i=1}^{\lfloor\frac{n}{2}\rfloor}z_i\cdot (X_{2i-1,2i-1}^*-X_{2i,2i}^*) \longmapsto (z_1,...,z_{\lfloor\frac{n}{2}\rfloor}),
\]
where the superscript $*$ refers to the dual basis. Then the set of dominant integral weights with respect to $\frak{b}^c$ consisting of $\underline{\kappa} = (\kappa_1,...,\kappa_{\lfloor\frac{n}{2}\rfloor}) \in \Z^{\lfloor\frac{n}{2}\rfloor}$ such that
\[
\kappa_1 \geq \cdots\geq  \kappa_{\lfloor\frac{n}{2}\rfloor-1} \geq 
\begin{cases}
|\kappa_{\lfloor\frac{n}{2}\rfloor}| & \mbox{ if $n$ is even},\\
\kappa_{\lfloor\frac{n}{2}\rfloor} \geq 0 & \mbox{ if $n$ is odd}.
\end{cases}
\]
For $\underline{\kappa} \in \Z^{\lfloor\frac{n}{2}\rfloor}$ dominant with respect to $\frak{b}^c$, we denote by $\tau_{\underline{\kappa}}$ the irreducible representation of ${\rm SO}_n(\R)$ with highest weight $\underline{\kappa}$.
For instance, $\tau_{\underline{e}}$ appears in $\extp^{b_n}(\frak{g}/\frak{k})$ with multiplicity one (cf.\,\cite[Proposition 6.1.3]{Mahnkopf2005}), where $\frak{k}=\frak{k}_n$ is the complexified Lie algebra of $K_n^\circ = \R_+\cdot {\rm SO}_n(\R)$ and
\[
\underline{e} = \underline{e}(n) = (n,n-2,...,n-2\lfloor\tfrac{n}{2}\rfloor+2).
\]
For $\mu \in X^+(T_n)$, let $\underline{\kappa} \in \Z^{\lfloor\frac{n}{2}\rfloor}$ be the dominant weight with respect to $\frak{b}^c$ defined in (\ref{E:parameter}). {\color{black}Then} $\underline{\kappa}-\underline{e}$ is dominant and $\tau_{\underline{\kappa}-\underline{e}}$ appears with multiplicity one in both $M_{\mu,\C}$ and $M_{\mu,\C}^\vee$.
%For $1 \leq i \leq \lfloor\frac{n}{2}\rfloor$, let $Y_i = -k_\infty\cdot E_{i,i+1}\cdot k_\infty^{-1} \in \frak{l}^c$.

\begin{prop}\label{L:cohomology}
Let $\pi \in \Omega_\mu \subset \Omega(n)$ with infinity type $(\underline{\kappa};\,\underline{{\sf w}})$. Then 
\[
{\rm Hom}_{{\rm SO}_n(\R)}(\tau_{\underline{e}(n)}\otimes\tau_{\underline{\kappa}-\underline{e}(n)},\pi)\subset {\rm Hom}_{K_n^\circ}\left(\extp^{b_n}(\frak{g}_n/\frak{k}_n)\otimes M_{\mu,\C}^\vee,\,\pi\right)
\]
is one-dimensional and contributes non-trivially to $H^{b_n}(\frak{g}_n,K_n^\circ;\pi\otimes M_{\mu,\C})$.
\end{prop}

\begin{proof}
Let $r = \lfloor\tfrac{n}{2}\rfloor$ and $P=P_{(2,...,2)}$ or $P=P_{(2,...,2,1)}$ be the standard parabolic subgroup of $\GL_n$ of type $(2,...,2)$ or $(2,...,2,1)$ depending on the parity of $n$. 
Then we have
\begin{align}\label{E:cohomology pf 1}
\pi\cong{\rm Ind}_{P(\R)}^{\GL_n(\R)}(\pi_1 \otimes \cdots\otimes \pi_k),
\end{align}
where $k = r$ (resp.\,$k=r+1$) if $n$ is even (resp.\,odd) and
\[
\pi_i = 
\begin{cases}
D_{\kappa_i}\otimes|\mbox{ }|^{{\sf w}_i/2} & \mbox{ if $1 \leq i \leq r$},\\ 
{\rm sgn}^\delta|\mbox{ }|^{{\sf w}_{r+1}/2} & \mbox{ if $n$ is odd, $i=r+1$, and $\pi = \pi_\mu\otimes{\rm sgn}^\delta$}.
\end{cases}
\]
In particular, we see that $\tau_{\underline{\kappa}}$ is a minimal ${\rm SO}_n(\R)$-type of $\pi$ appears with multiplicity one (cf.\,\cite[Lemma 6.1.1]{Mahnkopf2005}).
Also it is clear that $\tau_{\underline{\kappa}}$ appears with multiplicity one in $\tau_{\underline{e}} \otimes \tau_{\underline{\kappa} - \underline{e}}$.
Therefore, $\tau_{\underline{\kappa}}$ is the unique ${\rm SO}_n(\R)$-type that appears in both $\pi$ and $\tau_{\underline{e}} \otimes \tau_{\underline{\kappa} - \underline{e}}$, which implies that ${\rm Hom}_{{\rm SO}_n(\R)}(\tau_{\underline{e}}\otimes \tau_{\underline{\kappa}-\underline{e}},\pi)$ is one-dimensional.
To prove the assertion, by Proposition \ref{L:Delorme} it suffices to show that ${\rm Hom}_{{\rm SO}_n(\R)}(\tau_{\underline{e}}\otimes \tau_{\underline{\kappa}-\underline{e}},\pi)$ is contained in the image of the injective homomorphism $\iota_{\otimes_{i=1}^k{\pi_i}}$ therein. 
The $M_P(\C)$-representation 
\[
W := \left(\bigotimes_{i=1}^k \extp^{b_{n_i}}(\frak{g}_{n_i}/\frak{k}_{n_i})\right) \otimes \left(\bigotimes_{i=1}^kM_{\mu^{(i)},\C}^\vee\right)
\]
decomposes into characters of ${\rm SO}_2(\R)^r$.
The weight $\underline{\kappa}$ appears with multiplicity one with an eigenvector 
\[
v_1^+ \otimes v_2^+ :=
\begin{cases}
\left( \bigotimes_{i=1}^rY_+ \right) \otimes \left( \bigotimes_{i=1}^r(x+\sqrt{-1}\,y)^{\kappa_i-2}\right) & \mbox{ if $n$ is even},\\
\left( \bigotimes_{i=1}^rY_+ \otimes 1 \right) \otimes \left( \bigotimes_{i=1}^r(x+\sqrt{-1}\,y)^{\kappa_i-2} \otimes 1\right) & \mbox{ if $n$ is odd}.
\end{cases}
\]
We refer to \S\,\ref{SS:generator} for the notation.
Fix a non-zero vector $v_3^+$ in $\otimes_{i=1}^k (\pi_i\otimes\delta_i)$ of weight $\underline{\kappa}$. 
Let 
\[
f^+ \in {\rm Hom}_{{\rm SO}_2(\R)^r}\left(W,\,\otimes_{i=1}^k (\pi_i\otimes\delta_i)\right)
\]
defined as follows: For $v \in W$, there exists a unique $z \in \C$ such that $v-z\cdot(v_1^+\otimes v_2^+)$ is a sum of eigenvectors of weight other than $\underline{\kappa}$. Then $f^+(v) := z\cdot v_3^+$. 
%It is clear that $f^+$ belongs to the image of the natural isomorphism
%\[
%\bigotimes_{i=1}^k H^{b_{n_i}}(\frak{g}_{n_i},K_{n_i}^\circ;(\pi_i\otimes \delta_i) \otimes M_{\mu^{(i)},\C}) \longrightarrow {\rm Hom}_{{\rm SO}_2(\R)^r}\left(W, \,\otimes_{i=1}^k (\pi_i\otimes\delta_i)\right).
%\]
Let
\[
f^\circ:=\sum_{a \in \pi_0(K_n^{M_P})}a\cdot f^+ \in {\rm Hom}_{{\rm SO}_2(\R)^r}\left(W,\,\otimes_{i=1}^k (\pi_i\otimes\delta_i)\right)^{\pi_0(K_n^{M_P})}.
\]
Now we show that 
\[
\C\cdot \iota_{\otimes_{i=1}^k{\pi_i}}(f^\circ) = {\rm Hom}_{{\rm SO}_n(\R)}(\tau_{\underline{e}}\otimes \tau_{\underline{\kappa}-\underline{e}},\pi).
\]
Since the right-hand side is one-dimensional, it suffices to show that $\iota_{\otimes_{i=1}^k{\pi_i}}(f^\circ)$ is contained in the right-hand side.
By definition, $\iota_{\otimes_{i=1}^k{\pi_i}}$ is the restriction of the inverse of (\ref{E:Delorme proof 1}) to ${\rm Hom}_{{\rm SO}_2(\R)^r}\left(W,\,\otimes_{i=1}^k (\pi_i\otimes\delta_i)\right)^{\pi_0(K_n^{M_P})}$, where $\pi$ is realized as the induced representation (\ref{E:cohomology pf 1}).
We have 
\[
\iota_{\otimes_{i=1}^k{\pi_i}}(f^\circ)(v)(k) = f^\circ(k\cdot v),\quad v \in W,\,k\in{\rm SO}_n(\R).
\]
Therefore, the support of $\iota_{\otimes_{i=1}^k{\pi_i}}(f^\circ)$ is contained in the representation of ${\rm SO}_n(\R)$ generated by $v_1^+\otimes v_2^+$.
Thus we are reduced to prove that $v_1^+\otimes v_2^+ \in \tau_{\underline{e}} \otimes \tau_{\underline{\kappa} - \underline{e}}$.
In our case, the Kostant representative $w \in W_n^P$ in Lemma \ref{L:combinatorial} is given by
\[
w = \begin{cases}
\bp 1 &2& \cdots& r &r+1 &\cdots& 2r-1& 2r \\
1 & 3 & \cdots & 2r-1 & 2r & \cdots & 4 & 2
\ep & \mbox{ if $n$ is even},\\
\bp 1 &2& \cdots& r & r+1 &r+2&\cdots& 2r& 2r+1 \\
1 & 3 & \cdots & 2r-1 &2r+1 & 2r & \cdots & 4 & 2
\ep & \mbox{ if $n$ is odd}.
\end{cases}
\]
Note that 
\[
w(\mu+\rho_n)-\rho_n = (\mu^{(1)},...,\mu^{(k)}),\quad \ell(w) = b_n-r = \begin{cases}
r(r-1) & \mbox{ if $n$ is even},\\
r^2 & \mbox{ if $n$ is odd}.
\end{cases}
\]
In the definition of $\iota_{\otimes_{i=1}^k{\pi_i}}$, the embedding of $W$ into $\extp^{b_n}(\frak{g}/\frak{k})\otimes M_{\mu,\C}^\vee$ depends on an $M_P(\C)$-equivariant embedding from %$\otimes_{i=1}^k M_{\mu^{(i)},\C}$ to ${\rm Hom}\left(\extp^{\ell(w)}\frak{u}_P,M_{\mu,\C}\right)$.
%Since it appears with multiplicity one (cf.\,\cite[Corollary 5.7 and Theorem 5.14]{Kostant1961}), it is equivalent to fix an $M_P(\C)$-equivariant embedding from 
$\otimes_{i=1}^k M_{\mu^{(i)},\C}^\vee$ to $\extp^{\ell(w)}\frak{u}_P \otimes M_{\mu,\C}^\vee$. 
Up to $\C^\times$, such an embedding is determined by the following association between highest weight vectors for $\prod_{i=1}^kT_{n_i}(\C)$ with respect to $M_P(\C)\cap U_n(\C)$:
\begin{align*}
k_\infty^{-1}v_2^+ = 2^{\sum_{i=1}^r(\kappa_i-2)/2}\cdot\begin{cases}
\bigotimes_{i=1}^r x^{\kappa_i-2} & \mbox{ if $n$ is even}\\
\bigotimes_{i=1}^r x^{\kappa_i-2} \otimes 1 & \mbox{ if $n$ is odd}
\end{cases} \,\longmapsto\, \extp_{1 \leq i \leq r \atop 2i<j\leq n} e_{2i-1,j} \otimes w\cdot v_{\mu^\vee},
\end{align*}
where $v_{\mu^\vee}$ is a highest weight vector of $M_{\mu,\C}^\vee$.
Indeed, it is easy to verify that $\extp_{1 \leq i \leq r \atop 2i<j\leq n} e_{2i-1,j}$ has weight $(w\rho_n-\rho_n)^\vee$ and is invariant by $M_P(\C) \cap U_n(\C)$.
Also $w\cdot v_{\mu^\vee}$ has weight $w\mu^\vee$ and is invariant by $M_P(\C) \cap U_n(\C)$ since $w^{-1}\cdot\prod_{i=1}^kX^+(T_{n_i}) \subset X^+(T_n)$.
Therefore, $\extp_{1 \leq i \leq r \atop 2i<j\leq n} e_{2i-1,j} \otimes w\cdot v_{\mu^\vee}$ is a highest weight vector of weight
\[
(w\rho_n-\rho_n)^\vee + w\mu^\vee = \left((\mu^{(1)})^\vee,...,(\mu^{(k)})^\vee\right).
\]
Note that
\[
Y_+ = \frac{1}{2}\bp 1& \sqrt{-1} \\ \sqrt{-1} & 1\ep \bp 0&1\\0&0 \ep \bp 1& -\sqrt{-1} \\ -\sqrt{-1} & 1\ep.
\]
Therefore, under the inclusion $\oplus_{i=1}^k\frak{g}_{n_i} = \frak{m}_P \subset \frak{g}$, the wedge vector $v_1^+$ is represented by $\extp_{i=1}^{r} X_{2i-1,2i}$.
Also note that for $1 \leq i \leq r$ and $2i<j\leq n$, we have
\[
X_{2i-1,j} = \begin{cases}
{}^tX_{j-1,2i} & \mbox{ if $j$ is even},\\
{}^tX_{j+1,2i} & \mbox{ if $j$ is odd and $j<n$},\\
{}^tX_{n,2i} & \mbox{ if $j$ is odd and $j=n$}.
\end{cases}
\]
We thus conclude that the image of $v_1^+\otimes v_2^+$ in $\extp^{b_n}(\frak{g}/\frak{k})\otimes M_{\mu,\C}^\vee$ is equal to
\[
\left(\extp_{i=1}^{r} X_{2i-1,2i} \extp_{1 \leq i \leq r \atop 2i<j\leq n} {}^tX_{j,2i} \right) \otimes k_\infty w \cdot v_{\mu^\vee},
\]
where for $X \in \frak{g}$ we denote its image in $\frak{g}/\frak{k}$ also by $X$. It is clear that the above wedge vector in $\extp^{b_n}(\frak{g}/\frak{k})$ and $k_\infty w \cdot v_{\mu^\vee}\in M_{\mu,\C}^\vee$ have weights $\underline{e}$ and $\underline{\kappa}-\underline{e}$ respectively under the action of $\frak{l}^c$. Finally, they are also vanished by $\frak{n}^c$%(in fact, by $\frak{n}$),
, which can be verified using the fact that
\[
k_\infty w \cdot \frak{u} \cdot w^{-1}k_\infty^{-1} = \frak{n}.
\]
This completes the proof.
\end{proof}

\end{document}